\input amstex
\documentstyle{amsppt}
\magnification=\magstep1
\vsize =21 true cm
\hsize =16 true cm
\loadmsbm
\topmatter

\centerline{\bf Geometry and arithmetic associated to}
\centerline{\bf Appell hypergeometric partial differential
                equations}
\author{\smc Lei Yang}\endauthor
\endtopmatter
\document

\centerline{\bf Contents}
$$\aligned
 &\text{1. Introduction}\\
 &\text{2. Four derivatives and their properties}\\
 &\text{3. Nonlinear partial differential equations and modular
           functions associated to $U(2, 1)$}\\
 &\text{4. $\eta$-functions associated to $U(2, 1)$ and Picard curves}\\
 &\text{5. Transform problems and modular equations associated to algebraic surfaces}\\
 &\text{6. Triangular $s$-functions associated to $U(2, 1)$ and Picard curves}\\
 &\text{7. Four derivatives and nonlinear evolution equations}
\endaligned$$

\vskip 0.5 cm
\centerline{\bf 1. Introduction}
\vskip 0.5 cm

  Mathematicians often speak of the unity of their subject. Such a
unity from the 19th century origins in the history of linear
ordinary differential equations, especially those closely
connected with elliptic and modular functions. It is concerned
with the impact of group theoretical and geometrical ideas upon
the problem of understanding the nature of the solutions to a
differential equation (see \cite{Kl2}, \cite{Kl3} and \cite{G} as
excellent surveys).

  The theory of hypergeometric functions has a long history (see
\cite{Ho1}). The first occurrence of hypergeometric series goes
back to the 17th century. Euler was the first to consider the
series as a function of one variable $x$:
$$F(a, b; c; x)=1+\frac{a \cdot b}{1 \cdot 2} x+\frac{a(a+1) \cdot
  b(b+1)}{1 \cdot 2 \cdot c(c+1)} x^2+\cdots.$$
Euler found that $F$ is a solution of the hypergeometric
differential equation
$$x(1-x) y^{\prime \prime}+[c-(a+b+1)x] y^{\prime}-ab y=0$$
and up to a constant factor $F$ can be expressed as a
hypergeometric integral
$$\int_{0}^{1} t^{a-1} (1-t)^{c-a-1} (1-xt)^{-b} dt.$$

  The next developments are closely connected with the name of
Gauss, who considered the hypergeometric function as a function of
a complex variable. In a letter to Schumacher (September 17,
1808), Gauss outlined his new modern point of view (Gauss'
programme) as follows:
\roster
\item
  Proof of the existence of functions which satisfy our
  differential equation.
\item
  Investigation of the properties of these functions. We look
  upon every differential equation as defining a specific
  (possibly new) class of transcendental functions.
\endroster

  For this project the work of many outstanding mathematicians
was needed: Abel, Cauchy, Dedekind, Fuchs, Jacobi, Klein, Kummer,
Poincar\'{e}, Riemann, Schwarz, Weierstrass and others. The theory
of hypergeometric functions was a central theme in the function
theory of one complex variable. Now, we will give a survey as
follows:

  In 1812, Gauss studied the hypergeometric equation
$$x(1-x) \frac{d^2 y}{d x^2}+[c-(a+b+1)x] \frac{dy}{dx}-ab y=0.$$
This equation is important in its own right as a linear ordinary
differential equation for which explicit power-series solutions
can be given. Gauss' work in this direction was extended by
Riemann in 1857, who introduced the crucial idea of analytically
continuing the solutions around their singularities in the complex
domain. Riemann's monodromy approach gave a thorough global
understanding of the solutions of the hypergeometric equation. In
1872, Schwarz studied the hypergeometric equations all of whose
solutions were algebraic functions (see \cite{Sch}). It led to
Schwarz's discovery of a new class of transcendental functions
associated with the hypergeometric equation which turned out to be
of vital importance in the theory of automorphic functions. The
concept of of Schwarzian derivatives played a central role in his
study.

  It is known that the Schwarzian derivative
$$\{w; x\}:=\frac{w^{\prime \prime \prime}}{w^{\prime}}-
  \frac{3}{2} \left(\frac{w^{\prime \prime}}{w^{\prime}}\right)^{2}$$
satisfies that
$$\left\{ \frac{aw+b}{cw+d}; x \right\}=\{w; x\} \quad \text{for}
  \quad \left(\matrix a & b\\ c & d \endmatrix\right) \in GL(2,
  {\Bbb C}),$$
where $x$ is a complex number and $w$ is a complex-valued function
of $x$.

  The Schwarzian derivative enters into several branches of
complex analysis. It plays an important role in the theory of
linear second order differential equations, in conformal mapping,
and in the study of univalence and of uniformization (see
\cite{Hi}, \cite{Pa1} and \cite{Pa2}).

  The Schwarzian derivative also appears in the conformal field
theory (see \cite{Po}), which plays an important role in
mathematics, statistical mechanics and string theory. Some
infinite dimensional algebras are related to the conformal field
theory, such as Virasoro algebra, Kac-Moody Lie algebra and Hopf
algebra.

  Let us recall some basic facts about the Virasoro algebra (see
\cite{BMS} and \cite{Wi}). For $A={\Bbb C}[t, t^{-1}]$,
$\partial=\frac{d}{dt}$, the classical residue
$$\text{res}: \Omega \to {\Bbb C}, \quad
  \text{res} \left(\sum_{i} c_i t^i dt\right)=c_{-1}$$
defines the isomorphism $\Omega / dA={\Bbb C}$. The Virasoro
algebra $V$ is a central extension of the Lie algebra $A
\frac{d}{dt}$ by ${\Bbb C} z_0$ with the commutation rule
$$\left[f \frac{d}{dt}, g \frac{d}{dt}\right]_{z_0}=
  \left(f \frac{dg}{dt}-g \frac{df}{dt}\right) \frac{d}{dt}+
  \frac{1}{12} \text{res} g \frac{d^3 f}{dt^3} dt \cdot z_0.$$
The element $z_0$ or its proper value in a representation space of
$V$ is called the central charge.

  The list of fascinating aspects of the Korteweg-de Vries
equation is long. The well-known Miura transformations, Schwarzian
derivatives and the Hill operator are three examples from this
list. It is known that (see \cite{KN}, \cite{SS}, \cite{Do},
\cite{Wei1}, \cite{Wei2} and \cite{Wil}) the evolution equation
$$u_{t}=u_{xxx}-\frac{3}{2} u_{xx}^2 u_{x}^{-1}$$
for a function $u(x, t)$ is invariant with respect to the group
$PSL(2)$ of linear fractional transformations
$$u \mapsto \frac{a u+b}{c u+d},$$
and implies that the Korteweg-de Vries equation
$$v_{t}=v_{xxx}+6 v v_{x}$$
for the variable
$$v=\frac{1}{2}\left(u_{xxx} u_{x}^{-1}-\frac{3}{2} u_{xx}^2
    u_{x}^{-2}\right)=\frac{1}{2} \{u, x\}.$$
The above evolution equation is equivalent to
$$\frac{u_{t}}{u_{x}}=\{u, x\}.$$

  The hypergeometric equation is also important because it contains
many interesting equations as special cases. We are interested in
the Legendre's equation:
$$(1-\kappa^2) \frac{d^2 y}{d \kappa^2}+\frac{1-3 \kappa^2}{\kappa}
  \frac{dy}{d \kappa}-y=0.$$
There are two linearly independent solutions:
$$K=\int_{0}^{1} \frac{dx}{\sqrt{(1-x^2)(1-\kappa^2 x^2)}}, \quad
  K^{\prime}=\int_{1}^{\frac{1}{\kappa}}
  \frac{dx}{\sqrt{(1-x^2)(1-\kappa^2 x^2)}},$$
which are the periods of elliptic integrals considered as
functions of the modulus $\kappa^2$. This equation was studied by
Legendre, Gauss, Abel, Jacobi and Kummer. It is important in the
study of the arithmetical theory of modular functions.

  The arithmetical theory of modular functions is a result of the
synthesis of numerous ideas of many outstanding mathematicians
(see \cite{V} for an exciting survey). The complex multiplication
of elliptic functions is one of the main sources of this theory.
Complex multiplication theory is usually associated with the name
of Kronecker, who contributed much to it. Some of the ideas of
complex multiplication theory can be found in earlier papers by
Gauss, Abel and Eisenstein (see \cite{V} and \cite{Lem}).

  The transform problem for elliptic functions, which is related
to the study of isogenies of elliptic curves, appeared on the
historic stage in the following way. In 1750's, Fagnano and Euler
discovered a transform of $\frac{dx}{\sqrt{F(x)}}$ to
$\frac{dy}{\sqrt{G(y)}}$, where $F$ and $G$ are polynomials of
degree four. This transform is given by the equation $f(x, y)=0$,
where $f$ is a polynomial of second degree in each variable. In
1775, Landen found the Landen transform which is of order two. In
1824, Legendre discovered a transform of order three. However,
transforms of orders $3$, $5$ and $7$ had been known to Gauss
since 1808.

  The general transform problem for an arbitrary order was stated
by Jacobi (see \cite{W} for the details). He proved the existence
of a transform of any odd order $n$, namely $y=U(x)/V(x)$, where
$U$ and $V$ are mutually prime polynomials, $U$ is an odd
polynomial of degree $n$, and $V$ is an even polynomial of degree
$n-1$, and
$$\frac{dy}{\sqrt{(1-y^2)(1-\lambda^2 y^2)}}=\frac{1}{M}
  \frac{dx}{\sqrt{(1-x^2)(1-\kappa^2 x^2)}},$$
where $M$ is the so-called multiplicator and $\lambda$ is the
transformed modulus ($\kappa$ being given). This leads to the
study of modular equations which are of great importance in the
theory of elliptic functions (see \cite{V}). Let $n$ be an odd
prime. The roots of the modular equation are the $(n+1)$ values of
moduli $\lambda$ which are images of the modulus $\kappa$ via all
the maps of order $n$. For $n=3$, Legendre obtained (1824) that
the modular equation has the form: $\sqrt{\kappa
\lambda}+\sqrt{\kappa^{\prime} \lambda^{\prime}}=1$, where
$\kappa^2+{\kappa^{\prime}}^2=1$ and
$\lambda^2+{\lambda^{\prime}}^2=1$. Jacobi derived the modular
equation for $n=3, 5$ (see \cite{W}). In 1832, Galois claimed that
the Galois group of the modular equation is $PGL(2, {\Bbb
F}_{n})$, containing the subgroup $PSL(2, {\Bbb F}_{n})$ of index
$2$ and of order $\frac{1}{2} n(n^2-1)$. For and only for $n=5, 7,
11$, $PSL(2, {\Bbb F}_{n})$ contains a subgroup of index $n$.
Hence in these cases the modular equation can be reduced to an
equation of degree $n$. In 1858, Hermite (see \cite{He}, pp.5-12),
Kronecker and Brioschi used the modular equation to solve
algebraic equations of degree four and five. A general equation of
degree four is reduced to the modular equation of order three, a
general equation of degree five is reduced to the modular equation
of degree five. Klein (1879) studied the equation with the group
$PGL(2, {\Bbb F}_{7})$, the Galois group of the modular equation
of degree seven.

  Klein's work grew out of an attempt to reformulate the theory of
modular functions and modular transformations without using the
theory of elliptic functions. This had been almost completely
achieved by Dedekind in his celebrated paper of 1877 (see
\cite{De}). Klein (see \cite{Kl1}) enriched Dedekind's approach
with an explicit use of group theoretic ideas and an essentially
Galois theoretic approach to the fields of rational and modular
functions. Both Dedekind (see \cite{De}) and Klein (see
\cite{Kl1}) relied on the theory of the hypergeometric equation at
technical points in their arguments. In particular, for $0, 1,
\lambda, \infty \in {\Bbb P}^{1}({\Bbb C})$, $\lambda \neq 0, 1$,
Klein (see \cite{Kl1}) set
$$J=J(\lambda)=\frac{4}{27} \frac{(\lambda^2-\lambda+1)^{3}}{\lambda^2
    (\lambda-1)^2}.$$

  It is well-known that the Gauss hypergeometric function
$$F(a, b; c; x)=\frac{\Gamma(c)}{\Gamma(a) \Gamma(c-a)} \int_{0}^{1}
  t^{a-1} (1-t)^{c-a-1} (1-tx)^{-b} dt$$
satisfies the following equation:
$$x(1-x) \frac{d^2 y}{dx^2}+[c-(a+b+1)x] \frac{dy}{dx}-aby=0.$$

  For the group $GL(2)$, the elliptic integral
$$y=\int_{0}^{1} \frac{dt}{\sqrt{t(t-1)(t-x)}}$$
satisfies the following equation:
$$x(x-1) \frac{d^2 y}{dx^2}+(2x-1) \frac{dy}{dx}+\frac{1}{4} y=0.$$
Thus, a solution is given by
$$y=F\left(\frac{1}{2}, \frac{1}{2}; 1; x\right), \quad
  (a, b, c)=\left(\frac{1}{2}, \frac{1}{2}, 1\right).$$
The elliptic modular function is associated to the following
equation (see \cite{De} and \cite{Kl1}):
$$x(1-x) \frac{d^2 y}{dx^2}+\left(\frac{2}{3}-\frac{7}{6}x\right) \frac{dy}{dx}-\frac{1}{144} y=0.$$
A solution is given by
$$y=F\left(\frac{1}{12}, \frac{1}{12}; \frac{2}{3}; x\right), \quad
  (a, b, c)=\left(\frac{1}{12}, \frac{1}{12}, \frac{2}{3}\right).$$

  Poincar\'{e}'s papers (see \cite{P}) in the 1880's on Fuchsian and
Kleinian groups are of great interest. The papers are the climax
of the geometric theory of functions initiated by Riemann, and the
ideal representatives of the unity between analysis, geometry,
topology and algebra. In particular, Poincar\'{e} proved two
theorems:

{\smc Theorem 1}. {\it Between two Fuchsian functions with the
same group and without any essential singular point which is not a
consequence of their definition, there is an algebraic relation.}

{\smc Theorem 2}. {\it Any Fuchsian function $F(z)$ permits the
integration of a linear equation with algebraic coefficients in
the following way. If we put
$$x=F(z), \quad y_1=\sqrt{\frac{dF}{dz}}, \quad
  y_2=z \sqrt{\frac{dF}{dz}},$$
then $y_1$ and $y_2$ satisfy the differential equation
$$\frac{d^2 y}{dx^2}=y \phi(x),$$
where $\phi(x)$ is algebraic in $x$.}

  In fact, the hypergeometric functions are heavily involved in
special functions, representation theory, theory of automorphic
forms, elliptic curves and integrals and more algebraic geometry,
and number theory.

  In 1880 the two variable story began when Appell gave a
generalization of the hypergeometric equation as an equation in
two variables, and found four double series generalizing the
hypergeometric series which he called $F_1$, $F_2$, $F_3$ and
$F_4$. Appell showed that each of these new functions satisfies
two simultaneous partial differential equations (see \cite{AK}).
In particular, the function $F_1=F_1(a; b, b^{\prime}; c; x, y)$
satisfies the following partial differential equations:
$$\left\{\aligned
  x(1-x) \frac{\partial^2 z}{\partial x^2}+y(1-x) \frac{\partial^2
  z}{\partial x \partial y}+[c-(a+b+1)x] \frac{\partial
  z}{\partial x}-by \frac{\partial z}{\partial y}-abz &=0,\\
  y(1-y) \frac{\partial^2 z}{\partial y^2}+x(1-y) \frac{\partial^2
  z}{\partial x \partial y}+[c-(a+b^{\prime}+1)y] \frac{\partial
  z}{\partial y}-b^{\prime} x \frac{\partial z}{\partial x}-a
  b^{\prime} z &=0.
\endaligned\right.$$

  Picard (see \cite{Pi1}, \cite{Pi2}, \cite{Pi3} and \cite{Pi4})
studied the Appell hypergeometric partial differential equations
and the hyperfuchsian functions, which led him to investigate the
theory of algebraic surfaces. In particular, in his paper
\cite{Pi2}, Picard generalized Poincar\'{e}'s Theorem 1 and
Theorem 2 as follows:

{\smc Theorem 3}. {\it Suppose that three hyperfuchsian functions
$$x=F_1(\xi, \eta), \quad y=F_2(\xi, \eta), \quad z=F(\xi, \eta)$$
satisfy the algebraic relation
$$f(x, y, z)=0.$$
One can form a system of three partial differential equations
$$\left\{\aligned
  \frac{\partial^2 z}{\partial x^2} &=a \frac{\partial z}{\partial x}
  +b \frac{\partial z}{\partial y}+c z,\\
  \frac{\partial^2 z}{\partial x \partial y} &=a_1 \frac{\partial z}
  {\partial x}+b_1 \frac{\partial z}{\partial y}+c_1 z,\\
  \frac{\partial^2 z}{\partial y^2} &=a_2 \frac{\partial z}{\partial x}
  +b_2 \frac{\partial z}{\partial y}+c_2 z,
\endaligned\right.$$
where $a, b, c$ are rational functions of $x$, $y$ and $z$ ($z$ is
an algebraic function of $x$, $y$ defined by $f(x, y, z)=0$.)
These three equations have three linearly independent solutions
$z_1$, $z_2$, $z_3$. If one sets
$$\frac{z_2}{z_1}=\xi, \quad \frac{z_3}{z_1}=\eta,$$
these equations can be resolved by $x$, $y$ as follows:
$$x=F_1(\xi, \eta), \quad y=F_2(\xi, \eta).$$}

  In 1921, Giraud \cite{Gir} studied the two variable automorphic
functions and the complex hyperbolic geometry. According to
\cite{Gir} and \cite{Pi2}, set $\xi=\frac{z_2}{z_1}$ and
$\eta=\frac{z_3}{z_1}$, we have
$$\vmatrix \format \c \quad & \c \quad & \c\\
  \frac{\partial z_1}{\partial x} & \frac{\partial z_1}{\partial y} &z_1\\
  \frac{\partial z_2}{\partial x} & \frac{\partial z_2}{\partial y} &z_2\\
  \frac{\partial z_3}{\partial x} & \frac{\partial z_3}{\partial y} &z_3
  \endvmatrix
 =z_1^2 \vmatrix \format \c \quad & \c \quad & \c\\
  0 & 0 & z_1\\
  \frac{\partial \left(\frac{z_2}{z_1}\right)}{\partial x}
 &\frac{\partial \left(\frac{z_2}{z_1}\right)}{\partial y} & z_2\\
  \frac{\partial \left(\frac{z_3}{z_1}\right)}{\partial x}
 &\frac{\partial \left(\frac{z_3}{z_1}\right)}{\partial y} & z_3
 \endvmatrix=z_1^3 \frac{\partial(\xi, \eta)}{\partial(x, y)}
=1,$$ when $z_1$, $z_2$ and $z_3$ satisfy the following relations:
$$z_1=\root 3 \of{\frac{\partial(x, y)}{\partial(\xi, \eta)}}, \quad
  z_2=\xi \root 3 \of{\frac{\partial(x, y)}{\partial(\xi, \eta)}}, \quad
  z_3=\eta \root 3 \of{\frac{\partial(x, y)}{\partial(\xi, \eta)}}.$$

  In 1964, Shimura (see \cite{Sh1}) considered the classifying spaces
for the isomorphy classes of plane curves of the following types:
\roster
\item $y^3=p_4(x)$ (non-hyperelliptic),
\item $y^4=p_2(x) \cdot p_3(x)^2$ (hyperelliptic),
\endroster
where $p_k(x)$ is a polynomial of degree $k$ and without multiple
root. These are curves of genus $3$. By a result of Shimura (see
\cite{Sh1}) these curves are essentially parametrized by \roster
\item ${\Bbb B}^2/U(2, 1; {\Bbb Z}[\omega])$,
      $\omega=\frac{-1+\sqrt{-3}}{2}$,
\item ${\Bbb B}^2/U(2, 1; {\Bbb Z}[i])$, $i=\sqrt{-1}$,
\endroster
where ${\Bbb B}^2=\{(z_1, z_2) \in {\Bbb C}^2: |z_1|^2+|z_2|^2<1
\}$ is the complex ball. Here, $U(2, 1)=\{g \in GL(3, {\Bbb C}):
g^{*} I_{2, 1} g=I_{2, 1}\}$ with $I_{2, 1}=\text{diag}(1, 1,
-1)$. $U(2, 1; {\Bbb Z}[\omega])=\{ g=(a_{ij}) \in U(2, 1): a_{ij}
\in {\Bbb Z}[\omega] \}$, $U(2, 1; {\Bbb Z}[i])=\{ g=(a_{ij}) \in
U(2, 1): a_{ij} \in {\Bbb Z}[i] \}$.

  Moreover, Shimura (see \cite{Sh1}) proved that in these two cases,
the field of all automorphic functions on ${\Bbb B}^2$ with
respect to $\Gamma=U(2, 1; {\Bbb Z}[\omega])$ or $\Gamma=U(2, 1;
{\Bbb Z}[i])$ is a purely transcendental extension of ${\Bbb C}$
of dimension $2$.

  In fact, Picard (see \cite{Pi1}) investigated the curve $y^3=p_4(x)$
(Picard curve) and observed that moduli of such curves give
automorphic functions on ${\Bbb B}^2$. Heegner (see \cite{H})
studied the abelian integrals associated to these curves and the
transformations of automorphic functions.

  In his monographs \cite{Ho1}, \cite{Ho2}, Holzapfel studied a
concrete hypergeometric object of complex dimension two and the
complete solution of the associated system of Euler partial
differential equations. Holzapfel's results and methods involve
several parts of mathematics, notably algebraic geometry (theory
of surfaces, Gauss-Manin connection), theory of automorphic forms,
monodromy groups of differential equations and number theory.

  In their papers \cite{CW1} and \cite{CW2}, Cohen and Wolfart
classified the algebraic Appell-Lauricella functions $F_1$ in
several variables, which extended the famous 1872 list of Schwarz
of the algebraic Gauss hypergeometric functions (see \cite{Sch}).
The monodromy groups of the Appell-Lauricella functions $F_1$ in
several variables can be discontinuous, but this does not in
general entail that they be arithmetically defined (see
\cite{DM1}, \cite{M} and \cite{DM2}). Nonetheless, Cohen and
Wolfart associated to them families of abelian varieties for which
they play a role similar to that of modular groups.

  In fact, the geometric objects associated to Appell hypergeometric
partial differential equations are significant in automorphic
representation theory and arithmetic algebraic geometry (see the
monograph \cite{LR} as an excellent survey, Rogawski's work
\cite{R}, \cite{GRS} and Shimura's work \cite{Sh2}, \cite{Sh4}).
According to \cite{LR}, the Picard surfaces have an immediate
geometric appeal, their connected components being quotients of
two-dimensional balls by discrete groups were introduced by Picard
during his study of the monodromy of Appell hypergeometric partial
differential equations. We are interested in the groups which are
defined by discrete groups associated to unitary groups in three
variables over imaginary quadratic extensions of ${\Bbb Q}$. The
quotients are non-compact. The corresponding varieties are
investigated in \cite{LR}. As for many Shimura varieties the basic
arithmetic properties of the Picard surfaces are most easily
presented as aspects of the solution of a moduli problem.

  In fact, after Kronecker, the history of the development of
complex multiplication theory at the end of the 19th century is
mostly connected with the name of Weber. Takagi, Fueter and Hasse
devoted to the proof of Kronecker's Jugendtraum. Deuring
introduced the algebraic approach to complex multiplication (see
\cite{V}). The complex multiplication theory which realizes the
Jugendtraum is one-dimensional in the following sense: it
considers modular functions in one variable and one-dimensional
abelian varieties, that is, elliptic curves.

  In 1911, Hecke studied complex multiplications for Hilbert
modular functions in two variables. In 1955, some principal ideas
of multidimensional complex multiplication theory were introduced
by Weil and by Shimura and Taniyama. A thorough development of
these ideas was accomplished by Shimura and Taniyama (see
\cite{Sh3} and \cite{Sh5}). In a sense, the theory of Shimura
varieties can be considered as a far-reaching generalization of
complex multiplication theory.

  However, Ovseevi\v{c} (see \cite{O}) showed that if
$\text{dim} A \neq 1$, where $A$ is an abelian variety, that is,
if $[K: {\Bbb Q}] \neq 2$, then the analogue of the completeness
theorem is false: the maximal Shimura-Taniyama field does not
coincide with the maximal abelian extension of $K$. Thus if a
CM-field is not complex quadratic, abelian varieties give only a
partial solution of Hilbert's 12th problem.

  According to Langlands (see \cite{L}), the 12th problem of Hilbert
reminds us of the relationship of three subjects. The first is the
theory of class fields or abelian extensions of number fields. The
second is the algebraic theory of elliptic curves and more
generally, abelian varieties. The third is the theory of
automorphic functions which is connected with the study of abelian
varieties, especially of their moduli.

  In this paper, we will study the monodromy of Appell hypergeometric
partial differential equations, which lead us to find four
derivatives which are associated to the group $GL(3)$. Our four
derivatives have the remarkable properties, which play a role
similar to the Schwarzian derivative in the case of $GL(2)$ and
Gauss hypergeometric equation. We find that Appell hypergeometric
partial differential equations can be reduced to four nonlinear
partial differential equations by the use of our four derivatives.
We will investigate these four nonlinear partial differential
equations. We find that the Jacobian of the map which is defined
by the above nonlinear system satisfies a system of linear partial
differential equations.

  We are interested in some particular cases which grow out of
the maps on the complex hyperbolic space and the unitary group in
three variables $U(2, 1)=\{ g \in GL(3, {\Bbb C}): g^{*} J g=J \}$
where $J=\left(\matrix
    &    & 1\\
    & -1 &  \\
  1 &    &
   \endmatrix\right)$. The corresponding complex domain
${\frak S}_{2}=\{(z_1, z_2) \in {\Bbb C}^2: z_1+\overline{z_1}
-z_2 \overline{z_2}>0 \}$ is the simplest example of a bounded
symmetric domain which is not a tube domain. In these cases, the
system of linear partial differential equations has the remarkable
properties. Moreover, it can be solved by Appell hypergeometric
partial differential equations. This lead us to study the
properties of solutions. In particular, we find the
$\eta$-functions associated to the unitary group $U(2, 1)$, which
has the interesting properties. Unlike the Dedekind
$\eta$-functions in the case of $GL(2)$, our $\eta$-functions are
involved in the non-commutative translations.

  Our $\eta$-functions lead us to study the arithmetic of Picard
curves, especially the transform problems (moduli problems) and
the corresponding modular equations, which can be considered as a
natural generalization of complex multiplication. Hence, they are
related to Hilbert's 12th problem. Our approach is as follows:
instead of the Jacobian variety of Picard curves, which is an
abelian variety, we will investigate the Picard curves themselves.
We will give a rational transformation between two irrational
integrals which are associated to some algebraic surfaces of
degree seven. We find a modular equation which give the relation
between the moduli of two Picard curves. On the other hand, by the
investigation of the solutions, we find some interesting
transformations. The solutions have some remarkable properties
under these transformations. This leads us to obtain the
$J$-invariants associated to the Picard curves.

  Finally, we will study a system of nonlinear evolution equations
which are associated to our four derivatives. We find that they
are $GL(3)$-invariant. Moreover, our derivatives satisfy two first
order nonlinear partial differential equations, which play a role
similar to the Korteweg-de Vries equation in the case of $GL(2)$.

  Now, we will give a survey of our main results.

{\bf Main Results}.

  Now, we define the following four derivatives:

{\it Definition 1.1}. Four derivatives associated to the group
$GL(3)$ are given by:
$$\aligned
   \{u_1, u_2; x, y\}_{x}
&:=\frac{\vmatrix \format \c \quad & \c \\
   \frac{\partial u_1}{\partial x} &
   \frac{\partial u_2}{\partial x}\\
   \frac{\partial^2 u_1}{\partial x^2} &
   \frac{\partial^2 u_2}{\partial x^2}\endvmatrix}
  {\vmatrix \format \c \quad & \c\\
   \frac{\partial u_1}{\partial x} &
   \frac{\partial u_2}{\partial x}\\
   \frac{\partial u_1}{\partial y} &
   \frac{\partial u_2}{\partial y}
   \endvmatrix}, \quad \quad
   \{u_1, u_2; x, y\}_{y}
 :=\frac{\vmatrix \format \c \quad & \c \\
   \frac{\partial u_1}{\partial y} &
   \frac{\partial u_2}{\partial y}\\
   \frac{\partial^2 u_1}{\partial y^2} &
   \frac{\partial^2 u_2}{\partial y^2}\endvmatrix}
  {\vmatrix \format \c \quad & \c\\
   \frac{\partial u_1}{\partial y} &
   \frac{\partial u_2}{\partial y}\\
   \frac{\partial u_1}{\partial x} &
   \frac{\partial u_2}{\partial x}
   \endvmatrix},\\
   [u_1, u_2; x, y]_{x}
&:=\frac{\vmatrix \format \c \quad & \c \\
   \frac{\partial u_1}{\partial y} &
   \frac{\partial u_2}{\partial y}\\
   \frac{\partial^2 u_1}{\partial x^2} &
   \frac{\partial^2 u_2}{\partial x^2}\endvmatrix+2
   \vmatrix \format \c \quad & \c \\
   \frac{\partial u_1}{\partial x} &
   \frac{\partial u_2}{\partial x}\\
   \frac{\partial^2 u_1}{\partial x \partial y} &
   \frac{\partial^2 u_2}{\partial x \partial y}\endvmatrix}
  {\vmatrix \format \c \quad & \c\\
   \frac{\partial u_1}{\partial x} &
   \frac{\partial u_2}{\partial x}\\
   \frac{\partial u_1}{\partial y} &
   \frac{\partial u_2}{\partial y}
   \endvmatrix},\\
   [u_1, u_2; x, y]_{y}
&:=\frac{\vmatrix \format \c \quad & \c \\
   \frac{\partial u_1}{\partial x} &
   \frac{\partial u_2}{\partial x}\\
   \frac{\partial^2 u_1}{\partial y^2} &
   \frac{\partial^2 u_2}{\partial y^2}\endvmatrix+2
   \vmatrix \format \c \quad & \c \\
   \frac{\partial u_1}{\partial y} &
   \frac{\partial u_2}{\partial y}\\
   \frac{\partial^2 u_1}{\partial x \partial y} &
   \frac{\partial^2 u_2}{\partial x \partial y}\endvmatrix}
  {\vmatrix \format \c \quad & \c\\
   \frac{\partial u_1}{\partial y} &
   \frac{\partial u_2}{\partial y}\\
   \frac{\partial u_1}{\partial x} &
   \frac{\partial u_2}{\partial x}
   \endvmatrix},
\endaligned\tag 1.1$$
where $x$ and $y$ are complex numbers, $u_1$ and $u_2$ are
complex-valued functions of $x$ and $y$.

{\smc Proposition 1.2}. {\it Our derivatives are invariant under
linear fractional transformation acting on the first arguments:
$$\left\{\aligned
  \{v_1, v_2; x, y\}_{x} &=\{u_1, u_2; x, y\}_{x},\\
  \{v_1, v_2; x, y\}_{y} &=\{u_1, u_2; x, y\}_{y},\\
  [v_1, v_2; x, y]_{x} &=[u_1, u_2; x, y]_{x},\\
  [v_1, v_2; x, y]_{y} &=[u_1, u_2; x, y]_{y},
\endaligned\right.\tag 1.2$$
where
$$(v_1, v_2)=\gamma(u_1, u_2)=\left(\frac{a_1 u_1+a_2 u_2+a_3}
  {c_1 u_1+c_2 u_2+c_3}, \frac{b_1 u_1+b_2 u_2+b_3}{c_1 u_1+c_2
  u_2+c_3}\right)\tag 1.3$$
for $\gamma=\left(\matrix a_1 & a_2 & a_3\\ b_1 & b_2 & b_3\\
c_1 & c_2 & c_3 \endmatrix\right) \in GL(3, {\Bbb C})$.}

  Thus our derivatives are differential invariants of the general
linear group acting as a linear fractional transformation.

{\smc Corollary 1.3}. {\it Our derivatives vanish when $(w_1,
w_2)$ are linear fractional functions of $(z_1, z_2)$:
$$\left\{\aligned
  \{w_1, w_2; z_1, z_2\}_{z_1} &=0,\\
  \{w_1, w_2; z_1, z_2\}_{z_2} &=0,\\
    [w_1, w_2; z_1, z_2]_{z_1} &=0,\\
    [w_1, w_2; z_1, z_2]_{z_2} &=0,
\endaligned\right.\tag 1.4$$
where
$$(w_1, w_2)=\gamma(z_1, z_2)=
  \left(\frac{a_1 z_1+a_2 z_2+a_3}{c_1 z_1+c_2 z_2+c_3},
  \frac{b_1 z_1+b_2 z_2+b_3}{c_1 z_1+c_2 z_2+c_3}\right),$$
for $\gamma=\left(\matrix a_1 & a_2 & a_3\\ b_1 & b_2 & b_3\\
c_1 & c_2 & c_3 \endmatrix\right) \in GL(3, {\Bbb C})$.}

{\smc Proposition 1.4 (Main Proposition 1)}. {\it The change of
independent variables obeys the laws:
$$\aligned
 &\{u_1, u_2; x, y\}_{x}\\
=&\frac{(\frac{\partial w_1}{\partial x})^3}{\frac{\partial(w_1,
  w_2)}{\partial(x, y)}} \{u_1, u_2; w_1, w_2\}_{w_1}+
  \frac{(\frac{\partial w_1}{\partial x})^2 \frac{\partial w_2}
  {\partial x}}{\frac{\partial(w_1, w_2)}{\partial(x, y)}}
  [u_1, u_2; w_1, w_2]_{w_1}+\\
+&\frac{(\frac{\partial w_2}{\partial x})^3}{\frac{\partial(w_2,
  w_1)}{\partial(x, y)}} \{u_1, u_2; w_1, w_2\}_{w_2}+
  \frac{\frac{\partial w_1}{\partial x} (\frac{\partial w_2}
  {\partial x})^2}{\frac{\partial(w_2, w_1)}{\partial(x, y)}}
  [u_1, u_2; w_1, w_2]_{w_2}+\{w_1, w_2; x, y\}_{x}.
\endaligned\tag 1.5$$
$$\aligned
 &\{u_1, u_2; x, y\}_{y}\\
=&\frac{(\frac{\partial w_1}{\partial y})^3}{\frac{\partial(w_2,
  w_1)}{\partial(x, y)}} \{u_1, u_2; w_1, w_2\}_{w_1}+
  \frac{(\frac{\partial w_1}{\partial y})^2 \frac{\partial w_2}
  {\partial y}}{\frac{\partial(w_2, w_1)}{\partial(x, y)}}
  [u_1, u_2; w_1, w_2]_{w_1}+\\
+&\frac{(\frac{\partial w_2}{\partial y})^3}{\frac{\partial(w_1,
  w_2)}{\partial(x, y)}} \{u_1, u_2; w_1, w_2\}_{w_2}+
  \frac{\frac{\partial w_1}{\partial y} (\frac{\partial w_2}
  {\partial y})^2}{\frac{\partial(w_1, w_2)}{\partial(x, y)}}
  [u_1, u_2; w_1, w_2]_{w_2}+\{w_1, w_2; x, y\}_{y}.
\endaligned\tag 1.6$$
$$\aligned
 &[u_1, u_2; x, y]_{x}\\
=&\frac{3 (\frac{\partial w_1}{\partial x})^2 \frac{\partial
  w_1}{\partial y}}{\frac{\partial(w_1, w_2)}{\partial(x, y)}}
  \{u_1, u_2; w_1, w_2\}_{w_1}+\frac{(\frac{\partial w_1}{\partial
  x})^2 \frac{\partial w_2}{\partial y}+2 \frac{\partial w_1}{\partial x}
  \frac{\partial w_2}{\partial x} \frac{\partial w_1}{\partial y}}
  {\frac{\partial(w_1, w_2)}{\partial(x, y)}} [u_1, u_2; w_1,
  w_2]_{w_1}+\\
+&\frac{3 (\frac{\partial w_2}{\partial x})^2 \frac{\partial
  w_2}{\partial y}}{\frac{\partial(w_2, w_1)}{\partial(x, y)}}
  \{u_1, u_2; w_1, w_2\}_{w_2}+\frac{(\frac{\partial w_2}{\partial
  x})^2 \frac{\partial w_1}{\partial y}+2 \frac{\partial w_1}{\partial x}
  \frac{\partial w_2}{\partial x} \frac{\partial w_2}{\partial y}}
  {\frac{\partial(w_2, w_1)}{\partial(x, y)}} [u_1, u_2; w_1,
  w_2]_{w_2}+\\
+&[w_1, w_2; x, y]_{x}.
\endaligned\tag 1.7$$
$$\aligned
 &[u_1, u_2; x, y]_{y}\\
=&\frac{3 \frac{\partial w_1}{\partial x} (\frac{\partial
  w_1}{\partial y})^2}{\frac{\partial(w_2, w_1)}{\partial(x, y)}}
  \{u_1, u_2; w_1, w_2\}_{w_1}+\frac{\frac{\partial w_2}{\partial
  x} (\frac{\partial w_1}{\partial y})^2+2 \frac{\partial w_1}{\partial x}
  \frac{\partial w_1}{\partial y} \frac{\partial w_2}{\partial y}}
  {\frac{\partial(w_2, w_1)}{\partial(x, y)}} [u_1, u_2; w_1,
  w_2]_{w_1}+\\
+&\frac{3 \frac{\partial w_2}{\partial x} (\frac{\partial
  w_2}{\partial y})^2}{\frac{\partial(w_1, w_2)}{\partial(x, y)}}
  \{u_1, u_2; w_1, w_2\}_{w_2}+\frac{\frac{\partial w_1}{\partial
  x} (\frac{\partial w_2}{\partial y})^2+2 \frac{\partial w_2}{\partial x}
  \frac{\partial w_1}{\partial y} \frac{\partial w_2}{\partial y}}
  {\frac{\partial(w_1, w_2)}{\partial(x, y)}} [u_1, u_2; w_1,
  w_2]_{w_2}+\\
+&[w_1, w_2; x, y]_{y}.
\endaligned\tag 1.8$$}

{\smc Corollary 1.5}. {\it The exchange of the first and the
second arguments obeys the laws:
$$\aligned
 &\{w_1, w_2; x, y\}_{x}\\
=&-\frac{(\frac{\partial w_1}{\partial x})^3}{\frac{\partial(w_1,
  w_2)}{\partial(x, y)}} \{x, y; w_1, w_2\}_{w_1}-\frac{(\frac{\partial
  w_1}{\partial x})^2 \frac{\partial w_2}{\partial x}}{\frac{\partial
  (w_1, w_2)}{\partial(x, y)}} [x, y; w_1, w_2]_{w_1}+\\
 &+\frac{(\frac{\partial w_2}{\partial x})^3}{\frac{\partial(w_1, w_2)}
  {\partial(x, y)}} \{x, y; w_1, w_2\}_{w_2}+\frac{\frac{\partial w_1}
  {\partial x} (\frac{\partial w_2}{\partial x})^2}{\frac{\partial(w_1,
  w_2)}{\partial(x, y)}} [x, y; w_1, w_2]_{w_2},
\endaligned\tag 1.9$$
$$\aligned
 &\{w_1, w_2; x, y\}_{y}\\
=&\frac{(\frac{\partial w_1}{\partial y})^3}{\frac{\partial(w_1,
  w_2)}{\partial(x, y)}} \{x, y; w_1, w_2\}_{w_1}+\frac{(\frac{\partial
  w_1}{\partial y})^2 \frac{\partial w_2}{\partial y}}{\frac{\partial
  (w_1, w_2)}{\partial(x, y)}} [x, y; w_1, w_2]_{w_1}+\\
 &-\frac{(\frac{\partial w_2}{\partial y})^3}{\frac{\partial(w_1, w_2)}
  {\partial(x, y)}} \{x, y; w_1, w_2\}_{w_2}-\frac{\frac{\partial w_1}
  {\partial y} (\frac{\partial w_2}{\partial y})^2}{\frac{\partial(w_1,
  w_2)}{\partial(x, y)}} [x, y; w_1, w_2]_{w_2},
\endaligned\tag 1.10$$
$$\aligned
 &[w_1, w_2; x, y]_{x}\\
=&-\frac{3 (\frac{\partial w_1}{\partial x})^2 \frac{\partial
  w_1}{\partial y}}{\frac{\partial(w_1, w_2)}{\partial(x, y)}}
  \{x, y; w_1, w_2\}_{w_1}-\frac{(\frac{\partial w_1}{\partial
  x})^2 \frac{\partial w_2}{\partial y}+2 \frac{\partial w_1}
  {\partial x} \frac{\partial w_2}{\partial x} \frac{\partial
  w_1}{\partial y}}{\frac{\partial(w_1, w_2)}{\partial(x, y)}}
  [x, y; w_1, w_2]_{w_1}+\\
 &+\frac{3 (\frac{\partial w_2}{\partial x})^2 \frac{\partial w_2}
  {\partial y}}{\frac{\partial(w_1, w_2)}{\partial(x, y)}}
  \{x, y; w_1, w_2\}_{w_2}+\frac{(\frac{\partial w_2}{\partial x})^2
  \frac{\partial w_1}{\partial y}+2 \frac{\partial w_1}{\partial x}
  \frac{\partial w_2}{\partial x} \frac{\partial w_2}{\partial y}}
  {\frac{\partial(w_1, w_2)}{\partial(x, y)}} [x, y; w_1, w_2]_{w_2},
\endaligned\tag 1.11$$
$$\aligned
 &[w_1, w_2; x, y]_{y}\\
=&\frac{3 \frac{\partial w_1}{\partial x} (\frac{\partial
  w_1}{\partial y})^2}{\frac{\partial(w_1, w_2)}{\partial(x, y)}}
  \{x, y; w_1, w_2\}_{w_1}+\frac{\frac{\partial w_2}{\partial
  x} (\frac{\partial w_1}{\partial y})^2+2 \frac{\partial w_1}
  {\partial x} \frac{\partial w_1}{\partial y} \frac{\partial
  w_2}{\partial y}}{\frac{\partial(w_1, w_2)}{\partial(x, y)}}
  [x, y; w_1, w_2]_{w_1}+\\
 &-\frac{3 \frac{\partial w_2}{\partial x} (\frac{\partial w_2}
  {\partial y})^2}{\frac{\partial(w_1, w_2)}{\partial(x, y)}}
  \{x, y; w_1, w_2\}_{w_2}-\frac{\frac{\partial w_1}{\partial x}
  (\frac{\partial w_2}{\partial y})^2+2 \frac{\partial w_2}{\partial x}
  \frac{\partial w_1}{\partial y} \frac{\partial w_2}{\partial y}}
  {\frac{\partial(w_1, w_2)}{\partial(x, y)}} [x, y; w_1, w_2]_{w_2}.
\endaligned\tag 1.12$$}

  We may apply a linear fractional transformation to the second
arguments instead of the first. The result is as follows:

{\smc Corollary 1.6}. {\it For
$$(w_1, w_2)=\gamma(x, y)=\left(\frac{a_1 x+a_2 y+a_3}{c_1 x+c_2 y+c_3},
  \frac{b_1 x+b_2 y+b_3}{c_1 x+c_2 y+c_3}\right),\tag 1.13$$
where $\gamma=\left(\matrix a_1 & a_2 & a_3\\ b_1 & b_2 & b_3\\
c_1 & c_2 & c_3 \endmatrix\right) \in GL(3, {\Bbb C})$, the
following four identities hold:
$$\aligned
 &\{u_1, u_2; x, y\}_{x}\\
=&\frac{A_1^3}{\Delta (c_1 x+c_2 y+c_3)^3}
  \left\{u_1, u_2; \frac{a_1 x+a_2 y+a_3}{c_1 x+c_2 y+c_3},
  \frac{b_1 x+b_2 y+b_3}{c_1 x+c_2 y+c_3}\right\}_{\frac{a_1
  x+a_2 y+a_3}{c_1 x+c_2 y+c_3}}+\\
+&\frac{A_1^2 B_1}{\Delta (c_1 x+c_2 y+c_3)^3}
  \left[u_1, u_2; \frac{a_1 x+a_2 y+a_3}{c_1 x+c_2 y+c_3},
  \frac{b_1 x+b_2 y+b_3}{c_1 x+c_2 y+c_3}\right]_{\frac{a_1
  x+a_2 y+a_3}{c_1 x+c_2 y+c_3}}+\\
-&\frac{B_1^3}{\Delta (c_1 x+c_2 y+c_3)^3}
  \left\{u_1, u_2; \frac{a_1 x+a_2 y+a_3}{c_1 x+c_2 y+c_3},
  \frac{b_1 x+b_2 y+b_3}{c_1 x+c_2 y+c_3}\right\}_{\frac{b_1
  x+b_2 y+b_3}{c_1 x+c_2 y+c_3}}+\\
-&\frac{A_1 B_1^2}{\Delta (c_1 x+c_2 y+c_3)^3}
  \left[u_1, u_2; \frac{a_1 x+a_2 y+a_3}{c_1 x+c_2 y+c_3},
  \frac{b_1 x+b_2 y+b_3}{c_1 x+c_2 y+c_3}\right]_{\frac{b_1
  x+b_2 y+b_3}{c_1 x+c_2 y+c_3}},
\endaligned\tag 1.14$$
$$\aligned
 &\{u_1, u_2; x, y\}_{y}\\
=&-\frac{A_2^3}{\Delta (c_1 x+c_2 y+c_3)^3}
  \left\{u_1, u_2; \frac{a_1 x+a_2 y+a_3}{c_1 x+c_2 y+c_3},
  \frac{b_1 x+b_2 y+b_3}{c_1 x+c_2 y+c_3}\right\}_{\frac{a_1
  x+a_2 y+a_3}{c_1 x+c_2 y+c_3}}+\\
-&\frac{A_2^2 B_2}{\Delta (c_1 x+c_2 y+c_3)^3}
  \left[u_1, u_2; \frac{a_1 x+a_2 y+a_3}{c_1 x+c_2 y+c_3},
  \frac{b_1 x+b_2 y+b_3}{c_1 x+c_2 y+c_3}\right]_{\frac{a_1
  x+a_2 y+a_3}{c_1 x+c_2 y+c_3}}+\\
+&\frac{B_2^3}{\Delta (c_1 x+c_2 y+c_3)^3}
  \left\{u_1, u_2; \frac{a_1 x+a_2 y+a_3}{c_1 x+c_2 y+c_3},
  \frac{b_1 x+b_2 y+b_3}{c_1 x+c_2 y+c_3}\right\}_{\frac{b_1
  x+b_2 y+b_3}{c_1 x+c_2 y+c_3}}+\\
+&\frac{A_2 B_2^2}{\Delta (c_1 x+c_2 y+c_3)^3}
  \left[u_1, u_2; \frac{a_1 x+a_2 y+a_3}{c_1 x+c_2 y+c_3},
  \frac{b_1 x+b_2 y+b_3}{c_1 x+c_2 y+c_3}\right]_{\frac{b_1
  x+b_2 y+b_3}{c_1 x+c_2 y+c_3}},
\endaligned\tag 1.15$$
$$\aligned
 &[u_1, u_2; x, y]_{x}\\
=&\frac{3 A_1^2 A_2}{\Delta (c_1 x+c_2 y+c_3)^3}
  \left\{u_1, u_2; \frac{a_1 x+a_2 y+a_3}{c_1 x+c_2 y+c_3},
  \frac{b_1 x+b_2 y+b_3}{c_1 x+c_2 y+c_3}\right\}_{\frac{a_1
  x+a_2 y+a_3}{c_1 x+c_2 y+c_3}}+\\
+&\frac{A_1^2 B_2+2 A_1 B_1 A_2}{\Delta (c_1 x+c_2 y+c_3)^3}
  \left[u_1, u_2; \frac{a_1 x+a_2 y+a_3}{c_1 x+c_2 y+c_3},
  \frac{b_1 x+b_2 y+b_3}{c_1 x+c_2 y+c_3}\right]_{\frac{a_1
  x+a_2 y+a_3}{c_1 x+c_2 y+c_3}}+\\
-&\frac{3 B_1^2 B_2}{\Delta (c_1 x+c_2 y+c_3)^3}
  \left\{u_1, u_2; \frac{a_1 x+a_2 y+a_3}{c_1 x+c_2 y+c_3},
  \frac{b_1 x+b_2 y+b_3}{c_1 x+c_2 y+c_3}\right\}_{\frac{b_1
  x+b_2 y+b_3}{c_1 x+c_2 y+c_3}}+\\
-&\frac{B_1^2 A_2+2 A_1 B_1 B_2}{\Delta (c_1 x+c_2 y+c_3)^3}
  \left[u_1, u_2; \frac{a_1 x+a_2 y+a_3}{c_1 x+c_2 y+c_3},
  \frac{b_1 x+b_2 y+b_3}{c_1 x+c_2 y+c_3}\right]_{\frac{b_1
  x+b_2 y+b_3}{c_1 x+c_2 y+c_3}},
\endaligned\tag 1.16$$
$$\aligned
 &[u_1, u_2; x, y]_{y}\\
=&-\frac{3 A_1 A_2^2}{\Delta (c_1 x+c_2 y+c_3)^3}
  \left\{u_1, u_2; \frac{a_1 x+a_2 y+a_3}{c_1 x+c_2 y+c_3},
  \frac{b_1 x+b_2 y+b_3}{c_1 x+c_2 y+c_3}\right\}_{\frac{a_1
  x+a_2 y+a_3}{c_1 x+c_2 y+c_3}}+\\
-&\frac{B_1 A_2^2+2 A_1 A_2 B_2}{\Delta (c_1 x+c_2 y+c_3)^3}
  \left[u_1, u_2; \frac{a_1 x+a_2 y+a_3}{c_1 x+c_2 y+c_3},
  \frac{b_1 x+b_2 y+b_3}{c_1 x+c_2 y+c_3}\right]_{\frac{a_1
  x+a_2 y+a_3}{c_1 x+c_2 y+c_3}}+\\
+&\frac{3 B_1 B_2^2}{\Delta (c_1 x+c_2 y+c_3)^3}
  \left\{u_1, u_2; \frac{a_1 x+a_2 y+a_3}{c_1 x+c_2 y+c_3},
  \frac{b_1 x+b_2 y+b_3}{c_1 x+c_2 y+c_3}\right\}_{\frac{b_1
  x+b_2 y+b_3}{c_1 x+c_2 y+c_3}}+\\
+&\frac{A_1 B_2^2+2 B_1 A_2 B_2}{\Delta (c_1 x+c_2 y+c_3)^3}
  \left[u_1, u_2; \frac{a_1 x+a_2 y+a_3}{c_1 x+c_2 y+c_3},
  \frac{b_1 x+b_2 y+b_3}{c_1 x+c_2 y+c_3}\right]_{\frac{b_1
  x+b_2 y+b_3}{c_1 x+c_2 y+c_3}}.
\endaligned\tag 1.17$$
Here, $\Delta=\det(\gamma)$ and
$$\left\{\aligned
 A_1 &=(a_1 c_2-a_2 c_1) y+(a_1 c_3-a_3 c_1),\\
 A_2 &=(a_2 c_1-a_1 c_2) x+(a_2 c_3-a_3 c_2),\\
 B_1 &=(b_1 c_2-b_2 c_1) y+(b_1 c_3-b_3 c_1),\\
 B_2 &=(b_2 c_1-b_1 c_2) x+(b_2 c_3-b_3 c_2).
\endaligned\right.$$}

  The Proposition 1.4 implies that
$$\left(\matrix
  \{u_1, u_2; x, y\}_{x}\\
  \{u_1, u_2; x, y\}_{y}\\
   [u_1, u_2; x, y]_{x}\\
   [u_1, u_2; x, y]_{y}
  \endmatrix\right)
 =\frac{M(w_1, w_2; x, y)}{\frac{\partial(w_1, w_2)}{\partial(x, y)}}
  \left(\matrix
  \{u_1, u_2; w_1, w_2\}_{w_1}\\
  \{u_1, u_2; w_1, w_2\}_{w_2}\\
   [u_1, u_2; w_1, w_2]_{w_1}\\
   [u_1, u_2; w_1, w_2]_{w_2}
  \endmatrix\right)+\left(\matrix
  \{w_1, w_2; x, y\}_{x}\\
  \{w_1, w_2; x, y\}_{y}\\
   [w_1, w_2; x, y]_{x}\\
   [w_1, w_2; x, y]_{y}
  \endmatrix\right),\tag 1.18$$
where
$$\aligned
 &M(w_1, w_2; x, y):=\\
 &\left(\matrix
  (\frac{\partial w_1}{\partial x})^3 & -(\frac{\partial w_2}{\partial x})^3
 &(\frac{\partial w_1}{\partial x})^2 \frac{\partial w_2}{\partial x}
 &-\frac{\partial w_1}{\partial x} (\frac{\partial w_2}{\partial x})^2\\
  -(\frac{\partial w_1}{\partial y})^3 & (\frac{\partial w_2}{\partial y})^3
 &-(\frac{\partial w_1}{\partial y})^2 \frac{\partial w_2}{\partial y}
 &\frac{\partial w_1}{\partial y} (\frac{\partial w_2}{\partial y})^2\\
  3 (\frac{\partial w_1}{\partial x})^2 \frac{\partial w_1}{\partial y}
 &-3 (\frac{\partial w_2}{\partial x})^2 \frac{\partial w_2}{\partial y}
 &\matrix
  &(\frac{\partial w_1}{\partial x})^2 \frac{\partial w_2}{\partial y}+\\
  &+2 \frac{\partial w_1}{\partial x} \frac{\partial w_2}{\partial x}
   \frac{\partial w_1}{\partial y}
  \endmatrix
 &\matrix
  &-(\frac{\partial w_2}{\partial x})^2 \frac{\partial w_1}{\partial y}+\\
  &-2 \frac{\partial w_1}{\partial x} \frac{\partial w_2}{\partial x}
   \frac{\partial w_2}{\partial y}
  \endmatrix\\
  -3 \frac{\partial w_1}{\partial x} (\frac{\partial w_1}{\partial y})^2
 &3 \frac{\partial w_2}{\partial x} (\frac{\partial w_2}{\partial y})^2
 &\matrix
  &-\frac{\partial w_2}{\partial x} (\frac{\partial w_1}{\partial y})^2+\\
  &-2 \frac{\partial w_1}{\partial x} \frac{\partial w_1}{\partial y}
   \frac{\partial w_2}{\partial y}
  \endmatrix
 &\matrix
  &\frac{\partial w_1}{\partial x} (\frac{\partial w_2}{\partial y})^2+\\
  &+2 \frac{\partial w_2}{\partial x} \frac{\partial w_1}{\partial y}
   \frac{\partial w_2}{\partial y}
  \endmatrix
\endmatrix\right).
\endaligned\tag 1.19$$

{\smc Proposition 1.7}. {\it
$$M(w_1, w_2; x, y) \cdot M(u_1, u_2; w_1, w_2)
 =M(u_1, u_2; x, y).\tag 1.20$$}

  Set
$$U_{(w_1, w_2; x, y)}:=\left(\matrix
  M(w_1, w_2; x, y) & c \frac{\partial(w_1, w_2)}{\partial(x, y)}
                      \left(\matrix  \{w_1, w_2; x, y\}_{x}\\
                                     \{w_1, w_2; x, y\}_{y}\\
                                       [w_1, w_2; x, y]_{x}\\
                                       [w_1, w_2; x, y]_{y}
                      \endmatrix\right)\\
     0              & \frac{\partial(w_1, w_2)}{\partial(x, y)}
  \endmatrix\right).\tag 1.21$$
Then the above proposition implies that
$$U_{(w_1, w_2; x, y)} \cdot U_{(u_1, u_2; w_1, w_2)}
 =U_{(u_1, u_2; x, y)}.\tag 1.22$$

  The equation
$$U_{W, Z} U_{U, W}=U_{U, Z}$$
with $W=(w_1, w_2)$, $Z=(z_1, z_2)$ and $U=(u_1, u_2)$ is a
cocycle condition, which means that the $U$'s can be interpreted
as transition functions of a vector bundle. More precisely, when
we choose an open cover $\Sigma=\bigcup_{i} \Sigma_i$, with a
local coordinate $Z_i$ on $\Sigma_i$, then the above equation
permits us to interpret the $U_{Z_i, Z_j}$ as transition functions
on the intersection $\Sigma_i \cap \Sigma_j$. Let ${\Cal V}_{0}$
be this vector bundle. It is a five dimensional sub-bundle of the
infinite dimensional vector bundle ${\Cal V}$.

   The next Main Proposition 2 gives the relation between the
deformation of the Jacobian and our four derivatives:

{\smc Proposition 1.8 (Main Proposition 2)}. {\it The following
identity holds:
$$\aligned
 &\frac{\partial(f_1, f_2)}{\partial(z_1, z_2)}\\
=&\frac{\partial(\widehat{f_1}, \widehat{f_2})}{\partial(w_1,
  w_2)}-\widehat{f_1} \frac{\partial \widehat{f_1}}{\partial w_2}
  \{z_1, z_2; w_1, w_2\}_{w_1}-\widehat{f_2} \frac{\partial
  \widehat{f_2}}{\partial w_1} \{z_1, z_2; w_1, w_2\}_{w_2}+\\
 &-\widehat{f_1} \frac{\partial \widehat{f_2}}{\partial w_2}
  \frac{\vmatrix \format \c \quad & \c \\
  \frac{\partial z_1}{\partial w_2} &
  \frac{\partial z_2}{\partial w_2}\\
  \frac{\partial^2 z_1}{\partial w_1^2} &
  \frac{\partial^2 z_2}{\partial w_1^2}\endvmatrix}
 {\vmatrix \format \c \quad & \c\\
  \frac{\partial z_1}{\partial w_1} &
  \frac{\partial z_2}{\partial w_1}\\
  \frac{\partial z_1}{\partial w_2} &
  \frac{\partial z_2}{\partial w_2}
  \endvmatrix}
  -\left(\widehat{f_2} \frac{\partial \widehat{f_1}}{\partial w_2}-
  \widehat{f_1} \frac{\partial \widehat{f_1}}{\partial w_1}\right)
  \frac{\vmatrix \format \c \quad & \c \\
  \frac{\partial z_1}{\partial w_1} &
  \frac{\partial z_2}{\partial w_1}\\
  \frac{\partial^2 z_1}{\partial w_1 \partial w_2} &
  \frac{\partial^2 z_2}{\partial w_1 \partial w_2}\endvmatrix}
 {\vmatrix \format \c \quad & \c\\
  \frac{\partial z_1}{\partial w_1} &
  \frac{\partial z_2}{\partial w_1}\\
  \frac{\partial z_1}{\partial w_2} &
  \frac{\partial z_2}{\partial w_2}
  \endvmatrix}+\\
 &-\widehat{f_2} \frac{\partial \widehat{f_1}}{\partial w_1}
  \frac{\vmatrix \format \c \quad & \c \\
  \frac{\partial z_1}{\partial w_1} &
  \frac{\partial z_2}{\partial w_1}\\
  \frac{\partial^2 z_1}{\partial w_2^2} &
  \frac{\partial^2 z_2}{\partial w_2^2}\endvmatrix}
 {\vmatrix \format \c \quad & \c\\
  \frac{\partial z_1}{\partial w_2} &
  \frac{\partial z_2}{\partial w_2}\\
  \frac{\partial z_1}{\partial w_1} &
  \frac{\partial z_2}{\partial w_1}
  \endvmatrix}-\left(\widehat{f_1} \frac{\partial \widehat{f_2}}
  {\partial w_1}-\widehat{f_2} \frac{\partial \widehat{f_2}}
  {\partial w_2}\right) \frac{\vmatrix \format \c \quad & \c \\
  \frac{\partial z_1}{\partial w_2} &
  \frac{\partial z_2}{\partial w_2}\\
  \frac{\partial^2 z_1}{\partial w_1 \partial w_2} &
  \frac{\partial^2 z_2}{\partial w_1 \partial w_2}\endvmatrix}
 {\vmatrix \format \c \quad & \c\\
  \frac{\partial z_1}{\partial w_2} &
  \frac{\partial z_2}{\partial w_2}\\
  \frac{\partial z_1}{\partial w_1} &
  \frac{\partial z_2}{\partial w_1}
  \endvmatrix}+\\
 &+{\widehat{f_1}}^2 \frac{\vmatrix \format \c \quad & \c \\
  \frac{\partial^2 z_1}{\partial w_1^2} &
  \frac{\partial^2 z_2}{\partial w_1^2}\\
  \frac{\partial^2 z_1}{\partial w_1 \partial w_2} &
  \frac{\partial^2 z_2}{\partial w_1 \partial w_2}\endvmatrix}
 {\vmatrix \format \c \quad & \c\\
  \frac{\partial z_1}{\partial w_1} &
  \frac{\partial z_2}{\partial w_1}\\
  \frac{\partial z_1}{\partial w_2} &
  \frac{\partial z_2}{\partial w_2}
  \endvmatrix}+{\widehat{f_2}}^2
  \frac{\vmatrix \format \c \quad & \c \\
  \frac{\partial^2 z_1}{\partial w_2^2} &
  \frac{\partial^2 z_2}{\partial w_2^2}\\
  \frac{\partial^2 z_1}{\partial w_1 \partial w_2} &
  \frac{\partial^2 z_2}{\partial w_1 \partial w_2}\endvmatrix}
 {\vmatrix \format \c \quad & \c\\
  \frac{\partial z_1}{\partial w_2} &
  \frac{\partial z_2}{\partial w_2}\\
  \frac{\partial z_1}{\partial w_1} &
  \frac{\partial z_2}{\partial w_1}
  \endvmatrix}+\widehat{f_1} \widehat{f_2}
  \frac{\vmatrix \format \c \quad & \c \\
  \frac{\partial^2 z_1}{\partial w_1^2} &
  \frac{\partial^2 z_2}{\partial w_1^2}\\
  \frac{\partial^2 z_1}{\partial w_2^2} &
  \frac{\partial^2 z_2}{\partial w_2^2}\endvmatrix}
 {\vmatrix \format \c \quad & \c\\
  \frac{\partial z_1}{\partial w_1} &
  \frac{\partial z_2}{\partial w_1}\\
  \frac{\partial z_1}{\partial w_2} &
  \frac{\partial z_2}{\partial w_2}
  \endvmatrix},
\endaligned\tag 1.23$$
where
$$(f_1, f_2)=(\widehat{f_1}, \widehat{f_2})
  \left(\matrix
  \frac{\partial z_1}{\partial w_1} &
  \frac{\partial z_2}{\partial w_1}\\
  \frac{\partial z_1}{\partial w_2} &
  \frac{\partial z_2}{\partial w_2}
  \endmatrix\right).$$}

  Now, put
$$\phi(f_1, f_2):=\iint \frac{\partial(f_1, f_2)}{\partial(t_1, t_2)}
                  dt_1 \wedge dt_2.$$
We define an operation
$$\aligned
  &\left[f_1 \frac{\partial}{\partial t_1} \otimes \frac{\partial}{\partial t_2},
    f_2 \frac{\partial}{\partial t_1} \otimes \frac{\partial}{\partial t_2}\right]_{Z_0}^{*}\\
:=&c \cdot \phi(f_1, f_2)-f_1 \frac{\partial f_1}{\partial t_2}
   \frac{\partial}{\partial t_1}-f_1 \frac{\partial f_2}{\partial t_2}
   \frac{\partial}{\partial t_1}-\left(f_2 \frac{\partial f_1}{\partial t_2}
   -f_1 \frac{\partial f_1}{\partial t_1}\right) \frac{\partial}{\partial t_1}+\\
  &-f_2 \frac{\partial f_2}{\partial t_1} \frac{\partial}{\partial t_2}
   -f_2 \frac{\partial f_1}{\partial t_1} \frac{\partial}{\partial t_2}
   -\left(f_1 \frac{\partial f_2}{\partial t_1}-f_2 \frac{\partial f_2}
   {\partial t_2}\right) \frac{\partial}{\partial t_2}+f_1^2
   \frac{\partial^2}{\partial t_1^2}+f_2^2 \frac{\partial^2}{\partial t_2^2}
   +f_1 f_2 \frac{\partial^2}{\partial t_1 \partial t_2}.
\endaligned$$
Note that
$$f_1 \frac{\partial f_1}{\partial t_2}+f_1 \frac{\partial f_2}{\partial t_2}
  +\left(f_2 \frac{\partial f_1}{\partial t_2}-f_1 \frac{\partial f_1}{\partial
  t_1}\right)
 =(f_1+f_2) \frac{\partial f_1}{\partial t_2}+f_1 \left(\frac{\partial f_2}
  {\partial t_2}-\frac{\partial f_1}{\partial t_1}\right)$$
and
$$f_2 \frac{\partial f_2}{\partial t_1}+f_2 \frac{\partial f_1}{\partial t_1}
  +\left(f_1 \frac{\partial f_2}{\partial t_1}-f_2 \frac{\partial f_2}{\partial
  t_2}\right)
 =(f_1+f_2) \frac{\partial f_2}{\partial t_1}+f_2 \left(\frac{\partial f_1}
  {\partial t_1}-\frac{\partial f_2}{\partial t_2}\right).$$
Modulo the terms of second order partial derivatives, we define
the following operation:
$$\aligned
  &\left[f_1 \frac{\partial}{\partial t_1} \otimes \frac{\partial}{\partial t_2},
   f_2 \frac{\partial}{\partial t_1} \otimes \frac{\partial}{\partial
   t_2}\right]_{Z_0}
 :=\left[(f_1+f_2) \frac{\partial f_1}{\partial t_2}+f_1
   \left(\frac{\partial f_2}{\partial t_2}-\frac{\partial
   f_1}{\partial t_1}\right)\right] \frac{\partial}{\partial t_1}+\\
 +&\left[(f_1+f_2) \frac{\partial f_2}{\partial t_1}+f_2
   \left(\frac{\partial f_1}{\partial t_1}-\frac{\partial
   f_2}{\partial t_2}\right)\right] \frac{\partial}{\partial t_2}+
   c \cdot \phi(f_1, f_2).
\endaligned$$

{\smc Proposition 1.9}. {\it For the system of partial
differential equations:
$$\left\{\aligned
  r &=a_1 p+a_2 q+a_3 z,\\
  s &=b_1 p+b_2 q+b_3 z,\\
  t &=c_1 p+c_2 q+c_3 z,
\endaligned\right.$$
where $a_i=a_i(x, y)$, $b_i=b_i(x, y)$, $c_i=c_i(x, y)$ $(i=1, 2,
3)$ and $r=\frac{\partial^2 z}{\partial x^2}$, $s=\frac{\partial^2
z}{\partial x \partial y}$, $t=\frac{\partial^2 z}{\partial y^2}$,
$p=\frac{\partial z}{\partial x}$, $q=\frac{\partial z}{\partial
y}$, there exists three linearly independent solutions: $z_1$,
$z_2$ and $z_3$. Put $u_1=z_1/z_3$ and $u_2=z_2/z_3$. Then
$$\left\{\aligned
  \{ u_1, u_2; x, y \}_{x} &=a_2,\\
  \{ u_1, u_2; x, y \}_{y} &=c_1,\\
      [u_1, u_2; x, y]_{x} &=2 b_2-a_1,\\
      [u_1, u_2; x, y]_{y} &=2 b_1-c_2.
\endaligned\right.\tag 1.24$$}

{\smc Theorem 1.10 (Main Theorem 1)}. {\it For the over-determined
system of second order, nonlinear $($quasilinear$)$ complex
partial differential equations:
$$\left\{\aligned
  \{ w_1, w_2; v_1, v_2 \}_{v_1} &=F_1(v_1, v_2),\\
  \{ w_1, w_2; v_1, v_2 \}_{v_2} &=F_2(v_1, v_2),\\
    [ w_1, w_2; v_1, v_2 ]_{v_1} &=P_1(v_1, v_2),\\
    [ w_1, w_2; v_1, v_2 ]_{v_2} &=P_2(v_1, v_2),
  \endaligned\right.\tag 1.25$$
set the Jacobian
$$\Delta:=\frac{\partial(v_1, v_2)}{\partial(w_1, w_2)}
 =\vmatrix \format \c \quad & \c\\
  \frac{\partial v_1}{\partial w_1} &
  \frac{\partial v_2}{\partial w_1}\\
  \frac{\partial v_1}{\partial w_2} &
  \frac{\partial v_2}{\partial w_2}
  \endvmatrix.\tag 1.26$$
Then $z=\root 3 \of{\Delta}$ satisfies the following linear
differential equations:
$$\left\{\aligned
  \frac{\partial^2 z}{\partial v_1^2}+\frac{1}{3} P_1 \frac{\partial z}
  {\partial v_1}-F_1 \frac{\partial z}{\partial v_2}+\left(\frac{\partial
  F_1}{\partial v_2}-\frac{1}{3} \frac{\partial P_1}{\partial v_1}-
  \frac{2}{9} P_1^2-\frac{2}{3} F_1 P_2\right) z &=0,\\
  \frac{\partial^2 z}{\partial v_1 \partial v_2}-\frac{1}{3} P_2 \frac{\partial z}
  {\partial v_1}-\frac{1}{3} P_1 \frac{\partial z}{\partial v_2}+\left(\frac{1}{3}
  \frac{\partial P_2}{\partial v_1}+\frac{1}{3} \frac{\partial P_1}{\partial v_2}+
  \frac{1}{9} P_1 P_2-F_1 F_2\right) z &=0,\\
  \frac{\partial^2 z}{\partial v_2^2}+\frac{1}{3} P_2 \frac{\partial z}
  {\partial v_2}-F_2 \frac{\partial z}{\partial v_1}+\left(\frac{\partial
  F_2}{\partial v_1}-\frac{1}{3} \frac{\partial P_2}{\partial v_2}-
  \frac{2}{9} P_2^2-\frac{2}{3} F_2 P_1\right) z &=0.
  \endaligned\right.\tag 1.27$$}

{\smc Theorem 1.11 (Main Theorem 2)}. {\it For the system of
linear partial differential equations:
$$\left\{\aligned
  \frac{\partial^2 z}{\partial v_1^2}+\frac{1}{3} P_1 \frac{\partial z}
  {\partial v_1}-F_1 \frac{\partial z}{\partial v_2}+\left(\frac{\partial
  F_1}{\partial v_2}-\frac{1}{3} \frac{\partial P_1}{\partial v_1}-
  \frac{2}{9} P_1^2-\frac{2}{3} F_1 P_2\right) z &=0,\\
  \frac{\partial^2 z}{\partial v_1 \partial v_2}-\frac{1}{3} P_2 \frac{\partial z}
  {\partial v_1}-\frac{1}{3} P_1 \frac{\partial z}{\partial v_2}+\left(\frac{1}{3}
  \frac{\partial P_2}{\partial v_1}+\frac{1}{3} \frac{\partial P_1}{\partial v_2}+
  \frac{1}{9} P_1 P_2-F_1 F_2\right) z &=0,\\
  \frac{\partial^2 z}{\partial v_2^2}+\frac{1}{3} P_2 \frac{\partial z}
  {\partial v_2}-F_2 \frac{\partial z}{\partial v_1}+\left(\frac{\partial
  F_2}{\partial v_1}-\frac{1}{3} \frac{\partial P_2}{\partial v_2}-
  \frac{2}{9} P_2^2-\frac{2}{3} F_2 P_1\right) z &=0,
  \endaligned\right.\tag 1.28$$
put
$$P_1=\frac{\alpha}{v_1}+\frac{\beta}{v_1-1}+\frac{\gamma}{v_1-v_2}, \quad
  P_2=\frac{\alpha}{v_2}+\frac{\beta}{v_2-1}+\frac{\gamma}{v_2-v_1},\tag 1.29$$
$$F_1=\frac{-\gamma v_2(v_2-1)}{v_1(v_1-1)(v_1-v_2)}, \quad
  F_2=\frac{-\gamma v_1(v_1-1)}{v_2(v_2-1)(v_2-v_1)}.\tag 1.30$$
Set
$$w=v_1^{-\frac{1}{3} \alpha} v_2^{-\frac{1}{3} \alpha}
    (v_1-1)^{-\frac{1}{3} \beta} (v_2-1)^{-\frac{1}{3} \beta}
    (v_1-v_2)^{\frac{2}{3} \gamma} z,\tag 1.31$$
then the function $w=w(v_1, v_2)$ satisfies the following
equations:
$$\left\{\aligned
 &\frac{\partial^2 w}{\partial v_1^2}+\left(\frac{\alpha}{v_1}+
  \frac{\beta}{v_1-1}+\frac{-\gamma}{v_1-v_2}\right) \frac{\partial
  w}{\partial v_1}+\frac{\gamma v_2(v_2-1)}{v_1(v_1-1)(v_1-v_2)}
  \frac{\partial w}{\partial v_2}+\\
 &+\frac{(1-\alpha-\beta) \gamma}{v_1(v_1-1)} w=0,\\
 &\frac{\partial^2 w}{\partial v_2^2}+\left(\frac{\alpha}{v_2}+
  \frac{\beta}{v_2-1}+\frac{-\gamma}{v_2-v_1}\right) \frac{\partial
  w}{\partial v_2}+\frac{\gamma v_1(v_1-1)}{v_2(v_2-1)(v_2-v_1)}
  \frac{\partial w}{\partial v_1}+\\
 &+\frac{(1-\alpha-\beta) \gamma}{v_2(v_2-1)} w=0,\\
 &\frac{\partial^2 w}{\partial v_1 \partial v_2}=\frac{-\gamma}
  {v_1-v_2} \frac{\partial w}{\partial v_1}+\frac{\gamma}{v_1-v_2}
  \frac{\partial w}{\partial v_2}.
\endaligned\right.\tag 1.32$$
A solution is given by
$$\aligned
  w=&\text{const} \cdot F_1(\alpha+\beta-1; -\gamma, -\gamma;
     \alpha-\gamma; v_1, v_2)+\\
   +&\text{const} \cdot F_1(\alpha+\beta-1; -\gamma, -\gamma;
     \beta-\gamma; 1-v_1, 1-v_2),
\endaligned\tag 1.33$$
where $F_1=F_1(a; b, b^{\prime}; c; x, y)$ is the two variable
hypergeometric function $($see \cite{AK}$)$. Consequently,
$$\aligned
 z=&v_1^{\frac{1}{3} \alpha} v_2^{\frac{1}{3} \alpha}
    (v_1-1)^{\frac{1}{3} \beta} (v_2-1)^{\frac{1}{3} \beta}
    (v_1-v_2)^{-\frac{2}{3} \gamma} \times\\
 \times &[\text{const} \cdot F_1(\alpha+\beta-1; -\gamma, -\gamma;
    \alpha-\gamma; v_1, v_2)+\\
  +&\text{const} \cdot F_1(\alpha+\beta-1; -\gamma, -\gamma; \beta
   -\gamma; 1-v_1, 1-v_2)].
\endaligned\tag 1.34$$}

  In fact, when $(w_1, w_2)$ is a solution of the system of nonlinear
partial differential equations:
$$\left\{\aligned
  \{ w_1, w_2; v_1, v_2 \}_{v_1} &=F_1(v_1, v_2),\\
  \{ w_1, w_2; v_1, v_2 \}_{v_2} &=F_2(v_1, v_2),\\
      [w_1, w_2; v_1, v_2]_{v_1} &=P_1(v_1, v_2),\\
      [w_1, w_2; v_1, v_2]_{v_2} &=P_2(v_1, v_2),
\endaligned\right.$$
the general solution has the form
$$\left(\frac{a_1 w_1+a_2 w_2+a_3}{c_1 w_1+c_2 w_2+c_3},
  \frac{b_1 w_1+b_2 w_2+b_3}{c_1 w_1+c_2 w_2+c_3}\right),$$
where $\left(\matrix a_1 & a_2 & a_3\\ b_1 & b_2 & b_3\\ c_1 & c_2
& c_3 \endmatrix\right) \in GL(3, {\Bbb C})$.

  It is known that the two variable hypergeometric function
$$F_{1}(a; b, b^{\prime}; c; x, y)=\frac{\Gamma(c)}{\Gamma(a) \Gamma(c-a)}
  \int_{0}^{1} t^{a-1} (1-t)^{c-a-1} (1-tx)^{-b} (1-ty)^{-b^{\prime}} dt$$
satisfies the following equations:
$$\left\{\aligned
 x(1-x) \frac{\partial^2 z}{\partial x^2}+y(1-x) \frac{\partial^2 z}
 {\partial x \partial y}+[c-(a+b+1)x] \frac{\partial z}{\partial x}-
 by \frac{\partial z}{\partial y}-abz=0,\\
 y(1-y) \frac{\partial^2 z}{\partial y^2}+x(1-y) \frac{\partial^2 z}
 {\partial x \partial y}+[c-(a+b^{\prime}+1)y] \frac{\partial z}{\partial y}-
 b^{\prime} x \frac{\partial z}{\partial x}-a b^{\prime} z=0.
\endaligned\right.$$

  For the group $GL(3)$, the Picard integral
$$z=\int_{0}^{1} \frac{dt}{\root 3 \of{t(t-1)(t-x)(t-y)}}$$
satisfies the following equations:
$$\left\{\aligned
 x(1-x) \frac{\partial^2 z}{\partial x^2}+y(1-x) \frac{\partial^2 z}
 {\partial x \partial y}+\left(1-\frac{5}{3}x\right) \frac{\partial z}{\partial
 x}-\frac{1}{3} y \frac{\partial z}{\partial y}-\frac{1}{9}z=0,\\
 y(1-y) \frac{\partial^2 z}{\partial y^2}+x(1-y) \frac{\partial^2 z}
 {\partial x \partial y}+\left(1-\frac{5}{3}y\right) \frac{\partial z}{\partial
 y}-\frac{1}{3} x \frac{\partial z}{\partial x}-\frac{1}{9}z=0.
\endaligned\right.\tag 1.35$$
Thus, a solution is given by
$$z=F_1\left(\frac{1}{3}; \frac{1}{3}, \frac{1}{3}; 1; x, y\right), \quad
  (a, b, b^{\prime}, c)=\left(\frac{1}{3}, \frac{1}{3}, \frac{1}{3}, 1\right).\tag 1.36$$
The Picard modular function is associated to the following
equations:
$$\left\{\aligned
 x(1-x) \frac{\partial^2 z}{\partial x^2}+y(1-x) \frac{\partial^2 z}
 {\partial x \partial y}+\left(1-\frac{3}{2} x\right) \frac{\partial z}
 {\partial x}-\frac{1}{4} y \frac{\partial z}{\partial y}-\frac{1}{16} z=0,\\
 y(1-y) \frac{\partial^2 z}{\partial y^2}+x(1-y) \frac{\partial^2 z}
 {\partial x \partial y}+\left(1-\frac{3}{2} y\right) \frac{\partial z}
 {\partial y}-\frac{1}{4} x \frac{\partial z}{\partial x}-\frac{1}{16} z=0.
\endaligned\right.\tag 1.37$$
A solution is given by
$$z=F_1\left(\frac{1}{4}; \frac{1}{4}, \frac{1}{4}; 1; x, y\right), \quad
  (a, b, b^{\prime}, c)=\left(\frac{1}{4}, \frac{1}{4}, \frac{1}{4}, 1\right).\tag 1.38$$

  The $\eta$-function associated to the unitary group $U(2, 1)$ is given
as follows:
$$\eta(w_1, w_2):=v_1^{-\frac{1}{12}} v_2^{-\frac{1}{12}}
  (1-v_1)^{-\frac{1}{18}} (1-v_2)^{-\frac{1}{18}} (v_1-v_2)^{-\frac{1}{18}}
  \left[\frac{\partial(v_1, v_2)}{\partial(w_1, w_2)}\right]^{\frac{1}{9}}.\tag 1.39$$
It satisfies the following transformation:
$$\eta \left(\frac{a_1 w_1+a_2 w_2+a_3}{c_1 w_1+c_2 w_2+c_3},
  \frac{b_1 w_1+b_2 w_2+b_3}{c_1 w_1+c_2 w_2+c_3} \right)
 =\Delta^{-\frac{1}{9}} (c_1 w_1+c_2 w_2+c_3)^{\frac{1}{3}} \eta(w_1, w_2),\tag 1.40$$
where $\left(\matrix a_1 & a_2 & a_3\\ b_1 & b_2 & b_3\\ c_1 & c_2
& c_3 \endmatrix\right) \in GL(3, {\Bbb C})$ and
$\Delta=\vmatrix \format \c \quad & \c \quad & \c\\
a_1 & a_2 & a_3\\ b_1 & b_2 & b_3\\ c_1 & c_2 & c_3
\endvmatrix$.

  Let $U(2, 1)=\{ g \in GL(3, {\Bbb C}): g^{*} J g=J \}$
where $J=\left(\matrix
    &    & 1\\
    & -1 &  \\
  1 &    &
   \endmatrix\right)$. The corresponding complex domain
${\frak S}_{2}=\{(z_1, z_2) \in {\Bbb C}^2: z_1+\overline{z_1}
-z_2 \overline{z_2}>0 \}$ is the simplest example of a bounded
symmetric domain which is not a tube domain. Put
$\rho(Z):=z_1+\overline{z_1}-z_2 \overline{z_2}$ for $Z=(z_1, z_2)
\in {\frak S}_{2}$.

  The imaginary quadratic field $K={\Bbb Q}(\sqrt{-3})$ is
called the field of Eisenstein numbers. Its ring of integers
$${\Cal O}_K=\{ a/2+b \sqrt{-3}/2: a, b \in {\Bbb Z},
  a \equiv b (\text{mod} 2) \}={\Bbb Z}+{\Bbb Z} [\omega]$$
with $\omega=\frac{-1+\sqrt{-3}}{2}$ is called the ring of
Eisenstein integers.

  Let
$$T_1=\left(\matrix
      1 & 1 & -\omega\\
        & 1 &       1\\
        &   &       1
      \endmatrix\right), \quad
  T_2=\left(\matrix
      1 & \omega & -\omega\\
        &      1 & \overline{\omega}\\
        &        & 1
      \endmatrix\right), \quad
  S=-\overline{\omega} J
   =\left(\matrix
                       &                   & -\overline{\omega}\\
                       & \overline{\omega} &     \\
    -\overline{\omega} &                   &
    \endmatrix\right),\tag 1.41$$
$$U_1=\left(\matrix
      1 &         &  \\
        & -\omega &  \\
        &         & 1
      \endmatrix\right), \quad
  U_2=\left(\matrix
      -1 &         &   \\
         & -\omega &   \\
         &         & -1
      \endmatrix\right).\tag 1.42$$
The actions of the above five generators on $(w_1, w_2) \in {\frak
S}_{2}$ are given as follows:
$$\left\{\aligned
  T_1(w_1, w_2) &=(w_1+w_2-\omega, w_2+1),\\
  T_2(w_1, w_2) &=(w_1+\omega w_2-\omega, w_2+\overline{\omega}),\\
    S(w_1, w_2) &=\left(\frac{1}{w_1}, -\frac{w_2}{w_1}\right),\\
  U_1(w_1, w_2) &=(w_1, -\omega w_2),\\
  U_2(w_1, w_2) &=(w_1, \omega w_2).
  \endaligned\right.\tag 1.43$$

  Put
$$\eta_1(w_1, w_2):=\eta\left(\left(\matrix
  3 &          &  \\
    & 1-\omega &  \\
    &          & 1
  \endmatrix\right)(w_1, w_2)\right)
 =\eta(3 w_1, (1-\omega) w_2),\tag 1.44$$
$$\eta_2(w_1, w_2):=\eta\left(\left(\matrix
  1 &          &         \\
    & 1-\omega &         \\
    &          & 3
  \endmatrix\right)(w_1, w_2)\right)
 =\eta\left(\frac{w_1}{3}, \frac{w_2}{1-\overline{\omega}}\right).\tag 1.45$$
$$\eta_3(w_1, w_2):=\eta_2([T_1, T_2]^{-1}(w_1, w_2)).\tag 1.46$$
$$\eta_4(w_1, w_2):=\eta_1(S^3 [T_1, T_2]^{-1} S^3(w_1, w_2)).\tag 1.47$$
$$\eta_5(w_1, w_2):=\eta_2([T_1, T_2]^{-1} S^3 [T_1, T_2]^{-1} S^3(w_1, w_2)).\tag 1.48$$
Here,
$$[T_1, T_2]:=T_1 T_2 T_1^{-1} T_2^{-1}
 =\left(\matrix
  1 & 0 & \overline{\omega}-\omega\\
    & 1 & 0\\
    &   & 1
  \endmatrix\right).\tag 1.49$$

  Put
$$\left\{\aligned
  \varphi_0(w_1, w_2)=\frac{\eta_1(w_1, w_2)}{\eta_2(w_1, w_2)}
 &=:\kappa_0^{\frac{1}{9}}=:k_0^{\frac{1}{27}},\\
  \varphi_1(w_1, w_2)=\frac{\eta_1(w_1, w_2)}{\eta_3(w_1, w_2)}
 &=:\kappa_1^{\frac{1}{9}}=:k_1^{\frac{1}{27}},\\
  \varphi_2(w_1, w_2)=\frac{\eta_4(w_1, w_2)}{\eta_2(w_1, w_2)}
 &=:\kappa_2^{\frac{1}{9}}=:k_2^{\frac{1}{27}},\\
  \varphi_3(w_1, w_2)=\frac{\eta_4(w_1, w_2)}{\eta_5(w_1, w_2)}
 &=:\kappa_3^{\frac{1}{9}}=:k_3^{\frac{1}{27}}.
\endaligned\right.\tag 1.50$$
Here, there are four functions $\varphi_0$, $\varphi_1$,
$\varphi_2$ and $\varphi_3$, $C_4^2=6$. In fact, the dimension of
the moduli space of algebraic curves of genus three is:
$\text{dim} {\Cal M}_3=3 \times 3-3=6$. Namely, we define
$$\left\{\aligned
 K_1 &:=\int_{0}^{1} \frac{dx}{\root 3 \of{x^2 (1-x)(1-k_2 x)(1-k_3 x)}},\\
 K_2 &:=\int_{0}^{1} \frac{dx}{\root 3 \of{x^2 (1-x)(1-k_3 x)(1-k_1 x)}},\\
 K_3 &:=\int_{0}^{1} \frac{dx}{\root 3 \of{x^2 (1-x)(1-k_1 x)(1-k_2 x)}}.
\endaligned\right.\tag 1.51$$
$$\left\{\aligned
 K_4 &:=\int_{0}^{1} \frac{dx}{\root 3 \of{x^2 (1-x)(1-k_0 x)(1-k_1 x)}},\\
 K_5 &:=\int_{0}^{1} \frac{dx}{\root 3 \of{x^2 (1-x)(1-k_0 x)(1-k_2 x)}},\\
 K_6 &:=\int_{0}^{1} \frac{dx}{\root 3 \of{x^2 (1-x)(1-k_0 x)(1-k_3 x)}}.
\endaligned\right.\tag 1.52$$

  From the point of view of group representation theory, the
classical trigonometric integrals and elliptic integrals can be
interpreted as follows:
$$\text{$GL(1)$: trigonometric functions, cyclotomic units.}$$
$$\frac{dy}{\sqrt{1-y^2}}=\frac{1}{M} \frac{dx}{\sqrt{1-x^2}}.$$
$$\text{$GL(2)$: elliptic functions, elliptic units, complex
        multiplication.}$$
$$\frac{dy}{\sqrt{(1-y^2)(1-\lambda^2 y^2)}}=\frac{1}{M}
  \frac{dx}{\sqrt{(1-x^2)(1-\kappa^2 x^2)}}.$$
For the group $GL(3)$, we have the following theorem which can be
considered as a natural generalization of the above story from
$GL(1)$, $GL(2)$ to $GL(3)$. In fact, it extends the concept of
complex multiplication.

  Let us consider the following two homogeneous polynomials of
degree $9$:
$$\left\{\aligned
 f_1(x, y, z) &:=x^2 y^2 z^2 (z-x)(z-\lambda_1^3 x)(z-\lambda_2^3 x),\\
 f_2(x, y, z) &:=x^2 y^2 z^2 (z-y)(z-\lambda_0^3 y)(z-\lambda_3^3 y).
\endaligned\right.\tag 1.53$$
We will study the following two differential forms:
$$\omega_1:=\frac{x dy \wedge dz+y dz \wedge dx+z dx \wedge dy}
            {\root 3 \of{f_1(x, y, z)}}, \quad
  \omega_2:=\frac{x dy \wedge dz+y dz \wedge dx+z dx \wedge dy}
            {\root 3 \of{f_2(x, y, z)}}.\tag 1.54$$

{\smc Theorem 1.12 (Main Theorem 3)}. {\it The following identity
of differential forms holds:
$$\frac{dw_1 \wedge dw_2}{\root 3 \of{w_1^2 w_2^2 (1-w_1)
  (1-\lambda_1^3 w_1)(1-\lambda_2^3 w_1)}}
 =-\frac{5 t_1}{t_2^2} \frac{dt_1 \wedge dt_2}{\root 3
  \of{t_1^2 t_2^2 (1-t_1)(1-\kappa_1^3 t_1)(1-\kappa_2^3 t_1)}},\tag 1.55$$
where
$$w_1=\frac{(\beta+\gamma) t_1+\alpha}{\beta t_1+(\alpha+\gamma)}, \quad
  w_2=\frac{t_1 [(\beta+\gamma) t_1+\alpha]^2 [\beta t_1+(\alpha+\gamma)]}
      {t_2^5},\tag 1.56$$
which is a rational transformation of order five. Here,
$$\left\{\aligned
  \alpha &=\frac{(v_1-1)(v_2-1)(v_1+v_2-2)-(u_1-1)(u_2-1)(u_1+u_2-2)}
           {2(u_1-1)(u_2-1)(v_1-1)(v_2-1)},\\
  \beta  &=\frac{-(v_1-1)(v_2-1)(2 v_1 v_2-v_1-v_2)+(u_1-1)(u_2-1)(2
           u_1 u_2-u_1-u_2)}{2(u_1-1)(u_2-1)(v_1-1)(v_2-1)},\\
  \gamma &=\frac{(v_1-1)(v_2-1)}{(u_1-1)(u_2-1)},
\endaligned\right.\tag 1.57$$
with
$$u_1=\kappa_1^3, \quad u_2=\kappa_2^3, \quad
  v_1=\lambda_1^3, \quad v_2=\lambda_2^3.$$
Moreover, the moduli $(\kappa_1, \kappa_2)$ and $(\lambda_1,
\lambda_2)$ satisfy the modular equation:
$$(\kappa_1^3-1)(\kappa_2^3-1)(\kappa_1^3-\kappa_2^3)
 =(\lambda_1^3-1)(\lambda_2^3-1)(\lambda_1^3-\lambda_2^3),\tag 1.58$$
which is an algebraic variety of dimension three. The
corresponding two algebraic surfaces are:
$$\left\{\aligned
 w_3^3 &=w_1^2 w_2^2 (1-w_1)(1-\lambda_1^3 w_1)(1-\lambda_2^3 w_1),\\
 t_3^3 &=t_1^2 t_2^2 (1-t_1)(1-\kappa_1^3 t_1)(1-\kappa_2^3 t_1),
\endaligned\right.\tag 1.59$$
which give the homogeneous algebraic equations of degree seven:
$$\left\{\aligned
 w_3^3 w_4^4 &=w_1^2 w_2^2 (w_4-w_1)(w_4-\lambda_1^3 w_1)(w_4-\lambda_2^3 w_1),\\
 t_3^3 t_4^4 &=t_1^2 t_2^2 (t_4-t_1)(t_4-\kappa_1^3 t_1)(t_4-\kappa_2^3 t_1).
\endaligned\right.\tag 1.60$$}

{\smc Corollary 1.13}. {\it The following identity holds:
$$\frac{dy_1 \wedge dy_2}{\root 3 \of{(1-y_1^3)(1-\lambda_1^3
  y_1^3)(1-\lambda_2^3 y_1^3)}}=-\frac{5 x_1^3}{x_2^6} \frac{dx_1
  \wedge dx_2}{\root 3 \of{(1-x_1^3)(1-\kappa_1^3 x_1^3)(1-
  \kappa_2^3 x_1^3)}},\tag 1.61$$
where
$$y_1^3=\frac{(\beta+\gamma) x_1^3+\alpha}{\beta x_1^3+(\alpha+\gamma)}, \quad
  y_2^3=\frac{x_1^3 [(\beta+\gamma) x_1^3+\alpha]^2 [\beta x_1^3+(\alpha+\gamma)]}
      {x_2^{15}},\tag 1.62$$
and
$$\left\{\aligned
  \alpha &=\frac{(v_1-1)(v_2-1)(v_1+v_2-2)-(u_1-1)(u_2-1)(u_1+u_2-2)}
           {2(u_1-1)(u_2-1)(v_1-1)(v_2-1)},\\
  \beta  &=\frac{-(v_1-1)(v_2-1)(2 v_1 v_2-v_1-v_2)+(u_1-1)(u_2-1)(2
           u_1 u_2-u_1-u_2)}{2(u_1-1)(u_2-1)(v_1-1)(v_2-1)},\\
  \gamma &=\frac{(v_1-1)(v_2-1)}{(u_1-1)(u_2-1)},
\endaligned\right.$$
with
$$u_1=\kappa_1^3, \quad u_2=\kappa_2^3, \quad
  v_1=\lambda_1^3, \quad v_2=\lambda_2^3.$$
Moreover, the moduli $(\kappa_1, \kappa_2)$ and $(\lambda_1,
\lambda_2)$ satisfy the modular equation:
$$(\kappa_1^3-1)(\kappa_2^3-1)(\kappa_1^3-\kappa_2^3)
 =(\lambda_1^3-1)(\lambda_2^3-1)(\lambda_1^3-\lambda_2^3).$$}

  We will study the system of partial differential equations:
$$\left\{\aligned
 \{w_1, w_2; v_1, v_2\}_{v_1} &=F_1(\gamma; v_1, v_2),\\
 \{w_1, w_2; v_1, v_2\}_{v_2} &=F_2(\gamma; v_1, v_2),\\
 [w_1, w_2; v_1, v_2]_{v_1} &=P_1(\alpha, \beta, \gamma; v_1, v_2),\\
 [w_1, w_2; v_1, v_2]_{v_2} &=P_2(\alpha, \beta, \gamma; v_1, v_2),
\endaligned\right.\tag 1.63$$
where
$$\left\{\aligned
 F_1(\gamma; v_1, v_2) &=\frac{-\gamma v_2(v_2-1)}{v_1(v_1-1)(v_1-v_2)},\\
 F_2(\gamma; v_1, v_2) &=\frac{-\gamma v_1(v_1-1)}{v_2(v_2-1)(v_2-v_1)},\\
 P_1(\alpha, \beta, \gamma; v_1, v_2) &=\frac{\alpha}{v_1}+
 \frac{\beta}{v_1-1}+\frac{\gamma}{v_1-v_2},\\
 P_2(\alpha, \beta, \gamma; v_1, v_2) &=\frac{\alpha}{v_2}+
 \frac{\beta}{v_2-1}+\frac{\gamma}{v_2-v_1}.
\endaligned\right.\tag 1.64$$
The solutions are denoted by
$$w_1=s_1(\alpha, \beta, \gamma; v_1, v_2), \quad
  w_2=s_2(\alpha, \beta, \gamma; v_1, v_2).\tag 1.65$$
Consequently,
$$\left\{\aligned
 \{s_1(\alpha, \beta, \gamma; v_1, v_2), s_2(\alpha, \beta, \gamma;
 v_1, v_2); v_1, v_2\}_{v_1} &=F_1(\gamma; v_1, v_2),\\
 \{s_1(\alpha, \beta, \gamma; v_1, v_2), s_2(\alpha, \beta, \gamma;
 v_1, v_2); v_1, v_2\}_{v_2} &=F_2(\gamma; v_1, v_2),\\
 [s_1(\alpha, \beta, \gamma; v_1, v_2), s_2(\alpha, \beta, \gamma;
 v_1, v_2); v_1, v_2]_{v_1} &=P_1(\alpha, \beta, \gamma; v_1, v_2),\\
 [s_1(\alpha, \beta, \gamma; v_1, v_2), s_2(\alpha, \beta, \gamma;
 v_1, v_2); v_1, v_2]_{v_2} &=P_2(\alpha, \beta, \gamma; v_1, v_2),
\endaligned\right.\tag 1.66$$

  Let
$$T=\left(\matrix
    -1 &  0 & 1\\
       & -1 & 1\\
       &    & 1
    \endmatrix\right), \quad
  S_1=\left(\matrix
      1 & 0 & 0\\
      0 & 0 & 1\\
      0 & 1 & 0
      \endmatrix\right), \quad
  S_2=\left(\matrix
        &   & 1\\
        & 1 &  \\
      1 &   &
      \endmatrix\right).\tag 1.67$$
Then
$$T(x, y)=(1-x, 1-y).\tag 1.68$$
$$\left\{\aligned
  S_1(x, y) &=\left(\frac{x}{y}, \frac{1}{y}\right),\\
  S_1 T(x, y) &=\left(\frac{1-x}{1-y}, \frac{1}{1-y}\right),\\
  T S_1(x, y) &=\left(\frac{y-x}{y}, \frac{y-1}{y}\right),\\
  S_1 T S_1(x, y) &=\left(\frac{y-x}{y-1}, \frac{y}{y-1}\right).
\endaligned\right.\tag 1.69$$
$$\left\{\aligned
  S_2(x, y) &=\left(\frac{1}{x}, \frac{y}{x}\right),\\
  S_2 T(x, y) &=\left(\frac{1}{1-x}, \frac{1-y}{1-x}\right),\\
  T S_2(x, y) &=\left(\frac{x-1}{x}, \frac{x-y}{x}\right),\\
  S_2 T S_2(x, y) &=\left(\frac{x}{x-1}, \frac{x-y}{x-1}\right).
\endaligned\right.\tag 1.70$$

  We have
$$\left\{\aligned
 \{s_1(\beta, \alpha, \gamma; v_1, v_2), s_2(\beta, \alpha,
 \gamma; v_1, v_2); v_1, v_2 \}_{v_1} &=F_1(-\gamma; T(v_1, v_2)),\\
 \{s_1(\beta, \alpha, \gamma; v_1, v_2), s_2(\beta, \alpha,
 \gamma; v_1, v_2); v_1, v_2 \}_{v_2} &=F_2(-\gamma; T(v_1, v_2)),\\
 [s_1(\beta, \alpha, \gamma; v_1, v_2), s_2(\beta, \alpha,
 \gamma; v_1, v_2); v_1, v_2 ]_{v_1} &=-P_1(\alpha, \beta, \gamma; T(v_1, v_2)),\\
 [s_1(\beta, \alpha, \gamma; v_1, v_2), s_2(\beta, \alpha,
 \gamma; v_1, v_2); v_1, v_2 ]_{v_2} &=-P_2(\alpha, \beta, \gamma; T(v_1, v_2)).
\endaligned\right.\tag 1.71$$
$$\left\{\aligned
 \{s_1(\alpha, \gamma, \beta; v_1, v_2), s_2(\alpha, \gamma,
 \beta; v_1, v_2); v_1, v_2 \}_{v_1} &=F_1(-\beta; S_1(v_1, v_2)),\\
 \{s_1(\alpha, \gamma, \beta; v_1, v_2), s_2(\alpha, \gamma,
 \beta; v_1, v_2); v_1, v_2 \}_{v_2} &=F_2(-\beta; S_2(v_1, v_2)),\\
 [s_1(\alpha, \gamma, \beta; v_1, v_2), s_2(\alpha, \gamma,
 \beta; v_1, v_2); v_1, v_2 ]_{v_1} &=\frac{1}{v_2} P_1(\alpha, \beta, \gamma; S_1(v_1, v_2)),\\
 [s_1(\alpha, \gamma, \beta; v_1, v_2), s_2(\alpha, \gamma,
 \beta; v_1, v_2); v_1, v_2 ]_{v_2} &=\frac{1}{v_1} P_2(\alpha, \beta, \gamma; S_2(v_1, v_2)).
\endaligned\right.\tag 1.72$$
$$\left\{\aligned
 \{s_1(\gamma, \alpha, \beta; v_1, v_2), s_2(\gamma, \alpha,
 \beta; v_1, v_2); v_1, v_2 \}_{v_1} &=F_1(\beta; S_1 T(v_1, v_2)),\\
 \{s_1(\gamma, \alpha, \beta; v_1, v_2), s_2(\gamma, \alpha,
 \beta; v_1, v_2); v_1, v_2 \}_{v_2} &=F_2(\beta; S_2 T(v_1, v_2)),\\
 [s_1(\gamma, \alpha, \beta; v_1, v_2), s_2(\gamma, \alpha,
 \beta; v_1, v_2); v_1, v_2 ]_{v_1} &=\frac{1}{v_2-1} P_1(\alpha, \beta, \gamma; S_1 T(v_1, v_2)),\\
 [s_1(\gamma, \alpha, \beta; v_1, v_2), s_2(\gamma, \alpha,
 \beta; v_1, v_2); v_1, v_2 ]_{v_2} &=\frac{1}{v_1-1} P_2(\alpha, \beta, \gamma; S_2 T(v_1, v_2)).
\endaligned\right.\tag 1.73$$
$$\left\{\aligned
 \{s_1(\beta, \gamma, \alpha; v_1, v_2), s_2(\beta, \gamma,
 \alpha; v_1, v_2); v_1, v_2 \}_{v_1} &=F_1(\alpha; T S_1(v_1, v_2)),\\
 \{s_1(\beta, \gamma, \alpha; v_1, v_2), s_2(\beta, \gamma,
 \alpha; v_1, v_2); v_1, v_2 \}_{v_2} &=F_2(\alpha; T S_2(v_1, v_2)),\\
 [s_1(\beta, \gamma, \alpha; v_1, v_2), s_2(\beta, \gamma,
 \alpha; v_1, v_2); v_1, v_2 ]_{v_1} &=-\frac{1}{v_2} P_1(\alpha, \beta, \gamma; T S_1(v_1, v_2)),\\
 [s_1(\beta, \gamma, \alpha; v_1, v_2), s_2(\beta, \gamma,
 \alpha; v_1, v_2); v_1, v_2 ]_{v_2} &=-\frac{1}{v_1} P_2(\alpha, \beta, \gamma; T S_2(v_1, v_2)).
\endaligned\right.\tag 1.74$$
$$\left\{\aligned
 \{s_1(\gamma, \beta, \alpha; v_1, v_2), s_2(\gamma, \beta,
 \alpha; v_1, v_2); v_1, v_2 \}_{v_1} &=F_1(-\alpha; S_1 T S_1(v_1, v_2)),\\
 \{s_1(\gamma, \beta, \alpha; v_1, v_2), s_2(\gamma, \beta,
 \alpha; v_1, v_2); v_1, v_2 \}_{v_2} &=F_2(-\alpha; S_2 T S_2(v_1, v_2)),\\
 [s_1(\gamma, \beta, \alpha; v_1, v_2), s_2(\gamma, \beta,
 \alpha; v_1, v_2); v_1, v_2 ]_{v_1} &=\frac{1}{1-v_2} P_1(\alpha, \beta, \gamma; S_1 T S_1(v_1, v_2)),\\
 [s_1(\gamma, \beta, \alpha; v_1, v_2), s_2(\gamma, \beta,
 \alpha; v_1, v_2); v_1, v_2 ]_{v_2} &=\frac{1}{1-v_1} P_2(\alpha, \beta, \gamma; S_2 T S_2(v_1, v_2)).
\endaligned\right.\tag 1.75$$

{\smc Proposition 1.14 (Main Proposition 3)}. {\it The triangular
$s$-functions have the following properties of transformations:
\roster
\item
$$\left\{\aligned
  \{s_1(\alpha, \beta, \gamma; T(v_1, v_2)), s_2(\alpha, \beta, \gamma;
  T(v_1, v_2)); v_1, v_2\}_{v_1} &=F_1(\gamma; v_1, v_2),\\
  \{s_1(\alpha, \beta, \gamma; T(v_1, v_2)), s_2(\alpha, \beta, \gamma;
  T(v_1, v_2)); v_1, v_2\}_{v_2} &=F_2(\gamma; v_1, v_2),\\
  [s_1(\alpha, \beta, \gamma; T(v_1, v_2)), s_2(\alpha, \beta, \gamma;
  T(v_1, v_2)); v_1, v_2]_{v_1} &=P_1(\beta, \alpha, \gamma; v_1, v_2),\\
  [s_1(\alpha, \beta, \gamma; T(v_1, v_2)), s_2(\alpha, \beta, \gamma;
  T(v_1, v_2)); v_1, v_2]_{v_2} &=P_2(\beta, \alpha, \gamma; v_1, v_2).
\endaligned\right.\tag 1.76$$
\item
$$\left\{\aligned
 \{s_1(\alpha, \beta, \gamma; S_1(v_1, v_2)), s_2(\alpha, \beta, \gamma;
  S_1(v_1, v_2)); v_1, v_2\}_{v_1} =&F_1(\gamma; v_1, v_2),\\
 \{s_1(\alpha, \beta, \gamma; S_2(v_1, v_2)), s_2(\alpha, \beta, \gamma;
  S_2(v_1, v_2)); v_1, v_2\}_{v_2} =&F_2(\gamma; v_1, v_2),\\
  [s_1(\alpha, \beta, \gamma; S_1(v_1, v_2)), s_2(\alpha, \beta, \gamma;
  S_1(v_1, v_2)); v_1, v_2]_{v_1} =&P_1(\alpha, -2 \gamma, \beta+3 \gamma; v_1, v_2),\\
  [s_1(\alpha, \beta, \gamma; S_2(v_1, v_2)), s_2(\alpha, \beta, \gamma;
  S_2(v_1, v_2)); v_1, v_2]_{v_2} =&P_2(\alpha, -2 \gamma, \beta+3 \gamma; v_1, v_2).
\endaligned\right.\tag 1.77$$
\item
$$\left\{\aligned
 \{s_1(\alpha, \beta, \gamma; S_1 T(v_1, v_2)), s_2(\alpha, \beta, \gamma;
  S_1 T(v_1, v_2)); v_1, v_2\}_{v_1} =&F_1(\gamma; v_1, v_2),\\
 \{s_1(\alpha, \beta, \gamma; S_2 T(v_1, v_2)), s_2(\alpha, \beta, \gamma;
  S_2 T(v_1, v_2)); v_1, v_2\}_{v_2} =&F_2(\gamma; v_1, v_2),\\
  [s_1(\alpha, \beta, \gamma; S_1 T(v_1, v_2)), s_2(\alpha, \beta, \gamma;
  S_1 T(v_1, v_2)); v_1, v_2]_{v_1} =&P_1(-2 \gamma, \alpha, \beta+3 \gamma; v_1, v_2),\\
  [s_1(\alpha, \beta, \gamma; S_2 T(v_1, v_2)), s_2(\alpha, \beta, \gamma;
  S_2 T(v_1, v_2)); v_1, v_2]_{v_2} =&P_2(-2 \gamma, \alpha, \beta+3 \gamma; v_1, v_2).
\endaligned\right.\tag 1.78$$
\item
$$\left\{\aligned
 \{s_1(\alpha, \beta, \gamma; T S_1(v_1, v_2)), s_2(\alpha, \beta, \gamma;
  T S_1(v_1, v_2)); v_1, v_2\}_{v_1} =&F_1(\gamma; v_1, v_2),\\
 \{s_1(\alpha, \beta, \gamma; T S_2(v_1, v_2)), s_2(\alpha, \beta, \gamma;
  T S_2(v_1, v_2)); v_1, v_2\}_{v_2} =&F_2(\gamma; v_1, v_2),\\
  [s_1(\alpha, \beta, \gamma; T S_1(v_1, v_2)), s_2(\alpha, \beta, \gamma;
  T S_1(v_1, v_2)); v_1, v_2]_{v_1} =&P_1(\beta, -2 \gamma, \alpha+3 \gamma; v_1, v_2),\\
  [s_1(\alpha, \beta, \gamma; T S_2(v_1, v_2)), s_2(\alpha, \beta, \gamma;
  T S_2(v_1, v_2)); v_1, v_2]_{v_2} =&P_2(\beta, -2 \gamma, \alpha+3 \gamma; v_1, v_2).
\endaligned\right.\tag 1.79$$
\item
$$\left\{\aligned
 &\{s_1(\alpha, \beta, \gamma; S_1 T S_1(v_1, v_2)), s_2(\alpha, \beta, \gamma;
  S_1 T S_1(v_1, v_2)); v_1, v_2\}_{v_1}=F_1(\gamma; v_1, v_2),\\
 &\{s_1(\alpha, \beta, \gamma; S_2 T S_2(v_1, v_2)), s_2(\alpha, \beta, \gamma;
  S_2 T S_2(v_1, v_2)); v_1, v_2\}_{v_2}=F_2(\gamma; v_1, v_2),\\
 &[s_1(\alpha, \beta, \gamma; S_1 T S_1(v_1, v_2)), s_2(\alpha, \beta, \gamma;
  S_1 T S_1(v_1, v_2)); v_1, v_2]_{v_1}\\
=&P_1(-2 \gamma, \beta, \alpha+3 \gamma; v_1, v_2),\\
 &[s_1(\alpha, \beta, \gamma; S_2 T S_2(v_1, v_2)), s_2(\alpha, \beta, \gamma;
  S_2 T S_2(v_1, v_2)); v_1, v_2]_{v_2}\\
=&P_2(-2 \gamma, \beta, \alpha+3 \gamma; v_1, v_2).
\endaligned\right.\tag 1.80$$
\endroster}

{\smc Corollary 1.15}. {\it The triangular $s$-functions have the
following properties of transformations: \roster
\item
$$\left\{\aligned
 &\{s_1(\alpha, \beta, \gamma; T(v_1, v_2)), s_2(\alpha, \beta, \gamma;
  T(v_1, v_2)); v_1, v_2\}_{v_1}\\
=&\{s_1(\alpha, \beta, \gamma; v_1, v_2), s_2(\alpha, \beta,
  \gamma; v_1, v_2); v_1, v_2\}_{v_1},\\
 &\{s_1(\alpha, \beta, \gamma; T(v_1, v_2)), s_2(\alpha, \beta, \gamma;
  T(v_1, v_2)); v_1, v_2\}_{v_2}\\
=&\{s_1(\alpha, \beta, \gamma; v_1, v_2), s_2(\alpha, \beta,
  \gamma; v_1, v_2); v_1, v_2\}_{v_2},\\
 &[s_1(\alpha, \beta, \gamma; T(v_1, v_2)), s_2(\alpha, \beta, \gamma;
  T(v_1, v_2)); v_1, v_2]_{v_1}\\
=&[s_1(\beta, \alpha, \gamma; v_1, v_2), s_2(\beta, \alpha,
  \gamma; v_1, v_2); v_1, v_2]_{v_1},\\
 &[s_1(\alpha, \beta, \gamma; T(v_1, v_2)), s_2(\alpha, \beta, \gamma;
  T(v_1, v_2)); v_1, v_2]_{v_2}\\
=&[s_1(\beta, \alpha, \gamma; v_1, v_2), s_2(\beta, \alpha,
  \gamma; v_1, v_2); v_1, v_2]_{v_2}.
\endaligned\right.\tag 1.81$$
\item
$$\left\{\aligned
 &\{s_1(\alpha, \beta, \gamma; S_1(v_1, v_2)), s_2(\alpha, \beta, \gamma;
  S_1(v_1, v_2)); v_1, v_2\}_{v_1}\\
=&\{s_1(\alpha, \beta, \gamma; v_1, v_2), s_2(\alpha, \beta,
  \gamma; v_1, v_2); v_1, v_2\}_{v_1},\\
 &\{s_1(\alpha, \beta, \gamma; S_2(v_1, v_2)), s_2(\alpha, \beta, \gamma;
  S_2(v_1, v_2)); v_1, v_2\}_{v_2}\\
=&\{s_1(\alpha, \beta, \gamma; v_1, v_2), s_2(\alpha, \beta,
  \gamma; v_1, v_2); v_1, v_2\}_{v_2},\\
 &[s_1(\alpha, \beta, \gamma; S_1(v_1, v_2)), s_2(\alpha, \beta, \gamma;
  S_1(v_1, v_2)); v_1, v_2]_{v_1}\\
=&[s_1(\alpha, -2 \gamma, \beta+3 \gamma; v_1, v_2), s_2(\alpha,
  -2 \gamma, \beta+3 \gamma; v_1, v_2); v_1, v_2]_{v_1},\\
 &[s_1(\alpha, \beta, \gamma; S_2(v_1, v_2)), s_2(\alpha, \beta, \gamma;
  S_2(v_1, v_2)); v_1, v_2]_{v_2}\\
=&[s_1(\alpha, -2 \gamma, \beta+3 \gamma; v_1, v_2), s_2(\alpha,
  -2 \gamma, \beta+3 \gamma; v_1, v_2); v_1, v_2]_{v_2}.
\endaligned\right.\tag 1.82$$
\item
$$\left\{\aligned
 &\{s_1(\alpha, \beta, \gamma; S_1 T(v_1, v_2)), s_2(\alpha, \beta, \gamma;
  S_1 T(v_1, v_2)); v_1, v_2\}_{v_1}\\
=&\{s_1(\alpha, \beta, \gamma; v_1, v_2), s_2(\alpha, \beta,
  \gamma; v_1, v_2); v_1, v_2\}_{v_1},\\
 &\{s_1(\alpha, \beta, \gamma; S_2 T(v_1, v_2)), s_2(\alpha, \beta, \gamma;
  S_2 T(v_1, v_2)); v_1, v_2\}_{v_2}\\
=&\{s_1(\alpha, \beta, \gamma; v_1, v_2), s_2(\alpha, \beta,
  \gamma; v_1, v_2); v_1, v_2\}_{v_2},\\
 &[s_1(\alpha, \beta, \gamma; S_1 T(v_1, v_2)), s_2(\alpha, \beta, \gamma;
  S_1 T(v_1, v_2)); v_1, v_2]_{v_1}\\
=&[s_1(-2 \gamma, \alpha, \beta+3 \gamma; v_1, v_2), s_2(-2
  \gamma, \alpha, \beta+3 \gamma; v_1, v_2); v_1, v_2]_{v_1},\\
 &[s_1(\alpha, \beta, \gamma; S_2 T(v_1, v_2)), s_2(\alpha, \beta, \gamma;
  S_2 T(v_1, v_2)); v_1, v_2]_{v_2}\\
=&[s_1(-2 \gamma, \alpha, \beta+3 \gamma; v_1, v_2), s_2(-2
  \gamma, \alpha, \beta+3 \gamma; v_1, v_2); v_1, v_2]_{v_2}.
\endaligned\right.\tag 1.83$$
\item
$$\left\{\aligned
 &\{s_1(\alpha, \beta, \gamma; T S_1(v_1, v_2)), s_2(\alpha, \beta, \gamma;
  T S_1(v_1, v_2)); v_1, v_2\}_{v_1}\\
=&\{s_1(\alpha, \beta, \gamma; v_1, v_2), s_2(\alpha, \beta,
  \gamma; v_1, v_2); v_1, v_2\}_{v_1},\\
 &\{s_1(\alpha, \beta, \gamma; T S_2(v_1, v_2)), s_2(\alpha, \beta, \gamma;
  T S_2(v_1, v_2)); v_1, v_2\}_{v_2}\\
=&\{s_1(\alpha, \beta, \gamma; v_1, v_2), s_2(\alpha, \beta,
  \gamma; v_1, v_2); v_1, v_2\}_{v_2},\\
 &[s_1(\alpha, \beta, \gamma; T S_1(v_1, v_2)), s_2(\alpha, \beta, \gamma;
  T S_1(v_1, v_2)); v_1, v_2]_{v_1}\\
=&[s_1(\beta, -2 \gamma, \alpha+3 \gamma; v_1, v_2), s_2(\beta, -2
  \gamma, \alpha+3 \gamma; v_1, v_2); v_1, v_2]_{v_1},\\
 &[s_1(\alpha, \beta, \gamma; T S_2(v_1, v_2)), s_2(\alpha, \beta, \gamma;
  T S_2(v_1, v_2)); v_1, v_2]_{v_2}\\
=&[s_1(\beta, -2 \gamma, \alpha+3 \gamma; v_1, v_2), s_2(\beta, -2
  \gamma, \alpha+3 \gamma; v_1, v_2); v_1, v_2]_{v_2}.
\endaligned\right.\tag 1.84$$
\item
$$\left\{\aligned
 &\{s_1(\alpha, \beta, \gamma; S_1 T S_1(v_1, v_2)), s_2(\alpha, \beta, \gamma;
  S_1 T S_1(v_1, v_2)); v_1, v_2\}_{v_1}\\
=&\{s_1(\alpha, \beta, \gamma; v_1, v_2), s_2(\alpha, \beta,
  \gamma; v_1, v_2); v_1, v_2\}_{v_1},\\
 &\{s_1(\alpha, \beta, \gamma; S_2 T S_2(v_1, v_2)), s_2(\alpha, \beta, \gamma;
  S_2 T S_2(v_1, v_2)); v_1, v_2\}_{v_2}\\
=&\{s_1(\alpha, \beta, \gamma; v_1, v_2), s_2(\alpha, \beta,
  \gamma; v_1, v_2); v_1, v_2\}_{v_2},\\
 &[s_1(\alpha, \beta, \gamma; S_1 T S_1(v_1, v_2)), s_2(\alpha, \beta, \gamma;
  S_1 T S_1(v_1, v_2)); v_1, v_2]_{v_1}\\
=&[s_1(-2 \gamma, \beta, \alpha+3 \gamma; v_1, v_2), s_2(-2
  \gamma, \beta, \alpha+3 \gamma; v_1, v_2); v_1, v_2]_{v_1},\\
 &[s_1(\alpha, \beta, \gamma; S_2 T S_2(v_1, v_2)), s_2(\alpha, \beta, \gamma;
  S_2 T S_2(v_1, v_2)); v_1, v_2]_{v_2}\\
=&[s_1(-2 \gamma, \beta, \alpha+3 \gamma; v_1, v_2), s_2(-2
  \gamma, \beta, \alpha+3 \gamma; v_1, v_2); v_1, v_2]_{v_2}.
\endaligned\right.\tag 1.85$$
\endroster}

  The significant properties of Picard curves come from the
following facts:

{\it Definition}. Suppose that $C$ is an algebraic curve of genus
$g \geq 2$. Let $\{ \omega_1, \cdots, \omega_g \}$ be a basis of
$\Omega^{1}(C)$. Then the map
$$\aligned
  \varphi_{K}: C &\to {\Bbb P}^{g-1}\\
               p &\mapsto [\omega_1(p), \cdots, \omega_g(p)]
\endaligned$$
is called a canonical map over $C$.

  In fact, Suppose that $z$ is a local coordinate near the point $p$
and $\omega_{\alpha}=f_{\alpha}(z) dz$ ($\alpha=1, 2, \cdots, g$).
Then
$$\varphi_{K}(p)=[f_1(z(p)), \cdots, f_{g}(z(p))].$$

{\smc Theorem}. {\it Any algebraic curve of genus $g=2$ is a
hyperelliptic curve.}

{\it Definition}. Suppose that $C$ is a non-hyperelliptic
algebraic curve of genus $g \geq 2$. Let
$$\varphi_{K}: C \to {\Bbb P}^{g-1}$$
be a canonical map. Then
$$\varphi_{K}(C) \subset {\Bbb P}^{g-1}$$
is called a canonical curve.

  Let $C$ be a canonical curve of genus $g$. Then
$\text{deg}(C)=2g-2$.

{\smc Theorem}. {\it Any canonical curve of genus $g=3$ is a
smooth quartic plane algebraic curve.}

{\smc Theorem}. {\it Any canonical curve of genus $g=4$ is the
intersection of a quadratic surface and a cubic surface in ${\Bbb
P}^{3}$.}

  For $0, 1, \lambda_1, \lambda_2, \infty \in {\Bbb P}^{1}({\Bbb
C})$, $\lambda_1, \lambda_2 \notin \{0, 1\}$ and $\lambda_1 \neq
\lambda_2$, set
$$\left\{\aligned
  J_1=J_1(\lambda_1, \lambda_2):=\frac{\lambda_2^2 (\lambda_2-1)^2}
      {\lambda_1^2 (\lambda_1-1)^2 (\lambda_1-\lambda_2)^2},\\
  J_2=J_2(\lambda_1, \lambda_2):=\frac{\lambda_1^2 (\lambda_1-1)^2}
      {\lambda_2^2 (\lambda_2-1)^2 (\lambda_2-\lambda_1)^2}.
\endaligned\right.\tag 1.86$$
The corresponding Picard curve is
$$y^3=x(x-1)(x-\lambda_1)(x-\lambda_2).\tag 1.87$$
The functions $J_1$ and $J_2$ satisfy the following relations:
$$\aligned
 &J_1(\lambda_1, \lambda_2)=J_1(T(\lambda_1, \lambda_2))
 =J_1(S_1(\lambda_1, \lambda_2))=J_1(S_1 T(\lambda_1, \lambda_2))\\
=&J_1(T S_1(\lambda_1, \lambda_2))=J_1(S_1 T S_1(\lambda_1,
  \lambda_2)).
\endaligned\tag 1.88$$
$$\aligned
 &J_2(\lambda_1, \lambda_2)=J_2(T(\lambda_1, \lambda_2))
 =J_2(S_2(\lambda_1, \lambda_2))=J_2(S_2 T(\lambda_1, \lambda_2))\\
=&J_2(T S_2(\lambda_1, \lambda_2))=J_2(S_2 T S_2(\lambda_1,
  \lambda_2)).
\endaligned\tag 1.89$$

  We will study the following four evolution equations for $GL(3)$:
$$\frac{\vmatrix \format \c \quad & \c\\
  \frac{\partial u_1}{\partial t_1} &
  \frac{\partial u_2}{\partial t_1}\\
  \frac{\partial u_1}{\partial x} &
  \frac{\partial u_2}{\partial x}
  \endvmatrix}{\vmatrix \format \c \quad & \c\\
  \frac{\partial u_1}{\partial y} &
  \frac{\partial u_2}{\partial y}\\
  \frac{\partial u_1}{\partial x} &
  \frac{\partial u_2}{\partial x}
  \endvmatrix}
 =\{u_1, u_2; x, y\}_{x}, \quad
  \frac{\vmatrix \format \c \quad & \c\\
  \frac{\partial u_1}{\partial t_1} &
  \frac{\partial u_2}{\partial t_1}\\
  \frac{\partial u_1}{\partial y} &
  \frac{\partial u_2}{\partial y}
  \endvmatrix}{\vmatrix \format \c \quad & \c\\
  \frac{\partial u_1}{\partial y} &
  \frac{\partial u_2}{\partial y}\\
  \frac{\partial u_1}{\partial x} &
  \frac{\partial u_2}{\partial x}
  \endvmatrix}
 =[u_1, u_2; x, y]_{x},\tag 1.90$$
$$\frac{\vmatrix \format \c \quad & \c\\
  \frac{\partial u_1}{\partial t_2} &
  \frac{\partial u_2}{\partial t_2}\\
  \frac{\partial u_1}{\partial x} &
  \frac{\partial u_2}{\partial x}
  \endvmatrix}{\vmatrix \format \c \quad & \c\\
  \frac{\partial u_1}{\partial x} &
  \frac{\partial u_2}{\partial x}\\
  \frac{\partial u_1}{\partial y} &
  \frac{\partial u_2}{\partial y}
  \endvmatrix}
 =[u_1, u_2; x, y]_{y}, \quad
  \frac{\vmatrix \format \c \quad & \c\\
  \frac{\partial u_1}{\partial t_2} &
  \frac{\partial u_2}{\partial t_2}\\
  \frac{\partial u_1}{\partial y} &
  \frac{\partial u_2}{\partial y}
  \endvmatrix}{\vmatrix \format \c \quad & \c\\
  \frac{\partial u_1}{\partial x} &
  \frac{\partial u_2}{\partial x}\\
  \frac{\partial u_1}{\partial y} &
  \frac{\partial u_2}{\partial y}
  \endvmatrix}
 =\{u_1, u_2; x, y\}_{y},\tag 1.91$$
where the functions $u_1=u_1(x, y; t_1, t_2)$ and $u_2=u_2(x, y;
t_1, t_2)$. The most important property of the above equations is
that they are invariant with respect to the group $GL(3, {\Bbb
C})$ of linear fractional transformations:
$$(u_1, u_2) \mapsto \left(\frac{a_1 u_1+a_2 u_2+a_3}
  {c_1 u_1+c_2 u_2+c_3}, \frac{b_1 u_1+b_2 u_2+b_3}
  {c_1 u_1+c_2 u_2+c_3}\right)$$
with $\gamma=\left(\matrix a_1 & a_2 & a_3\\ b_1 & b_2 & b_3\\
c_1 & c_2 & c_3 \endmatrix\right) \in GL(3, {\Bbb C})$.

{\smc Theorem 1.16 (Main Theorem 4)}. {\it For the system of
partial differential equations:
$$\frac{\vmatrix \format \c \quad & \c\\
  \frac{\partial u_1}{\partial t_1} &
  \frac{\partial u_2}{\partial t_1}\\
  \frac{\partial u_1}{\partial x} &
  \frac{\partial u_2}{\partial x}
  \endvmatrix}{\vmatrix \format \c \quad & \c\\
  \frac{\partial u_1}{\partial y} &
  \frac{\partial u_2}{\partial y}\\
  \frac{\partial u_1}{\partial x} &
  \frac{\partial u_2}{\partial x}
  \endvmatrix}
 =\{u_1, u_2; x, y\}_{x}, \quad
  \frac{\vmatrix \format \c \quad & \c\\
  \frac{\partial u_1}{\partial t_1} &
  \frac{\partial u_2}{\partial t_1}\\
  \frac{\partial u_1}{\partial y} &
  \frac{\partial u_2}{\partial y}
  \endvmatrix}{\vmatrix \format \c \quad & \c\\
  \frac{\partial u_1}{\partial y} &
  \frac{\partial u_2}{\partial y}\\
  \frac{\partial u_1}{\partial x} &
  \frac{\partial u_2}{\partial x}
  \endvmatrix}
 =[u_1, u_2; x, y]_{x},\tag 1.92$$
$$\frac{\vmatrix \format \c \quad & \c\\
  \frac{\partial u_1}{\partial t_2} &
  \frac{\partial u_2}{\partial t_2}\\
  \frac{\partial u_1}{\partial x} &
  \frac{\partial u_2}{\partial x}
  \endvmatrix}{\vmatrix \format \c \quad & \c\\
  \frac{\partial u_1}{\partial x} &
  \frac{\partial u_2}{\partial x}\\
  \frac{\partial u_1}{\partial y} &
  \frac{\partial u_2}{\partial y}
  \endvmatrix}
 =[u_1, u_2; x, y]_{y}, \quad
  \frac{\vmatrix \format \c \quad & \c\\
  \frac{\partial u_1}{\partial t_2} &
  \frac{\partial u_2}{\partial t_2}\\
  \frac{\partial u_1}{\partial y} &
  \frac{\partial u_2}{\partial y}
  \endvmatrix}{\vmatrix \format \c \quad & \c\\
  \frac{\partial u_1}{\partial x} &
  \frac{\partial u_2}{\partial x}\\
  \frac{\partial u_1}{\partial y} &
  \frac{\partial u_2}{\partial y}
  \endvmatrix}
 =\{u_1, u_2; x, y\}_{y},\tag 1.93$$
set
$$\left\{\aligned
  v_1 &:=\{u_1, u_2; x, y\}_{x},\\
  w_1 &:=[u_1, u_2; x, y]_{x},\\
  v_2 &:=\{u_1, u_2; x, y\}_{y},\\
  w_2 &:=[u_1, u_2; x, y]_{y}.
\endaligned\right.\tag 1.94$$
Then $v_1$, $v_2$, $w_1$ and $w_2$ satisfy the following two
equations:
$$\left\{\aligned
  \frac{\partial w_1}{\partial t_2}+\frac{\partial
  v_2}{\partial t_1}-v_1 \frac{\partial v_2}{\partial y}-
  v_2 \frac{\partial w_1}{\partial x}+w_1 \frac{\partial
  v_2}{\partial x}+w_2 \frac{\partial w_1}{\partial y} &=0,\\
  \frac{\partial w_2}{\partial t_1}+\frac{\partial
  v_1}{\partial t_2}-v_1 \frac{\partial w_2}{\partial y}-
  v_2 \frac{\partial v_1}{\partial x}+w_1 \frac{\partial
  w_2}{\partial x}+w_2 \frac{\partial v_1}{\partial y} &=0.
\endaligned\right.\tag 1.95$$}

  In fact, The above two equations are invariant under the
action of
$$\left\{\aligned
  v_1(x, y; t_1, t_2) &\mapsto v_1(x-a_1 t_1, y-b_1 t_2; t_1, t_2),\\
  v_2(x, y; t_1, t_2) &\mapsto v_2(x-a_2 t_1, y-b_2 t_2; t_1, t_2),\\
  w_1(x, y; t_1, t_2) &\mapsto w_1(x-a_1 t_1, y-b_1 t_2; t_1, t_2)+a_2,\\
  w_2(x, y; t_1, t_2) &\mapsto w_2(x-a_2 t_1, y-b_2 t_2; t_1,
  t_2)+b_1.
\endaligned\right.\tag 1.96$$

  This paper consists of seven sections. In section two, we will
study the monodromy of Appell hypergeometric partial differential
equations and give the construction of our four derivatives. We
will investigate the properties of our four derivatives. In
section three, we will study the system of nonlinear partial
differential equations associated to our four derivatives. We will
give a solution which is associated to the Picard modular
function. In section four, we will introduce the $\eta$-functions
associated to the unitary group $U(2, 1)$ and study their
properties, especially the relation with Picard curves. In section
five, we will study the transform problems (moduli problems) and
the associated modular equations which give the connection between
the moduli of two Picard curves. In section six, we will study the
properties of solutions of the above system of nonlinear partial
differential equations, which we call the triangular
$s$-functions. We obtain the $J$-invariants associated to the
Picard curves. In section seven, we will investigate a system of
nonlinear evolution equations associated to our four derivatives.

\vskip 0.5 cm
\centerline{\bf 2. Four derivatives and their properties}
\vskip 0.5 cm

  Let us recall some basic facts about algebraic curves, which
show the significant properties of Picard curves.

{\it Definition}. Suppose that $C$ is an algebraic curve of genus
$g \geq 2$. Let $\{ \omega_1, \cdots, \omega_g \}$ be a basis of
$\Omega^{1}(C)$. Then the map
$$\aligned
  \varphi_{K}: C &\to {\Bbb P}^{g-1}\\
               p &\mapsto [\omega_1(p), \cdots, \omega_g(p)]
\endaligned$$
is called a canonical map over $C$.

  In fact, Suppose that $z$ is a local coordinate near the point $p$
and $\omega_{\alpha}=f_{\alpha}(z) dz$ ($\alpha=1, 2, \cdots, g$).
Then
$$\varphi_{K}(p)=[f_1(z(p)), \cdots, f_{g}(z(p))].$$

{\smc Theorem}. {\it Any algebraic curve of genus $g=2$ is a
hyperelliptic curve.}

{\it Definition}. Suppose that $C$ is a non-hyperelliptic
algebraic curve of genus $g \geq 2$. Let
$$\varphi_{K}: C \to {\Bbb P}^{g-1}$$
be a canonical map. Then
$$\varphi_{K}(C) \subset {\Bbb P}^{g-1}$$
is called a canonical curve.

  Let $C$ be a canonical curve of genus $g$. Then
$\text{deg}(C)=2g-2$.

{\smc Theorem}. {\it Any canonical curve of genus $g=3$ is a
smooth quartic plane algebraic curve.}

{\smc Theorem}. {\it Any canonical curve of genus $g=4$ is the
intersection of a quadratic surface and a cubic surface in ${\Bbb
P}^{3}$.}

  Let us consider the Picard integral
$$z(x, y)=\int_{0}^{1} \frac{dt}{\root 3 \of{t(t-1)(t-x)(t-y)}}.$$
The corresponding Picard curve is given by
$$w^3=t(t-1)(t-x)(t-y).$$
In fact,
$$z(x, y)=-x^{-\frac{1}{3}} y^{-\frac{1}{3}} \int_{0}^{1}
  t^{-\frac{1}{3}} (1-t)^{-\frac{1}{3}} (1-t x^{-1})^{-\frac{1}{3}}
  (1-t y^{-1})^{-\frac{1}{3}} dt.$$

  It is known that the two variable hypergeometric function
(Appell hypergeometric function)
$$F_{1}(a; b, b^{\prime}; c; x, y)
 =\frac{\Gamma(c)}{\Gamma(a) \Gamma(c-a)} \int_{0}^{1}
  t^{a-1} (1-t)^{c-a-1} (1-tx)^{-b} (1-ty)^{-b^{\prime}} dt$$
satisfies the following equations:
$$\left\{\aligned
 x(1-x)r+y(1-x)s+[c-(a+b+1)x]p-byq-abz=0,\\
 y(1-y)t+x(1-y)s+[c-(a+b^{\prime}+1)y]q-b^{\prime}xp-ab^{\prime}z=0,
\endaligned\right.$$
where
$$r=\frac{\partial^2 z}{\partial x^2}, \quad
  s=\frac{\partial^2 z}{\partial x \partial y}, \quad
  t=\frac{\partial^2 z}{\partial y^2}, \quad
  p=\frac{\partial z}{\partial x}, \quad
  q=\frac{\partial z}{\partial y}.$$
Thus,
$$z(x, y)=-\frac{\Gamma\left(\frac{2}{3}\right)^2}{\Gamma\left(\frac{4}{3}
  \right)} x^{-\frac{1}{3}} y^{-\frac{1}{3}} F_{1}\left(\frac{2}{3};
  \frac{1}{3}, \frac{1}{3}; \frac{4}{3}; \frac{1}{x}, \frac{1}{y}\right).$$
Set $u=\frac{1}{x}$ and $v=\frac{1}{y}$. $f(u,
v):=F_{1}\left(\frac{2}{3}; \frac{1}{3}, \frac{1}{3}; \frac{4}{3};
u, v\right)$ satisfies the following equations:
$$\left\{\aligned
 u(1-u) \frac{\partial^2 f}{\partial u^2}+v(1-u) \frac{\partial^2 f}
 {\partial u \partial v}+\left(\frac{4}{3}-2u\right) \frac{\partial f}{\partial u}
 -\frac{1}{3} v \frac{\partial f}{\partial v}-\frac{2}{9} f=0,\\
 v(1-v) \frac{\partial^2 f}{\partial v^2}+u(1-v) \frac{\partial^2 f}
 {\partial u \partial v}+\left(\frac{4}{3}-2v\right) \frac{\partial f}{\partial v}
 -\frac{1}{3} u \frac{\partial f}{\partial u}-\frac{2}{9} f=0.
\endaligned\right.$$
By a straightforward calculation, we find that $z=z(x, y)$
satisfies the following equations:
$$\left\{\aligned
 x(1-x)r+y(1-x)s+\left(1-\frac{5}{3}x\right)p-\frac{1}{3}yq-\frac{1}{9}z=0,\\
 y(1-y)t+x(1-y)s+\left(1-\frac{5}{3}y\right)q-\frac{1}{3}xp-\frac{1}{9}z=0.
\endaligned\right.$$
Let $\omega_1$, $\omega_2$ and $\omega_3$ be three linearly
independent solutions of the above partial differential equations.
Set
$$\frac{\omega_2}{\omega_1}=u(x, y), \quad
  \frac{\omega_3}{\omega_1}=v(x, y).$$
In his paper \cite{Pi1}, Picard proved that
$$2 \text{Re}(v)+|u|^2<0.$$

  In general, according to \cite{AK}, p.54, Th\'{e}or\`{e}me b and p.81, we have
the following Theorem:

{\smc Theorem}. {\it Suppose that there are four functions $z_1$,
$z_2$, $z_3$ and $z_4$ which satisfy the equations:
$$\left\{\aligned
  r &=a_1 p+a_2 q+a_3 z,\\
  s &=b_1 p+b_2 q+b_3 z,\\
  t &=c_1 p+c_2 q+c_3 z,
\endaligned\right.$$
where $a_i=a_i(x, y)$, $b_i=b_i(x, y)$ and $c_i=c_i(x, y)$ $(i=1,
2, 3)$ are functions of $x$ and $y$. There exists a linear
relation of constant coefficients for all these functions such
that
$$C_1 z_1+C_2 z_2+C_3 z_3+C_4 z_4=0.$$
The above partial differential equations can admit three linearly
independent solutions.}

  Now, let us consider the following partial differential equations:
$$\left\{\aligned
  x(1-x) \frac{\partial^2 z}{\partial x^2}+y(1-x) \frac{\partial^2
  z}{\partial x \partial y}+[c-(a+b+1)x] \frac{\partial
  z}{\partial x}-by \frac{\partial z}{\partial y}-abz &=0,\\
  y(1-y) \frac{\partial^2 z}{\partial y^2}+x(1-y) \frac{\partial^2
  z}{\partial x \partial y}+[c-(a+b^{\prime}+1)y] \frac{\partial
  z}{\partial y}-b^{\prime} x \frac{\partial z}{\partial x}-a
  b^{\prime} z &=0,
\endaligned\right.\tag 2.1$$
which are equivalent to the following three equations:
$$\left\{\aligned
 &\frac{\partial^2 z}{\partial x^2}+\left(\frac{c-b^{\prime}}{x}+
  \frac{a+b-c+1}{x-1}+\frac{b^{\prime}}{x-y}\right) \frac{\partial
  z}{\partial x}-\frac{by(y-1)}{x(x-1)(x-y)} \frac{\partial z}{\partial
  y}+\frac{abz}{x(x-1)}=0,\\
 &\frac{\partial^2 z}{\partial y^2}+\left(\frac{c-b}{y}+
  \frac{a+b^{\prime}-c+1}{y-1}+\frac{b}{y-x}\right)
  \frac{\partial z}{\partial y}-\frac{b^{\prime}x(x-1)}{y(y-1)(y-x)}
  \frac{\partial z}{\partial x}+\frac{ab^{\prime}z}{y(y-1)}=0,\\
 &\frac{\partial^2 z}{\partial x \partial y}=\frac{b^{\prime}}{x-y}
  \frac{\partial z}{\partial x}-\frac{b}{x-y} \frac{\partial z}
  {\partial y}.
\endaligned\right.\tag 2.2$$
There are three linearly independent solutions $z_1$, $z_2$ and
$z_3$. Consider the monodromy:
$$\left\{\aligned
  w_1 &=a_1 z_1+a_2 z_2+a_3 z_3,\\
  w_2 &=b_1 z_1+b_2 z_2+b_3 z_3,\\
  w_3 &=c_1 z_1+c_2 z_2+c_3 z_3.
\endaligned\right.$$
We have
$$\left\{\aligned
  \frac{w_1}{w_3} &=\frac{a_1 z_1+a_2 z_2+a_3 z_3}{c_1 z_1+c_2
  z_2+c_3 z_3}=\frac{a_1 \frac{z_1}{z_3}+a_2
  \frac{z_2}{z_3}+a_3}{c_1 \frac{z_1}{z_3}+c_2
  \frac{z_2}{z_3}+c_3},\\
  \frac{w_2}{w_3} &=\frac{b_1 z_1+b_2 z_2+b_3 z_3}{c_1 z_1+c_2
  z_2+c_3 z_3}=\frac{b_1 \frac{z_1}{z_3}+b_2
  \frac{z_2}{z_3}+b_3}{c_1 \frac{z_1}{z_3}+c_2
  \frac{z_2}{z_3}+c_3}.
\endaligned\right.$$
Set $u_1=\frac{z_1}{z_3}$, $u_2=\frac{z_2}{z_3}$ and
$v_1=\frac{w_1}{w_3}$, $v_2=\frac{w_2}{w_3}$, then
$$v_1=\frac{a_1 u_1+a_2 u_2+a_3}{c_1 u_1+c_2 u_2+c_3}, \quad
  v_2=\frac{b_1 u_1+b_2 u_2+b_3}{c_1 u_1+c_2 u_2+c_3}.$$

  Note that
$$\left\{\aligned
  \frac{\partial z_1}{\partial x} z_3-z_1 \frac{\partial z_3}{\partial x}
&=z_3^2 \frac{\partial u_1}{\partial x},\\
  \frac{\partial z_1}{\partial y} z_3-z_1 \frac{\partial z_3}{\partial y}
&=z_3^2 \frac{\partial u_1}{\partial y},\\
  \frac{\partial z_2}{\partial x} z_3-z_2 \frac{\partial z_3}{\partial x}
&=z_3^2 \frac{\partial u_2}{\partial x},\\
  \frac{\partial z_2}{\partial y} z_3-z_2 \frac{\partial z_3}{\partial y}
&=z_3^2 \frac{\partial u_2}{\partial y}.
\endaligned\right.$$
We find that
$$\left\{\aligned
  \frac{\partial^2 z_1}{\partial x^2} z_3-z_1 \frac{\partial^2 z_3}{\partial x^2}
&=z_3^2 \frac{\partial^2 u_1}{\partial x^2}+2 z_3 \frac{\partial
  z_3}{\partial x} \frac{\partial u_1}{\partial x},\\
  \frac{\partial^2 z_1}{\partial y^2} z_3-z_1 \frac{\partial^2 z_3}{\partial y^2}
&=z_3^2 \frac{\partial^2 u_1}{\partial y^2}+2 z_3 \frac{\partial
  z_3}{\partial y} \frac{\partial u_1}{\partial y},\\
  \frac{\partial^2 z_2}{\partial x^2} z_3-z_2 \frac{\partial^2 z_3}{\partial x^2}
&=z_3^2 \frac{\partial^2 u_2}{\partial x^2}+2 z_3 \frac{\partial
  z_3}{\partial x} \frac{\partial u_2}{\partial x},\\
  \frac{\partial^2 z_2}{\partial y^2} z_3-z_2 \frac{\partial^2 z_3}{\partial y^2}
&=z_3^2 \frac{\partial^2 u_2}{\partial y^2}+2 z_3 \frac{\partial
  z_3}{\partial y} \frac{\partial u_2}{\partial y}.
\endaligned\right.$$
$$\left\{\aligned
  \frac{\partial^2 z_1}{\partial x \partial y} z_3-
  z_1 \frac{\partial^2 z_3}{\partial x \partial y}
&=z_3^2 \frac{\partial^2 u_1}{\partial x \partial y}+
  z_3 \frac{\partial z_3}{\partial x} \frac{\partial u_1}{\partial y}+
  z_3 \frac{\partial z_3}{\partial y} \frac{\partial u_1}{\partial x},\\
  \frac{\partial^2 z_2}{\partial x \partial y} z_3-
  z_2 \frac{\partial^2 z_3}{\partial x \partial y}
&=z_3^2 \frac{\partial^2 u_2}{\partial x \partial y}+
  z_3 \frac{\partial z_3}{\partial x} \frac{\partial u_2}{\partial y}+
  z_3 \frac{\partial z_3}{\partial y} \frac{\partial u_2}{\partial x}.
\endaligned\right.$$
Substituting these expressions in the above three partial
differential equations, we get the next six equations:
$$\left\{\aligned
  \frac{\partial^2 u_1}{\partial x^2}+\left(2 \frac{\frac{\partial
  z_3}{\partial x}}{z_3}+\frac{c-b^{\prime}}{x}+\frac{a+b-c+1}
  {x-1}+\frac{b^{\prime}}{x-y}\right) \frac{\partial u_1}{\partial x}-
  \frac{by(y-1)}{x(x-1)(x-y)} \frac{\partial u_1}{\partial y} &=0,\\
  \frac{\partial^2 u_2}{\partial x^2}+\left(2 \frac{\frac{\partial
  z_3}{\partial x}}{z_3}+\frac{c-b^{\prime}}{x}+\frac{a+b-c+1}
  {x-1}+\frac{b^{\prime}}{x-y}\right) \frac{\partial u_2}{\partial x}-
  \frac{by(y-1)}{x(x-1)(x-y)} \frac{\partial u_2}{\partial y} &=0,\\
  \frac{\partial^2 u_1}{\partial y^2}+\left(2 \frac{\frac{\partial
  z_3}{\partial y}}{z_3}+\frac{c-b}{y}+\frac{a+b^{\prime}-c+1}
  {y-1}+\frac{b}{y-x}\right) \frac{\partial u_1}{\partial y}-
  \frac{b^{\prime} x(x-1)}{y(y-1)(y-x)} \frac{\partial u_1}{\partial x} &=0,\\
  \frac{\partial^2 u_2}{\partial y^2}+\left(2 \frac{\frac{\partial
  z_3}{\partial y}}{z_3}+\frac{c-b}{y}+\frac{a+b^{\prime}-c+1}
  {y-1}+\frac{b}{y-x}\right) \frac{\partial u_2}{\partial y}-
  \frac{b^{\prime} x(x-1)}{y(y-1)(y-x)} \frac{\partial u_2}{\partial x} &=0.
\endaligned\right.$$
$$\left\{\aligned
  \frac{\partial^2 u_1}{\partial x \partial y}+\left(\frac{\frac{\partial
  z_3}{\partial y}}{z_3}-\frac{b^{\prime}}{x-y}\right) \frac{\partial u_1}
  {\partial x}+\left(\frac{\frac{\partial z_3}{\partial x}}{z_3}+\frac{b}
  {x-y}\right) \frac{\partial u_1}{\partial y} &=0,\\
  \frac{\partial^2 u_2}{\partial x \partial y}+\left(\frac{\frac{\partial
  z_3}{\partial y}}{z_3}-\frac{b^{\prime}}{x-y}\right) \frac{\partial u_2}
  {\partial x}+\left(\frac{\frac{\partial z_3}{\partial x}}{z_3}+\frac{b}
  {x-y}\right) \frac{\partial u_2}{\partial y} &=0.
\endaligned\right.$$

  It is known that the Schwarzian derivative
$$\{w; x\}:=\frac{w^{\prime \prime \prime}}{w^{\prime}}-
  \frac{3}{2} \left(\frac{w^{\prime \prime}}{w^{\prime}}\right)^{2}$$
satisfies that
$$\left\{ \frac{aw+b}{cw+d}; x \right\}=\{w; x\} \quad \text{for}
  \quad \left(\matrix a & b\\ c & d \endmatrix\right) \in GL(2,
  {\Bbb C}).$$
Now, we define the following four derivatives:

{\it Definition} 2.1. Four derivatives associated to the group
$GL(3)$ are given by:
$$\{u_1, u_2; x, y\}_{x}
:=\frac{\vmatrix \format \c \quad & \c \\
  \frac{\partial u_1}{\partial x} &
  \frac{\partial u_2}{\partial x}\\
  \frac{\partial^2 u_1}{\partial x^2} &
  \frac{\partial^2 u_2}{\partial x^2}\endvmatrix}
 {\vmatrix \format \c \quad & \c\\
  \frac{\partial u_1}{\partial x} &
  \frac{\partial u_2}{\partial x}\\
  \frac{\partial u_1}{\partial y} &
  \frac{\partial u_2}{\partial y}
  \endvmatrix}, \quad
  \{u_1, u_2; x, y\}_{y}
:=\frac{\vmatrix \format \c \quad & \c \\
  \frac{\partial u_1}{\partial y} &
  \frac{\partial u_2}{\partial y}\\
  \frac{\partial^2 u_1}{\partial y^2} &
  \frac{\partial^2 u_2}{\partial y^2}\endvmatrix}
 {\vmatrix \format \c \quad & \c\\
  \frac{\partial u_1}{\partial y} &
  \frac{\partial u_2}{\partial y}\\
  \frac{\partial u_1}{\partial x} &
  \frac{\partial u_2}{\partial x}
  \endvmatrix},\tag 2.3$$
$$[u_1, u_2; x, y]_{x}
:=\frac{\vmatrix \format \c \quad & \c \\
  \frac{\partial u_1}{\partial y} &
  \frac{\partial u_2}{\partial y}\\
  \frac{\partial^2 u_1}{\partial x^2} &
  \frac{\partial^2 u_2}{\partial x^2}\endvmatrix+2
  \vmatrix \format \c \quad & \c \\
  \frac{\partial u_1}{\partial x} &
  \frac{\partial u_2}{\partial x}\\
  \frac{\partial^2 u_1}{\partial x \partial y} &
  \frac{\partial^2 u_2}{\partial x \partial y}\endvmatrix}
 {\vmatrix \format \c \quad & \c\\
  \frac{\partial u_1}{\partial x} &
  \frac{\partial u_2}{\partial x}\\
  \frac{\partial u_1}{\partial y} &
  \frac{\partial u_2}{\partial y}
  \endvmatrix},\tag 2.4$$
$$[u_1, u_2; x, y]_{y}
:=\frac{\vmatrix \format \c \quad & \c \\
  \frac{\partial u_1}{\partial x} &
  \frac{\partial u_2}{\partial x}\\
  \frac{\partial^2 u_1}{\partial y^2} &
  \frac{\partial^2 u_2}{\partial y^2}\endvmatrix+2
  \vmatrix \format \c \quad & \c \\
  \frac{\partial u_1}{\partial y} &
  \frac{\partial u_2}{\partial y}\\
  \frac{\partial^2 u_1}{\partial x \partial y} &
  \frac{\partial^2 u_2}{\partial x \partial y}\endvmatrix}
 {\vmatrix \format \c \quad & \c\\
  \frac{\partial u_1}{\partial y} &
  \frac{\partial u_2}{\partial y}\\
  \frac{\partial u_1}{\partial x} &
  \frac{\partial u_2}{\partial x}
  \endvmatrix},\tag 2.5$$
where $x$ and $y$ are complex numbers, $u_1$ and $u_2$ are
complex-valued functions of $x$ and $y$.

  The above six equations give that
$$\left\{\aligned
  \{u_1, u_2; x, y\}_{x} &=\frac{by(y-1)}{x(x-1)(x-y)},\\
  \{u_1, u_2; x, y\}_{y} &=\frac{b^{\prime} x(x-1)}{y(y-1)(y-x)},\\
  [u_1, u_2; x, y]_{x} &=\frac{c-b^{\prime}}{x}+\frac{a+b-c+1}{x-1}
                        +\frac{b^{\prime}-2b}{x-y},\\
  [u_1, u_2; x, y]_{y} &=\frac{c-b}{y}+\frac{a+b^{\prime}-c+1}{y-1}
                        +\frac{b-2b^{\prime}}{y-x}.
\endaligned\right.\tag 2.6$$

  In general, we have the following result:

{\smc Proposition 2.2}. {\it For the system of partial
differential equations:
$$\left\{\aligned
  r &=a_1 p+a_2 q+a_3 z,\\
  s &=b_1 p+b_2 q+b_3 z,\\
  t &=c_1 p+c_2 q+c_3 z,
\endaligned\right.\tag 2.7$$
where $a_i=a_i(x, y)$, $b_i=b_i(x, y)$, $c_i=c_i(x, y)$ $(i=1, 2,
3)$ and $r=\frac{\partial^2 z}{\partial x^2}$, $s=\frac{\partial^2
z}{\partial x \partial y}$, $t=\frac{\partial^2 z}{\partial y^2}$,
$p=\frac{\partial z}{\partial x}$, $q=\frac{\partial z}{\partial
y}$, there exists three linearly independent solutions: $z_1$,
$z_2$ and $z_3$. Put $u_1=z_1/z_3$ and $u_2=z_2/z_3$. Then
$$\left\{\aligned
  \{ u_1, u_2; x, y \}_{x} &=a_2,\\
  \{ u_1, u_2; x, y \}_{y} &=c_1,\\
      [u_1, u_2; x, y]_{x} &=2 b_2-a_1,\\
      [u_1, u_2; x, y]_{y} &=2 b_1-c_2.
\endaligned\right.\tag 2.8$$}

{\smc Proposition 2.3}. {\it Our derivatives are invariant under
linear fractional transformation acting on the first arguments
$$\left\{\aligned
  \{v_1, v_2; x, y\}_{x} &=\{u_1, u_2; x, y\}_{x},\\
  \{v_1, v_2; x, y\}_{y} &=\{u_1, u_2; x, y\}_{y},\\
  [v_1, v_2; x, y]_{x} &=[u_1, u_2; x, y]_{x},\\
  [v_1, v_2; x, y]_{y} &=[u_1, u_2; x, y]_{y},
\endaligned\right.\tag 2.9$$
where
$$(v_1, v_2)=\gamma(u_1, u_2)=\left(\frac{a_1 u_1+a_2 u_2+a_3}
  {c_1 u_1+c_2 u_2+c_3}, \frac{b_1 u_1+b_2 u_2+b_3}{c_1 u_1+c_2
  u_2+c_3}\right)\tag 2.10$$
for $\gamma=\left(\matrix a_1 & a_2 & a_3\\ b_1 & b_2 & b_3\\
c_1 & c_2 & c_3 \endmatrix\right) \in GL(3, {\Bbb C})$.}

{\it Proof}. We have
$$\vmatrix \format \c \quad & \c \\
  \frac{\partial v_1}{\partial x} & \frac{\partial v_2}{\partial x}\\
  \frac{\partial v_1}{\partial y} & \frac{\partial v_2}{\partial y}
  \endvmatrix
 =\Delta (c_1 u_1+c_2 u_2+c_3)^{-3}
  \vmatrix \format \c \quad & \c \\
  \frac{\partial u_1}{\partial x} & \frac{\partial u_2}{\partial x}\\
  \frac{\partial u_1}{\partial y} & \frac{\partial u_2}{\partial y}
  \endvmatrix,$$
where
$$\Delta=\vmatrix \format \c \quad & \c \quad & \c \\
  a_1 & a_2 & a_3\\ b_1 & b_2 & b_3\\ c_1 & c_2 & c_3
  \endvmatrix.$$
Note that
$$\vmatrix \format \c \quad & \c \\ a & b\\ c+c^{\prime} &
  d+d^{\prime} \endvmatrix=\vmatrix \format \c \quad & \c \\
  a & b\\ c & d \endvmatrix+\vmatrix \format \c \quad & \c \\
  a & b\\ c^{\prime} & d^{\prime} \endvmatrix.$$
A straightforward calculation gives that
$$\vmatrix \format \c \quad & \c \\
  \frac{\partial v_1}{\partial x} &
  \frac{\partial v_2}{\partial x}\\
  \frac{\partial^2 v_1}{\partial x^2} &
  \frac{\partial^2 v_2}{\partial x^2}
  \endvmatrix
 =\Delta (c_1 u_1+c_2 u_2+c_3)^{-3}
  \vmatrix \format \c \quad & \c \\
  \frac{\partial u_1}{\partial x} &
  \frac{\partial u_2}{\partial x}\\
  \frac{\partial^2 u_1}{\partial x^2} &
  \frac{\partial^2 u_2}{\partial x^2}
  \endvmatrix,$$
$$\aligned
  \vmatrix \format \c \quad & \c \\
  \frac{\partial v_1}{\partial y} &
  \frac{\partial v_2}{\partial y}\\
  \frac{\partial^2 v_1}{\partial x^2} &
  \frac{\partial^2 v_2}{\partial x^2}
  \endvmatrix
=&2 \Delta (c_1 u_1+c_2 u_2+c_3)^{-4}
  (c_1 \frac{\partial u_1}{\partial x}
  +c_2 \frac{\partial u_2}{\partial x})
  \vmatrix \format \c \quad & \c \\
  \frac{\partial u_1}{\partial x} &
  \frac{\partial u_2}{\partial x}\\
  \frac{\partial u_1}{\partial y} &
  \frac{\partial u_2}{\partial y}
  \endvmatrix+\\
+&\Delta (c_1 u_1+c_2 u_2+c_3)^{-3}
  \vmatrix \format \c \quad & \c \\
  \frac{\partial u_1}{\partial y} &
  \frac{\partial u_2}{\partial y}\\
  \frac{\partial^2 u_1}{\partial x^2} &
  \frac{\partial^2 u_2}{\partial x^2}
  \endvmatrix.
\endaligned$$
$$\aligned
  \vmatrix \format \c \quad & \c \\
  \frac{\partial v_1}{\partial x} &
  \frac{\partial v_2}{\partial x}\\
  \frac{\partial^2 v_1}{\partial x \partial y} &
  \frac{\partial^2 v_2}{\partial x \partial y}
  \endvmatrix
=&-\Delta (c_1 u_1+c_2 u_2+c_3)^{-4}
  (c_1 \frac{\partial u_1}{\partial x}
  +c_2 \frac{\partial u_2}{\partial x})
  \vmatrix \format \c \quad & \c \\
  \frac{\partial u_1}{\partial x} &
  \frac{\partial u_2}{\partial x}\\
  \frac{\partial u_1}{\partial y} &
  \frac{\partial u_2}{\partial y}
  \endvmatrix+\\
+&\Delta (c_1 u_1+c_2 u_2+c_3)^{-3}
  \vmatrix \format \c \quad & \c \\
  \frac{\partial u_1}{\partial x} &
  \frac{\partial u_2}{\partial x}\\
  \frac{\partial^2 u_1}{\partial x \partial y} &
  \frac{\partial^2 u_2}{\partial x \partial y}
  \endvmatrix.
\endaligned$$
Hence,
$$\{v_1, v_2; x, y\}_{x}=\{u_1, u_2; x, y\}_{x}, \quad
    [v_1, v_2; x, y]_{x}=[u_1, u_2; x, y]_{x}.$$
Exchange the position of $x$ and $y$, we get the other identities.
\flushpar
$\qquad \qquad \qquad \qquad \qquad \qquad \qquad \qquad \qquad
 \qquad \qquad \qquad \qquad \qquad \qquad \qquad \qquad \qquad
 \quad \boxed{}$

  Thus our derivatives are differential invariants of the general
linear group acting as a linear fractional transformation.

{\smc Corollary 2.4}. {\it Our derivatives vanish when $(w_1,
w_2)$ are linear fractional functions of $(z_1, z_2)$:
$$\left\{\aligned
  \{w_1, w_2; z_1, z_2\}_{z_1} &=0,\\
  \{w_1, w_2; z_1, z_2\}_{z_2} &=0,\\
    [w_1, w_2; z_1, z_2]_{z_1} &=0,\\
    [w_1, w_2; z_1, z_2]_{z_2} &=0,
\endaligned\right.\tag 2.11$$
where
$$(w_1, w_2)=\gamma(z_1, z_2)=
  \left(\frac{a_1 z_1+a_2 z_2+a_3}{c_1 z_1+c_2 z_2+c_3},
  \frac{b_1 z_1+b_2 z_2+b_3}{c_1 z_1+c_2 z_2+c_3}\right),\tag 2.12$$
for $\gamma=\left(\matrix a_1 & a_2 & a_3\\ b_1 & b_2 & b_3\\
c_1 & c_2 & c_3 \endmatrix\right) \in GL(3, {\Bbb C})$.}

{\it Proof}. Note that
$$\left\{\aligned
  \{z_1, z_2; z_1, z_2\}_{z_1} &=0,\\
  \{z_1, z_2; z_1, z_2\}_{z_2} &=0,\\
    [z_1, z_2; z_1, z_2]_{z_1} &=0,\\
    [z_1, z_2; z_1, z_2]_{z_2} &=0.
\endaligned\right.$$
The above Proposition gives the desired result.
\flushpar
$\qquad \qquad \qquad \qquad \qquad \qquad \qquad \qquad \qquad
 \qquad \qquad \qquad \qquad \qquad \qquad \qquad \qquad \qquad
 \quad \boxed{}$

{\smc Proposition 2.5 (Main Proposition 1)}. {\it The change of
independent variables obeys the laws:
$$\aligned
 &\{u_1, u_2; x, y\}_{x}\\
=&\frac{(\frac{\partial w_1}{\partial x})^3}{\frac{\partial(w_1,
  w_2)}{\partial(x, y)}} \{u_1, u_2; w_1, w_2\}_{w_1}+
  \frac{(\frac{\partial w_1}{\partial x})^2 \frac{\partial w_2}
  {\partial x}}{\frac{\partial(w_1, w_2)}{\partial(x, y)}}
  [u_1, u_2; w_1, w_2]_{w_1}+\\
+&\frac{(\frac{\partial w_2}{\partial x})^3}{\frac{\partial(w_2,
  w_1)}{\partial(x, y)}} \{u_1, u_2; w_1, w_2\}_{w_2}+
  \frac{\frac{\partial w_1}{\partial x} (\frac{\partial w_2}
  {\partial x})^2}{\frac{\partial(w_2, w_1)}{\partial(x, y)}}
  [u_1, u_2; w_1, w_2]_{w_2}+\{w_1, w_2; x, y\}_{x}.
\endaligned\tag 2.13$$
$$\aligned
 &\{u_1, u_2; x, y\}_{y}\\
=&\frac{(\frac{\partial w_1}{\partial y})^3}{\frac{\partial(w_2,
  w_1)}{\partial(x, y)}} \{u_1, u_2; w_1, w_2\}_{w_1}+
  \frac{(\frac{\partial w_1}{\partial y})^2 \frac{\partial w_2}
  {\partial y}}{\frac{\partial(w_2, w_1)}{\partial(x, y)}}
  [u_1, u_2; w_1, w_2]_{w_1}+\\
+&\frac{(\frac{\partial w_2}{\partial y})^3}{\frac{\partial(w_1,
  w_2)}{\partial(x, y)}} \{u_1, u_2; w_1, w_2\}_{w_2}+
  \frac{\frac{\partial w_1}{\partial y} (\frac{\partial w_2}
  {\partial y})^2}{\frac{\partial(w_1, w_2)}{\partial(x, y)}}
  [u_1, u_2; w_1, w_2]_{w_2}+\{w_1, w_2; x, y\}_{y}.
\endaligned\tag 2.14$$
$$\aligned
 &[u_1, u_2; x, y]_{x}\\
=&\frac{3 (\frac{\partial w_1}{\partial x})^2 \frac{\partial
  w_1}{\partial y}}{\frac{\partial(w_1, w_2)}{\partial(x, y)}}
  \{u_1, u_2; w_1, w_2\}_{w_1}+\frac{(\frac{\partial w_1}{\partial
  x})^2 \frac{\partial w_2}{\partial y}+2 \frac{\partial w_1}{\partial x}
  \frac{\partial w_2}{\partial x} \frac{\partial w_1}{\partial y}}
  {\frac{\partial(w_1, w_2)}{\partial(x, y)}} [u_1, u_2; w_1,
  w_2]_{w_1}+\\
+&\frac{3 (\frac{\partial w_2}{\partial x})^2 \frac{\partial
  w_2}{\partial y}}{\frac{\partial(w_2, w_1)}{\partial(x, y)}}
  \{u_1, u_2; w_1, w_2\}_{w_2}+\frac{(\frac{\partial w_2}{\partial
  x})^2 \frac{\partial w_1}{\partial y}+2 \frac{\partial w_1}{\partial x}
  \frac{\partial w_2}{\partial x} \frac{\partial w_2}{\partial y}}
  {\frac{\partial(w_2, w_1)}{\partial(x, y)}} [u_1, u_2; w_1,
  w_2]_{w_2}+\\
+&[w_1, w_2; x, y]_{x}.
\endaligned\tag 2.15$$
$$\aligned
 &[u_1, u_2; x, y]_{y}\\
=&\frac{3 \frac{\partial w_1}{\partial x} (\frac{\partial
  w_1}{\partial y})^2}{\frac{\partial(w_2, w_1)}{\partial(x, y)}}
  \{u_1, u_2; w_1, w_2\}_{w_1}+\frac{\frac{\partial w_2}{\partial
  x} (\frac{\partial w_1}{\partial y})^2+2 \frac{\partial w_1}{\partial x}
  \frac{\partial w_1}{\partial y} \frac{\partial w_2}{\partial y}}
  {\frac{\partial(w_2, w_1)}{\partial(x, y)}} [u_1, u_2; w_1,
  w_2]_{w_1}+\\
+&\frac{3 \frac{\partial w_2}{\partial x} (\frac{\partial
  w_2}{\partial y})^2}{\frac{\partial(w_1, w_2)}{\partial(x, y)}}
  \{u_1, u_2; w_1, w_2\}_{w_2}+\frac{\frac{\partial w_1}{\partial
  x} (\frac{\partial w_2}{\partial y})^2+2 \frac{\partial w_2}{\partial x}
  \frac{\partial w_1}{\partial y} \frac{\partial w_2}{\partial y}}
  {\frac{\partial(w_1, w_2)}{\partial(x, y)}} [u_1, u_2; w_1,
  w_2]_{w_2}+\\
+&[w_1, w_2; x, y]_{y}.
\endaligned\tag 2.16$$}

{\it Proof}. We have
$$\frac{\partial u_1}{\partial x}
 =\frac{\partial u_1}{\partial w_1} \frac{\partial w_1}{\partial x}
 +\frac{\partial u_1}{\partial w_2} \frac{\partial w_2}{\partial x},
  \quad
  \frac{\partial u_2}{\partial x}
 =\frac{\partial u_2}{\partial w_1} \frac{\partial w_1}{\partial x}
 +\frac{\partial u_2}{\partial w_2} \frac{\partial w_2}{\partial x}.$$
$$\frac{\partial u_1}{\partial y}
 =\frac{\partial u_1}{\partial w_1} \frac{\partial w_1}{\partial y}
 +\frac{\partial u_1}{\partial w_2} \frac{\partial w_2}{\partial y},
  \quad
  \frac{\partial u_2}{\partial y}
 =\frac{\partial u_2}{\partial w_1} \frac{\partial w_1}{\partial y}
 +\frac{\partial u_2}{\partial w_2} \frac{\partial w_2}{\partial y}.$$
$$\frac{\partial^2 u_1}{\partial x^2}
 =\frac{\partial^2 u_1}{\partial w_1^2} \left(\frac{\partial w_1}
  {\partial x}\right)^2+2 \frac{\partial^2 u_1}{\partial w_1 \partial
  w_2} \frac{\partial w_1}{\partial x} \frac{\partial w_2}{\partial x}
  +\frac{\partial^2 u_1}{\partial w_2^2} \left(\frac{\partial w_2}
  {\partial x}\right)^2+\frac{\partial u_1}{\partial w_1}
  \frac{\partial^2 w_1}{\partial x^2}+\frac{\partial u_1}{\partial
  w_2} \frac{\partial^2 w_2}{\partial x^2}.$$
$$\frac{\partial^2 u_2}{\partial x^2}
 =\frac{\partial^2 u_2}{\partial w_1^2} \left(\frac{\partial w_1}
  {\partial x}\right)^2+2 \frac{\partial^2 u_2}{\partial w_1 \partial
  w_2} \frac{\partial w_1}{\partial x} \frac{\partial w_2}{\partial x}
  +\frac{\partial^2 u_2}{\partial w_2^2} \left(\frac{\partial w_2}
  {\partial x}\right)^2+\frac{\partial u_2}{\partial w_1}
  \frac{\partial^2 w_1}{\partial x^2}+\frac{\partial u_2}{\partial
  w_2} \frac{\partial^2 w_2}{\partial x^2}.$$
$$\frac{\partial^2 u_1}{\partial y^2}
 =\frac{\partial^2 u_1}{\partial w_1^2} \left(\frac{\partial w_1}
  {\partial y}\right)^2+2 \frac{\partial^2 u_1}{\partial w_1 \partial
  w_2} \frac{\partial w_1}{\partial y} \frac{\partial w_2}{\partial y}
  +\frac{\partial^2 u_1}{\partial w_2^2} \left(\frac{\partial w_2}
  {\partial y}\right)^2+\frac{\partial u_1}{\partial w_1}
  \frac{\partial^2 w_1}{\partial y^2}+\frac{\partial u_1}{\partial
  w_2} \frac{\partial^2 w_2}{\partial y^2}.$$
$$\frac{\partial^2 u_2}{\partial y^2}
 =\frac{\partial^2 u_2}{\partial w_1^2} \left(\frac{\partial w_1}
  {\partial y}\right)^2+2 \frac{\partial^2 u_2}{\partial w_1 \partial
  w_2} \frac{\partial w_1}{\partial y} \frac{\partial w_2}{\partial y}
  +\frac{\partial^2 u_2}{\partial w_2^2} \left(\frac{\partial w_2}
  {\partial y}\right)^2+\frac{\partial u_2}{\partial w_1}
  \frac{\partial^2 w_1}{\partial y^2}+\frac{\partial u_2}{\partial
  w_2} \frac{\partial^2 w_2}{\partial y^2}.$$
$$\aligned
  \frac{\partial^2 u_1}{\partial x \partial y}
=&\frac{\partial^2 u_1}{\partial w_1^2} \frac{\partial
  w_1}{\partial x} \frac{\partial w_1}{\partial y}+\frac{\partial^2
  u_1}{\partial w_1 \partial w_2} \left(\frac{\partial w_1}{\partial x}
  \frac{\partial w_2}{\partial y}+\frac{\partial w_2}{\partial x}
  \frac{\partial w_1}{\partial y}\right)+\frac{\partial^2 u_1}{\partial w_2^2}
  \frac{\partial w_2}{\partial x} \frac{\partial w_2}{\partial y}+\\
+&\frac{\partial u_1}{\partial w_1} \frac{\partial^2 w_1}{\partial
  x \partial y}+\frac{\partial u_1}{\partial w_2} \frac{\partial^2
  w_2}{\partial x \partial y}.
\endaligned$$
$$\aligned
  \frac{\partial^2 u_2}{\partial x \partial y}
=&\frac{\partial^2 u_2}{\partial w_1^2} \frac{\partial
  w_1}{\partial x} \frac{\partial w_1}{\partial y}+\frac{\partial^2
  u_2}{\partial w_1 \partial w_2} \left(\frac{\partial w_1}{\partial x}
  \frac{\partial w_2}{\partial y}+\frac{\partial w_2}{\partial x}
  \frac{\partial w_1}{\partial y}\right)+\frac{\partial^2 u_2}{\partial w_2^2}
  \frac{\partial w_2}{\partial x} \frac{\partial w_2}{\partial y}+\\
+&\frac{\partial u_2}{\partial w_1} \frac{\partial^2 w_1}{\partial
  x \partial y}+\frac{\partial u_2}{\partial w_2} \frac{\partial^2
  w_2}{\partial x \partial y}.
\endaligned$$
Hence,
$$\aligned
 &\vmatrix \format \c \quad & \c \\
  \frac{\partial u_1}{\partial x} &
  \frac{\partial u_2}{\partial x}\\
  \frac{\partial^2 u_1}{\partial x^2} &
  \frac{\partial^2 u_2}{\partial x^2}
  \endvmatrix\\
=&\left(\frac{\partial w_1}{\partial x}\right)^3
  \vmatrix \format \c \quad & \c \\
  \frac{\partial u_1}{\partial w_1} &
  \frac{\partial u_2}{\partial w_1}\\
  \frac{\partial^2 u_1}{\partial w_1^2} &
  \frac{\partial^2 u_2}{\partial w_1^2}
  \endvmatrix
 +\left(\frac{\partial w_1}{\partial x}\right)^2 \frac{\partial
  w_2}{\partial x} \left[
  \vmatrix \format \c \quad & \c \\
  \frac{\partial u_1}{\partial w_2} &
  \frac{\partial u_2}{\partial w_2}\\
  \frac{\partial^2 u_1}{\partial w_1^2} &
  \frac{\partial^2 u_2}{\partial w_1^2}
  \endvmatrix+2
  \vmatrix \format \c \quad & \c \\
  \frac{\partial u_1}{\partial w_1} &
  \frac{\partial u_2}{\partial w_1}\\
  \frac{\partial^2 u_1}{\partial w_1 \partial w_2} &
  \frac{\partial^2 u_2}{\partial w_1 \partial w_2}
  \endvmatrix\right]+\\
+&\left(\frac{\partial w_2}{\partial x}\right)^3
  \vmatrix \format \c \quad & \c \\
  \frac{\partial u_1}{\partial w_2} &
  \frac{\partial u_2}{\partial w_2}\\
  \frac{\partial^2 u_1}{\partial w_2^2} &
  \frac{\partial^2 u_2}{\partial w_2^2}
  \endvmatrix
 +\frac{\partial w_1}{\partial x} \left(\frac{\partial
  w_2}{\partial x}\right)^2 \left[
  \vmatrix \format \c \quad & \c \\
  \frac{\partial u_1}{\partial w_1} &
  \frac{\partial u_2}{\partial w_1}\\
  \frac{\partial^2 u_1}{\partial w_2^2} &
  \frac{\partial^2 u_2}{\partial w_2^2}
  \endvmatrix+2
  \vmatrix \format \c \quad & \c \\
  \frac{\partial u_1}{\partial w_2} &
  \frac{\partial u_2}{\partial w_2}\\
  \frac{\partial^2 u_1}{\partial w_1 \partial w_2} &
  \frac{\partial^2 u_2}{\partial w_1 \partial w_2}
  \endvmatrix\right]+\\
+&\vmatrix \format \c \quad & \c \\
  \frac{\partial u_1}{\partial w_1} &
  \frac{\partial u_2}{\partial w_1}\\
  \frac{\partial u_1}{\partial w_2} &
  \frac{\partial u_2}{\partial w_2}
  \endvmatrix
  \vmatrix \format \c \quad & \c \\
  \frac{\partial w_1}{\partial x} &
  \frac{\partial w_2}{\partial x}\\
  \frac{\partial^2 w_1}{\partial x^2} &
  \frac{\partial^2 w_2}{\partial x^2}
  \endvmatrix.
\endaligned$$
Note that
$$\frac{\partial(u_1, u_2)}{\partial(x, y)}
 =\frac{\partial(u_1, u_2)}{\partial(w_1, w_2)}
  \frac{\partial(w_1, w_2)}{\partial(x, y)}.$$
We find that
$$\aligned
 &\{u_1, u_2; x, y\}_{x}\\
=&\frac{(\frac{\partial w_1}{\partial x})^3}{\frac{\partial(w_1,
  w_2)}{\partial(x, y)}} \{u_1, u_2; w_1, w_2\}_{w_1}+
  \frac{(\frac{\partial w_1}{\partial x})^2 \frac{\partial w_2}
  {\partial x}}{\frac{\partial(w_1, w_2)}{\partial(x, y)}}
  [u_1, u_2; w_1, w_2]_{w_1}+\\
+&\frac{(\frac{\partial w_2}{\partial x})^3}{\frac{\partial(w_2,
  w_1)}{\partial(x, y)}} \{u_1, u_2; w_1, w_2\}_{w_2}+
  \frac{\frac{\partial w_1}{\partial x} (\frac{\partial w_2}
  {\partial x})^2}{\frac{\partial(w_2, w_1)}{\partial(x, y)}}
  [u_1, u_2; w_1, w_2]_{w_2}+\{w_1, w_2; x, y\}_{x}.
\endaligned$$

  On the other hand,
$$\aligned
 &\vmatrix \format \c \quad & \c \\
  \frac{\partial u_1}{\partial y} &
  \frac{\partial u_2}{\partial y}\\
  \frac{\partial^2 u_1}{\partial x^2} &
  \frac{\partial^2 u_2}{\partial x^2}
  \endvmatrix\\
=&\left(\frac{\partial w_1}{\partial x}\right)^2 \frac{\partial
  w_1}{\partial y}
  \vmatrix \format \c \quad & \c \\
  \frac{\partial u_1}{\partial w_1} &
  \frac{\partial u_2}{\partial w_1}\\
  \frac{\partial^2 u_1}{\partial w_1^2} &
  \frac{\partial^2 u_2}{\partial w_1^2}
  \endvmatrix+
  \left(\frac{\partial w_1}{\partial x}\right)^2 \frac{\partial
  w_2}{\partial y}
  \vmatrix \format \c \quad & \c \\
  \frac{\partial u_1}{\partial w_2} &
  \frac{\partial u_2}{\partial w_2}\\
  \frac{\partial^2 u_1}{\partial w_1^2} &
  \frac{\partial^2 u_2}{\partial w_1^2}
  \endvmatrix+\\
+&2 \frac{\partial w_1}{\partial x} \frac{\partial w_2}{\partial
  x} \frac{\partial w_1}{\partial y}
  \vmatrix \format \c \quad & \c \\
  \frac{\partial u_1}{\partial w_1} &
  \frac{\partial u_2}{\partial w_1}\\
  \frac{\partial^2 u_1}{\partial w_1 \partial w_2} &
  \frac{\partial^2 u_2}{\partial w_1 \partial w_2}
  \endvmatrix+
  2 \frac{\partial w_1}{\partial x} \frac{\partial w_2}{\partial
  x} \frac{\partial w_2}{\partial y}
  \vmatrix \format \c \quad & \c \\
  \frac{\partial u_1}{\partial w_2} &
  \frac{\partial u_2}{\partial w_2}\\
  \frac{\partial^2 u_1}{\partial w_1 \partial w_2} &
  \frac{\partial^2 u_2}{\partial w_1 \partial w_2}
  \endvmatrix+\\
+&\left(\frac{\partial w_2}{\partial x}\right)^2 \frac{\partial
  w_1}{\partial y}
  \vmatrix \format \c \quad & \c \\
  \frac{\partial u_1}{\partial w_1} &
  \frac{\partial u_2}{\partial w_1}\\
  \frac{\partial^2 u_1}{\partial w_2^2} &
  \frac{\partial^2 u_2}{\partial w_2^2}
  \endvmatrix+
  \left(\frac{\partial w_2}{\partial x}\right)^2 \frac{\partial
  w_2}{\partial y}
  \vmatrix \format \c \quad & \c \\
  \frac{\partial u_1}{\partial w_2} &
  \frac{\partial u_2}{\partial w_2}\\
  \frac{\partial^2 u_1}{\partial w_2^2} &
  \frac{\partial^2 u_2}{\partial w_2^2}
  \endvmatrix+\\
&+\vmatrix \format \c \quad & \c \\
  \frac{\partial w_1}{\partial y} &
  \frac{\partial w_2}{\partial y}\\
  \frac{\partial^2 w_1}{\partial x^2} &
  \frac{\partial^2 w_2}{\partial x^2}
  \endvmatrix
  \vmatrix \format \c \quad & \c \\
  \frac{\partial u_1}{\partial w_1} &
  \frac{\partial u_2}{\partial w_1}\\
  \frac{\partial u_1}{\partial w_2} &
  \frac{\partial u_2}{\partial w_2}
  \endvmatrix.
\endaligned$$
$$\aligned
 &\vmatrix \format \c \quad & \c \\
  \frac{\partial u_1}{\partial x} &
  \frac{\partial u_2}{\partial x}\\
  \frac{\partial^2 u_1}{\partial x \partial y} &
  \frac{\partial^2 u_2}{\partial x \partial y}
  \endvmatrix\\
=&\left(\frac{\partial w_1}{\partial x}\right)^2 \frac{\partial
  w_1}{\partial y}
  \vmatrix \format \c \quad & \c \\
  \frac{\partial u_1}{\partial w_1} &
  \frac{\partial u_2}{\partial w_1}\\
  \frac{\partial^2 u_1}{\partial w_1^2} &
  \frac{\partial^2 u_2}{\partial w_1^2}
  \endvmatrix+
  \frac{\partial w_1}{\partial x} \frac{\partial w_1}{\partial y}
  \frac{\partial w_2}{\partial x}
  \vmatrix \format \c \quad & \c \\
  \frac{\partial u_1}{\partial w_2} &
  \frac{\partial u_2}{\partial w_2}\\
  \frac{\partial^2 u_1}{\partial w_1^2} &
  \frac{\partial^2 u_2}{\partial w_1^2}
  \endvmatrix+\\
+&\frac{\partial w_1}{\partial x} \left(\frac{\partial w_1}
  {\partial x} \frac{\partial w_2}{\partial y}+\frac{\partial
  w_2}{\partial x} \frac{\partial w_1}{\partial y}\right)
  \vmatrix \format \c \quad & \c \\
  \frac{\partial u_1}{\partial w_1} &
  \frac{\partial u_2}{\partial w_1}\\
  \frac{\partial^2 u_1}{\partial w_1 \partial w_2} &
  \frac{\partial^2 u_2}{\partial w_1 \partial w_2}
  \endvmatrix+\\
+&\frac{\partial w_2}{\partial x} \left(\frac{\partial w_1}
  {\partial x} \frac{\partial w_2}{\partial y}+\frac{\partial
  w_2}{\partial x} \frac{\partial w_1}{\partial y}\right)
  \vmatrix \format \c \quad & \c \\
  \frac{\partial u_1}{\partial w_2} &
  \frac{\partial u_2}{\partial w_2}\\
  \frac{\partial^2 u_1}{\partial w_1 \partial w_2} &
  \frac{\partial^2 u_2}{\partial w_1 \partial w_2}
  \endvmatrix+\\
+&\frac{\partial w_1}{\partial x} \frac{\partial w_2}{\partial x}
  \frac{\partial w_2}{\partial y}
  \vmatrix \format \c \quad & \c \\
  \frac{\partial u_1}{\partial w_1} &
  \frac{\partial u_2}{\partial w_1}\\
  \frac{\partial^2 u_1}{\partial w_2^2} &
  \frac{\partial^2 u_2}{\partial w_2^2}
  \endvmatrix+
  \left(\frac{\partial w_2}{\partial x}\right)^2 \frac{\partial
  w_2}{\partial y}
  \vmatrix \format \c \quad & \c \\
  \frac{\partial u_1}{\partial w_2} &
  \frac{\partial u_2}{\partial w_2}\\
  \frac{\partial^2 u_1}{\partial w_2^2} &
  \frac{\partial^2 u_2}{\partial w_2^2}
  \endvmatrix+\\
&+\vmatrix \format \c \quad & \c \\
  \frac{\partial w_1}{\partial x} &
  \frac{\partial w_2}{\partial x}\\
  \frac{\partial^2 w_1}{\partial x \partial y} &
  \frac{\partial^2 w_2}{\partial x \partial y}
  \endvmatrix
  \vmatrix \format \c \quad & \c \\
  \frac{\partial u_1}{\partial w_1} &
  \frac{\partial u_2}{\partial w_1}\\
  \frac{\partial u_1}{\partial w_2} &
  \frac{\partial u_2}{\partial w_2}
  \endvmatrix.
\endaligned$$
We find that
$$\aligned
 &[u_1, u_2; x, y]_{x}\\
=&\frac{3 (\frac{\partial w_1}{\partial x})^2 \frac{\partial
  w_1}{\partial y}}{\frac{\partial(w_1, w_2)}{\partial(x, y)}}
  \{u_1, u_2; w_1, w_2\}_{w_1}+\frac{(\frac{\partial w_1}{\partial
  x})^2 \frac{\partial w_2}{\partial y}+2 \frac{\partial w_1}{\partial x}
  \frac{\partial w_2}{\partial x} \frac{\partial w_1}{\partial y}}
  {\frac{\partial(w_1, w_2)}{\partial(x, y)}} [u_1, u_2; w_1,
  w_2]_{w_1}+\\
+&\frac{3 (\frac{\partial w_2}{\partial x})^2 \frac{\partial
  w_2}{\partial y}}{\frac{\partial(w_2, w_1)}{\partial(x, y)}}
  \{u_1, u_2; w_1, w_2\}_{w_2}+\frac{(\frac{\partial w_2}{\partial
  x})^2 \frac{\partial w_1}{\partial y}+2 \frac{\partial w_1}{\partial x}
  \frac{\partial w_2}{\partial x} \frac{\partial w_2}{\partial y}}
  {\frac{\partial(w_2, w_1)}{\partial(x, y)}} [u_1, u_2; w_1,
  w_2]_{w_2}+\\
+&[w_1, w_2; x, y]_{x}.
\endaligned$$

  Exchange the position of $x$ and $y$, we get the other identities.
\flushpar
$\qquad \qquad \qquad \qquad \qquad \qquad \qquad \qquad \qquad
 \qquad \qquad \qquad \qquad \qquad \qquad \qquad \qquad \qquad
 \quad \boxed{}$

{\smc Corollary 2.6}. {\it The exchange of the first and the
second arguments obeys the laws:
$$\aligned
 &\{w_1, w_2; x, y\}_{x}\\
=&-\frac{(\frac{\partial w_1}{\partial x})^3}{\frac{\partial(w_1,
  w_2)}{\partial(x, y)}} \{x, y; w_1, w_2\}_{w_1}-\frac{(\frac{\partial
  w_1}{\partial x})^2 \frac{\partial w_2}{\partial x}}{\frac{\partial
  (w_1, w_2)}{\partial(x, y)}} [x, y; w_1, w_2]_{w_1}+\\
 &+\frac{(\frac{\partial w_2}{\partial x})^3}{\frac{\partial(w_1, w_2)}
  {\partial(x, y)}} \{x, y; w_1, w_2\}_{w_2}+\frac{\frac{\partial w_1}
  {\partial x} (\frac{\partial w_2}{\partial x})^2}{\frac{\partial(w_1,
  w_2)}{\partial(x, y)}} [x, y; w_1, w_2]_{w_2},
\endaligned\tag 2.17$$
$$\aligned
 &\{w_1, w_2; x, y\}_{y}\\
=&\frac{(\frac{\partial w_1}{\partial y})^3}{\frac{\partial(w_1,
  w_2)}{\partial(x, y)}} \{x, y; w_1, w_2\}_{w_1}+\frac{(\frac{\partial
  w_1}{\partial y})^2 \frac{\partial w_2}{\partial y}}{\frac{\partial
  (w_1, w_2)}{\partial(x, y)}} [x, y; w_1, w_2]_{w_1}+\\
 &-\frac{(\frac{\partial w_2}{\partial y})^3}{\frac{\partial(w_1, w_2)}
  {\partial(x, y)}} \{x, y; w_1, w_2\}_{w_2}-\frac{\frac{\partial w_1}
  {\partial y} (\frac{\partial w_2}{\partial y})^2}{\frac{\partial(w_1,
  w_2)}{\partial(x, y)}} [x, y; w_1, w_2]_{w_2},
\endaligned\tag 2.18$$
$$\aligned
 &[w_1, w_2; x, y]_{x}\\
=&-\frac{3 (\frac{\partial w_1}{\partial x})^2 \frac{\partial
  w_1}{\partial y}}{\frac{\partial(w_1, w_2)}{\partial(x, y)}}
  \{x, y; w_1, w_2\}_{w_1}-\frac{(\frac{\partial w_1}{\partial
  x})^2 \frac{\partial w_2}{\partial y}+2 \frac{\partial w_1}
  {\partial x} \frac{\partial w_2}{\partial x} \frac{\partial
  w_1}{\partial y}}{\frac{\partial(w_1, w_2)}{\partial(x, y)}}
  [x, y; w_1, w_2]_{w_1}+\\
 &+\frac{3 (\frac{\partial w_2}{\partial x})^2 \frac{\partial w_2}
  {\partial y}}{\frac{\partial(w_1, w_2)}{\partial(x, y)}}
  \{x, y; w_1, w_2\}_{w_2}+\frac{(\frac{\partial w_2}{\partial x})^2
  \frac{\partial w_1}{\partial y}+2 \frac{\partial w_1}{\partial x}
  \frac{\partial w_2}{\partial x} \frac{\partial w_2}{\partial y}}
  {\frac{\partial(w_1, w_2)}{\partial(x, y)}} [x, y; w_1, w_2]_{w_2},
\endaligned\tag 2.19$$
$$\aligned
 &[w_1, w_2; x, y]_{y}\\
=&\frac{3 \frac{\partial w_1}{\partial x} (\frac{\partial
  w_1}{\partial y})^2}{\frac{\partial(w_1, w_2)}{\partial(x, y)}}
  \{x, y; w_1, w_2\}_{w_1}+\frac{\frac{\partial w_2}{\partial
  x} (\frac{\partial w_1}{\partial y})^2+2 \frac{\partial w_1}
  {\partial x} \frac{\partial w_1}{\partial y} \frac{\partial
  w_2}{\partial y}}{\frac{\partial(w_1, w_2)}{\partial(x, y)}}
  [x, y; w_1, w_2]_{w_1}+\\
 &-\frac{3 \frac{\partial w_2}{\partial x} (\frac{\partial w_2}
  {\partial y})^2}{\frac{\partial(w_1, w_2)}{\partial(x, y)}}
  \{x, y; w_1, w_2\}_{w_2}-\frac{\frac{\partial w_1}{\partial x}
  (\frac{\partial w_2}{\partial y})^2+2 \frac{\partial w_2}{\partial x}
  \frac{\partial w_1}{\partial y} \frac{\partial w_2}{\partial y}}
  {\frac{\partial(w_1, w_2)}{\partial(x, y)}} [x, y; w_1, w_2]_{w_2}.
\endaligned\tag 2.20$$}

{\it Proof}. Put $u_1=x$ and $u_2=y$ in the Main Proposition 1.
Note that
$$\left\{\aligned
  \{x, y; x, y\}_{x} &=0,\\
  \{x, y; x, y\}_{y} &=0,\\
    [x, y; x, y]_{x} &=0,\\
    [x, y; x, y]_{y} &=0.
\endaligned\right.$$
We get the desired result.
\flushpar
$\qquad \qquad \qquad \qquad \qquad \qquad \qquad \qquad \qquad
 \qquad \qquad \qquad \qquad \qquad \qquad \qquad \qquad \qquad
 \quad \boxed{}$

  We may apply a linear fractional transformation to the second
arguments instead of the first. The result is as follows:

{\smc Corollary 2.7}. {\it For
$$(w_1, w_2)=\gamma(x, y)=\left(\frac{a_1 x+a_2 y+a_3}{c_1 x+c_2 y+c_3},
  \frac{b_1 x+b_2 y+b_3}{c_1 x+c_2 y+c_3}\right),$$
where $\gamma=\left(\matrix a_1 & a_2 & a_3\\ b_1 & b_2 & b_3\\
c_1 & c_2 & c_3 \endmatrix\right) \in GL(3, {\Bbb C})$, the
following four identities hold:
$$\aligned
 &\{u_1, u_2; x, y\}_{x}\\
=&\frac{A_1^3}{\Delta (c_1 x+c_2 y+c_3)^3}
  \left\{u_1, u_2; \frac{a_1 x+a_2 y+a_3}{c_1 x+c_2 y+c_3},
  \frac{b_1 x+b_2 y+b_3}{c_1 x+c_2 y+c_3}\right\}_{\frac{a_1
  x+a_2 y+a_3}{c_1 x+c_2 y+c_3}}+\\
+&\frac{A_1^2 B_1}{\Delta (c_1 x+c_2 y+c_3)^3}
  \left[u_1, u_2; \frac{a_1 x+a_2 y+a_3}{c_1 x+c_2 y+c_3},
  \frac{b_1 x+b_2 y+b_3}{c_1 x+c_2 y+c_3}\right]_{\frac{a_1
  x+a_2 y+a_3}{c_1 x+c_2 y+c_3}}+\\
-&\frac{B_1^3}{\Delta (c_1 x+c_2 y+c_3)^3}
  \left\{u_1, u_2; \frac{a_1 x+a_2 y+a_3}{c_1 x+c_2 y+c_3},
  \frac{b_1 x+b_2 y+b_3}{c_1 x+c_2 y+c_3}\right\}_{\frac{b_1
  x+b_2 y+b_3}{c_1 x+c_2 y+c_3}}+\\
-&\frac{A_1 B_1^2}{\Delta (c_1 x+c_2 y+c_3)^3}
  \left[u_1, u_2; \frac{a_1 x+a_2 y+a_3}{c_1 x+c_2 y+c_3},
  \frac{b_1 x+b_2 y+b_3}{c_1 x+c_2 y+c_3}\right]_{\frac{b_1
  x+b_2 y+b_3}{c_1 x+c_2 y+c_3}},
\endaligned\tag 2.21$$
$$\aligned
 &\{u_1, u_2; x, y\}_{y}\\
=&-\frac{A_2^3}{\Delta (c_1 x+c_2 y+c_3)^3}
  \left\{u_1, u_2; \frac{a_1 x+a_2 y+a_3}{c_1 x+c_2 y+c_3},
  \frac{b_1 x+b_2 y+b_3}{c_1 x+c_2 y+c_3}\right\}_{\frac{a_1
  x+a_2 y+a_3}{c_1 x+c_2 y+c_3}}+\\
-&\frac{A_2^2 B_2}{\Delta (c_1 x+c_2 y+c_3)^3}
  \left[u_1, u_2; \frac{a_1 x+a_2 y+a_3}{c_1 x+c_2 y+c_3},
  \frac{b_1 x+b_2 y+b_3}{c_1 x+c_2 y+c_3}\right]_{\frac{a_1
  x+a_2 y+a_3}{c_1 x+c_2 y+c_3}}+\\
+&\frac{B_2^3}{\Delta (c_1 x+c_2 y+c_3)^3}
  \left\{u_1, u_2; \frac{a_1 x+a_2 y+a_3}{c_1 x+c_2 y+c_3},
  \frac{b_1 x+b_2 y+b_3}{c_1 x+c_2 y+c_3}\right\}_{\frac{b_1
  x+b_2 y+b_3}{c_1 x+c_2 y+c_3}}+\\
+&\frac{A_2 B_2^2}{\Delta (c_1 x+c_2 y+c_3)^3}
  \left[u_1, u_2; \frac{a_1 x+a_2 y+a_3}{c_1 x+c_2 y+c_3},
  \frac{b_1 x+b_2 y+b_3}{c_1 x+c_2 y+c_3}\right]_{\frac{b_1
  x+b_2 y+b_3}{c_1 x+c_2 y+c_3}},
\endaligned\tag 2.22$$
$$\aligned
 &[u_1, u_2; x, y]_{x}\\
=&\frac{3 A_1^2 A_2}{\Delta (c_1 x+c_2 y+c_3)^3}
  \left\{u_1, u_2; \frac{a_1 x+a_2 y+a_3}{c_1 x+c_2 y+c_3},
  \frac{b_1 x+b_2 y+b_3}{c_1 x+c_2 y+c_3}\right\}_{\frac{a_1
  x+a_2 y+a_3}{c_1 x+c_2 y+c_3}}+\\
+&\frac{A_1^2 B_2+2 A_1 B_1 A_2}{\Delta (c_1 x+c_2 y+c_3)^3}
  \left[u_1, u_2; \frac{a_1 x+a_2 y+a_3}{c_1 x+c_2 y+c_3},
  \frac{b_1 x+b_2 y+b_3}{c_1 x+c_2 y+c_3}\right]_{\frac{a_1
  x+a_2 y+a_3}{c_1 x+c_2 y+c_3}}+\\
-&\frac{3 B_1^2 B_2}{\Delta (c_1 x+c_2 y+c_3)^3}
  \left\{u_1, u_2; \frac{a_1 x+a_2 y+a_3}{c_1 x+c_2 y+c_3},
  \frac{b_1 x+b_2 y+b_3}{c_1 x+c_2 y+c_3}\right\}_{\frac{b_1
  x+b_2 y+b_3}{c_1 x+c_2 y+c_3}}+\\
-&\frac{B_1^2 A_2+2 A_1 B_1 B_2}{\Delta (c_1 x+c_2 y+c_3)^3}
  \left[u_1, u_2; \frac{a_1 x+a_2 y+a_3}{c_1 x+c_2 y+c_3},
  \frac{b_1 x+b_2 y+b_3}{c_1 x+c_2 y+c_3}\right]_{\frac{b_1
  x+b_2 y+b_3}{c_1 x+c_2 y+c_3}},
\endaligned\tag 2.23$$
$$\aligned
 &[u_1, u_2; x, y]_{y}\\
=&-\frac{3 A_1 A_2^2}{\Delta (c_1 x+c_2 y+c_3)^3}
  \left\{u_1, u_2; \frac{a_1 x+a_2 y+a_3}{c_1 x+c_2 y+c_3},
  \frac{b_1 x+b_2 y+b_3}{c_1 x+c_2 y+c_3}\right\}_{\frac{a_1
  x+a_2 y+a_3}{c_1 x+c_2 y+c_3}}+\\
-&\frac{B_1 A_2^2+2 A_1 A_2 B_2}{\Delta (c_1 x+c_2 y+c_3)^3}
  \left[u_1, u_2; \frac{a_1 x+a_2 y+a_3}{c_1 x+c_2 y+c_3},
  \frac{b_1 x+b_2 y+b_3}{c_1 x+c_2 y+c_3}\right]_{\frac{a_1
  x+a_2 y+a_3}{c_1 x+c_2 y+c_3}}+\\
+&\frac{3 B_1 B_2^2}{\Delta (c_1 x+c_2 y+c_3)^3}
  \left\{u_1, u_2; \frac{a_1 x+a_2 y+a_3}{c_1 x+c_2 y+c_3},
  \frac{b_1 x+b_2 y+b_3}{c_1 x+c_2 y+c_3}\right\}_{\frac{b_1
  x+b_2 y+b_3}{c_1 x+c_2 y+c_3}}+\\
+&\frac{A_1 B_2^2+2 B_1 A_2 B_2}{\Delta (c_1 x+c_2 y+c_3)^3}
  \left[u_1, u_2; \frac{a_1 x+a_2 y+a_3}{c_1 x+c_2 y+c_3},
  \frac{b_1 x+b_2 y+b_3}{c_1 x+c_2 y+c_3}\right]_{\frac{b_1
  x+b_2 y+b_3}{c_1 x+c_2 y+c_3}}.
\endaligned\tag 2.24$$
Here, $\Delta=\det(\gamma)$ and
$$\left\{\aligned
 A_1 &=(a_1 c_2-a_2 c_1) y+(a_1 c_3-a_3 c_1),\\
 A_2 &=(a_2 c_1-a_1 c_2) x+(a_2 c_3-a_3 c_2),\\
 B_1 &=(b_1 c_2-b_2 c_1) y+(b_1 c_3-b_3 c_1),\\
 B_2 &=(b_2 c_1-b_1 c_2) x+(b_2 c_3-b_3 c_2).
\endaligned\right.$$}

{\it Proof}. We have
$$\left\{\aligned
  \frac{\partial w_1}{\partial x} &=(c_1 x+c_2 y+c_3)^{-2} A_1,\\
  \frac{\partial w_2}{\partial x} &=(c_1 x+c_2 y+c_3)^{-2} B_1,\\
  \frac{\partial w_1}{\partial y} &=(c_1 x+c_2 y+c_3)^{-2} A_2,\\
  \frac{\partial w_2}{\partial y} &=(c_1 x+c_2 y+c_3)^{-2} B_2.\\
\endaligned\right.$$
Note that
$$\left\{\aligned
  \{w_1, w_2; x, y\}_{x} &=0,\\
  \{w_1, w_2; x, y\}_{y} &=0,\\
    [w_1, w_2; x, y]_{x} &=0,\\
    [w_1, w_2; x, y]_{y} &=0.
\endaligned\right.$$
A straightforward calculation gives the result.
\flushpar
$\qquad \qquad \qquad \qquad \qquad \qquad \qquad \qquad \qquad
 \qquad \qquad \qquad \qquad \qquad \qquad \qquad \qquad \qquad
 \quad \boxed{}$

  The Main Proposition 1 implies that
$$\left(\matrix
  \{u_1, u_2; x, y\}_{x}\\
  \{u_1, u_2; x, y\}_{y}\\
   [u_1, u_2; x, y]_{x}\\
   [u_1, u_2; x, y]_{y}
  \endmatrix\right)
 =\frac{M(w_1, w_2; x, y)}{\frac{\partial(w_1, w_2)}{\partial(x, y)}}
  \left(\matrix
  \{u_1, u_2; w_1, w_2\}_{w_1}\\
  \{u_1, u_2; w_1, w_2\}_{w_2}\\
   [u_1, u_2; w_1, w_2]_{w_1}\\
   [u_1, u_2; w_1, w_2]_{w_2}
  \endmatrix\right)+\left(\matrix
  \{w_1, w_2; x, y\}_{x}\\
  \{w_1, w_2; x, y\}_{y}\\
   [w_1, w_2; x, y]_{x}\\
   [w_1, w_2; x, y]_{y}
  \endmatrix\right),\tag 2.25$$
where
$$\aligned
 &M(w_1, w_2; x, y):=\\
 &\left(\matrix
  (\frac{\partial w_1}{\partial x})^3 & -(\frac{\partial w_2}{\partial x})^3
 &(\frac{\partial w_1}{\partial x})^2 \frac{\partial w_2}{\partial x}
 &-\frac{\partial w_1}{\partial x} (\frac{\partial w_2}{\partial x})^2\\
  -(\frac{\partial w_1}{\partial y})^3 & (\frac{\partial w_2}{\partial y})^3
 &-(\frac{\partial w_1}{\partial y})^2 \frac{\partial w_2}{\partial y}
 &\frac{\partial w_1}{\partial y} (\frac{\partial w_2}{\partial y})^2\\
  3 (\frac{\partial w_1}{\partial x})^2 \frac{\partial w_1}{\partial y}
 &-3 (\frac{\partial w_2}{\partial x})^2 \frac{\partial w_2}{\partial y}
 &\matrix
  &(\frac{\partial w_1}{\partial x})^2 \frac{\partial w_2}{\partial y}+\\
  &+2 \frac{\partial w_1}{\partial x} \frac{\partial w_2}{\partial x}
   \frac{\partial w_1}{\partial y}
  \endmatrix
 &\matrix
  &-(\frac{\partial w_2}{\partial x})^2 \frac{\partial w_1}{\partial y}+\\
  &-2 \frac{\partial w_1}{\partial x} \frac{\partial w_2}{\partial x}
   \frac{\partial w_2}{\partial y}
  \endmatrix\\
  -3 \frac{\partial w_1}{\partial x} (\frac{\partial w_1}{\partial y})^2
 &3 \frac{\partial w_2}{\partial x} (\frac{\partial w_2}{\partial y})^2
 &\matrix
  &-\frac{\partial w_2}{\partial x} (\frac{\partial w_1}{\partial y})^2+\\
  &-2 \frac{\partial w_1}{\partial x} \frac{\partial w_1}{\partial y}
   \frac{\partial w_2}{\partial y}
  \endmatrix
 &\matrix
  &\frac{\partial w_1}{\partial x} (\frac{\partial w_2}{\partial y})^2+\\
  &+2 \frac{\partial w_2}{\partial x} \frac{\partial w_1}{\partial y}
   \frac{\partial w_2}{\partial y}
  \endmatrix
\endmatrix\right).
\endaligned\tag 2.26$$

{\smc Proposition 2.8}. {\it
$$M(w_1, w_2; x, y) \cdot M(u_1, u_2; w_1, w_2)
 =M(u_1, u_2; x, y).\tag 2.27$$}

{\it Proof}. Put
$$(a_{ij})=M(w_1, w_2; x, y) \cdot M(u_1, u_2; w_1, w_2).$$
We find that
$$\aligned
  a_{1, 1}
=&\left(\frac{\partial w_1}{\partial x}\right)^3
  \left(\frac{\partial u_1}{\partial w_1}\right)^3+
  \left(\frac{\partial w_2}{\partial x}\right)^3
  \left(\frac{\partial u_1}{\partial w_2}\right)^3+\\
+&\left(\frac{\partial w_1}{\partial x}\right)^2
  \frac{\partial w_2}{\partial x} \cdot 3
  \left(\frac{\partial u_1}{\partial w_1}\right)^2
  \frac{\partial u_1}{\partial w_2}+
  \frac{\partial w_1}{\partial x}
  \left(\frac{\partial w_2}{\partial x}\right)^2 \cdot 3
  \frac{\partial u_1}{\partial w_1}
  \left(\frac{\partial u_1}{\partial w_2}\right)^2\\
=&\left(\frac{\partial u_1}{\partial w_1} \frac{\partial w_1}
  {\partial x}+\frac{\partial u_1}{\partial w_2} \frac{\partial
  w_2}{\partial x}\right)^3\\
=&\left(\frac{\partial u_1}{\partial x}\right)^3.
\endaligned$$
$$\aligned
  a_{1, 3}
=&\left(\frac{\partial w_1}{\partial x}\right)^3
  \left(\frac{\partial u_1}{\partial w_1}\right)^2
  \frac{\partial u_2}{\partial w_1}+
  \left(\frac{\partial w_2}{\partial x}\right)^3
  \left(\frac{\partial u_1}{\partial w_2}\right)^2
  \frac{\partial u_2}{\partial w_2}+\\
+&\left(\frac{\partial w_1}{\partial x}\right)^2
  \frac{\partial w_2}{\partial x}
  \left[\left(\frac{\partial u_1}{\partial w_1}\right)^2
  \frac{\partial u_2}{\partial w_2}+2 \frac{\partial u_1}
  {\partial w_1} \frac{\partial u_2}{\partial w_1}
  \frac{\partial u_1}{\partial w_2}\right]+\\
+&\frac{\partial w_1}{\partial x}
  \left(\frac{\partial w_2}{\partial x}\right)^2
  \left[\frac{\partial u_2}{\partial w_1}
  \left(\frac{\partial u_1}{\partial w_2}\right)^2+2
  \frac{\partial u_1}{\partial w_1} \frac{\partial u_1}
  {\partial w_2} \frac{\partial u_2}{\partial w_2}\right]\\
=&\left(\frac{\partial u_1}{\partial w_1} \frac{\partial w_1}
  {\partial x}+\frac{\partial u_1}{\partial w_2} \frac{\partial
  w_2}{\partial x}\right)^2 \left(\frac{\partial u_2}{\partial
  w_1} \frac{\partial w_1}{\partial x}+\frac{\partial u_2}
  {\partial w_2} \frac{\partial w_2}{\partial x}\right)\\
=&\left(\frac{\partial u_1}{\partial x}\right)^2
  \frac{\partial u_2}{\partial x}.
\endaligned$$
$$\aligned
  a_{3, 1}
=&3 \left(\frac{\partial w_1}{\partial x}\right)^2
  \frac{\partial w_1}{\partial y} \left(\frac{\partial u_1}
  {\partial w_1}\right)^3+3 \left(\frac{\partial w_2}{\partial
  x}\right)^2 \frac{\partial w_2}{\partial y} \left(\frac{\partial
  u_1}{\partial w_2}\right)^3+\\
+&\left[\left(\frac{\partial w_1}{\partial x}\right)^2
  \frac{\partial w_2}{\partial y}+2 \frac{\partial w_1}
  {\partial x} \frac{\partial w_2}{\partial x}
  \frac{\partial w_1}{\partial y}\right] \cdot 3
  \left(\frac{\partial u_1}{\partial w_1}\right)^2
  \frac{\partial u_1}{\partial w_2}+\\
+&\left[\left(\frac{\partial w_2}{\partial x}\right)^2
  \frac{\partial w_1}{\partial y}+2 \frac{\partial w_1}
  {\partial x} \frac{\partial w_2}{\partial x}
  \frac{\partial w_2}{\partial y}\right] \cdot 3
  \frac{\partial u_1}{\partial w_1}
  \left(\frac{\partial u_1}{\partial w_2}\right)^2\\
=&3 \left(\frac{\partial u_1}{\partial w_1} \frac{\partial w_1}
  {\partial x}+\frac{\partial u_1}{\partial w_2} \frac{\partial
  w_2}{\partial x}\right)^2 \left(\frac{\partial u_1}{\partial
  w_1} \frac{\partial w_1}{\partial y}+\frac{\partial u_1}
  {\partial w_2} \frac{\partial w_2}{\partial y}\right)\\
=&3 \left(\frac{\partial u_1}{\partial x}\right)^2 \frac{\partial
  u_1}{\partial y}.
\endaligned$$
$$\aligned
  a_{3, 3}
=&3 \left(\frac{\partial w_1}{\partial x}\right)^2
  \frac{\partial w_1}{\partial y} \left(\frac{\partial u_1}
  {\partial w_1}\right)^2 \frac{\partial u_2}{\partial w_1}+
  3 \left(\frac{\partial w_2}{\partial x}\right)^2
  \frac{\partial w_2}{\partial y} \left(\frac{\partial u_1}
  {\partial w_2}\right)^2 \frac{\partial u_2}{\partial w_2}+\\
+&\left[\left(\frac{\partial w_1}{\partial x}\right)^2
  \frac{\partial w_2}{\partial y}+2 \frac{\partial w_1}{\partial x}
  \frac{\partial w_2}{\partial x} \frac{\partial w_1}{\partial y}\right]
  \left[\left(\frac{\partial u_1}{\partial w_1}\right)^2
  \frac{\partial u_2}{\partial w_2}+2 \frac{\partial u_1}{\partial w_1}
  \frac{\partial u_2}{\partial w_1} \frac{\partial u_1}{\partial w_2}\right]+\\
+&\left[\left(\frac{\partial w_2}{\partial x}\right)^2
  \frac{\partial w_1}{\partial y}+2 \frac{\partial w_1}{\partial x}
  \frac{\partial w_2}{\partial x} \frac{\partial w_2}{\partial y}\right]
  \left[\frac{\partial u_2}{\partial w_1}\left(\frac{\partial u_1}{\partial
  w_2}\right)^2+2 \frac{\partial u_1}{\partial w_1} \frac{\partial u_1}
  {\partial w_2} \frac{\partial u_2}{\partial w_2}\right]\\
=&\left(\frac{\partial u_1}{\partial w_1} \frac{\partial w_1}
  {\partial x}+\frac{\partial u_1}{\partial w_2} \frac{\partial
  w_2}{\partial x}\right)^2 \left(\frac{\partial u_2}{\partial
  w_1} \frac{\partial w_1}{\partial y}+\frac{\partial u_2}
  {\partial w_2} \frac{\partial w_2}{\partial y}\right)+\\
+&2 \left(\frac{\partial u_1}{\partial w_1} \frac{\partial w_1}
  {\partial x}+\frac{\partial u_1}{\partial w_2} \frac{\partial
  w_2}{\partial x}\right) \left(\frac{\partial u_2}{\partial w_1}
  \frac{\partial w_1}{\partial x}+\frac{\partial u_2}{\partial
  w_2} \frac{\partial w_2}{\partial x}\right) \left(\frac{\partial
  u_1}{\partial w_1} \frac{\partial w_1}{\partial y}+\frac{\partial
  u_1}{\partial w_2} \frac{\partial w_2}{\partial y}\right)\\
=&\left(\frac{\partial u_1}{\partial x}\right)^2 \frac{\partial
  u_2}{\partial y}+2 \frac{\partial u_1}{\partial x} \frac{\partial
  u_2}{\partial x} \frac{\partial u_1}{\partial y}.
\endaligned$$

  By the symmetry property, we can get the other terms in the matrix
$(a_{ij})$.
\flushpar
$\qquad \qquad \qquad \qquad \qquad \qquad \qquad \qquad \qquad
 \qquad \qquad \qquad \qquad \qquad \qquad \qquad \qquad \qquad
 \quad \boxed{}$

  In particular,
$$M(w_1, w_2; x, y) \cdot M(x, y; w_1, w_2)=M(x, y; x, y)=I.$$
Hence,
$$M(w_1, w_2; x, y)^{-1}=M(x, y; w_1, w_2).\tag 2.28$$

  Set
$$U_{(w_1, w_2; x, y)}:=\left(\matrix
  M(w_1, w_2; x, y) & c \frac{\partial(w_1, w_2)}{\partial(x, y)}
                 \left(\matrix  \{w_1, w_2; x, y\}_{x}\\
                                \{w_1, w_2; x, y\}_{y}\\
                                  [w_1, w_2; x, y]_{x}\\
                                  [w_1, w_2; x, y]_{y}
                      \endmatrix\right)\\
     0              & \frac{\partial(w_1, w_2)}{\partial(x, y)}
  \endmatrix\right).\tag 2.29$$
Then the above proposition implies that
$$U_{(w_1, w_2; x, y)} \cdot U_{(u_1, u_2; w_1, w_2)}
 =U_{(u_1, u_2; x, y)}.\tag 2.30$$

  The equation
$$U_{W, Z} U_{U, W}=U_{U, Z}$$
with $W=(w_1, w_2)$, $Z=(z_1, z_2)$ and $U=(u_1, u_2)$ is a
cocycle condition, which means that the $U$'s can be interpreted
as transition functions of a vector bundle. More precisely, when
we choose an open cover $\Sigma=\bigcup_{i} \Sigma_i$, with a
local coordinate $Z_i$ on $\Sigma_i$, then the above equation
permits us to interpret the $U_{Z_i, Z_j}$ as transition functions
on the intersection $\Sigma_i \cap \Sigma_j$. Let ${\Cal V}_{0}$
be this vector bundle. It is a five dimensional sub-bundle of the
infinite dimensional vector bundle ${\Cal V}$.

  We have
$$U_{(w_1, w_2; z_1, z_2)} \left(\matrix \Lambda\\ f \endmatrix\right)
 =\left(\matrix
  M(w_1, w_2; z_1, z_2) \Lambda+c \frac{\partial(w_1, w_2)}
  {\partial(z_1, z_2)}f \left(\matrix
  \{w_1, w_2; z_1, z_2\}_{z_1}\\
  \{w_1, w_2; z_1, z_2\}_{z_2}\\
    [w_1, w_2; z_1, z_2]_{z_1}\\
    [w_1, w_2; z_1, z_2]_{z_2}
  \endmatrix\right)\\
  \frac{\partial(w_1, w_2)}{\partial(z_1, z_2)} f
  \endmatrix\right),$$
where $\Lambda=\left(\matrix \lambda_1\\ \lambda_2\\ \lambda_3\\
\lambda_4 \endmatrix\right)$. In particular, when $c=0$, one has
$$U_{(w_1, w_2; z_1, z_2)} \left(\matrix \Lambda\\ f \endmatrix\right)
 =\left(\matrix
  M(w_1, w_2; z_1, z_2) \Lambda\\
  \frac{\partial(w_1, w_2)}{\partial(z_1, z_2)} f
  \endmatrix\right).$$

  At $c=0$, the pair $(\Lambda, f)$ is a section of the bundle $N
\oplus T$, where $N$ is the rank four vector bundle and $T$ is the
line bundle. At $c \neq 0$, the pair $(\Lambda, f)$ is a section
of a deformation ${\Cal V}_{0}$ of $N \oplus T$.

  The transition matrices $U_{(w_1, w_2; z_1, z_2)}$ are
triangular, this means that there is an exact sequence
$$0 \rightarrow N \rightarrow {\Cal V}_{0} \rightarrow T
    \rightarrow 0.$$
There is a natural map $\alpha$ from sections of $N$ to sections
of ${\Cal V}_{0}$, namely $\Lambda \mapsto (\Lambda, 0)$; and a
natural map $\beta$ from sections of ${\Cal V}_{0}$ to sections of
$T$, namely $(\Lambda, f) \mapsto f$. Moreover, $\beta \circ
\alpha=0$ and the image of $\alpha$ is the kernel of $\beta$.

  Let us consider the compact complex surfaces with $c_1^2=3c_2$,
where $c_i$ is the $i$-th Chern class, i.e., the ball quotient
surfaces. A uniformization of a ball quotient surface $\Sigma$
gives a covering by open sets $\Sigma_i$ with local parameters
$Z_i=(z_{1, i}, z_{2, i})$ such that on $\Sigma_i \cap \Sigma_j$
the transformations from $Z_i$ to $Z_j$ are $U(2, 1)$
transformations
$$(z_{1, j}, z_{2, j})=\gamma(z_{1, i}, z_{2, i})
 =\left(\frac{a_1 z_{1, i}+a_2 z_{2, i}+a_3}{c_1 z_{1, i}+
  c_2 z_{2, i}+c_3}, \frac{b_1 z_{1, i}+b_2 z_{2, i}+b_3}{
  c_1 z_{1, i}+c_2 z_{2, i}+c_3}\right)$$
with $\gamma=\left(\matrix a_1 & a_2 & a_3\\ b_1 & b_2 & b_3\\
c_1 & c_2 & c_3 \endmatrix\right) \in U(2, 1)$. For such
transformations, our four derivatives
$$\{z_{1, i}, z_{2, i}; z_{1, j}, z_{2, j}\}_{z_{1, j}}, \quad
  \{z_{1, i}, z_{2, i}; z_{1, j}, z_{2, j}\}_{z_{2, j}}$$
and
$$[z_{1, i}, z_{2, i}; z_{1, j}, z_{2, j}]_{z_{1, j}}, \quad
  [z_{1, i}, z_{2, i}; z_{1, j}, z_{2, j}]_{z_{2, j}}$$
vanish, so the transition matrices become diagonal, and the
extension ${\Cal V}_{0}$ is split as $N \oplus T$.

  Recall that
$$\{z_1, z_2; w_1, w_2\}_{w_1}
:=\frac{\vmatrix \format \c \quad & \c \\
  \frac{\partial z_1}{\partial w_1} &
  \frac{\partial z_2}{\partial w_1}\\
  \frac{\partial^2 z_1}{\partial w_1^2} &
  \frac{\partial^2 z_2}{\partial w_1^2}\endvmatrix}
 {\vmatrix \format \c \quad & \c\\
  \frac{\partial z_1}{\partial w_1} &
  \frac{\partial z_2}{\partial w_1}\\
  \frac{\partial z_1}{\partial w_2} &
  \frac{\partial z_2}{\partial w_2}
  \endvmatrix}, \quad
  \{z_1, z_2; w_1, w_2\}_{w_2}
:=\frac{\vmatrix \format \c \quad & \c \\
  \frac{\partial z_1}{\partial w_2} &
  \frac{\partial z_2}{\partial w_2}\\
  \frac{\partial^2 z_1}{\partial w_2^2} &
  \frac{\partial^2 z_2}{\partial w_2^2}\endvmatrix}
 {\vmatrix \format \c \quad & \c\\
  \frac{\partial z_1}{\partial w_2} &
  \frac{\partial z_2}{\partial w_2}\\
  \frac{\partial z_1}{\partial w_1} &
  \frac{\partial z_2}{\partial w_1}
  \endvmatrix},$$
$$[z_1, z_2; w_1, w_2]_{w_1}
:=\frac{\vmatrix \format \c \quad & \c \\
  \frac{\partial z_1}{\partial w_2} &
  \frac{\partial z_2}{\partial w_2}\\
  \frac{\partial^2 z_1}{\partial w_1^2} &
  \frac{\partial^2 z_2}{\partial w_1^2}\endvmatrix+2
  \vmatrix \format \c \quad & \c \\
  \frac{\partial z_1}{\partial w_1} &
  \frac{\partial z_2}{\partial w_1}\\
  \frac{\partial^2 z_1}{\partial w_1 \partial w_2} &
  \frac{\partial^2 z_2}{\partial w_1 \partial w_2}\endvmatrix}
 {\vmatrix \format \c \quad & \c\\
  \frac{\partial z_1}{\partial w_1} &
  \frac{\partial z_2}{\partial w_1}\\
  \frac{\partial z_1}{\partial w_2} &
  \frac{\partial z_2}{\partial w_2}
  \endvmatrix},$$
$$[z_1, z_2; w_1, w_2]_{w_2}
:=\frac{\vmatrix \format \c \quad & \c \\
  \frac{\partial z_1}{\partial w_1} &
  \frac{\partial z_2}{\partial w_1}\\
  \frac{\partial^2 z_1}{\partial w_2^2} &
  \frac{\partial^2 z_2}{\partial w_2^2}\endvmatrix+2
  \vmatrix \format \c \quad & \c \\
  \frac{\partial z_1}{\partial w_2} &
  \frac{\partial z_2}{\partial w_2}\\
  \frac{\partial^2 z_1}{\partial w_1 \partial w_2} &
  \frac{\partial^2 z_2}{\partial w_1 \partial w_2}\endvmatrix}
 {\vmatrix \format \c \quad & \c\\
  \frac{\partial z_1}{\partial w_2} &
  \frac{\partial z_2}{\partial w_2}\\
  \frac{\partial z_1}{\partial w_1} &
  \frac{\partial z_2}{\partial w_1}
  \endvmatrix}.$$

  The next Main Proposition 2 gives the relation between the
deformation of the Jacobian and our four derivatives:

{\smc Proposition 2.9 (Main Proposition 2)}. {\it The following
identity holds:
$$\aligned
 &\frac{\partial(f_1, f_2)}{\partial(z_1, z_2)}\\
=&\frac{\partial(\widehat{f_1}, \widehat{f_2})}{\partial(w_1,
  w_2)}-\widehat{f_1} \frac{\partial \widehat{f_1}}{\partial w_2}
  \{z_1, z_2; w_1, w_2\}_{w_1}-\widehat{f_2} \frac{\partial
  \widehat{f_2}}{\partial w_1} \{z_1, z_2; w_1, w_2\}_{w_2}+\\
 &-\widehat{f_1} \frac{\partial \widehat{f_2}}{\partial w_2}
  \frac{\vmatrix \format \c \quad & \c \\
  \frac{\partial z_1}{\partial w_2} &
  \frac{\partial z_2}{\partial w_2}\\
  \frac{\partial^2 z_1}{\partial w_1^2} &
  \frac{\partial^2 z_2}{\partial w_1^2}\endvmatrix}
 {\vmatrix \format \c \quad & \c\\
  \frac{\partial z_1}{\partial w_1} &
  \frac{\partial z_2}{\partial w_1}\\
  \frac{\partial z_1}{\partial w_2} &
  \frac{\partial z_2}{\partial w_2}
  \endvmatrix}
  -\left(\widehat{f_2} \frac{\partial \widehat{f_1}}{\partial w_2}-
  \widehat{f_1} \frac{\partial \widehat{f_1}}{\partial w_1}\right)
  \frac{\vmatrix \format \c \quad & \c \\
  \frac{\partial z_1}{\partial w_1} &
  \frac{\partial z_2}{\partial w_1}\\
  \frac{\partial^2 z_1}{\partial w_1 \partial w_2} &
  \frac{\partial^2 z_2}{\partial w_1 \partial w_2}\endvmatrix}
 {\vmatrix \format \c \quad & \c\\
  \frac{\partial z_1}{\partial w_1} &
  \frac{\partial z_2}{\partial w_1}\\
  \frac{\partial z_1}{\partial w_2} &
  \frac{\partial z_2}{\partial w_2}
  \endvmatrix}+\\
 &-\widehat{f_2} \frac{\partial \widehat{f_1}}{\partial w_1}
  \frac{\vmatrix \format \c \quad & \c \\
  \frac{\partial z_1}{\partial w_1} &
  \frac{\partial z_2}{\partial w_1}\\
  \frac{\partial^2 z_1}{\partial w_2^2} &
  \frac{\partial^2 z_2}{\partial w_2^2}\endvmatrix}
 {\vmatrix \format \c \quad & \c\\
  \frac{\partial z_1}{\partial w_2} &
  \frac{\partial z_2}{\partial w_2}\\
  \frac{\partial z_1}{\partial w_1} &
  \frac{\partial z_2}{\partial w_1}
  \endvmatrix}-\left(\widehat{f_1} \frac{\partial \widehat{f_2}}
  {\partial w_1}-\widehat{f_2} \frac{\partial \widehat{f_2}}
  {\partial w_2}\right) \frac{\vmatrix \format \c \quad & \c \\
  \frac{\partial z_1}{\partial w_2} &
  \frac{\partial z_2}{\partial w_2}\\
  \frac{\partial^2 z_1}{\partial w_1 \partial w_2} &
  \frac{\partial^2 z_2}{\partial w_1 \partial w_2}\endvmatrix}
 {\vmatrix \format \c \quad & \c\\
  \frac{\partial z_1}{\partial w_2} &
  \frac{\partial z_2}{\partial w_2}\\
  \frac{\partial z_1}{\partial w_1} &
  \frac{\partial z_2}{\partial w_1}
  \endvmatrix}+\\
 &+{\widehat{f_1}}^2 \frac{\vmatrix \format \c \quad & \c \\
  \frac{\partial^2 z_1}{\partial w_1^2} &
  \frac{\partial^2 z_2}{\partial w_1^2}\\
  \frac{\partial^2 z_1}{\partial w_1 \partial w_2} &
  \frac{\partial^2 z_2}{\partial w_1 \partial w_2}\endvmatrix}
 {\vmatrix \format \c \quad & \c\\
  \frac{\partial z_1}{\partial w_1} &
  \frac{\partial z_2}{\partial w_1}\\
  \frac{\partial z_1}{\partial w_2} &
  \frac{\partial z_2}{\partial w_2}
  \endvmatrix}+{\widehat{f_2}}^2
  \frac{\vmatrix \format \c \quad & \c \\
  \frac{\partial^2 z_1}{\partial w_2^2} &
  \frac{\partial^2 z_2}{\partial w_2^2}\\
  \frac{\partial^2 z_1}{\partial w_1 \partial w_2} &
  \frac{\partial^2 z_2}{\partial w_1 \partial w_2}\endvmatrix}
 {\vmatrix \format \c \quad & \c\\
  \frac{\partial z_1}{\partial w_2} &
  \frac{\partial z_2}{\partial w_2}\\
  \frac{\partial z_1}{\partial w_1} &
  \frac{\partial z_2}{\partial w_1}
  \endvmatrix}+\widehat{f_1} \widehat{f_2}
  \frac{\vmatrix \format \c \quad & \c \\
  \frac{\partial^2 z_1}{\partial w_1^2} &
  \frac{\partial^2 z_2}{\partial w_1^2}\\
  \frac{\partial^2 z_1}{\partial w_2^2} &
  \frac{\partial^2 z_2}{\partial w_2^2}\endvmatrix}
 {\vmatrix \format \c \quad & \c\\
  \frac{\partial z_1}{\partial w_1} &
  \frac{\partial z_2}{\partial w_1}\\
  \frac{\partial z_1}{\partial w_2} &
  \frac{\partial z_2}{\partial w_2}
  \endvmatrix},
\endaligned\tag 2.31$$
where
$$(f_1, f_2)=(\widehat{f_1}, \widehat{f_2})
  \left(\matrix
  \frac{\partial z_1}{\partial w_1} &
  \frac{\partial z_2}{\partial w_1}\\
  \frac{\partial z_1}{\partial w_2} &
  \frac{\partial z_2}{\partial w_2}
  \endmatrix\right).$$}

{\it Proof}. Note that
$$\frac{\partial(f_1, f_2)}{\partial(z_1, z_2)}
 =\frac{\partial(f_1, f_2)}{\partial(w_1, w_2)}
  \frac{\partial(w_1, w_2)}{\partial(z_1, z_2)}.$$
Here,
$$f_1=\widehat{f_1} \frac{\partial z_1}{\partial w_1}+
      \widehat{f_2} \frac{\partial z_1}{\partial w_2}, \quad
  f_2=\widehat{f_1} \frac{\partial z_2}{\partial w_1}+
      \widehat{f_2} \frac{\partial z_2}{\partial w_2}.$$
We find that
$$\aligned
 &\frac{\partial(f_1, f_2)}{\partial(w_1, w_2)}\\
=&\frac{\partial(\widehat{f_1}, \widehat{f_2})}{\partial(w_1,
  w_2)} \frac{\partial(z_1, z_2)}{\partial(w_1, w_2)}-
  \widehat{f_1} \frac{\partial \widehat{f_1}}{\partial w_2}
  \vmatrix \format \c \quad & \c \\
  \frac{\partial z_1}{\partial w_1} &
  \frac{\partial z_2}{\partial w_1}\\
  \frac{\partial^2 z_1}{\partial w_1^2} &
  \frac{\partial^2 z_2}{\partial w_1^2}\endvmatrix+
  \widehat{f_2} \frac{\partial \widehat{f_2}}{\partial w_1}
  \vmatrix \format \c \quad & \c \\
  \frac{\partial z_1}{\partial w_2} &
  \frac{\partial z_2}{\partial w_2}\\
  \frac{\partial^2 z_1}{\partial w_2^2} &
  \frac{\partial^2 z_2}{\partial w_2^2}\endvmatrix+\\
-&\widehat{f_1} \frac{\partial \widehat{f_2}}{\partial w_2}
  \vmatrix \format \c \quad & \c \\
  \frac{\partial z_1}{\partial w_2} &
  \frac{\partial z_2}{\partial w_2}\\
  \frac{\partial^2 z_1}{\partial w_1^2} &
  \frac{\partial^2 z_2}{\partial w_1^2}\endvmatrix-
  \left(\widehat{f_2} \frac{\partial \widehat{f_1}}{\partial w_2}
  -\widehat{f_1} \frac{\partial \widehat{f_1}}{\partial w_1}\right)
  \vmatrix \format \c \quad & \c \\
  \frac{\partial z_1}{\partial w_1} &
  \frac{\partial z_2}{\partial w_1}\\
  \frac{\partial^2 z_1}{\partial w_1 \partial w_2} &
  \frac{\partial^2 z_2}{\partial w_1 \partial w_2}\endvmatrix+\\
+&\widehat{f_2} \frac{\partial \widehat{f_1}}{\partial w_1}
  \vmatrix \format \c \quad & \c \\
  \frac{\partial z_1}{\partial w_1} &
  \frac{\partial z_2}{\partial w_1}\\
  \frac{\partial^2 z_1}{\partial w_2^2} &
  \frac{\partial^2 z_2}{\partial w_2^2}\endvmatrix+
  \left(\widehat{f_1} \frac{\partial \widehat{f_2}}{\partial w_1}
  -\widehat{f_2} \frac{\partial \widehat{f_2}}{\partial w_2}\right)
  \vmatrix \format \c \quad & \c \\
  \frac{\partial z_1}{\partial w_2} &
  \frac{\partial z_2}{\partial w_2}\\
  \frac{\partial^2 z_1}{\partial w_1 \partial w_2} &
  \frac{\partial^2 z_2}{\partial w_1 \partial w_2}\endvmatrix+\\
+&\widehat{f_1}^2
  \vmatrix \format \c \quad & \c \\
  \frac{\partial^2 z_1}{\partial w_1^2} &
  \frac{\partial^2 z_2}{\partial w_1^2}\\
  \frac{\partial^2 z_1}{\partial w_1 \partial w_2} &
  \frac{\partial^2 z_2}{\partial w_1 \partial w_2}\endvmatrix-
  \widehat{f_2}^2
  \vmatrix \format \c \quad & \c \\
  \frac{\partial^2 z_1}{\partial w_2^2} &
  \frac{\partial^2 z_2}{\partial w_2^2}\\
  \frac{\partial^2 z_1}{\partial w_1 \partial w_2} &
  \frac{\partial^2 z_2}{\partial w_1 \partial w_2}\endvmatrix+
  \widehat{f_1} \widehat{f_2}
  \vmatrix \format \c \quad & \c \\
  \frac{\partial^2 z_1}{\partial w_1^2} &
  \frac{\partial^2 z_2}{\partial w_1^2}\\
  \frac{\partial^2 z_1}{\partial w_2^2} &
  \frac{\partial^2 z_2}{\partial w_2^2}\endvmatrix.
\endaligned$$
Thus, we can obtain the desired result.
\flushpar
$\qquad \qquad \qquad \qquad \qquad \qquad \qquad \qquad \qquad
 \qquad \qquad \qquad \qquad \qquad \qquad \qquad \qquad \qquad
 \quad \boxed{}$

  Now, put
$$\phi(f_1, f_2):=\iint \frac{\partial(f_1, f_2)}{\partial(t_1, t_2)}
                  dt_1 \wedge dt_2.$$
We define an operation
$$\aligned
  &\left[f_1 \frac{\partial}{\partial t_1} \otimes \frac{\partial}{\partial t_2},
    f_2 \frac{\partial}{\partial t_1} \otimes \frac{\partial}{\partial t_2}\right]_{Z_0}^{*}\\
:=&c \cdot \phi(f_1, f_2)-f_1 \frac{\partial f_1}{\partial t_2}
   \frac{\partial}{\partial t_1}-f_1 \frac{\partial f_2}{\partial t_2}
   \frac{\partial}{\partial t_1}-\left(f_2 \frac{\partial f_1}{\partial t_2}
   -f_1 \frac{\partial f_1}{\partial t_1}\right) \frac{\partial}{\partial t_1}+\\
  &-f_2 \frac{\partial f_2}{\partial t_1} \frac{\partial}{\partial t_2}
   -f_2 \frac{\partial f_1}{\partial t_1} \frac{\partial}{\partial t_2}
   -\left(f_1 \frac{\partial f_2}{\partial t_1}-f_2 \frac{\partial f_2}
   {\partial t_2}\right) \frac{\partial}{\partial t_2}+f_1^2
   \frac{\partial^2}{\partial t_1^2}+f_2^2 \frac{\partial^2}{\partial t_2^2}
   +f_1 f_2 \frac{\partial^2}{\partial t_1 \partial t_2}.
\endaligned$$
Note that
$$f_1 \frac{\partial f_1}{\partial t_2}+f_1 \frac{\partial f_2}{\partial t_2}
  +\left(f_2 \frac{\partial f_1}{\partial t_2}-f_1 \frac{\partial f_1}{\partial
  t_1}\right)
 =(f_1+f_2) \frac{\partial f_1}{\partial t_2}+f_1 \left(\frac{\partial f_2}
  {\partial t_2}-\frac{\partial f_1}{\partial t_1}\right)$$
and
$$f_2 \frac{\partial f_2}{\partial t_1}+f_2 \frac{\partial f_1}{\partial t_1}
  +\left(f_1 \frac{\partial f_2}{\partial t_1}-f_2 \frac{\partial f_2}{\partial
  t_2}\right)
 =(f_1+f_2) \frac{\partial f_2}{\partial t_1}+f_2 \left(\frac{\partial f_1}
  {\partial t_1}-\frac{\partial f_2}{\partial t_2}\right).$$
Modulo the terms of second order partial derivatives, we define
the following operation:
$$\aligned
  &\left[f_1 \frac{\partial}{\partial t_1} \otimes \frac{\partial}{\partial t_2},
   f_2 \frac{\partial}{\partial t_1} \otimes \frac{\partial}{\partial
   t_2}\right]_{Z_0}
 :=\left[(f_1+f_2) \frac{\partial f_1}{\partial t_2}+f_1
   \left(\frac{\partial f_2}{\partial t_2}-\frac{\partial
   f_1}{\partial t_1}\right)\right] \frac{\partial}{\partial t_1}+\\
 +&\left[(f_1+f_2) \frac{\partial f_2}{\partial t_1}+f_2
   \left(\frac{\partial f_1}{\partial t_1}-\frac{\partial
   f_2}{\partial t_2}\right)\right] \frac{\partial}{\partial t_2}+
   c \cdot \phi(f_1, f_2).
\endaligned\tag 2.32$$

\vskip 0.5 cm
\centerline{\bf 3. Nonlinear partial differential equations and
                   modular functions}
\centerline{\bf associated to $U(2, 1)$}
\vskip 0.5 cm

  Let $U(2, 1)=\{ g \in GL(3, {\Bbb C}): g^{*} J g=J \}$
where $J=\left(\matrix
    &    & 1\\
    & -1 &  \\
  1 &    &
   \endmatrix\right)$. The corresponding complex domain
${\frak S}_{2}=\{(z_1, z_2) \in {\Bbb C}^2: z_1+\overline{z_1}
-z_2 \overline{z_2}>0 \}$ is the simplest example of a bounded
symmetric domain which is not a tube domain. Put
$\rho(Z):=z_1+\overline{z_1}-z_2 \overline{z_2}$ for $Z=(z_1, z_2)
\in {\frak S}_{2}$.

  The imaginary quadratic field $K={\Bbb Q}(\sqrt{-3})$ is
called the field of Eisenstein numbers. Its ring of integers
$${\Cal O}_K=\{ a/2+b \sqrt{-3}/2: a, b \in {\Bbb Z},
  a \equiv b (\text{mod} 2) \}={\Bbb Z}+{\Bbb Z} [\omega]$$
with $\omega=\frac{-1+\sqrt{-3}}{2}$ is called the ring of
Eisenstein integers.

  Let $P$ be the parabolic subgroup of $U(2, 1)$ which consists
of upper triangular matrices. This group $P$ has the Langlands
decomposition $P=NAM$, where
$$N=\left\{ n=[\alpha, \beta]:=\left(\matrix 1 & \alpha & \beta\\
    & 1 & \overline{\alpha}\\ & & 1 \endmatrix\right): \alpha,
    \beta \in {\Bbb C}, \beta+\overline{\beta}=\alpha \overline{\alpha}
    \right\},$$
$$A=\left\{ a=\left(\matrix \rho^{-1} &  & \\ & 1 & \\ & & \rho
    \endmatrix\right): \rho \in {\Bbb R}, \rho>0 \right\},$$
$$M=\left\{ m=\left(\matrix \beta &  & \\ & \beta^{-2} &  \\
    &  & \beta \endmatrix\right): \beta \in {\Bbb C}, |\beta|^2=1
    \right\}.$$
Therefore,
$$N({\Cal O}_{K})=\{ n=[\alpha, \beta]: \alpha, \beta \in
  {\Cal O}_{K}, \beta+\overline{\beta}=\alpha \overline{\alpha}\}.\tag 3.1$$
Put $\Gamma=U(2, 1; {\Cal O}_{K})$. Then $\Gamma \cap N=N({\Cal
O}_{K})$.

  Let
$$T_1=\left(\matrix
      1 & 1 & -\omega\\
        & 1 &       1\\
        &   &       1
      \endmatrix\right), \quad
  T_2=\left(\matrix
      1 & \omega & -\omega\\
        &      1 & \overline{\omega}\\
        &        & 1
      \endmatrix\right), \quad
  S=-\overline{\omega} J
   =\left(\matrix
                       &                   & -\overline{\omega}\\
                       & \overline{\omega} &     \\
    -\overline{\omega} &                   &
    \endmatrix\right).$$

{\smc Theorem 3.1}. {\it The subgroup $N({\Cal O}_{K})$ is
generated by $T_1$ and $T_2$.}

{\it Proof}. Suppose that $\beta+\overline{\beta}=\alpha
\overline{\alpha}$ with $\alpha, \beta \in {\Cal O}_{K}$. Set
$\alpha=m+n \omega$, $\beta=b_1+b_2 \omega$, where $m, n, b_1, b_2
\in {\Bbb Z}$. Then $2b_1-b_2=m^2-mn+n^2$. The general solution is
given by
$$b_2=m+n+mn+2l, \quad b_1=\frac{1}{2}(m^2+n^2+m+n)+l, \quad (l
  \in {\Bbb Z}).$$
Hence,
$$\beta=\frac{1}{2}(m^2-mn+n^2)+\left[\frac{1}{2}(m+n+mn)+l\right]
        (\omega-\overline{\omega}).$$
Set
$$[\alpha, \beta]:=\left(\matrix 1 & \alpha & \beta\\
  & 1 & \overline{\alpha}\\ &  & 1 \endmatrix\right),
  \quad \beta+\overline{\beta}=\alpha \overline{\alpha},
  \quad \alpha, \beta \in {\Cal O}_{K}.$$
Then
$$[\alpha_1, \beta_1] [\alpha_2, \beta_2]=[\alpha_1+\alpha_2,
  \beta_1+\beta_2+\alpha_1 \overline{\alpha_2}].$$
$$[\alpha, \beta]^{-1}=[-\alpha, \overline{\beta}], \quad
  [\alpha_1, \beta_1] [\alpha_2, \beta_2] [\alpha_1, \beta_1]^{-1}
  [\alpha_2, \beta_2]^{-1}=[0, \alpha_1 \overline{\alpha_2}-
  \overline{\alpha_1} \alpha_2].$$
We have
$$[1, -\omega] [\omega, -\omega] [1, -\omega]^{-1}
  [\omega, -\omega]^{-1}=[0, \overline{\omega}-\omega].$$
$$\aligned
 [1, -\omega]^{m}=\left[m, \frac{1}{2}m(m-1)-m \omega\right], \quad
 &(m \in {\Bbb Z}),\\
 [\omega, -\omega]^{n}=\left[n \omega, \frac{1}{2}n(n-1)-n \omega\right],
 \quad &(n \in {\Bbb Z}),\\
 [0, \overline{\omega}-\omega]^{l}=[0, l(\overline{\omega}-\omega)],
 \quad &(l \in {\Bbb Z}).
\endaligned$$
This gives that
$$\aligned
 &[1, -\omega]^{m} [\omega, -\omega]^{n} [0, \overline{\omega}-\omega]^{l}\\
=&\left[m+n \omega, \frac{1}{2}(m^2-mn+n^2)+
  \left[\frac{1}{2}(m+n+mn)+l\right](\overline{\omega}-\omega)\right].
\endaligned$$
Thus,
$$\aligned
 &[1, -\omega]^{m} [\omega, -\omega]^{n} [0, \overline{\omega}-
  \omega]^{-l-m-n-mn}\\
=&\left[m+n \omega, \frac{1}{2}(m^2-mn+n^2)+
  \left[\frac{1}{2}(m+n+mn)+l\right](\omega-\overline{\omega})\right].
\endaligned$$
Note that $T_1=[1, -\omega]$ and $T_2=[\omega, -\omega]$. This
completes the proof of Theorem 3.1.

\flushpar
$\qquad \qquad \qquad \qquad \qquad \qquad \qquad \qquad\qquad
 \qquad \qquad \qquad \qquad \qquad \qquad \qquad \qquad \qquad
 \quad \boxed{}$

  In fact, $T_1$ and $T_2$ are two non-commutative translations,
their commutator $[T_1, T_2]:=T_1 T_2 T_1^{-1} T_2^{-1}=[0,
\overline{\omega}-\omega]$ is a commutative translation. We have
$$T_{1}(z_1, z_2)=(z_1+z_2-\omega, z_2+1), \quad
  T_{2}(z_1, z_2)=(z_1+\omega z_2-\omega,
  z_2+\overline{\omega}),$$
$$S(z_1, z_2)=\left(\frac{1}{z_1}, -\frac{z_2}{z_1}\right),$$
where $(z_1, z_2) \in {\frak S}_{2}$.

  When $T_1(z_1, z_2)=S(z_1, z_2)$, then
$$z_1+z_2-\omega=\frac{1}{z_1}, \quad z_2+1=-\frac{z_2}{z_1}.\tag 3.2$$
The solutions are given by
$$(z_1, z_2)=(\omega, \overline{\omega}), \quad
  \left(\omega i, \frac{-\omega i}{\omega i+1}\right), \quad
  \left(-\omega i, \frac{\omega i}{-\omega i+1}\right).$$
Here, $\rho(\omega, \overline{\omega})=-2<0$, $\rho\left(\omega i,
\frac{-\omega i}{\omega i+1}\right)=-2-2 \sqrt{3}<0$ and
$\rho\left(-\omega i, \frac{\omega i}{1-\omega
i}\right)=2(\sqrt{3}-1)>0$. Thus, $\left(-\omega i, \frac{\omega
i}{1-\omega i}\right) \in {\frak S}_{2}$.

  Let
$$(-\omega, 0) \mapsto 1, \quad (-\overline{\omega}, 0) \mapsto z_2,
  \quad \infty \mapsto \infty, \quad \left(-\omega i, \frac{\omega i}{1-\omega i}\right)
  \mapsto 0.\tag 3.3$$
The point $0$ is of order $3$, $1$ is of order $1$, $z_2$ is of
order $1$ and $\infty$ is a logarithmic branch point.

  By
$$S\left(-\omega i, \frac{\omega i}{1-\omega i}\right)
 =T_1\left(-\omega i, \frac{\omega i}{1-\omega i}\right),$$
we have
$$\left(-\omega i, \frac{\omega i}{1-\omega i}\right)
 =S^2 \left(-\omega i, \frac{\omega i}{1-\omega i}\right)
 =S T_1\left(-\omega i, \frac{\omega i}{1-\omega i}\right).$$
On the other hand, $(S T_1)^4=\omega I$. This implies that $S T_1$
is of order four. In fact, $S$ is of order two and
$$S(-\omega, 0)=(-\overline{\omega}, 0), \quad
  S(-\overline{\omega}, 0)=(-\omega, 0).$$
Note that
$$\frac{\partial (v_i^{\frac{1}{4}})}{\partial w_j}
 =\frac{1}{4} v_i^{-\frac{3}{4}} \frac{\partial v_i}{\partial w_j}, \quad
  \frac{\partial (v_i^{\frac{1}{2}})}{\partial w_j}
 =\frac{1}{2} v_i^{-\frac{1}{2}} \frac{\partial v_i}{\partial w_j}, \quad
  \frac{\partial (\log v_i)}{\partial w_j}
 =\frac{1}{v_i} \frac{\partial v_i}{\partial w_j}, \quad
  i, j=1, 2.$$

  This leads us to consider
$$v_1^{-\frac{3}{4}} (1-v_1)^{-\frac{1}{2}} (v_2-v_1)^{-\frac{1}{2}}
  \vmatrix \format \c \quad & \c\\
  \frac{\partial v_1}{\partial w_1} & \frac{\partial v_2}{\partial w_1}\\
  \frac{\partial v_1}{\partial w_2} & \frac{\partial v_2}{\partial w_2}
  \endvmatrix.\tag 3.4$$

  When $T_2(z_1, z_2)=S(z_1, z_2)$, then
$$z_1+\omega z_2-\omega=\frac{1}{z_1}, \quad
  z_2+\overline{\omega}=-\frac{z_2}{z_1}.\tag 3.5$$
The solutions are given by
$$(z_1, z_2)=(\omega, \omega), \quad \left(\omega i, \frac{-i}{\omega
  i+1}\right), \quad \left(-\omega i, \frac{i}{-\omega i+1}\right).$$
Here, $\rho(\omega, \omega)=-2<0$, $\rho\left(\omega i,
\frac{-i}{\omega i+1}\right)=-2-2 \sqrt{3}<0$ and
$\rho\left(-\omega i, \frac{i}{1-\omega
i}\right)=2(\sqrt{3}-1)>0$. Thus, $\left(-\omega i,
\frac{i}{1-\omega i}\right) \in {\frak S}_{2}$.

  Let
$$(-\omega, 0) \mapsto 1, \quad (-\overline{\omega}, 0) \mapsto
  z_1, \quad \infty \mapsto \infty, \quad \left(-\omega i, \frac{i}{1
  -\omega i}\right) \mapsto 0.\tag 3.6$$
The point $0$ is of order $3$, $1$ is of order $1$, $z_1$ is of
order $1$ and $\infty$ is a logarithmic branch point.

  By
$$S\left(-\omega i, \frac{i}{1-\omega i}\right)
 =T_2\left(-\omega i, \frac{i}{1-\omega i}\right),$$
we have
$$\left(-\omega i, \frac{i}{1-\omega i}\right)
 =S^2 \left(-\omega i, \frac{i}{1-\omega i}\right)
 =S T_2\left(-\omega i, \frac{i}{1-\omega i}\right).$$
On the other hand, $(S T_2)^4=\omega I$. This implies that $S T_2$
is of order four. By the same argument as above, this leads us to
consider
$$v_2^{-\frac{3}{4}} (1-v_2)^{-\frac{1}{2}} (v_1-v_2)^{-\frac{1}{2}}
  \vmatrix \format \c \quad & \c\\
  \frac{\partial v_2}{\partial w_1} & \frac{\partial v_1}{\partial w_1}\\
  \frac{\partial v_2}{\partial w_2} & \frac{\partial v_1}{\partial w_2}
  \endvmatrix.\tag 3.7$$

  In fact, for
$$u_1=\frac{a_1 w_1+a_2 w_2+a_3}{c_1 w_1+c_2 w_2+c_3}, \quad
  u_2=\frac{b_1 w_1+b_2 w_2+b_3}{c_1 w_1+c_2 w_2+c_3}$$
with $\Delta=\vmatrix \format \c \quad & \c \quad & \c\\
a_1 & a_2 & a_3\\ b_1 & b_2 & b_3\\ c_1 & c_2 & c_3
\endvmatrix$, we have
$$\vmatrix \format \c \quad & \c\\
  \frac{\partial v_1}{\partial w_1} & \frac{\partial v_2}{\partial w_1}\\
  \frac{\partial v_1}{\partial w_2} & \frac{\partial v_2}{\partial w_2}
  \endvmatrix
 =\Delta (c_1 w_1+c_2 w_2+c_3)^{-3}
  \vmatrix \format \c \quad & \c\\
  \frac{\partial v_1}{\partial u_1} & \frac{\partial v_2}{\partial u_1}\\
  \frac{\partial v_1}{\partial u_2} & \frac{\partial v_2}{\partial u_2}
  \endvmatrix.\tag 3.8$$

{\smc Theorem 3.2 (Main Theorem 1)}. {\it For the over-determined
system of second order, nonlinear $($quasilinear$)$ complex
partial differential equations:
$$\left\{\aligned
  \{ w_1, w_2; v_1, v_2 \}_{v_1} &=F_1(v_1, v_2),\\
  \{ w_1, w_2; v_1, v_2 \}_{v_2} &=F_2(v_1, v_2),\\
    [ w_1, w_2; v_1, v_2 ]_{v_1} &=P_1(v_1, v_2),\\
    [ w_1, w_2; v_1, v_2 ]_{v_2} &=P_2(v_1, v_2),
  \endaligned\right.\tag 3.9$$
set the Jacobian
$$\Delta:=\frac{\partial(v_1, v_2)}{\partial(w_1, w_2)}
 =\vmatrix \format \c \quad & \c\\
  \frac{\partial v_1}{\partial w_1} &
  \frac{\partial v_2}{\partial w_1}\\
  \frac{\partial v_1}{\partial w_2} &
  \frac{\partial v_2}{\partial w_2}
  \endvmatrix.\tag 3.10$$
Then $z=\root 3 \of{\Delta}$ satisfies the following linear
differential equations:
$$\left\{\aligned
  \frac{\partial^2 z}{\partial v_1^2}+\frac{1}{3} P_1 \frac{\partial z}
  {\partial v_1}-F_1 \frac{\partial z}{\partial v_2}+\left(\frac{\partial
  F_1}{\partial v_2}-\frac{1}{3} \frac{\partial P_1}{\partial v_1}-
  \frac{2}{9} P_1^2-\frac{2}{3} F_1 P_2\right) z &=0,\\
  \frac{\partial^2 z}{\partial v_1 \partial v_2}-\frac{1}{3} P_2 \frac{\partial z}
  {\partial v_1}-\frac{1}{3} P_1 \frac{\partial z}{\partial v_2}+\left(\frac{1}{3}
  \frac{\partial P_2}{\partial v_1}+\frac{1}{3} \frac{\partial P_1}{\partial v_2}+
  \frac{1}{9} P_1 P_2-F_1 F_2\right) z &=0,\\
  \frac{\partial^2 z}{\partial v_2^2}+\frac{1}{3} P_2 \frac{\partial z}
  {\partial v_2}-F_2 \frac{\partial z}{\partial v_1}+\left(\frac{\partial
  F_2}{\partial v_1}-\frac{1}{3} \frac{\partial P_2}{\partial v_2}-
  \frac{2}{9} P_2^2-\frac{2}{3} F_2 P_1\right) z &=0.
  \endaligned\right.\tag 3.11$$}

{\it Proof}. Note that
$$\left\{\aligned
  1 &=\frac{\partial v_1}{\partial v_1}
     =\frac{\partial v_1}{\partial w_1} \frac{\partial w_1}{\partial v_1}
     +\frac{\partial v_1}{\partial w_2} \frac{\partial w_2}{\partial v_1},\\
  0 &=\frac{\partial v_1}{\partial v_2}
     =\frac{\partial v_1}{\partial w_1} \frac{\partial w_1}{\partial v_2}
     +\frac{\partial v_1}{\partial w_2} \frac{\partial w_2}{\partial v_2},\\
  0 &=\frac{\partial v_2}{\partial v_1}
     =\frac{\partial v_2}{\partial w_1} \frac{\partial w_1}{\partial v_1}
     +\frac{\partial v_2}{\partial w_2} \frac{\partial w_2}{\partial v_1},\\
  1 &=\frac{\partial v_2}{\partial v_2}
     =\frac{\partial v_2}{\partial w_1} \frac{\partial w_1}{\partial v_2}
     +\frac{\partial v_2}{\partial w_2} \frac{\partial w_2}{\partial v_2}.
\endaligned\right.$$
Thus,
$$\frac{\partial v_1}{\partial w_1}=\frac{\frac{\partial
  w_2}{\partial v_2}}{\vmatrix \format \c \quad & \c\\
  \frac{\partial w_1}{\partial v_1} & \frac{\partial w_2}{\partial v_1}\\
  \frac{\partial w_1}{\partial v_2} & \frac{\partial w_2}{\partial v_2}
  \endvmatrix}, \quad
  \frac{\partial v_1}{\partial w_2}=\frac{-\frac{\partial
  w_1}{\partial v_2}}{\vmatrix \format \c \quad & \c\\
  \frac{\partial w_1}{\partial v_1} & \frac{\partial w_2}{\partial v_1}\\
  \frac{\partial w_1}{\partial v_2} & \frac{\partial w_2}{\partial v_2}
  \endvmatrix},$$
$$\frac{\partial v_2}{\partial w_1}=\frac{-\frac{\partial
  w_2}{\partial v_1}}{\vmatrix \format \c \quad & \c\\
  \frac{\partial w_1}{\partial v_1} & \frac{\partial w_2}{\partial v_1}\\
  \frac{\partial w_1}{\partial v_2} & \frac{\partial w_2}{\partial v_2}
  \endvmatrix}, \quad
  \frac{\partial v_2}{\partial w_2}=\frac{\frac{\partial
  w_1}{\partial v_1}}{\vmatrix \format \c \quad & \c\\
  \frac{\partial w_1}{\partial v_1} & \frac{\partial w_2}{\partial v_1}\\
  \frac{\partial w_1}{\partial v_2} & \frac{\partial w_2}{\partial v_2}
  \endvmatrix}.$$
Hence,
$$\vmatrix \format \c \quad & \c\\
  \frac{\partial v_1}{\partial w_1} & \frac{\partial v_2}{\partial w_1}\\
  \frac{\partial v_1}{\partial w_2} & \frac{\partial v_2}{\partial w_2}
  \endvmatrix
 =\frac{1}{\vmatrix \format \c \quad & \c\\
  \frac{\partial w_1}{\partial v_1} & \frac{\partial w_2}{\partial v_1}\\
  \frac{\partial w_1}{\partial v_2} & \frac{\partial w_2}{\partial v_2}
  \endvmatrix}.$$
We have
$$\aligned
 &\frac{\partial}{\partial v_1} \left(\frac{\partial v_2}{\partial w_2}\right)=
  \frac{\partial}{\partial v_1} \left(\frac{\frac{\partial w_1}
  {\partial v_1}}{\vmatrix \format \c \quad & \c\\
  \frac{\partial w_1}{\partial v_1} & \frac{\partial w_2}{\partial v_1}\\
  \frac{\partial w_1}{\partial v_2} & \frac{\partial w_2}{\partial v_2}
  \endvmatrix}\right)\\
=&\vmatrix \format \c \quad & \c\\
  \frac{\partial w_1}{\partial v_1} & \frac{\partial w_2}{\partial v_1}\\
  \frac{\partial w_1}{\partial v_2} & \frac{\partial w_2}{\partial v_2}
  \endvmatrix^{-2}
  \left[-\left(\frac{\partial w_1}{\partial v_1}\right)^{3}
  \frac{\partial}{\partial v_2}\left(\frac{\frac{\partial
  w_2}{\partial v_1}}{\frac{\partial w_1}{\partial v_1}}\right)+
  \frac{\partial w_1}{\partial v_2} \vmatrix \format \c \quad & \c\\
  \frac{\partial w_1}{\partial v_1} & \frac{\partial w_2}{\partial v_1}\\
  \frac{\partial^2 w_1}{\partial v_1^2} &
  \frac{\partial^2 w_2}{\partial v_1^2}
  \endvmatrix\right]\\
=&\left(\frac{\partial v_2}{\partial w_2}\right)^{3} \vmatrix \format \c \quad & \c\\
  \frac{\partial v_1}{\partial w_1} & \frac{\partial v_2}{\partial w_1}\\
  \frac{\partial v_1}{\partial w_2} & \frac{\partial v_2}{\partial w_2}
  \endvmatrix^{-1}
  \frac{\partial}{\partial v_2} \left(\frac{\frac{\partial v_2}{\partial
  w_1}}{\frac{\partial v_2}{\partial w_2}}\right)-\frac{\partial
  v_1}{\partial w_2} F_1(v_1, v_2).
\endaligned$$
Similarly,
$$\frac{\partial}{\partial v_2} \left(\frac{\partial v_1}{\partial w_1}\right)
 =\left(\frac{\partial v_1}{\partial w_1}\right)^{3} \vmatrix \format \c \quad & \c\\
  \frac{\partial v_1}{\partial w_1} & \frac{\partial v_2}{\partial w_1}\\
  \frac{\partial v_1}{\partial w_2} & \frac{\partial v_2}{\partial w_2}
  \endvmatrix^{-1}
  \frac{\partial}{\partial v_1} \left(\frac{\frac{\partial v_1}{\partial
  w_2}}{\frac{\partial v_1}{\partial w_1}}\right)-\frac{\partial
  v_2}{\partial w_1} F_2(v_1, v_2).$$
$$\frac{\partial}{\partial v_1} \left(\frac{\partial v_2}{\partial w_1}\right)
 =-\frac{\partial v_1}{\partial w_1} F_1(v_1, v_2)+\vmatrix \format
  \c \quad & \c\\
  \frac{\partial v_1}{\partial w_1} & \frac{\partial v_2}{\partial w_1}\\
  \frac{\partial v_1}{\partial w_2} & \frac{\partial v_2}{\partial w_2}
  \endvmatrix^{-1}
  \frac{\partial v_2}{\partial w_1} \left(\frac{\partial v_2}{\partial w_2}\right)^2
  \frac{\partial}{\partial v_2} \left(\frac{\frac{\partial v_2}{\partial w_1}}
  {\frac{\partial v_2}{\partial w_2}}\right).$$
$$\frac{\partial}{\partial v_2} \left(\frac{\partial v_1}{\partial w_2}\right)
 =-\frac{\partial v_2}{\partial w_2} F_2(v_1, v_2)+\vmatrix \format
  \c \quad & \c\\
  \frac{\partial v_1}{\partial w_1} & \frac{\partial v_2}{\partial w_1}\\
  \frac{\partial v_1}{\partial w_2} & \frac{\partial v_2}{\partial w_2}
  \endvmatrix^{-1}
  \frac{\partial v_1}{\partial w_2} \left(\frac{\partial v_1}{\partial w_1}\right)^2
  \frac{\partial}{\partial v_1} \left(\frac{\frac{\partial v_1}{\partial w_2}}
  {\frac{\partial v_1}{\partial w_1}}\right).$$
Put
$$\frac{\partial v_1}{\partial w_1}=f_{1, 1}(v_1, v_2), \quad
  \frac{\partial v_1}{\partial w_2}=f_{1, 2}(v_1, v_2), \quad
  \frac{\partial v_2}{\partial w_1}=f_{2, 1}(v_1, v_2), \quad
  \frac{\partial v_2}{\partial w_2}=f_{2, 2}(v_1, v_2).$$
Set $\Delta=f_{1, 1} f_{2, 2}-f_{1, 2} f_{2, 1}$, then
$$\left\{\aligned
 f_{1, 2} F_1 &=\frac{f_{2, 2}^{2}}{\Delta} \frac{\partial f_{2, 1}}
 {\partial v_2}-\frac{f_{2, 1} f_{2, 2}}{\Delta} \frac{\partial f_{2, 2}}
 {\partial v_2}-\frac{\partial f_{2, 2}}{\partial v_1},\\
 f_{2, 1} F_2 &=\frac{f_{1, 1}^{2}}{\Delta} \frac{\partial f_{1, 2}}
 {\partial v_1}-\frac{f_{1, 1} f_{1, 2}}{\Delta} \frac{\partial f_{1, 1}}
 {\partial v_1}-\frac{\partial f_{1, 1}}{\partial v_2},\\
 f_{1, 1} F_1 &=\frac{f_{2, 1} f_{2, 2}}{\Delta} \frac{\partial f_{2, 1}}
 {\partial v_2}-\frac{f_{2, 1}^{2}}{\Delta} \frac{\partial f_{2, 2}}
 {\partial v_2}-\frac{\partial f_{2, 1}}{\partial v_1},\\
 f_{2, 2} F_2 &=\frac{f_{1, 1} f_{1, 2}}{\Delta} \frac{\partial f_{1, 2}}
 {\partial v_1}-\frac{f_{1, 2}^{2}}{\Delta} \frac{\partial f_{1, 1}}
 {\partial v_1}-\frac{\partial f_{1, 2}}{\partial v_2},
\endaligned\right.$$
which are equivalent to
$$\left\{\aligned
  \Delta F_1 &=f_{2, 1} \frac{\partial f_{2, 2}}{\partial v_1}-f_{2,
  2} \frac{\partial f_{2, 1}}{\partial v_1},\\
  \Delta F_2 &=f_{1, 2} \frac{\partial f_{1, 1}}{\partial v_2}-f_{1,
  1} \frac{\partial f_{1, 2}}{\partial v_2},\\
  f_{1, 1} \frac{\partial f_{2, 2}}{\partial v_1}-f_{1, 2} \frac{\partial
  f_{2, 1}}{\partial v_1} &=f_{2, 2} \frac{\partial f_{2, 1}}{\partial
  v_2}-f_{2, 1} \frac{\partial f_{2, 2}}{\partial v_2},\\
  f_{1, 1} \frac{\partial f_{1, 2}}{\partial v_1}-f_{1, 2} \frac{\partial
  f_{1, 1}}{\partial v_1} &=f_{2, 2} \frac{\partial f_{1, 1}}{\partial
  v_2}-f_{2, 1} \frac{\partial f_{1, 2}}{\partial v_2}.
\endaligned\right.$$

  We have
$$\aligned
 &\{w_1, w_2; v_1, v_2\}_{v_1}=\frac{\partial v_2}{\partial w_2}
  \frac{\partial}{\partial v_1} \left(\frac{\partial w_2}{\partial v_1}\right)+
  \frac{\partial v_2}{\partial w_1} \frac{\partial}{\partial v_1} \left(
  \frac{\partial w_1}{\partial v_1}\right)\\
=&f_{2, 2} \frac{\partial}{\partial v_1} \left(-\frac{f_{2, 1}}
  {\Delta}\right)+f_{2, 1} \frac{\partial}{\partial v_1} \left(\frac{f_{2, 2}}
  {\Delta}\right)
 =\frac{1}{\Delta} \left(f_{2, 1} \frac{\partial f_{2, 2}}{\partial v_1}-
  f_{2, 2} \frac{\partial f_{2, 1}}{\partial v_1}\right).
\endaligned$$
Similarly,
$$\{w_1, w_2; v_1, v_2\}_{v_2}=\frac{1}{\Delta} \left(f_{1, 2} \frac{\partial
  f_{1, 1}}{\partial v_2}-f_{1, 1} \frac{\partial f_{1, 2}}{\partial v_2}\right).$$
$$[w_1, w_2; v_1, v_2]_{v_1}
 =\frac{1}{\Delta} \left[ \left(f_{1, 2} \frac{\partial f_{2, 1}}{\partial v_1}-
  f_{1, 1} \frac{\partial f_{2, 2}}{\partial v_1}\right)+2 \left(f_{2, 2}
  \frac{\partial f_{1, 1}}{\partial v_1}-f_{2, 1} \frac{\partial
  f_{1, 2}}{\partial v_1}\right)-\frac{\partial \Delta}{\partial v_1} \right].$$
$$[w_1, w_2; v_1, v_2]_{v_2}
 =\frac{1}{\Delta} \left[ \left(f_{2, 1} \frac{\partial f_{1, 2}}{\partial v_2}-
  f_{2, 2} \frac{\partial f_{1, 1}}{\partial v_2}\right)+2 \left(f_{1, 1}
  \frac{\partial f_{2, 2}}{\partial v_2}-f_{1, 2} \frac{\partial
  f_{2, 1}}{\partial v_2}\right)-\frac{\partial \Delta}{\partial v_2} \right].$$
By
$$\frac{\partial \Delta}{\partial v_1}=\frac{\partial}{\partial v_1}
  \left(f_{1, 1} f_{2, 2}-f_{1, 2} f_{2, 1}\right)
 =\frac{\partial f_{1, 1}}{\partial v_1} f_{2, 2}+
  f_{1, 1} \frac{\partial f_{2, 2}}{\partial v_1}-
  \frac{\partial f_{1, 2}}{\partial v_1} f_{2, 1}-
  f_{1, 2} \frac{\partial f_{2, 1}}{\partial v_1},$$
$$\frac{\partial \Delta}{\partial v_2}=\frac{\partial}{\partial v_2}
  \left(f_{1, 1} f_{2, 2}-f_{1, 2} f_{2, 1}\right)
 =\frac{\partial f_{1, 1}}{\partial v_2} f_{2, 2}+
  f_{1, 1} \frac{\partial f_{2, 2}}{\partial v_2}-
  \frac{\partial f_{1, 2}}{\partial v_2} f_{2, 1}-
  f_{1, 2} \frac{\partial f_{2, 1}}{\partial v_2},$$
we get
$$[w_1, w_2; v_1, v_2]_{v_1}
 =\frac{1}{\Delta} \left[2 \left(f_{1, 2} \frac{\partial f_{2, 1}}
  {\partial v_1}-f_{1, 1} \frac{\partial f_{2, 2}}{\partial v_1}\right)+
  \left(f_{2, 2} \frac{\partial f_{1, 1}}{\partial v_1}-f_{2, 1}
  \frac{\partial f_{1, 2}}{\partial v_1}\right) \right],$$
$$[w_1, w_2; v_1, v_2]_{v_2}
 =\frac{1}{\Delta} \left[2 \left(f_{2, 1} \frac{\partial f_{1, 2}}{\partial
  v_2}-f_{2, 2} \frac{\partial f_{1, 1}}{\partial v_2}\right)+\left(f_{1, 1}
  \frac{\partial f_{2, 2}}{\partial v_2}-f_{1, 2} \frac{\partial
  f_{2, 1}}{\partial v_2}\right) \right].$$
The equations
$$\left\{\aligned
  \{ w_1, w_2; v_1, v_2 \}_{v_1} &=F_1(v_1, v_2),\\
  \{ w_1, w_2; v_1, v_2 \}_{v_2} &=F_2(v_1, v_2),\\
    [ w_1, w_2; v_1, v_2 ]_{v_1} &=P_1(v_1, v_2),\\
    [ w_1, w_2; v_1, v_2 ]_{v_2} &=P_2(v_1, v_2)
  \endaligned\right.$$
give that
$$\left\{\aligned
 \Delta F_1 &=f_{2, 1} \frac{\partial f_{2, 2}}{\partial v_1}-
              f_{2, 2} \frac{\partial f_{2, 1}}{\partial v_1},\\
 \Delta F_2 &=f_{1, 2} \frac{\partial f_{1, 1}}{\partial v_2}-
              f_{1, 1} \frac{\partial f_{1, 2}}{\partial v_2},\\
 \Delta P_1 &=2 \left(f_{1, 2} \frac{\partial f_{2, 1}}{\partial v_1}-
              f_{1, 1} \frac{\partial f_{2, 2}}{\partial v_1}\right)+
              \left(f_{2, 2} \frac{\partial f_{1, 1}}{\partial v_1}-
              f_{2, 1} \frac{\partial f_{1, 2}}{\partial v_1}\right),\\
 \Delta P_2 &=2 \left(f_{2, 1} \frac{\partial f_{1, 2}}{\partial v_2}-
              f_{2, 2} \frac{\partial f_{1, 1}}{\partial v_2}\right)+
              \left(f_{1, 1} \frac{\partial f_{2, 2}}{\partial v_2}-
              f_{1, 2} \frac{\partial f_{2, 1}}{\partial v_2}\right).
\endaligned\right.$$

  By
$$\left\{\aligned
  \frac{\partial \Delta}{\partial v_1}
=&\left(f_{2, 2} \frac{\partial f_{1, 1}}{\partial v_1}
  -f_{2, 1} \frac{\partial f_{1, 2}}{\partial v_1}\right)-
  \left(f_{1, 2} \frac{\partial f_{2, 1}}{\partial v_1}
  -f_{1, 1} \frac{\partial f_{2, 2}}{\partial v_1}\right),\\
  \Delta P_1
=&\left(f_{2, 2} \frac{\partial f_{1, 1}}{\partial v_1}
  -f_{2, 1} \frac{\partial f_{1, 2}}{\partial v_1}\right)+2
  \left(f_{1, 2} \frac{\partial f_{2, 1}}{\partial v_1}
  -f_{1, 1} \frac{\partial f_{2, 2}}{\partial v_1}\right),
\endaligned\right.$$
we have
$$\left\{\aligned
  f_{2, 2} \frac{\partial f_{1, 1}}{\partial v_1}
 -f_{2, 1} \frac{\partial f_{1, 2}}{\partial v_1}
&=\frac{2}{3} \frac{\partial \Delta}{\partial v_1}
 +\frac{1}{3} \Delta P_1,\\
  f_{1, 2} \frac{\partial f_{2, 1}}{\partial v_1}
 -f_{1, 1} \frac{\partial f_{2, 2}}{\partial v_1}
&=\frac{1}{3} \Delta P_1-\frac{1}{3} \frac{\partial
  \Delta}{\partial v_1}.
\endaligned\right.$$

  By
$$\left\{\aligned
  \frac{\partial \Delta}{\partial v_2}
=&\left(f_{1, 1} \frac{\partial f_{2, 2}}{\partial v_2}
  -f_{1, 2} \frac{\partial f_{2, 1}}{\partial v_2}\right)-
  \left(f_{2, 1} \frac{\partial f_{1, 2}}{\partial v_2}
  -f_{2, 2} \frac{\partial f_{1, 1}}{\partial v_2}\right),\\
  \Delta P_2
=&\left(f_{1, 1} \frac{\partial f_{2, 2}}{\partial v_2}
  -f_{1, 2} \frac{\partial f_{2, 1}}{\partial v_2}\right)+2
  \left(f_{2, 1} \frac{\partial f_{1, 2}}{\partial v_2}
  -f_{2, 2} \frac{\partial f_{1, 1}}{\partial v_2}\right),
\endaligned\right.$$
we have
$$\left\{\aligned
  f_{1, 1} \frac{\partial f_{2, 2}}{\partial v_2}
 -f_{1, 2} \frac{\partial f_{2, 1}}{\partial v_2}
&=\frac{2}{3} \frac{\partial \Delta}{\partial v_2}
 +\frac{1}{3} \Delta P_2,\\
  f_{2, 1} \frac{\partial f_{1, 2}}{\partial v_2}
 -f_{2, 2} \frac{\partial f_{1, 1}}{\partial v_2}
&=\frac{1}{3} \Delta P_2-\frac{1}{3} \frac{\partial
  \Delta}{\partial v_2}.
\endaligned\right.$$

  The above two equations
$$\left\{\aligned
  f_{1, 1} \frac{\partial f_{2, 2}}{\partial v_1}-f_{1, 2} \frac{\partial
  f_{2, 1}}{\partial v_1} &=f_{2, 2} \frac{\partial f_{2, 1}}{\partial
  v_2}-f_{2, 1} \frac{\partial f_{2, 2}}{\partial v_2},\\
  f_{1, 1} \frac{\partial f_{1, 2}}{\partial v_1}-f_{1, 2} \frac{\partial
  f_{1, 1}}{\partial v_1} &=f_{2, 2} \frac{\partial f_{1, 1}}{\partial
  v_2}-f_{2, 1} \frac{\partial f_{1, 2}}{\partial v_2}
\endaligned\right.$$
imply that
$$\left\{\aligned
  f_{2, 1} \frac{\partial f_{2, 2}}{\partial v_2}-
  f_{2, 2} \frac{\partial f_{2, 1}}{\partial v_2}
&=\frac{1}{3} \Delta P_1-\frac{1}{3} \frac{\partial
  \Delta}{\partial v_1},\\
  f_{1, 2} \frac{\partial f_{1, 1}}{\partial v_1}-
  f_{1, 1} \frac{\partial f_{1, 2}}{\partial v_1}
&=\frac{1}{3} \Delta P_2-\frac{1}{3} \frac{\partial
  \Delta}{\partial v_2}.
\endaligned\right.$$

  By
$$\left\{\aligned
  f_{1, 1} \frac{\partial f_{2, 2}}{\partial v_2}-
  f_{1, 2} \frac{\partial f_{2, 1}}{\partial v_2}
=&\frac{2}{3} \frac{\partial \Delta}{\partial v_2}+
  \frac{1}{3} \Delta P_2,\\
  f_{2, 1} \frac{\partial f_{2, 2}}{\partial v_2}-
  f_{2, 2} \frac{\partial f_{2, 1}}{\partial v_2}
=&-\frac{1}{3} \frac{\partial \Delta}{\partial v_1}
  +\frac{1}{3} \Delta P_1,
\endaligned\right.$$
we have
$$\left\{\aligned
  \frac{\partial f_{2, 2}}{\partial v_2}
=&\frac{1}{3} \frac{f_{1, 2}}{\Delta} \frac{\partial \Delta}
 {\partial v_1}+\frac{2}{3} \frac{f_{2, 2}}{\Delta}
  \frac{\partial \Delta}{\partial v_2}-\frac{1}{3} f_{1, 2} P_1+
  \frac{1}{3} f_{2, 2} P_2,\\
  \frac{\partial f_{2, 1}}{\partial v_2}
=&\frac{1}{3} \frac{f_{1, 1}}{\Delta} \frac{\partial \Delta}
 {\partial v_1}+\frac{2}{3} \frac{f_{2, 1}}{\Delta}
  \frac{\partial \Delta}{\partial v_2}-\frac{1}{3} f_{1, 1} P_1+
  \frac{1}{3} f_{2, 1} P_2.
\endaligned\right.$$
By
$$\left\{\aligned
  f_{2, 1} \frac{\partial f_{2, 2}}{\partial v_1}-
  f_{2, 2} \frac{\partial f_{2, 1}}{\partial v_1} &=\Delta F_1,\\
  f_{1, 1} \frac{\partial f_{2, 2}}{\partial v_1}-
  f_{1, 2} \frac{\partial f_{2, 1}}{\partial v_1} &=
  \frac{1}{3} \frac{\partial \Delta}{\partial v_1}-\frac{1}{3}
  \Delta P_1,
\endaligned\right.$$
we have
$$\left\{\aligned
  \frac{\partial f_{2, 2}}{\partial v_1}
&=\frac{1}{3} \frac{f_{2, 2}}{\Delta} \frac{\partial \Delta}
 {\partial v_1}-f_{1, 2} F_1-\frac{1}{3} f_{2, 2} P_1,\\
  \frac{\partial f_{2, 1}}{\partial v_1}
&=\frac{1}{3} \frac{f_{2, 1}}{\Delta} \frac{\partial \Delta}
 {\partial v_1}-f_{1, 1} F_1-\frac{1}{3} f_{2, 1} P_1.
\endaligned\right.$$
Differentiating the equation
$$\frac{1}{3} \left(\frac{\partial \Delta}{\partial v_1}-\Delta P_1\right)
 =f_{2, 2} \frac{\partial f_{2, 1}}{\partial v_2}-
  f_{2, 1} \frac{\partial f_{2, 2}}{\partial v_2}$$
with respect to $v_1$, we get
$$\frac{1}{3} \left(\frac{\partial^2 \Delta}{\partial v_1^2}-
  \frac{\partial \Delta}{\partial v_1} P_1-\Delta
  \frac{\partial P_1}{\partial v_1}\right)
 =\frac{\partial f_{2, 2}}{\partial v_1} \frac{\partial f_{2, 1}}
  {\partial v_2}-\frac{\partial f_{2, 1}}{\partial v_1}
  \frac{\partial f_{2, 2}}{\partial v_2}+f_{2, 2}
  \frac{\partial^2 f_{2, 1}}{\partial v_1 \partial v_2}-f_{2, 1}
  \frac{\partial^2 f_{2, 2}}{\partial v_1 \partial v_2}.$$
Differentiating the equation
$$\Delta F_1=f_{2, 1} \frac{\partial f_{2, 2}}{\partial v_1}-
             f_{2, 2} \frac{\partial f_{2, 1}}{\partial v_1}$$
with respect to $v_2$, we get
$$\frac{\partial \Delta}{\partial v_2} F_1+\Delta \frac{\partial
  F_1}{\partial v_2}=\frac{\partial f_{2, 1}}{\partial v_2}
  \frac{\partial f_{2, 2}}{\partial v_1}-\frac{\partial f_{2, 2}}
  {\partial v_2} \frac{\partial f_{2, 1}}{\partial v_1}+f_{2, 1}
  \frac{\partial^2 f_{2, 2}}{\partial v_1 \partial v_2}-f_{2, 2}
  \frac{\partial^2 f_{2, 1}}{\partial v_1 \partial v_2}.$$
Combine the above two equalities, we find that
$$\frac{1}{3} \left(\frac{\partial^2 \Delta}{\partial v_1^2}-
  \frac{\partial \Delta}{\partial v_1} P_1-\Delta \frac{\partial
  P_1}{\partial v_1}\right)+\frac{\partial \Delta}{\partial v_2}
  F_1+\Delta \frac{\partial F_1}{\partial v_2}
 =2 \left(\frac{\partial f_{2, 2}}{\partial v_1} \frac{\partial
  f_{2, 1}}{\partial v_2}-\frac{\partial f_{2, 1}}{\partial v_1}
  \frac{\partial f_{2, 2}}{\partial v_2}\right).$$
A straightforward calculation gives that
$$\frac{\partial f_{2, 2}}{\partial v_1} \frac{\partial f_{2, 1}}
  {\partial v_2}-\frac{\partial f_{2, 1}}{\partial v_1}
  \frac{\partial f_{2, 2}}{\partial v_2}
 =\frac{1}{9} \frac{1}{\Delta} \left(\frac{\partial \Delta}{\partial
  v_1}\right)^2-\frac{2}{9} P_1 \frac{\partial \Delta}{\partial v_1}
  +\frac{2}{3} F_1 \frac{\partial \Delta}{\partial v_2}+\frac{1}{9}
  \Delta P_1^2+\frac{1}{3} \Delta F_1 P_2.$$
This gives that the following equality:
$$\frac{1}{3} \frac{\partial^2 \Delta}{\partial v_1^2}-\frac{2}{9}
  \frac{1}{\Delta} \left(\frac{\partial \Delta}{\partial v_1}\right)^2
  +\frac{1}{9} P_1 \frac{\partial \Delta}{\partial v_1}-\frac{1}{3}
  F_1 \frac{\partial \Delta}{\partial v_2}+\left(\frac{\partial
  F_1}{\partial v_2}-\frac{1}{3} \frac{\partial P_1}{\partial v_1}-
  \frac{2}{9} P_1^2-\frac{2}{3} F_1 P_2\right) \Delta=0,$$
which is equivalent to
$$\frac{\partial^2 (\root 3 \of{\Delta})}{\partial v_1^2}+\frac{1}{3} P_1
  \frac{\partial (\root 3 \of{\Delta})}{\partial v_1}-F_1 \frac{\partial
  (\root 3 \of{\Delta})}{\partial v_2}+\left(\frac{\partial F_1}{\partial v_2}-\frac{1}{3}
  \frac{\partial P_1}{\partial v_1}-\frac{2}{9} P_1^2-\frac{2}{3} F_1 P_2
  \right) \root 3 \of{\Delta}=0.$$

  By the symmetry property, we get the other equation:
$$\frac{\partial^2 (\root 3 \of{\Delta})}{\partial v_2^2}+\frac{1}{3} P_2
  \frac{\partial (\root 3 \of{\Delta})}{\partial v_2}-F_2 \frac{\partial
  (\root 3 \of{\Delta})}{\partial v_1}+\left(\frac{\partial F_2}{\partial v_1}-\frac{1}{3}
  \frac{\partial P_2}{\partial v_2}-\frac{2}{9} P_2^2-\frac{2}{3} F_2 P_1
  \right) \root 3 \of{\Delta}=0.$$
Here,
$$\left\{\aligned
  \frac{\partial f_{1, 1}}{\partial v_1}
=&\frac{1}{3} \frac{f_{2, 1}}{\Delta} \frac{\partial \Delta}
 {\partial v_2}+\frac{2}{3} \frac{f_{1, 1}}{\Delta}
  \frac{\partial \Delta}{\partial v_1}-\frac{1}{3} f_{2, 1} P_2+
  \frac{1}{3} f_{1, 1} P_1,\\
  \frac{\partial f_{1, 2}}{\partial v_1}
=&\frac{1}{3} \frac{f_{2, 2}}{\Delta} \frac{\partial \Delta}
 {\partial v_2}+\frac{2}{3} \frac{f_{1, 2}}{\Delta}
  \frac{\partial \Delta}{\partial v_1}-\frac{1}{3} f_{2, 2} P_2+
  \frac{1}{3} f_{1, 2} P_1.
\endaligned\right.$$
$$\left\{\aligned
  \frac{\partial f_{1, 1}}{\partial v_2}
&=\frac{1}{3} \frac{f_{1, 1}}{\Delta} \frac{\partial \Delta}
 {\partial v_2}-f_{2, 1} F_2-\frac{1}{3} f_{1, 1} P_2,\\
  \frac{\partial f_{1, 2}}{\partial v_2}
&=\frac{1}{3} \frac{f_{1, 2}}{\Delta} \frac{\partial \Delta}
 {\partial v_2}-f_{2, 2} F_2-\frac{1}{3} f_{1, 2} P_2.
\endaligned\right.$$

  Differentiating the equation
$$f_{2, 2} \frac{\partial f_{2, 1}}{\partial v_2}-
  f_{2, 1} \frac{\partial f_{2, 2}}{\partial v_2}=
  \frac{1}{3} \left(\frac{\partial \Delta}{\partial v_1}-
  \Delta P_1\right)$$
with respect to $v_2$, we get
$$f_{2, 2} \frac{\partial^2 f_{2, 1}}{\partial v_2^2}-
  f_{2, 1} \frac{\partial^2 f_{2, 2}}{\partial v_2^2}=
  \frac{1}{3} \left(\frac{\partial^2 \Delta}{\partial v_1
  \partial v_2}-\frac{\partial \Delta}{\partial v_2} P_1-
  \Delta \frac{\partial P_1}{\partial v_2}\right).$$
Differentiating the equation
$$f_{1, 1} \frac{\partial f_{2, 2}}{\partial v_2}-
  f_{1, 2} \frac{\partial f_{2, 1}}{\partial v_2}=
  \frac{2}{3} \frac{\partial \Delta}{\partial v_2}+
  \frac{1}{3} \Delta P_2$$
with respect to $v_1$, we get
$$\frac{\partial f_{1, 1}}{\partial v_1} \frac{\partial f_{2, 2}}
 {\partial v_2}-\frac{\partial f_{1, 2}}{\partial v_1} \frac{\partial
  f_{2, 1}}{\partial v_2}+f_{1, 1} \frac{\partial^2 f_{2, 2}}{\partial
  v_1 \partial v_2}-f_{1, 2} \frac{\partial^2 f_{2, 1}}{\partial
  v_1 \partial v_2}
 =\frac{2}{3} \frac{\partial^2 \Delta}{\partial v_1 \partial v_2}+
  \frac{1}{3} \frac{\partial \Delta}{\partial v_1} P_2+
  \frac{1}{3} \Delta \frac{\partial P_2}{\partial v_1}.$$
Differentiating the equation
$$f_{1, 1} \frac{\partial f_{2, 2}}{\partial v_1}-
  f_{1, 2} \frac{\partial f_{2, 1}}{\partial v_1}=
  f_{2, 2} \frac{\partial f_{2, 1}}{\partial v_2}-
  f_{2, 1} \frac{\partial f_{2, 2}}{\partial v_2}$$
with respect to $v_2$, we get
$$f_{2, 2} \frac{\partial^2 f_{2, 1}}{\partial v_2^2}-
  f_{2, 1} \frac{\partial^2 f_{2, 2}}{\partial v_2^2}=
  \frac{\partial f_{1, 1}}{\partial v_2} \frac{\partial f_{2, 2}}
  {\partial v_1}-\frac{\partial f_{1, 2}}{\partial v_2}
  \frac{\partial f_{2, 1}}{\partial v_1}+f_{1, 1} \frac{\partial^2
  f_{2, 2}}{\partial v_1 \partial v_2}-f_{1, 2} \frac{\partial^2
  f_{2, 1}}{\partial v_1 \partial v_2}.$$

  Combining the above three equalities, we find that
$$\aligned
 &\frac{1}{3} \frac{\partial^2 \Delta}{\partial v_1 \partial v_2}+
  \frac{1}{3} \frac{\partial \Delta}{\partial v_1} P_2+
  \frac{1}{3} \Delta \frac{\partial P_2}{\partial v_1}+
  \frac{1}{3} \frac{\partial \Delta}{\partial v_2} P_1+
  \frac{1}{3} \Delta \frac{\partial P_1}{\partial v_2}\\
=&\frac{\partial f_{1, 2}}{\partial v_2} \frac{\partial f_{2, 1}}
  {\partial v_1}-\frac{\partial f_{1, 1}}{\partial v_2}
  \frac{\partial f_{2, 2}}{\partial v_1}+\frac{\partial f_{1, 1}}
  {\partial v_1} \frac{\partial f_{2, 2}}{\partial v_2}-
  \frac{\partial f_{1, 2}}{\partial v_1} \frac{\partial f_{2, 1}}
  {\partial v_2}.
\endaligned$$
A straightforward calculation gives that the right hand side of
the above equation is equal to
$$\frac{2}{9} \frac{1}{\Delta} \frac{\partial \Delta}{\partial v_1}
  \frac{\partial \Delta}{\partial v_2}+\frac{4}{9} P_2
  \frac{\partial \Delta}{\partial v_1}+\frac{4}{9} P_1
  \frac{\partial \Delta}{\partial v_2}+\Delta F_1 F_2-
  \frac{1}{9} \Delta P_1 P_2.$$
This implies that
$$\aligned
 &\frac{1}{3} \frac{\partial^2 \Delta}{\partial v_1 \partial v_2}-
  \frac{2}{9} \frac{1}{\Delta} \frac{\partial \Delta}{\partial v_1}
  \frac{\partial \Delta}{\partial v_2}-\frac{1}{9} P_2 \frac{\partial
  \Delta}{\partial v_1}-\frac{1}{9} P_1 \frac{\partial \Delta}{\partial
  v_2}+\\
 &+\left(\frac{1}{3} \frac{\partial P_2}{\partial v_1}+\frac{1}{3}
  \frac{\partial P_1}{\partial v_2}+\frac{1}{9} P_1 P_2-F_1 F_2\right)
  \Delta=0,
\endaligned$$
which is equivalent to
$$\frac{\partial^2 (\root 3 \of{\Delta})}{\partial v_1 \partial v_2}-
  \frac{1}{3} P_2 \frac{\partial (\root 3 \of{\Delta})}{\partial v_1}-
  \frac{1}{3} P_1 \frac{\partial (\root 3 \of{\Delta})}{\partial v_2}+
  \left(\frac{1}{3} \frac{\partial P_2}{\partial v_1}+\frac{1}{3}
  \frac{\partial P_1}{\partial v_2}+\frac{1}{9} P_1 P_2-F_1 F_2\right)
  \root 3 \of{\Delta}=0.$$
$\qquad \qquad \qquad \qquad \qquad \qquad \qquad \qquad \qquad
 \qquad \qquad \qquad \qquad \qquad \qquad \qquad \qquad \qquad
 \quad \boxed{}$

  Now, we will study the equations:
$$\left\{\aligned
  \{ w_1, w_2; v_1, v_2 \}_{v_1} &=F_1(v_1, v_2),\\
  \{ w_1, w_2; v_1, v_2 \}_{v_2} &=F_2(v_1, v_2),\\
    [ w_1, w_2; v_1, v_2 ]_{v_1} &=P_1(v_1, v_2),\\
    [ w_1, w_2; v_1, v_2 ]_{v_2} &=P_2(v_1, v_2).
\endaligned\right.$$
Here $F_1(v_1, v_2)$ is finite except at $v_1=0$, $1$, $v_2$,
$\infty$, where respectively, the products
$$v_1^{-\frac{3}{4}} \vmatrix \format \c \quad & \c\\
  \frac{\partial v_1}{\partial w_1} & \frac{\partial v_2}{\partial w_1}\\
  \frac{\partial v_1}{\partial w_2} & \frac{\partial v_2}{\partial w_2}
  \endvmatrix,
  (1-v_1)^{-\frac{1}{2}} \vmatrix \format \c \quad & \c\\
  \frac{\partial v_1}{\partial w_1} & \frac{\partial v_2}{\partial w_1}\\
  \frac{\partial v_1}{\partial w_2} & \frac{\partial v_2}{\partial w_2}
  \endvmatrix,
  (v_2-v_1)^{-\frac{1}{2}} \vmatrix \format \c \quad & \c\\
  \frac{\partial v_1}{\partial w_1} & \frac{\partial v_2}{\partial w_1}\\
  \frac{\partial v_1}{\partial w_2} & \frac{\partial v_2}{\partial w_2}
  \endvmatrix,
  v_1^{-1} \vmatrix \format \c \quad & \c\\
  \frac{\partial v_1}{\partial w_1} & \frac{\partial v_2}{\partial w_1}\\
  \frac{\partial v_1}{\partial w_2} & \frac{\partial v_2}{\partial w_2}
  \endvmatrix$$
are finite and non-zero. Similarly, $F_2(v_1, v_2)$ is finite
except at $v_2=0$, $1$, $v_1$, $\infty$, where respectively, the
products
$$v_2^{-\frac{3}{4}} \vmatrix \format \c \quad & \c\\
  \frac{\partial v_2}{\partial w_1} & \frac{\partial v_1}{\partial w_1}\\
  \frac{\partial v_2}{\partial w_2} & \frac{\partial v_1}{\partial w_2}
  \endvmatrix,
  (1-v_2)^{-\frac{1}{2}} \vmatrix \format \c \quad & \c\\
  \frac{\partial v_2}{\partial w_1} & \frac{\partial v_1}{\partial w_1}\\
  \frac{\partial v_2}{\partial w_2} & \frac{\partial v_1}{\partial w_2}
  \endvmatrix,
  (v_1-v_2)^{-\frac{1}{2}} \vmatrix \format \c \quad & \c\\
  \frac{\partial v_2}{\partial w_1} & \frac{\partial v_1}{\partial w_1}\\
  \frac{\partial v_2}{\partial w_2} & \frac{\partial v_1}{\partial w_2}
  \endvmatrix,
  v_2^{-1} \vmatrix \format \c \quad & \c\\
  \frac{\partial v_2}{\partial w_1} & \frac{\partial v_1}{\partial w_1}\\
  \frac{\partial v_2}{\partial w_2} & \frac{\partial v_1}{\partial w_2}
  \endvmatrix$$
are finite and non-zero.

  Set $v_1^{-\frac{1}{4}} \root 3 \of{\Delta}=y$, which is finite and
non-zero at $v_1=0$. Then
$$\aligned
 &v_1^2 \frac{\partial^2 y}{\partial v_1^2}+\left(\frac{1}{3} P_1
  v_1^2+\frac{1}{2} v_1\right) \frac{\partial y}{\partial v_1}-
  F_1 v_1^2 \frac{\partial y}{\partial v_2}+\\
 &+\left[ \left(\frac{\partial F_1}{\partial v_2}-\frac{1}{3}
  \frac{\partial P_1}{\partial v_1}-\frac{2}{9} P_1^2-\frac{2}{3}
  F_1 P_2 \right) v_1^2+\frac{1}{12} P_1 v_1-\frac{3}{16}\right] y=0
\endaligned$$
and
$$\aligned
 &v_1 \frac{\partial^2 y}{\partial v_1 \partial v_2}-\frac{1}{3}
  P_2 v_1 \frac{\partial y}{\partial v_1}+\left(-\frac{1}{3} P_1 v_1+
  \frac{1}{4}\right) \frac{\partial y}{\partial v_2}+\\
 &+\left[ \left(\frac{1}{3} \frac{\partial P_2}{\partial v_1}+
  \frac{1}{3} \frac{\partial P_1}{\partial v_2}+\frac{1}{9} P_1
  P_2-F_1 F_2\right) v_1-\frac{1}{12} P_2 \right] y=0
\endaligned$$
at $v_1=0$.

   Set $(1-v_1)^{-\frac{1}{6}} \root 3 \of{\Delta}=y$, which is finite and
non-zero at $v_1=1$. Then
$$\aligned
 &(1-v_1)^2 \frac{\partial^2 y}{\partial v_1^2}+\left[\frac{1}{3} P_1
  (1-v_1)^2-\frac{1}{3} (1-v_1)\right] \frac{\partial y}{\partial v_1}
  -F_1 (1-v_1)^2 \frac{\partial y}{\partial v_2}+\\
 &+\left[ \left(\frac{\partial F_1}{\partial v_2}-\frac{1}{3}
  \frac{\partial P_1}{\partial v_1}-\frac{2}{9} P_1^2-\frac{2}{3} F_1
  P_2\right) (1-v_1)^2-\frac{1}{18} P_1 (1-v_1)-\frac{5}{36} \right] y=0
\endaligned$$
and
$$\aligned
 &(1-v_1) \frac{\partial^2 y}{\partial v_1 \partial v_2}-\frac{1}{3}
  P_2 (1-v_1) \frac{\partial y}{\partial v_1}+\left[-\frac{1}{3} P_1
  (1-v_1)-\frac{1}{6}\right] \frac{\partial y}{\partial v_2}+\\
 &+\left[ \left(\frac{1}{3} \frac{\partial P_2}{\partial v_1}+\frac{1}{3}
  \frac{\partial P_1}{\partial v_2}+\frac{1}{9} P_1 P_2-F_1 F_2\right)
  (1-v_1)+\frac{1}{18} P_2 \right] y=0
\endaligned$$
at $v_1=1$.

   Set $(v_2-v_1)^{-\frac{1}{6}} \root 3 \of{\Delta}=y$, which is finite and
non-zero at $v_1=v_2$. Then
$$\aligned
 &(v_2-v_1)^2 \frac{\partial^2 y}{\partial v_1^2}+\left[\frac{1}{3} P_1
  (v_2-v_1)^2-\frac{1}{3} (v_2-v_1)\right] \frac{\partial y}{\partial v_1}
  -F_1 (v_2-v_1)^2 \frac{\partial y}{\partial v_2}+\\
 &+\left[\left(\frac{\partial F_1}{\partial v_2}-\frac{1}{3} \frac{\partial P_1}
  {\partial v_1}-\frac{2}{9} P_1^2-\frac{2}{3} F_1 P_2\right) (v_2-v_1)^2-
  \frac{1}{18} P_1 (v_2-v_1)+\right.\\
 &\left.-\frac{1}{6} F_1 (v_2-v_1)-\frac{5}{36}\right] y=0
\endaligned$$
at $v_1=v_2$.

   Set $v_1^{-\frac{1}{3}} \root 3 \of{\Delta}=y$, which is finite and non-zero at
$v_1=\infty$. Then
$$\aligned
 &v_1^2 \frac{\partial^2 y}{\partial v_1^2}+\left(\frac{1}{3} P_1
  v_1^2+\frac{2}{3} v_1\right) \frac{\partial y}{\partial v_1}-F_1
  v_1^2 \frac{\partial y}{\partial v_2}+\\
 &+\left[ \left(\frac{\partial F_1}{\partial v_2}-\frac{1}{3}
  \frac{\partial P_1}{\partial v_1}-\frac{2}{9} P_1^2-\frac{2}{3}
  F_1 P_2\right) v_1^2+\frac{1}{9} P_1 v_1-\frac{2}{9} \right] y=0
\endaligned$$
and
$$\aligned
 &v_1 \frac{\partial^2 y}{\partial v_1 \partial v_2}-\frac{1}{3} P_2
  v_1 \frac{\partial y}{\partial v_1}+\left(-\frac{1}{3} P_1 v_1+
  \frac{1}{3}\right) \frac{\partial y}{\partial v_2}+\\
 &+\left[ \left(\frac{1}{3} \frac{\partial P_2}{\partial v_1}+\frac{1}{3}
  \frac{\partial P_1}{\partial v_2}+\frac{1}{9} P_1 P_2-F_1 F_2\right)
  v_1-\frac{1}{9} P_2 \right] y=0
\endaligned$$
at $v_1=\infty$. Put $u_1=\frac{1}{v_1}$, then
$$\aligned
 &u_1^2 \frac{\partial^2 y}{\partial u_1^2}+\left(\frac{4}{3} u_1-
  \frac{1}{3} P_1\right) \frac{\partial y}{\partial u_1}-\frac{F_1}
  {u_1^2} \frac{\partial y}{\partial v_2}+\\
 &+\left[ \left(\frac{\partial F_1}{\partial v_2}+\frac{1}{3} u_1^2
  \frac{\partial P_1}{\partial u_1}-\frac{2}{9} P_1^2-\frac{2}{3}
  F_1 P_2\right) \frac{1}{u_1^2}+\frac{1}{9} \frac{P_1}{u_1}-
  \frac{2}{9}\right] y=0
\endaligned$$
and
$$\aligned
 &-u_1 \frac{\partial^2 y}{\partial u_1 \partial v_2}+\frac{1}{3} P_2
  u_1 \frac{\partial y}{\partial u_1}+\left(-\frac{1}{3} \frac{P_1}{u_1}
  +\frac{1}{3}\right) \frac{\partial y}{\partial v_2}+\\
 &+\left[ \left(-\frac{1}{3} u_1^2 \frac{\partial P_2}{\partial u_1}+
  \frac{1}{3} \frac{\partial P_1}{\partial v_2}+\frac{1}{9} P_1 P_2-
  F_1 F_2 \right) \frac{1}{u_1}-\frac{1}{9} P_2 \right] y=0
\endaligned$$
at $u_1=0$. Here, $F_1=F_1(\frac{1}{u_1}, v_2)$,
$F_2=F_2(\frac{1}{u_1}, v_2)$, $P_1=P_1(\frac{1}{u_1}, v_2)$ and
$P_2=P_2(\frac{1}{u_1}, v_2)$.

  Since
$$\left\{\aligned
  \frac{1}{3} P_1 v_1^2+\frac{1}{2} v_1 &=0 \quad \text{at} \quad v_1=0,\\
  \frac{1}{3} P_1 (1-v_1)^2-\frac{1}{3} (1-v_1) &=0 \quad \text{at} \quad v_1=1,\\
  \frac{1}{3} P_1 (v_2-v_1)^2-\frac{1}{3} (v_2-v_1) &=0 \quad \text{at} \quad v_1=v_2,\\
  \frac{4}{3} u_1-\frac{1}{3} P_1 &=0  \quad \text{at} \quad u_1=0,
\endaligned\right.$$
we have
$$P_1(v_1)=\frac{p_1(v_1)}{v_1 (1-v_1) (v_2-v_1)}$$
and
$$\frac{4}{3} u_1-\frac{1}{3} \frac{u_1^3 p_1\left(\frac{1}{u_1}\right)}
  {(u_1-1)(u_1 v_2-1)}=0 \quad \text{at} \quad u_1=0.$$
Hence, $p_1\left(\frac{1}{u_1}\right)=\frac{a_1}{u_1^2}+
\frac{b_1}{u_1}+c_1$ and
$$P_1(v_1)=\frac{a_1 v_1^2+b_1 v_1+c_1}{v_1 (1-v_1) (v_2-v_1)}.$$
Moreover,
$$\left\{\aligned
  -\frac{1}{3} P_1 v_1+\frac{1}{4} &=0 \quad \text{at} \quad v_1=0,\\
  -\frac{1}{3} P_1 (1-v_1)-\frac{1}{6} &=0 \quad \text{at} \quad v_1=1,\\
  -\frac{1}{3} \frac{P_1}{u_1}+\frac{1}{3} &=0 \quad \text{at} \quad u_1=0.
\endaligned\right.$$
This implies that
$$c_1=\frac{3}{4} v_2, \quad a_1+b_1+c_1=-\frac{1}{2}(v_2-1),
  \quad a_1=1.$$
Thus,
$$a_1=1, \quad b_1=-\frac{5}{4} v_2-\frac{1}{2}, \quad
  c_1=\frac{3}{4} v_2.$$
Hence,
$$P_1(v_1, v_2)=\frac{\frac{3}{4}}{v_1}+\frac{\frac{1}{2}}{v_1-1}+
           \frac{-\frac{1}{4}}{v_1-v_2}.$$

  For
$$\left\{\aligned
  F_1 v_1^2 &=0 \quad \text{at} \quad v_1=0,\\
  F_1 (1-v_1)^2 &=0 \quad \text{at} \quad v_1=1,\\
  F_1 (v_2-v_1)^2 &=0 \quad \text{at} \quad v_1=v_2,\\
  F_1 \frac{1}{u_1^2} &=0 \quad \text{at} \quad u_1=0,
\endaligned\right.$$
we have
$$F_1(v_1)=\frac{f_1(v_1)}{v_1 (1-v_1) (v_2-v_1)}$$
and
$$\frac{u_1 f_1\left(\frac{1}{u_1}\right)}{(u_1-1)(u_1 v_2-1)}=0 \quad
  \text{at} \quad u_1=0.$$
Hence,
$$F_1(v_1, v_2)=\frac{A(v_2)}{v_1 (v_1-1) (v_1-v_2)}.$$

  Similarly, we have
$$P_2(v_2)=\frac{a_2 v_2^2+b_2 v_2+c_2}{v_2(v_2-1)(v_2-v_1)}.$$
Here, $a_2=1$, $b_2=-\frac{5}{4} v_1-\frac{1}{2}$ and
$c_2=\frac{3}{4} v_1$. Hence,
$$P_2(v_1, v_2)=\frac{\frac{3}{4}}{v_2}+\frac{\frac{1}{2}}{v_2-1}+
           \frac{-\frac{1}{4}}{v_2-v_1}.$$
On the other hand,
$$F_2(v_1, v_2)=\frac{B(v_1)}{v_2 (v_2-1) (v_2-v_1)}.$$

  The equation
$$\left(\frac{\partial F_1}{\partial v_2}-\frac{1}{3} \frac{\partial P_1}
  {\partial v_1}-\frac{2}{9} P_1^2-\frac{2}{3} F_1 P_2\right) (v_2-v_1)^2
  -\frac{1}{18} P_1 (v_2-v_1)-\frac{1}{6} F_1 (v_2-v_1)-\frac{5}{36}=0$$
at $v_1=v_2$ gives that
$$A(v_2)=\frac{1}{4} v_2(v_2-1).$$
Hence,
$$F_1(v_1, v_2)=\frac{\frac{1}{4} v_2(v_2-1)}{v_1(v_1-1)(v_1-v_2)}.$$
Similarly,
$$F_2(v_1, v_2)=\frac{\frac{1}{4} v_1(v_1-1)}{v_2(v_2-1)(v_2-v_1)}.$$

  This leads us to study the following equations:

{\smc Theorem 3.3 (Main Theorem 2)}. {\it For the system of linear
partial differential equations:
$$\left\{\aligned
  \frac{\partial^2 z}{\partial v_1^2}+\frac{1}{3} P_1 \frac{\partial z}
  {\partial v_1}-F_1 \frac{\partial z}{\partial v_2}+\left(\frac{\partial
  F_1}{\partial v_2}-\frac{1}{3} \frac{\partial P_1}{\partial v_1}-
  \frac{2}{9} P_1^2-\frac{2}{3} F_1 P_2\right) z &=0,\\
  \frac{\partial^2 z}{\partial v_1 \partial v_2}-\frac{1}{3} P_2 \frac{\partial z}
  {\partial v_1}-\frac{1}{3} P_1 \frac{\partial z}{\partial v_2}+\left(\frac{1}{3}
  \frac{\partial P_2}{\partial v_1}+\frac{1}{3} \frac{\partial P_1}{\partial v_2}+
  \frac{1}{9} P_1 P_2-F_1 F_2\right) z &=0,\\
  \frac{\partial^2 z}{\partial v_2^2}+\frac{1}{3} P_2 \frac{\partial z}
  {\partial v_2}-F_2 \frac{\partial z}{\partial v_1}+\left(\frac{\partial
  F_2}{\partial v_1}-\frac{1}{3} \frac{\partial P_2}{\partial v_2}-
  \frac{2}{9} P_2^2-\frac{2}{3} F_2 P_1\right) z &=0,
  \endaligned\right.\tag 3.12$$
put
$$P_1=\frac{\alpha}{v_1}+\frac{\beta}{v_1-1}+\frac{\gamma}{v_1-v_2}, \quad
  P_2=\frac{\alpha}{v_2}+\frac{\beta}{v_2-1}+\frac{\gamma}{v_2-v_1},\tag 3.13$$
$$F_1=\frac{-\gamma v_2(v_2-1)}{v_1(v_1-1)(v_1-v_2)}, \quad
  F_2=\frac{-\gamma v_1(v_1-1)}{v_2(v_2-1)(v_2-v_1)}.\tag 3.14$$
Set
$$w=v_1^{-\frac{1}{3} \alpha} v_2^{-\frac{1}{3} \alpha}
    (v_1-1)^{-\frac{1}{3} \beta} (v_2-1)^{-\frac{1}{3} \beta}
    (v_1-v_2)^{\frac{2}{3} \gamma} z,\tag 3.15$$
then the function $w=w(v_1, v_2)$ satisfies the following
equations:
$$\left\{\aligned
 &\frac{\partial^2 w}{\partial v_1^2}+\left(\frac{\alpha}{v_1}+
  \frac{\beta}{v_1-1}+\frac{-\gamma}{v_1-v_2}\right) \frac{\partial
  w}{\partial v_1}+\frac{\gamma v_2(v_2-1)}{v_1(v_1-1)(v_1-v_2)}
  \frac{\partial w}{\partial v_2}+\\
 &+\frac{(1-\alpha-\beta) \gamma}{v_1(v_1-1)} w=0,\\
 &\frac{\partial^2 w}{\partial v_2^2}+\left(\frac{\alpha}{v_2}+
  \frac{\beta}{v_2-1}+\frac{-\gamma}{v_2-v_1}\right) \frac{\partial
  w}{\partial v_2}+\frac{\gamma v_1(v_1-1)}{v_2(v_2-1)(v_2-v_1)}
  \frac{\partial w}{\partial v_1}+\\
 &+\frac{(1-\alpha-\beta) \gamma}{v_2(v_2-1)} w=0,\\
 &\frac{\partial^2 w}{\partial v_1 \partial v_2}=\frac{-\gamma}
  {v_1-v_2} \frac{\partial w}{\partial v_1}+\frac{\gamma}{v_1-v_2}
  \frac{\partial w}{\partial v_2}.
\endaligned\right.\tag 3.16$$
A solution is given by
$$\aligned
  w=&\text{const} \cdot F_1(\alpha+\beta-1; -\gamma, -\gamma;
     \alpha-\gamma; v_1, v_2)+\\
   +&\text{const} \cdot F_1(\alpha+\beta-1; -\gamma, -\gamma;
     \beta-\gamma; 1-v_1, 1-v_2),
\endaligned\tag 3.17$$
where $F_1=F_1(a; b, b^{\prime}; c; x, y)$ is the two variable
hypergeometric function $($see \cite{AK}$)$. Consequently,
$$\aligned
 z=&v_1^{\frac{1}{3} \alpha} v_2^{\frac{1}{3} \alpha}
    (v_1-1)^{\frac{1}{3} \beta} (v_2-1)^{\frac{1}{3} \beta}
    (v_1-v_2)^{-\frac{2}{3} \gamma} \times\\
 \times &[\text{const} \cdot F_1(\alpha+\beta-1; -\gamma, -\gamma;
    \alpha-\gamma; v_1, v_2)+\\
  +&\text{const} \cdot F_1(\alpha+\beta-1; -\gamma, -\gamma; \beta
   -\gamma; 1-v_1, 1-v_2)].
\endaligned\tag 3.18$$}

{\it Proof}. Note that
$$z=v_1^{\frac{1}{3} \alpha} v_2^{\frac{1}{3} \alpha}
    (v_1-1)^{\frac{1}{3} \beta} (v_2-1)^{\frac{1}{3} \beta}
    (v_1-v_2)^{-\frac{2}{3} \gamma} w.$$
We have
$$\frac{\partial z}{\partial v_1}
 =v_1^{\frac{1}{3} \alpha} v_2^{\frac{1}{3} \alpha}
  (v_1-1)^{\frac{1}{3} \beta} (v_2-1)^{\frac{1}{3} \beta}
  (v_1-v_2)^{-\frac{2}{3} \gamma} \left[\frac{\partial w}
  {\partial v_1}+\left(\frac{1}{3} P_1-\frac{\gamma}{v_1-
  v_2}\right) w\right].$$
$$\frac{\partial z}{\partial v_2}
 =v_1^{\frac{1}{3} \alpha} v_2^{\frac{1}{3} \alpha}
  (v_1-1)^{\frac{1}{3} \beta} (v_2-1)^{\frac{1}{3} \beta}
  (v_1-v_2)^{-\frac{2}{3} \gamma} \left[\frac{\partial w}
  {\partial v_2}+\left(\frac{1}{3} P_2-\frac{\gamma}{v_2-
  v_1}\right) w\right].$$
$$\aligned
 &\frac{\partial^2 z}{\partial v_1^2}
=v_1^{\frac{1}{3} \alpha} v_2^{\frac{1}{3} \alpha}
  (v_1-1)^{\frac{1}{3} \beta} (v_2-1)^{\frac{1}{3} \beta}
  (v_1-v_2)^{-\frac{2}{3} \gamma} \times \\
 &\times \left[\frac{\partial^2 w}{\partial v_1^2}+2
  (\frac{1}{3} P_1-\frac{\gamma}{v_1-v_2}) \frac{\partial w}
  {\partial v_1}+((\frac{1}{3} P_1-\frac{\gamma}{v_1-v_2})^2+
  \frac{1}{3} \frac{\partial P_1}{\partial v_1}+\frac{\gamma}
  {(v_1-v_2)^2}) w\right].
\endaligned$$
$$\aligned
 &\frac{\partial^2 z}{\partial v_2^2}
=v_1^{\frac{1}{3} \alpha} v_2^{\frac{1}{3} \alpha}
  (v_1-1)^{\frac{1}{3} \beta} (v_2-1)^{\frac{1}{3} \beta}
  (v_1-v_2)^{-\frac{2}{3} \gamma} \times \\
 &\times \left[\frac{\partial^2 w}{\partial v_2^2}+2
  (\frac{1}{3} P_2-\frac{\gamma}{v_2-v_1}) \frac{\partial w}
  {\partial v_2}+((\frac{1}{3} P_2-\frac{\gamma}{v_2-v_1})^2+
  \frac{1}{3} \frac{\partial P_2}{\partial v_2}+\frac{\gamma}
  {(v_2-v_1)^2}) w\right].
\endaligned$$
$$\aligned
  \frac{\partial^2 z}{\partial v_1 \partial v_2}
=&v_1^{\frac{1}{3} \alpha} v_2^{\frac{1}{3} \alpha}
  (v_1-1)^{\frac{1}{3} \beta} (v_2-1)^{\frac{1}{3} \beta}
  (v_1-v_2)^{-\frac{2}{3} \gamma} \times \\
 &\times [\frac{\partial^2 w}{\partial v_1 \partial v_2}+
  (\frac{1}{3} P_2-\frac{\gamma}{v_2-v_1}) \frac{\partial w}
  {\partial v_1}+(\frac{1}{3} P_1-\frac{\gamma}{v_1-v_2})
  \frac{\partial w}{\partial v_2}+\\
 &+((\frac{1}{3} P_1-\frac{\gamma}{v_1-v_2})(\frac{1}{3} P_2
  -\frac{\gamma}{v_2-v_1})+\frac{1}{3} \frac{\partial P_1}
  {\partial v_2}-\frac{\gamma}{(v_1-v_2)^2}) w].
\endaligned$$
The equation
$$\frac{\partial^2 z}{\partial v_1^2}+\frac{1}{3} P_1 \frac{\partial z}
  {\partial v_1}-F_1 \frac{\partial z}{\partial v_2}+\left(\frac{\partial
  F_1}{\partial v_2}-\frac{1}{3} \frac{\partial P_1}{\partial v_1}-
  \frac{2}{9} P_1^2-\frac{2}{3} F_1 P_2\right) z=0$$
gives that
$$\aligned
 &\frac{\partial^2 w}{\partial v_1^2}+\left(P_1-\frac{2 \gamma}{v_1-v_2}
  \right) \frac{\partial w}{\partial v_1}-F_1 \frac{\partial w}{\partial
  v_2}+\\
 &+\left[-\frac{\gamma}{v_1-v_2} P_1+\frac{\gamma^2+\gamma}{(v_1-v_2)^2}
  -F_1 P_2+\frac{\gamma}{v_2-v_1} F_1+\frac{\partial F_1}{\partial
  v_2}\right] w=0.
\endaligned$$
Here,
$$\frac{\partial F_1}{\partial v_2}=\frac{\gamma}{v_1(v_1-1)}-
  \frac{\gamma}{(v_1-v_2)^2}.$$
Thus,
$$\aligned
 &-\frac{\gamma}{v_1-v_2} P_1+\frac{\gamma^2+\gamma}{(v_1-v_2)^2}-
  F_1 P_2+\frac{\gamma}{v_2-v_1} F_1+\frac{\partial F_1}{\partial v_2}\\
=&\frac{\gamma}{v_1-v_2} \left(-P_1+\frac{\gamma}{v_1-v_2}\right)+
  F_1 \left(-P_2+\frac{\gamma}{v_2-v_1}\right)+\frac{\gamma}{(v_1-v_2)^2}+
  \frac{\partial F_1}{\partial v_2}\\
=&\frac{\gamma}{v_1-v_2} \left(-\frac{\alpha}{v_1}-\frac{\beta}
  {v_1-1}\right)+\frac{\gamma v_2(v_2-1)}{v_1(v_1-1)(v_1-v_2)}
  \left(\frac{\alpha}{v_2}+\frac{\beta}{v_2-1}\right)+\frac{\gamma}
  {v_1(v_1-1)}\\
=&\frac{(1-\alpha-\beta) \gamma}{v_1(v_1-1)}.
\endaligned$$
Hence, we get the following equation:
$$\frac{\partial^2 w}{\partial v_1^2}+\left(\frac{\alpha}{v_1}+
  \frac{\beta}{v_1-1}+\frac{-\gamma}{v_1-v_2}\right) \frac{\partial
  w}{\partial v_1}+\frac{\gamma v_2(v_2-1)}{v_1(v_1-1)(v_1-v_2)}
  \frac{\partial w}{\partial v_2}+\frac{(1-\alpha-\beta) \gamma}
  {v_1(v_1-1)} w=0.$$
Exchange the position of $v_1$ and $v_2$, we get the other
equation.

  The equation
$$\frac{\partial^2 z}{\partial v_1 \partial v_2}-\frac{1}{3} P_2 \frac{\partial z}
  {\partial v_1}-\frac{1}{3} P_1 \frac{\partial z}{\partial v_2}+\left(\frac{1}{3}
  \frac{\partial P_2}{\partial v_1}+\frac{1}{3} \frac{\partial P_1}{\partial v_2}+
  \frac{1}{9} P_1 P_2-F_1 F_2\right) z=0$$
gives that
$$\frac{\partial^2 w}{\partial v_1 \partial v_2}+\frac{\gamma}
  {v_1-v_2} \frac{\partial w}{\partial v_1}+\frac{\gamma}{v_2-v_1}
  \frac{\partial w}{\partial v_2}+\left[-\frac{\gamma^2+\gamma}
  {(v_1-v_2)^2}+\frac{2}{3} \frac{\partial P_1}{\partial v_2}+
  \frac{1}{3} \frac{\partial P_2}{\partial v_1}-F_1 F_2\right]
  w=0.$$
Note that
$$\frac{\partial P_1}{\partial v_2}=\frac{\partial P_2}{\partial
  v_1}=\frac{\gamma}{(v_1-v_2)^2}, \quad
  F_1 F_2=\frac{-\gamma^2}{(v_1-v_2)^2},$$
we find that
$$-\frac{\gamma^2+\gamma}{(v_1-v_2)^2}+\frac{2}{3} \frac{\partial
  P_1}{\partial v_2}+\frac{1}{3} \frac{\partial P_2}{\partial v_1}
  -F_1 F_2=0.$$
Thus,
$$\frac{\partial^2 w}{\partial v_1 \partial v_2}+\frac{\gamma}
  {v_1-v_2} \frac{\partial w}{\partial v_1}+\frac{\gamma}{v_2-v_1}
  \frac{\partial w}{\partial v_2}=0.$$

  It is known that the equations
$$F_1: \left\{\aligned
  x(1-x) \frac{\partial^2 z}{\partial x^2}+y(1-x) \frac{\partial^2
  z}{\partial x \partial y}+[c-(a+b+1)x] \frac{\partial
  z}{\partial x}-by \frac{\partial z}{\partial y}-abz &=0,\\
  y(1-y) \frac{\partial^2 z}{\partial y^2}+x(1-y) \frac{\partial^2
  z}{\partial x \partial y}+[c-(a+b^{\prime}+1)y] \frac{\partial
  z}{\partial y}-b^{\prime} x \frac{\partial z}{\partial x}-a
  b^{\prime} z &=0
\endaligned\right.$$
are equivalent to the following three equations:
$$\left\{\aligned
 &\frac{\partial^2 z}{\partial x^2}+\left(\frac{c-b^{\prime}}{x}+
  \frac{a+b-c+1}{x-1}+\frac{b^{\prime}}{x-y}\right) \frac{\partial
  z}{\partial x}-\frac{by(y-1)}{x(x-1)(x-y)} \frac{\partial z}{\partial
  y}+\frac{abz}{x(x-1)}=0,\\
 &\frac{\partial^2 z}{\partial y^2}+\left(\frac{c-b}{y}+
  \frac{a+b^{\prime}-c+1}{y-1}+\frac{b}{y-x}\right)
  \frac{\partial z}{\partial y}-\frac{b^{\prime}x(x-1)}{y(y-1)(y-x)}
  \frac{\partial z}{\partial x}+\frac{ab^{\prime}z}{y(y-1)}=0,\\
 &\frac{\partial^2 z}{\partial x \partial y}=\frac{b^{\prime}}{x-y}
  \frac{\partial z}{\partial x}-\frac{b}{x-y} \frac{\partial z}
  {\partial y},
\endaligned\right.$$
where $F_1=F_1(a; b, b^{\prime}; c; x, y)$ is the two variable
hypergeometric function (see \cite{AK}). Here,
$$c-b^{\prime}=c-b=\alpha, \quad a+b-c+1=a+b^{\prime}-c+1=\beta,$$
$$b=b^{\prime}=-\gamma, \quad ab=ab^{\prime}=(1-\alpha-\beta) \gamma.$$
Thus,
$$a=\alpha+\beta-1, \quad b=-\gamma, \quad b^{\prime}=-\gamma,
  \quad c=\alpha-\gamma.$$
A solution is given by
$$\aligned
  w=&\text{const} \cdot F_1(\alpha+\beta-1; -\gamma, -\gamma; \alpha-\gamma; v_1, v_2)+\\
   +&\text{const} \cdot F_1(\alpha+\beta-1; -\gamma, -\gamma; \beta-\gamma; 1-v_1, 1-v_2).
\endaligned$$
Consequently,
$$\aligned
 z=&v_1^{\frac{1}{3} \alpha} v_2^{\frac{1}{3} \alpha}
    (v_1-1)^{\frac{1}{3} \beta} (v_2-1)^{\frac{1}{3} \beta}
    (v_1-v_2)^{-\frac{2}{3} \gamma} \times\\
 \times &[\text{const} \cdot F_1(\alpha+\beta-1; -\gamma, -\gamma;
    \alpha-\gamma; v_1, v_2)+\\
  +&\text{const} \cdot F_1(\alpha+\beta-1; -\gamma, -\gamma;
    \beta-\gamma; 1-v_1, 1-v_2)].
\endaligned$$
$\qquad \qquad \qquad \qquad \qquad \qquad \qquad \qquad \qquad
 \qquad \qquad \qquad \qquad \qquad \qquad \qquad \qquad \qquad
 \quad \boxed{}$

  In our case, $\alpha=\frac{3}{4}$, $\beta=\frac{1}{2}$,
$\gamma=-\frac{1}{4}$.
$$w=v_1^{-\frac{1}{4}} v_2^{-\frac{1}{4}} (v_1-1)^{-\frac{1}{6}}
    (v_2-1)^{-\frac{1}{6}} (v_1-v_2)^{-\frac{1}{6}} \root 3 \of{\Delta} \tag 3.19$$
satisfies the following equations:
$$\left\{\aligned
 &\frac{\partial^2 w}{\partial v_1^2}+\left(\frac{\frac{3}{4}}{v_1}+
  \frac{\frac{1}{2}}{v_1-1}+\frac{\frac{1}{4}}{v_1-v_2}\right) \frac{\partial
  w}{\partial v_1}+\frac{-\frac{1}{4} v_2(v_2-1)}{v_1(v_1-1)(v_1-v_2)}
  \frac{\partial w}{\partial v_2}+\frac{\frac{1}{16}}{v_1(v_1-1)} w=0,\\
 &\frac{\partial^2 w}{\partial v_2^2}+\left(\frac{\frac{3}{4}}{v_2}+
  \frac{\frac{1}{2}}{v_2-1}+\frac{\frac{1}{4}}{v_2-v_1}\right) \frac{\partial
  w}{\partial v_2}+\frac{-\frac{1}{4} v_1(v_1-1)}{v_2(v_2-1)(v_2-v_1)}
  \frac{\partial w}{\partial v_1}+\frac{\frac{1}{16}}{v_2(v_2-1)} w=0,\\
 &\frac{\partial^2 w}{\partial v_1 \partial v_2}=\frac{\frac{1}{4}}
  {v_1-v_2} \frac{\partial w}{\partial v_1}+\frac{-\frac{1}{4}}{v_1-v_2}
  \frac{\partial w}{\partial v_2}.
\endaligned\right.\tag 3.20$$
A solution is given by
$$w=\text{const} \cdot F_1\left(\frac{1}{4}; \frac{1}{4}, \frac{1}{4}; 1;
    v_1, v_2\right)+\text{const} \cdot F_1\left(\frac{1}{4}; \frac{1}{4},
    \frac{1}{4}; \frac{3}{4}; 1-v_1, 1-v_2\right).\tag 3.21$$
Consequently,
$$\aligned
 &\root 3 \of{\Delta}
 =v_1^{\frac{1}{4}} v_2^{\frac{1}{4}} (v_1-1)^{\frac{1}{6}}
  (v_2-1)^{\frac{1}{6}} (v_1-v_2)^{\frac{1}{6}} \times\\
 &\left[\text{const} \cdot F_1\left(\frac{1}{4}; \frac{1}{4},
  \frac{1}{4}; 1; v_1, v_2\right)+\text{const} \cdot
  F_1\left(\frac{1}{4}; \frac{1}{4}, \frac{1}{4}; \frac{3}{4};
  1-v_1, 1-v_2\right)\right].
\endaligned\tag 3.22$$

  In fact, when $(w_1, w_2)$ is a solution of the system of nonlinear
partial differential equations:
$$\left\{\aligned
  \{ w_1, w_2; v_1, v_2 \}_{v_1} &=F_1(v_1, v_2),\\
  \{ w_1, w_2; v_1, v_2 \}_{v_2} &=F_2(v_1, v_2),\\
      [w_1, w_2; v_1, v_2]_{v_1} &=P_1(v_1, v_2),\\
      [w_1, w_2; v_1, v_2]_{v_2} &=P_2(v_1, v_2),
\endaligned\right.$$
the general solution has the form
$$\left(\frac{a_1 w_1+a_2 w_2+a_3}{c_1 w_1+c_2 w_2+c_3},
  \frac{b_1 w_1+b_2 w_2+b_3}{c_1 w_1+c_2 w_2+c_3}\right),$$
where $\left(\matrix a_1 & a_2 & a_3\\ b_1 & b_2 & b_3\\ c_1 & c_2
& c_3 \endmatrix\right) \in GL(3, {\Bbb C})$.

\vskip 0.5 cm
\centerline{\bf 4. $\eta$-functions associated to
                   $U(2, 1)$ and Picard curves}
\vskip 0.5 cm

  The $\eta$-function associated to the unitary group $U(2, 1)$ is given
as follows:
$$\eta(w_1, w_2):=v_1^{-\frac{1}{12}} v_2^{-\frac{1}{12}}
  (1-v_1)^{-\frac{1}{18}} (1-v_2)^{-\frac{1}{18}} (v_1-v_2)^{-\frac{1}{18}}
  \left[\frac{\partial(v_1, v_2)}{\partial(w_1, w_2)}\right]^{\frac{1}{9}}.
  \tag 4.1$$
It satisfies the following transformation:
$$\eta \left(\frac{a_1 w_1+a_2 w_2+a_3}{c_1 w_1+c_2 w_2+c_3},
  \frac{b_1 w_1+b_2 w_2+b_3}{c_1 w_1+c_2 w_2+c_3} \right)
 =\Delta^{-\frac{1}{9}} (c_1 w_1+c_2 w_2+c_3)^{\frac{1}{3}} \eta(w_1, w_2),
  \tag 4.2$$
where $\left(\matrix a_1 & a_2 & a_3\\ b_1 & b_2 & b_3\\ c_1 & c_2
& c_3 \endmatrix\right) \in GL(3, {\Bbb C})$ and
$\Delta=\vmatrix \format \c \quad & \c \quad & \c\\
a_1 & a_2 & a_3\\ b_1 & b_2 & b_3\\ c_1 & c_2 & c_3
\endvmatrix$.

  Let
$$T_1=\left(\matrix
      1 & 1 & -\omega\\
        & 1 &       1\\
        &   &       1
      \endmatrix\right), \quad
  T_2=\left(\matrix
      1 & \omega & -\omega\\
        &      1 & \overline{\omega}\\
        &        & 1
      \endmatrix\right), \quad
  S=-\overline{\omega} J
   =\left(\matrix
                       &                   & -\overline{\omega}\\
                       & \overline{\omega} &     \\
    -\overline{\omega} &                   &
    \endmatrix\right),$$
$$U_1=\left(\matrix
      1 &         &  \\
        & -\omega &  \\
        &         & 1
      \endmatrix\right), \quad
  U_2=\left(\matrix
      -1 &         &   \\
         & -\omega &   \\
         &         & -1
      \endmatrix\right).$$
The actions of the above five generators on $(w_1, w_2) \in {\frak
S}_{2}$ are given as follows:
$$\left\{\aligned
  T_1(w_1, w_2) &=(w_1+w_2-\omega, w_2+1),\\
  T_2(w_1, w_2) &=(w_1+\omega w_2-\omega, w_2+\overline{\omega}),\\
    S(w_1, w_2) &=\left(\frac{1}{w_1}, -\frac{w_2}{w_1}\right),\\
  U_1(w_1, w_2) &=(w_1, -\omega w_2),\\
  U_2(w_1, w_2) &=(w_1, \omega w_2).
  \endaligned\right.\tag 4.3$$
We have
$$\eta(T_1(w_1, w_2))
 =1^{-\frac{1}{9}} 1^{\frac{1}{3}} \eta(w_1, w_2)
 =1^{\frac{2}{9}} \eta(w_1, w_2)
 =e^{\frac{4}{9} \pi i} \eta(w_1, w_2).$$
$$\eta(T_2(w_1, w_2))
 =1^{-\frac{1}{9}} 1^{\frac{1}{3}} \eta(w_1, w_2)
 =1^{\frac{2}{9}} \eta(w_1, w_2)
 =e^{\frac{4}{9} \pi i} \eta(w_1, w_2).$$
$$\aligned
  \eta(U_1(w_1, w_2))
&=(-\omega)^{-\frac{1}{9}} 1^{\frac{1}{3}} \eta(w_1, w_2)
 =(e^{\frac{5}{3} \pi i})^{-\frac{1}{9}} e^{\frac{2}{3} \pi i} \eta(w_1, w_2)\\
&=e^{\frac{13}{27} \pi i} \eta(w_1, w_2)
 =1^{\frac{13}{54}} \eta(w_1, w_2).
\endaligned$$
$$\aligned
  \eta(U_2(w_1, w_2))
&=(-\omega)^{-\frac{1}{9}} (-1)^{\frac{1}{3}} \eta(w_1, w_2)
 =(e^{\frac{5}{3} \pi i})^{-\frac{1}{9}} e^{\frac{1}{3} \pi i} \eta(w_1, w_2)\\
&=e^{\frac{4}{27} \pi i} \eta(w_1, w_2)
 =1^{\frac{2}{27}} \eta(w_1, w_2).
\endaligned$$
$$\aligned
  \eta(S(w_1, w_2))
&=(-1)^{-\frac{1}{9}} (-\overline{\omega} w_1)^{\frac{1}{3}}
  \eta(w_1, w_2)
 =(e^{\pi i})^{-\frac{1}{9}} (e^{\frac{1}{3} \pi i} w_1)^{\frac{1}{3}} \eta(w_1, w_2)\\
&=w_1^{\frac{1}{3}} \eta(w_1, w_2).
\endaligned$$

  Note that
$$[T_1, T_2]:=T_1 T_2 T_1^{-1} T_2^{-1}
 =\left(\matrix
  1 & 0 & \overline{\omega}-\omega\\
    & 1 & 0\\
    &   & 1
  \endmatrix\right).$$
In fact,
$$\eta(T_1^{-1}(w_1, w_2))=1^{-\frac{2}{9}} \eta(w_1, w_2), \quad
  \eta(T_2^{-1}(w_1, w_2))=1^{-\frac{2}{9}} \eta(w_1, w_2).$$
Hence,
$$\eta([T_1, T_2](w_1, w_2))=\eta(w_1, w_2).$$

  In fact, Holzapfel (see \cite{Ho1}) proved the following theorem:

{\smc Theorem (\cite{Ho1})}. {\it The monodromy group of the
hypergeometric differential equation system considered by Picard
is the principal congruence subgroup $U(2, 1; {\Cal
O}_{K})(1-\omega)$ of the full Eisenstein lattice of the ball
$U(2, 1; {\Cal O}_{K})$ with respect to the principal ideal ${\Cal
O}_{K} \cdot (1-\omega)$, $\omega=\exp(2 \pi i/3)$, $K={\Bbb
Q}(\sqrt{-3})$.}

  According to Picard (see \cite{Pi2} or \cite{Ho1}), the five
generators of the principal congruence subgroup $U(2, 1; {\Cal
O}_{K})(1-\omega)$ are given as follows:
$$g_1=\left(\matrix
  1 & \overline{\omega}-\omega & 1-\overline{\omega}\\
  0 & \overline{\omega}        & 1-\omega\\
  0 & 0                        & 1
  \endmatrix\right), \quad
  g_2=\left(\matrix
  1 & \overline{\omega}-1 & 1-\overline{\omega}\\
  0 & \overline{\omega}   & 1-\overline{\omega}\\
  0 & 0                   & 1
  \endmatrix\right), \quad
  g_3=\left(\matrix
  1 & 0      & 0\\
  0 & \omega & 0\\
  0 & 0      & 1
  \endmatrix\right),\tag 4.4$$
$$g_4=\left(\matrix
  -\omega             & 0  & \overline{\omega}-1\\
  0                   & -1 & 0\\
  \overline{\omega}-1 & 0  & 2 \omega
  \endmatrix\right), \quad
  g_5=\left(\matrix
  1 & 0 & 0\\
  \overline{\omega}-\omega & \overline{\omega} & 0\\
  1-\overline{\omega} & 1-\omega & 1
  \endmatrix\right).\tag 4.5$$
In fact,
$$\left\{\aligned
  g_1 &=U_1^{-4} T_1^{-1} T_2^{-2}, \\
  g_2 &=U_1^{-4} T_1^{-2} T_2^{-1}, \\
  g_3 &=U_1^4,\\
  g_4 &=S^3 [T_1, T_2] S^3 (S^4 U_2)^{-1} [T_1, T_2],\\
  g_5 &=S^3 U_1^{-4} T_1^{-1} T_2 S^3.
\endaligned\right.\tag 4.6$$

  For $g=\left(\matrix \alpha & &\\ & \beta & \\ & & \gamma
\endmatrix\right) \in U(2, 1; {\Cal O}_{K})$, we have
$\overline{\alpha} \gamma=\beta \overline{\beta}=\alpha
\overline{\gamma}=1$. Hence,
$$g=\left(\matrix \alpha &  & \\ & \beta & \\ & &
    \overline{\alpha}^{-1} \endmatrix\right), \quad
    \alpha, \beta \in \{ \pm 1, \pm \omega,
    \pm \overline{\omega} \}.$$
Note that
$$U_1^{3}=\left(\matrix 1 & & \\ & -1 & \\ & & 1
          \endmatrix\right), \quad
  (U_1 U_2)^{3}=\left(\matrix -1 & & \\ & 1 & \\ & & -1
                \endmatrix\right).$$
We only need to consider the following nine cases:
$$I=U_1^{6}=U_2^{6}=S^{6}, \quad
  \left(\matrix \omega & & \\ & 1 & \\ & & \omega
  \endmatrix\right)=S^{2} U_1^{2}, \quad
  \left(\matrix \overline{\omega} & & \\ & 1 & \\
  & & \overline{\omega} \endmatrix\right)=(S^{2} U_1^{2})^{2}.$$
$$\left(\matrix 1 & & \\ & \omega & \\ & & 1
  \endmatrix\right)=U_1^{4}=U_2^{4}, \quad
  \left(\matrix \omega & & \\ & \omega & \\
  & & \omega \endmatrix\right)=S^{2}, \quad
  \left(\matrix \overline{\omega} & & \\ & \omega & \\
  & & \overline{\omega} \endmatrix\right)=S^{4} U_1^{2}.$$
$$\left(\matrix 1 & & \\ & \overline{\omega} & \\ & & 1
  \endmatrix\right)=U_1^{2}=U_2^{2}, \quad
  \left(\matrix \omega & & \\ & \overline{\omega} & \\
  & & \omega \endmatrix\right)=(S^{4} U_1^{2})^{2}, \quad
  \left(\matrix \overline{\omega} & & \\ & \overline{\omega} & \\
  & & \overline{\omega} \endmatrix\right)=S^{4}.$$

  Put
$$\eta_1(w_1, w_2):=\eta\left(\left(\matrix
  3 &          &  \\
    & 1-\omega &  \\
    &          & 1
  \endmatrix\right)(w_1, w_2)\right)
 =\eta(3 w_1, (1-\omega) w_2),\tag 4.7$$
$$\eta_2(w_1, w_2):=\eta\left(\left(\matrix
  1 &          &         \\
    & 1-\omega &         \\
    &          & 3
  \endmatrix\right)(w_1, w_2)\right)
 =\eta\left(\frac{w_1}{3}, \frac{w_2}{1-\overline{\omega}}\right).\tag 4.8$$

{\smc Proposition 4.1}. {\it
$$\left\{\aligned
 \eta_1(T_1(w_1, w_2)) &=1^{\frac{2}{3}} \eta_1(w_1, w_2),\\
 \eta_2(T_1^3 (w_1, w_2)) &=1^{\frac{2}{3}} \eta_2(w_1, w_2),\\
 \eta_1(T_2(w_1, w_2)) &=\eta_1(w_1, w_2),\\
 \eta_2(T_2^3 (w_1, w_2)) &=\eta_2(w_1, w_2),\\
 \eta_1(S(w_1, w_2)) &=\root 3 \of{\frac{w_1}{3}} \eta_2(w_1, w_2).
\endaligned\right.\tag 4.9$$}

{\it Proof}. We have
$$\aligned
 &\eta_1(T_1(w_1, w_2))\\
=&\eta\left(\left(\matrix
  3 &          &  \\
    & 1-\omega &  \\
    &          & 1
  \endmatrix\right) T_1
  \left(\matrix
  \frac{1}{3} &                    &  \\
              & \frac{1}{1-\omega} &  \\
              &                    & 1
  \endmatrix\right)(3 w_1, (1-\omega) w_2)\right)\\
=&\eta\left(\left(\matrix
  1 & 1-\overline{\omega} & -3 \omega\\
  0 &                   1 &  1-\omega\\
  0 &                   0 &         1
  \endmatrix\right)(3 w_1, (1-\omega) w_2)\right)\\
=&\eta(T_1^2 T_2 [T_1, T_2]^{-1}(3 w_1, (1-\omega) w_2))\\
=&1^{\frac{2}{3}} \eta(3 w_1, (1-\omega) w_2)\\
=&1^{\frac{2}{3}} \eta_1(w_1, w_2).
\endaligned$$
$$\aligned
 &\eta_2(T_1^3(w_1, w_2))\\
=&\eta\left(\left(\matrix
  1 &          &  \\
    & 1-\omega &  \\
    &          & 3
  \endmatrix\right) T_1^3
  \left(\matrix
  1 &                    &            \\
    & \frac{1}{1-\omega} &            \\
    &                    & \frac{1}{3}
  \endmatrix\right)\left(\frac{w_1}{3},
  \frac{w_2}{1-\overline{\omega}}\right)\right)\\
=&\eta\left(\left(\matrix
  1 & 1-\overline{\omega} & 1-\omega\\
  0 &                   1 & 1-\omega\\
  0 &                   0 &        1
  \endmatrix\right) \left(\frac{w_1}{3},
  \frac{w_2}{1-\overline{\omega}}\right)\right)\\
=&\eta\left(T_1^2 T_2 [T_1, T_2]^{-2} \left(\frac{w_1}{3},
  \frac{w_2}{1-\overline{\omega}}\right)\right)\\
=&1^{\frac{2}{3}} \eta\left(\frac{w_1}{3}, \frac{w_2}{1-\overline{\omega}}\right)\\
=&1^{\frac{2}{3}} \eta_2(w_1, w_2).
\endaligned$$
$$\aligned
 &\eta_1(T_2(w_1, w_2))\\
=&\eta\left(\left(\matrix
  3 &          &  \\
    & 1-\omega &  \\
    &          & 1
  \endmatrix\right) T_2
  \left(\matrix
  \frac{1}{3} &                    &   \\
              & \frac{1}{1-\omega} &   \\
              &                    & 1
  \endmatrix\right) (3 w_1, (1-\omega) w_2)\right)\\
=&\eta\left(\left(\matrix
  1 & \omega-1 &           -3 \omega\\
  0 &        1 & \overline{\omega}-1\\
  0 &        0 &                   1
  \endmatrix\right)(3 w_1, (1-\omega) w_2)\right)\\
=&\eta(T_1^{-1} T_2 [T_1, T_2]^2 (3 w_1, (1-\omega) w_2))\\
=&\eta(3 w_1, (1-\omega) w_2)\\
=&\eta_1(w_1, w_2).
\endaligned$$
$$\aligned
 &\eta_2(T_2^3(w_1, w_2))\\
=&\eta\left(\left(\matrix
  1 &          &  \\
    & 1-\omega &  \\
    &          & 3
  \endmatrix\right) T_2^3
  \left(\matrix
  1 &                    &            \\
    & \frac{1}{1-\omega} &            \\
    &                    & \frac{1}{3}
  \endmatrix\right) \left(\frac{w_1}{3},
  \frac{w_2}{1-\overline{\omega}}\right)\right)\\
=&\eta\left(\left(\matrix
  1 & \omega-1 &            1-\omega\\
  0 &        1 & \overline{\omega}-1\\
  0 &        0 &                   1
  \endmatrix\right) \left(\frac{w_1}{3},
  \frac{w_2}{1-\overline{\omega}}\right)\right)\\
=&\eta\left(T_1^{-1} T_2 [T_1, T_2] \left(\frac{w_1}{3},
  \frac{w_2}{1-\overline{\omega}}\right)\right)\\
=&\eta\left(\frac{w_1}{3}, \frac{w_2}{1-\overline{\omega}}\right)\\
=&\eta_2(w_1, w_2).
\endaligned$$
$$\aligned
 &\eta_1(S(w_1, w_2))=\eta\left(\frac{3}{w_1}, (1-\omega)
  \left(-\frac{w_2}{w_1}\right)\right)\\
=&\eta\left(S\left(\frac{w_1}{3},
  \frac{w_2}{1-\overline{\omega}}\right)\right)
 =\root 3 \of{\frac{w_1}{3}} \eta\left(\frac{w_1}{3},
  \frac{w_2}{1-\overline{\omega}}\right)\\
=&\root 3 \of{\frac{w_1}{3}} \eta_2(w_1, w_2).
\endaligned$$
$\qquad \qquad \qquad \qquad \qquad \qquad \qquad \qquad \qquad
 \qquad \qquad \qquad \qquad \qquad \qquad \qquad \qquad \qquad
 \quad \boxed{}$

  Put
$$\eta_3(w_1, w_2):=\eta_2([T_1, T_2]^{-1}(w_1, w_2)),\tag 4.10$$
$$\eta_4(w_1, w_2):=\eta_1(S^3 [T_1, T_2]^{-1} S^3(w_1, w_2)),\tag 4.11$$
$$\eta_5(w_1, w_2):=\eta_2([T_1, T_2]^{-1} S^3 [T_1, T_2]^{-1} S^3(w_1, w_2)).\tag 4.12$$

{\smc Proposition 4.2}. {\it
$$\eta_3(g_1(w_1, w_2))=1^{\frac{7}{27}} \eta_2(w_1, w_2), \quad
  \eta_1(g_1(w_1, w_2))=1^{-\frac{5}{27}} \eta_1(w_1, w_2).\tag 4.13$$
$$\eta_3(g_2(w_1, w_2))=1^{\frac{1}{27}} \eta_2(w_1, w_2), \quad
  \eta_1(g_2(w_1, w_2))=1^{-\frac{23}{27}} \eta_1(w_1, w_2).\tag 4.14$$
$$\eta_1(g_3(w_1, w_2))=1^{\frac{26}{27}} \eta_1(w_1, w_2), \quad
  \eta_2(g_3(w_1, w_2))=1^{\frac{26}{27}} \eta_2(w_1, w_2).\tag 4.15$$
$$\eta_4(g_4(w_1, w_2))=1^{\frac{19}{27}} \eta_1(w_1, w_2), \quad
  \eta_5(g_4(w_1, w_2))=1^{\frac{19}{27}} \eta_2(w_1, w_2).\tag 4.16$$
$$\aligned
  \eta_4(g_5(w_1, w_2)) &=1^{\frac{7}{27}} \root 3 \of{-3 \overline{\omega}
  w_1+(1-\omega) w_2+1} \eta_1(w_1, w_2),\\
  \eta_2(g_5(w_1, w_2)) &=1^{-\frac{5}{27}} \root 3 \of{(1-\overline{\omega})
  w_1+(1-\omega) w_2+1} \eta_2(w_1, w_2).
\endaligned\tag 4.17$$}

{\it Proof}.
$$\aligned
 &\eta_3(g_1(w_1, w_2))\\
=&\eta\left(\left(\matrix 1 & & \\ & 1-\omega & \\ & & 3
  \endmatrix\right) [T_1, T_2]^{-1} g_1
  \left(\matrix 1 & & \\ & \frac{1}{1-\omega} & \\ & & \frac{1}{3}
  \endmatrix\right) \left(\frac{w_1}{3}, \frac{w_2}{1-\overline{\omega}}\right)\right)\\
=&\eta\left(\left(\matrix
  1 & -\omega & -\overline{\omega}\\
  0 & \overline{\omega} & -\omega\\
  0 & 0                 & 1
  \endmatrix\right)\left(\frac{w_1}{3}, \frac{w_2}{1-\overline{\omega}}\right)\right)\\
=&\eta\left(U_1^2 T_2^{-1}\left(\frac{w_1}{3}, \frac{w_2}{1-\overline{\omega}}\right)\right)\\
=&1^{\frac{7}{27}} \eta_2(w_1, w_2).
\endaligned$$
$$\aligned
 &\eta_1(g_1(w_1, w_2))\\
=&\eta\left(\left(\matrix 3 & & \\ & 1-\omega & \\ & & 1
  \endmatrix\right) g_1
  \left(\matrix \frac{1}{3} & & \\ & \frac{1}{1-\omega} & \\ & & 1
  \endmatrix\right)(3 w_1, (1-\omega) w_2)\right)\\
=&\eta\left(\left(\matrix
  1 & -3 \omega & 3(1-\overline{\omega})\\
  0 & \overline{\omega} & -3 \omega\\
  0 &                 0 &         1
  \endmatrix\right)(3 w_1, (1-\omega) w_2)\right)\\
=&\eta(U_1^2 T_2^{-3}(3 w_1, (1-\omega) w_2))\\
=&1^{-\frac{5}{27}} \eta_1(w_1, w_2).
\endaligned$$
$$\aligned
 &\eta_3(g_2(w_1, w_2))\\
=&\eta\left(\left(\matrix 1 & & \\ & 1-\omega & \\ & & 3
  \endmatrix\right) [T_1, T_2]^{-1} g_2
  \left(\matrix 1 & & \\ & \frac{1}{1-\omega} & \\ & & \frac{1}{3}
  \endmatrix\right)\left(\frac{w_1}{3}, \frac{w_2}{1-\overline{\omega}}\right)\right)\\
=&\eta\left(\left(\matrix
  1 & \overline{\omega} & -\overline{\omega}\\
  0 & \overline{\omega} & 1\\
  0 &                 0 & 1
  \endmatrix\right)\left(\frac{w_1}{3}, \frac{w_2}{1-\overline{\omega}}\right)\right)\\
=&\eta\left(U_1^2 T_1^{-1} T_2^{-1}\left(\frac{w_1}{3}, \frac{w_2}{1-\overline{\omega}}\right)\right)\\
=&1^{\frac{1}{27}} \eta_2(w_1, w_2).
\endaligned$$
$$\aligned
 &\eta_1(g_2(w_1, w_2))\\
=&\eta\left(\left(\matrix 3 & & \\ & 1-\omega & \\ & & 1
  \endmatrix\right) g_2
  \left(\matrix \frac{1}{3} & & \\ & \frac{1}{1-\omega} & \\ & & 1
  \endmatrix\right)(3 w_1, (1-\omega) w_2)\right)\\
=&\eta\left(\left(\matrix
  1 & 3 \overline{\omega} & 3(1-\overline{\omega})\\
  0 & \overline{\omega}   & 3\\
  0 &                 0   & 1
  \endmatrix\right)(3 w_1, (1-\omega) w_2)\right)\\
=&\eta(U_1^2 T_1^{-3} T_2^{-3} [T_1, T_2]^{-3}(3 w_1, (1-\omega) w_2))\\
=&1^{-\frac{23}{27}} \eta_1(w_1, w_2).
\endaligned$$
$$\aligned
 &\eta_1(g_3(w_1, w_2))\\
=&\eta\left(\left(\matrix 3 & & \\ & 1-\omega & \\ & & 1
  \endmatrix\right) g_3
  \left(\matrix \frac{1}{3} & & \\ & \frac{1}{1-\omega} & \\ & & 1
  \endmatrix\right)(3 w_1, (1-\omega) w_2)\right)\\
=&\eta(U_1^4(3 w_1, (1-\omega) w_2))\\
=&1^{\frac{26}{27}} \eta_1(w_1, w_2).
\endaligned$$
$$\aligned
 &\eta_2(g_3(w_1, w_2))\\
=&\eta\left(\left(\matrix 1 & & \\ & 1-\omega & \\ & & 3
  \endmatrix\right) g_3
  \left(\matrix 1 & & \\ & \frac{1}{1-\omega} & \\ & & \frac{1}{3}
  \endmatrix\right)\left(\frac{w_1}{3}, \frac{w_2}{1-\overline{\omega}}\right)\right)\\
=&\eta\left(U_1^4\left(\frac{w_1}{3}, \frac{w_2}{1-\overline{\omega}}\right)\right)\\
=&1^{\frac{26}{27}} \eta_2(w_1, w_2).
\endaligned$$
$$\aligned
 &\eta_4(g_4(w_1, w_2))\\
=&\eta\left(\left(\matrix 3 & & \\ & 1-\omega & \\ & & 1
  \endmatrix\right) S^3 [T_1, T_2]^{-1} S^3 g_4
  \left(\matrix \frac{1}{3} & & \\ & \frac{1}{1-\omega} & \\ & & 1
  \endmatrix\right)(3 w_1, (1-\omega) w_2)\right)\\
=&\eta\left(\left(\matrix
 -\omega &  0 & 3(\overline{\omega}-1)\\
       0 & -1 & 0\\
       0 &  0 & -\omega
  \endmatrix\right)(3 w_1, (1-\omega) w_2)\right)\\
=&\eta(U_1^2 U_2^3 S^2 [T_1, T_2]^3 (3 w_1, (1-\omega) w_2))\\
=&1^{\frac{19}{27}} \eta_1(w_1, w_2).
\endaligned$$
$$\aligned
 &\eta_5(g_4(w_1, w_2))\\
=&\eta\left(\left(\matrix 1 & & \\ & 1-\omega & \\ & & 3
  \endmatrix\right) [T_1, T_2]^{-1} S^3 [T_1, T_2]^{-1} S^3 g_4
  \left(\matrix 1 & & \\ & \frac{1}{1-\omega} & \\ & & \frac{1}{3}
  \endmatrix\right)\left(\frac{w_1}{3}, \frac{w_2}{1-\overline{\omega}}\right)\right)\\
=&\eta\left(\left(\matrix -\omega & & \\ & -1 & \\ & & -\omega
  \endmatrix\right)\left(\frac{w_1}{3}, \frac{w_2}{1-\overline{\omega}}\right)\right)\\
=&\eta\left(U_1^2 U_2^3 S^2\left(\frac{w_1}{3}, \frac{w_2}{1-\overline{\omega}}\right)\right)\\
=&1^{\frac{19}{27}} \eta_2(w_1, w_2).
\endaligned$$
$$\aligned
 &\eta_4(g_5(w_1, w_2))\\
=&\eta\left(\left(\matrix 3 & & \\ & 1-\omega & \\ & & 1
  \endmatrix\right) S^3 [T_1, T_2]^{-1} S^3 g_5
  \left(\matrix \frac{1}{3} & & \\ & \frac{1}{1-\omega} & \\ & & 1
  \endmatrix\right)(3 w_1, (1-\omega) w_2)\right)\\
=&\eta\left(\left(\matrix
  1 & 0 & 0\\
  \overline{\omega} & \overline{\omega} & 0\\
  -\overline{\omega} & 1 & 1
  \endmatrix\right)(3 w_1, (1-\omega) w_2)\right)\\
=&\eta(U_1^2 S^3 T_1^{-1} S^3(3 w_1, (1-\omega) w_2))\\
=&1^{\frac{7}{27}} \root 3 \of{-3 \overline{\omega} w_1+(1-\omega)
  w_2+1} \eta_1(w_1, w_2).
\endaligned$$
$$\aligned
 &\eta_2(g_5(w_1, w_2))\\
=&\eta\left(\left(\matrix 1 & & \\ & 1-\omega & \\ & & 3
  \endmatrix\right) g_5
  \left(\matrix 1 & & \\ & \frac{1}{1-\omega} & \\ & & \frac{1}{3}
  \endmatrix\right)\left(\frac{w_1}{3}, \frac{w_2}{1-\overline{\omega}}\right)\right)\\
=&\eta\left(\left(\matrix
  1 & 0 & 0\\
  3 \overline{\omega} & \overline{\omega} & 0\\
  3(1-\overline{\omega}) & 3 & 1
  \endmatrix\right)\left(\frac{w_1}{3}, \frac{w_2}{1-\overline{\omega}}\right)\right)\\
=&\eta\left(U_1^2 S^3 T_1^{-3} S^3\left(\frac{w_1}{3}, \frac{w_2}{1-\overline{\omega}}\right)\right)\\
=&1^{-\frac{5}{27}} \root 3 \of{(1-\overline{\omega})
  w_1+(1-\omega) w_2+1} \eta_2(w_1, w_2).
\endaligned$$
$\qquad \qquad \qquad \qquad \qquad \qquad \qquad \qquad \qquad
 \qquad \qquad \qquad \qquad \qquad \qquad \qquad \qquad \qquad
 \quad \boxed{}$

  In fact,
$$\left\{\aligned
  \eta_3([T_1, T_2](w_1, w_2)) &=\eta_2(w_1, w_2),\\
  \eta_4(S^3 [T_1, T_2] S^3(w_1, w_2)) &=\eta_1(w_1, w_2),\\
  \eta_5(S^3 [T_1, T_2] S^3 [T_1, T_2](w_1, w_2)) &=\eta_2(w_1, w_2).
\endaligned\right.\tag 4.18$$
Writing explicitly,
$$\eta_3(w_1, w_2)=\eta\left(\left(\matrix
  1 &        0 & \omega-\overline{\omega}\\
    & 1-\omega &                        0\\
    &          &                        3
  \endmatrix\right)(w_1, w_2)\right)=\eta \left(\frac{w_1+
  (\omega-\overline{\omega})}{3}, \frac{w_2}{1-
  \overline{\omega}}\right),\tag 4.19$$
$$\aligned
  \eta_4(w_1, w_2)
=&\eta\left(\left(\matrix
  3 &        0 & 0\\
  0 & 1-\omega & 0\\
  \omega-\overline{\omega} & 0 & 1
  \endmatrix\right)(w_1, w_2)\right)\\
=&\eta \left(\frac{3 w_1}{(\omega-\overline{\omega}) w_1+1},
  \frac{(1-\omega) w_2}{(\omega-\overline{\omega}) w_1
  +1}\right),
\endaligned\tag 4.20$$
$$\aligned
  \eta_5(w_1, w_2)
=&\eta\left(\left(\matrix
  -2 & 0 & \omega-\overline{\omega}\\
   0 & 1-\omega & 0\\
  3(\omega-\overline{\omega}) & 0 & 3
  \endmatrix\right)(w_1, w_2)\right)\\
=&\eta \left(\frac{-2 w_1+(\omega-\overline{\omega})}
  {3(\omega-\overline{\omega}) w_1+3}, \frac{(1-\omega)
  w_2}{3 (\omega-\overline{\omega}) w_1+3}\right).
\endaligned\tag 4.21$$

{\smc Proposition 4.3}. {\it
$$\eta_1([T_1, T_2]^{-1}(w_1, w_2))=\eta_1(w_1, w_2).\tag 4.22$$
$$\eta_2([T_1, T_2]^{-3}(w_1, w_2))=\eta_2(w_1, w_2).\tag 4.23$$
$$\eta_3([T_1, T_2]^{-3}(w_1, w_2))=\eta_3(w_1, w_2).\tag 4.24$$
$$\eta_4(S^3 [T_1, T_2]^{-3} S^3(w_1, w_2))=\root 3 \of{\frac{4
 (\omega-\overline{\omega})w_1+1}{(\omega-\overline{\omega})w_1+1}}
 \eta_4(w_1, w_2).\tag 4.25$$
$$\eta_5(S^3 [T_1, T_2] S^3 [T_1, T_2]^{-3} S^3 [T_1, T_2]^{-1}
  S^3(w_1, w_2))=\eta_5(w_1, w_2).\tag 4.26$$}

{\it Proof}. We have
$$\aligned
 &\eta_1([T_1, T_2]^{-1}(w_1, w_2))\\
=&\eta\left(\left(\matrix
  3 &          &  \\
    & 1-\omega &  \\
    &          & 1
  \endmatrix\right) [T_1, T_2]^{-1}(w_1, w_2)\right)\\
=&\eta\left(\left(\matrix
  3 & 0        & 3(\omega-\overline{\omega})\\
    & 1-\omega & 0\\
    &          & 1
  \endmatrix\right)(w_1, w_2)\right)\\
=&\eta\left(\left(\matrix
  1 & 0 & 3(\omega-\overline{\omega})\\
    & 1 & 0\\
    &   & 1
  \endmatrix\right)
  \left(\matrix
  3 &          &  \\
    & 1-\omega &  \\
    &          & 1
  \endmatrix\right)(w_1, w_2)\right)\\
=&\eta([T_1, T_2]^{-3}(3w_1, (1-\omega)w_2))\\
=&\eta(3w_1, (1-\omega)w_2)\\
=&\eta_1(w_1, w_2).
\endaligned$$
$$\aligned
 &\eta_2([T_1, T_2]^{-3}(w_1, w_2))\\
=&\eta\left(\left(\matrix
  1 &          &  \\
    & 1-\omega &  \\
    &          & 3
  \endmatrix\right) [T_1, T_2]^{-3}(w_1, w_2)\right)\\
=&\eta\left(\left(\matrix
  1 & 0        & 3(\omega-\overline{\omega})\\
  0 & 1-\omega & 0\\
  0 & 0        & 3
  \endmatrix\right)(w_1, w_2)\right)\\
=&\eta\left(\left(\matrix
  1 & 0 & \omega-\overline{\omega}\\
    & 1 & 0\\
    &   & 1
 \endmatrix\right)
 \left(\matrix
 1 &          &  \\
   & 1-\omega &  \\
   &          & 3
 \endmatrix\right)(w_1, w_2)\right)\\
=&\eta\left([T_1, T_2]^{-1}\left(\frac{w_1}{3}, \frac{w_2}{1-
  \overline{\omega}}\right)\right)\\
=&\eta\left(\frac{w_1}{3}, \frac{w_2}{1-\overline{\omega}}\right)\\
=&\eta_2(w_1, w_2).
\endaligned$$
$$\aligned
 &\eta_3([T_1, T_2]^{-3}(w_1, w_2))\\
=&\eta_2([T_1, T_2]^{-1} [T_1, T_2]^{-3}(w_1, w_2))\\
=&\eta_2([T_1, T_2]^{-1}(w_1, w_2))\\
=&\eta_3(w_1, w_2).
\endaligned$$
$$\aligned
 &\eta_4(S^3 [T_1, T_2]^{-3} S^3(w_1, w_2))\\
=&\eta\left(\left(\matrix
  3 &          &  \\
    & 1-\omega &  \\
    &          & 1
  \endmatrix\right) S^3 [T_1, T_2]^{-4} S^3(w_1, w_2)\right)\\
=&\eta\left(\left(\matrix
  3 & 0 & 0\\
  0 & 1-\omega & 0\\
  4(\omega-\overline{\omega}) & 0 & 1
  \endmatrix\right)(w_1, w_2)\right)\\
=&\eta\left(\left(\matrix
  1 & 0 & 0\\
  0 & 1 & 0\\
  \omega-\overline{\omega} & 0 & 1
  \endmatrix\right) \left(\matrix
  3 & 0 & 0\\
  0 & 1-\omega & 0\\
  \omega-\overline{\omega} & 0 & 1
  \endmatrix\right)(w_1, w_2)\right)\\
=&\eta\left(S^3 [T_1, T_2]^{-1} S^3
  \left(\matrix
  3 & 0        & 0\\
  0 & 1-\omega & 0\\
  \omega-\overline{\omega} & 0 & 1
  \endmatrix\right)(w_1, w_2)\right)\\
=&\root 3 \of{\frac{4(\omega-\overline{\omega}) w_1+1}{(\omega-
  \overline{\omega}) w_1+1}} \eta\left(\left(\matrix
  3 & 0 & 0\\
  0 & 1-\omega & 0\\
  \omega-\overline{\omega} & 0 & 1
  \endmatrix\right)(w_1, w_2)\right)\\
=&\root 3 \of{\frac{4(\omega-\overline{\omega}) w_1+1}{(\omega-
  \overline{\omega}) w_1+1}} \eta_4(w_1, w_2).
\endaligned$$
$$\aligned
 &\eta_5(S^3 [T_1, T_2] S^3 [T_1, T_2]^{-3} S^3 [T_1, T_2]^{-1}
  S^3(w_1, w_2))\\
=&\eta_2([T_1, T_2]^{-4} S^3 [T_1, T_2]^{-1} S^3(w_1, w_2))\\
=&\eta_2([T_1, T_2]^{-1} S^3 [T_1, T_2]^{-1} S^3(w_1, w_2))\\
=&\eta_5(w_1, w_2).
\endaligned$$
$\qquad \qquad \qquad \qquad \qquad \qquad \qquad \qquad \qquad
 \qquad \qquad \qquad \qquad \qquad \qquad \qquad \qquad \qquad
 \quad \boxed{}$

  Set
$$\aligned
 &\varphi_0(w_1, w_2):=\frac{\eta_1(w_1, w_2)}{\eta_2(w_1, w_2)}, \quad
  \varphi_1(w_1, w_2):=\frac{\eta_1(w_1, w_2)}{\eta_3(w_1, w_2)},\\
 &\varphi_2(w_1, w_2):=\frac{\eta_4(w_1, w_2)}{\eta_2(w_1, w_2)}, \quad
  \varphi_3(w_1, w_2):=\frac{\eta_4(w_1, w_2)}{\eta_5(w_1, w_2)}.
\endaligned\tag 4.27$$

{\smc Proposition 4.4}. {\it
$$\left\{\aligned
 &\varphi_1(g_1(w_1, w_2))=1^{-\frac{4}{9}} \varphi_0(w_1, w_2),\\
 &\varphi_1(g_2(w_1, w_2))=1^{-\frac{8}{9}} \varphi_0(w_1, w_2),\\
 &\varphi_0(g_3(w_1, w_2))=\varphi_0(w_1, w_2),\\
 &\varphi_3(g_4(w_1, w_2))=\varphi_0(w_1, w_2),\\
 &\varphi_2(g_5(w_1, w_2))=1^{\frac{4}{9}} \root 3 \of{\frac{-3
  \overline{\omega} w_1+(1-\omega) w_2+1}{(1-\overline{\omega})
  w_1+(1-\omega) w_2+1}} \varphi_0(w_1, w_2).
\endaligned\right.\tag 4.28$$}

{\it Proof}. This is the corollary of Proposition 4.2.
\flushpar
$\qquad \qquad \qquad \qquad \qquad \qquad \qquad \qquad \qquad
 \qquad \qquad \qquad \qquad \qquad \qquad \qquad \qquad \qquad
 \quad \boxed{}$

{\smc Proposition 4.5}. {\it
$$\left\{\aligned
 &\varphi_1([T_1, T_2](w_1, w_2))=\varphi_0(w_1, w_2),\\
 &\varphi_2(S^3 [T_1, T_2] S^3(w_1, w_2))=\frac{1}{\root 3 \of{(\overline{\omega}-
  \omega) w_1+1}} \varphi_0(w_1, w_2),\\
 &\varphi_3(S^3 [T_1, T_2] S^3 [T_1, T_2](w_1, w_2))=\varphi_0(w_1, w_2).
\endaligned\right.\tag 4.29$$}

{\it Proof}. We have
$$\varphi_1([T_1, T_2](w_1, w_2))
 =\frac{\eta_1([T_1, T_2](w_1, w_2))}{\eta_3([T_1, T_2](w_1, w_2)}
 =\frac{\eta_1(w_1, w_2)}{\eta_2(w_1, w_2)}=\varphi_0(w_1, w_2).$$
$$\varphi_2(S^3[T_1, T_2]S^3(w_1, w_2))
 =\frac{\eta_4(S^3[T_1, T_2]S^3(w_1, w_2))}{\eta_2(S^3[T_1, T_2]S^3(w_1, w_2))}.$$
Here,
$$\eta_4(S^3[T_1, T_2]S^3(w_1, w_2))=\eta_1(w_1, w_2).$$
$$\aligned
 &\eta_2(S^3[T_1, T_2]S^3(w_1, w_2))\\
=&\eta\left(\left(\matrix
  1 &          &  \\
    & 1-\omega &  \\
    &          & 3
  \endmatrix\right) S^3 [T_1, T_2] S^3
  \left(\matrix
  1 &                    &            \\
    & \frac{1}{1-\omega} &            \\
    &                    & \frac{1}{3}
  \endmatrix\right)
  \left(\frac{w_1}{3}, \frac{w_2}{1-\overline{\omega}}\right)\right)\\
=&\eta\left(\left(\matrix
  1                           & 0 & 0\\
  0                           & 1 & 0\\
  3(\overline{\omega}-\omega) & 0 & 1
  \endmatrix\right)
  \left(\frac{w_1}{3}, \frac{w_2}{1-\overline{\omega}}\right)\right)\\
=&\eta\left(S^3[T_1, T_2]^3 S^3\left(\frac{w_1}{3},
  \frac{w_2}{1-\overline{\omega}}\right)\right)\\
=&\root 3 \of{(\overline{\omega}-\omega) w_1+1}
  \eta\left(\frac{w_1}{3}, \frac{w_2}{1-\overline{\omega}}\right)\\
=&\root 3 \of{(\overline{\omega}-\omega) w_1+1} \eta_2(w_1, w_2).
\endaligned$$
Thus,
$$\varphi_2(S^3[T_1, T_2]S^3(w_1, w_2))=\frac{1}{\root 3
  \of{(\overline{\omega}-\omega) w_1+1}} \varphi_0(w_1, w_2).$$

  Note that
$$\eta_4(S^3[T_1, T_2]S^3[T_1, T_2](w_1, w_2))
 =\eta_1([T_1, T_2](w_1, w_2))
 =\eta_1(w_1, w_2).$$
$$\eta_5(S^3[T_1, T_2]S^3[T_1, T_2](w_1, w_2))
 =\eta_2(w_1, w_2).$$
This gives that
$$\varphi_3(S^3[T_1, T_2]S^3[T_1, T_2](w_1, w_2))
 =\frac{\eta_1(w_1, w_2)}{\eta_2(w_1, w_2)}
 =\varphi_0(w_1, w_2).$$
$\qquad \qquad \qquad \qquad \qquad \qquad \qquad \qquad \qquad
 \qquad \qquad \qquad \qquad \qquad \qquad \qquad \qquad \qquad
 \quad \boxed{}$

{\smc Proposition 4.6}. {\it
$$\varphi_0([T_1, T_2]^{-3}(w_1, w_2))=\varphi_0(w_1, w_2).\tag 4.30$$
$$\varphi_1([T_1, T_2]^{-3}(w_1, w_2))=\varphi_1(w_1, w_2).\tag 4.31$$
$$\varphi_2(S^3 [T_1, T_2]^{-3} S^3(w_1, w_2))
 =\root 3 \of{\frac{4(\omega-\overline{\omega})w_1+1}{[(\omega-
  \overline{\omega})w_1+1][3(\omega-\overline{\omega})w_1+1]}}
  \varphi_2(w_1, w_2).\tag 4.32$$
$$\varphi_3(S^3 [T_1, T_2] S^3 [T_1, T_2]^{-3} S^3 [T_1, T_2]^{-1}
  S^3(w_1, w_2))=\varphi_3(w_1, w_2).\tag 4.33$$}

{\it Proof}. We have
$$\varphi_0([T_1, T_2]^{-3}(w_1, w_2))
 =\frac{\eta_1([T_1, T_2]^{-3}(w_1, w_2))}{\eta_2([T_1, T_2]^{-3}(w_1, w_2))}
 =\frac{\eta_1(w_1, w_2)}{\eta_2(w_1, w_2)}
 =\varphi_0(w_1, w_2).$$
$$\varphi_1([T_1, T_2]^{-3}(w_1, w_2))
 =\frac{\eta_1([T_1, T_2]^{-3}(w_1, w_2))}{\eta_3([T_1, T_2]^{-3}(w_1, w_2))}
 =\frac{\eta_1(w_1, w_2)}{\eta_3(w_1, w_2)}
 =\varphi_1(w_1, w_2).$$
$$\varphi_2(S^3 [T_1, T_2]^{-3} S^3(w_1, w_2))
 =\frac{\eta_4(S^3 [T_1, T_2]^{-3} S^3(w_1, w_2))}
  {\eta_2(S^3 [T_1, T_2]^{-3} S^3(w_1, w_2))}.$$
Here,
$$\eta_4(S^3 [T_1, T_2]^{-3} S^3(w_1, w_2))
 =\root 3 \of{\frac{4(\omega-\overline{\omega}) w_1+1}{(\omega-
  \overline{\omega}) w_1+1}} \eta_4(w_1, w_2).$$
$$\aligned
 &\eta_2(S^3 [T_1, T_2]^{-3} S^3(w_1, w_2))\\
=&\eta\left(\left(\matrix
  1 &          &  \\
    & 1-\omega &  \\
    &          & 3
  \endmatrix\right) S^3 [T_1, T_2]^{-3} S^3
  \left(\matrix
  1 &                    &            \\
    & \frac{1}{1-\omega} &            \\
    &                    & \frac{1}{3}
  \endmatrix\right)
  \left(\frac{w_1}{3}, \frac{w_2}{1-\overline{\omega}}\right)\right)\\
=&\eta\left(\left(\matrix
  1                           & 0 & 0\\
  0                           & 1 & 0\\
  9(\omega-\overline{\omega}) & 0 & 1
  \endmatrix\right)
  \left(\frac{w_1}{3}, \frac{w_2}{1-\overline{\omega}}\right)\right)\\
=&\eta\left(S^3 [T_1, T_2]^{-9} S^3 \left(\frac{w_1}{3},
  \frac{w_2}{1-\overline{\omega}}\right)\right)\\
=&\root 3 \of{3(\omega-\overline{\omega}) w_1+1} \eta\left(\frac{w_1}{3},
  \frac{w_2}{1-\overline{\omega}}\right)\\
=&\root 3 \of{3(\omega-\overline{\omega}) w_1+1} \eta_2(w_1, w_2).
\endaligned$$
Thus,
$$\varphi_2(S^3 [T_1, T_2]^{-3} S^3(w_1, w_2))
 =\root 3 \of{\frac{4(\omega-\overline{\omega}) w_1+1}{[(\omega-
  \overline{\omega}) w_1+1][3(\omega-\overline{\omega}) w_1+1]}}
  \varphi_2(w_1, w_2).$$
Note that
$$\aligned
 &\eta_4(S^3 [T_1, T_2] S^3 [T_1, T_2]^{-3} S^3 [T_1, T_2]^{-1} S^3(w_1, w_2))\\
=&\eta_1([T_1, T_2]^{-3} S^3 [T_1, T_2]^{-1} S^3(w_1, w_2))\\
=&\eta_1(S^3 [T_1, T_2]^{-1} S^3(w_1, w_2))\\
=&\eta_4(w_1, w_2).
\endaligned$$
$$\eta_5(S^3 [T_1, T_2] S^3 [T_1, T_2]^{-3} S^3 [T_1, T_2]^{-1} S^3(w_1, w_2))
 =\eta_5(w_1, w_2).$$
This implies that
$$\varphi_3(S^3 [T_1, T_2] S^3 [T_1, T_2]^{-3} S^3 [T_1, T_2]^{-1} S^3(w_1, w_2))
 =\frac{\eta_4(w_1, w_2)}{\eta_5(w_1, w_2)}
 =\varphi_3(w_1, w_2).$$
$\qquad \qquad \qquad \qquad \qquad \qquad \qquad \qquad \qquad
 \qquad \qquad \qquad \qquad \qquad \qquad \qquad \qquad \qquad
 \quad \boxed{}$

  Note that
$$\frac{\varphi_1(w_1, w_2)}{\varphi_2(w_1, w_2)}=\frac{\eta_1(w_1, w_2)
  \eta_2(w_1, w_2)}{\eta_3(w_1, w_2) \eta_4(w_1, w_2)}.$$

  Put
$$\left\{\aligned
  \varphi_0(w_1, w_2)=\frac{\eta_1(w_1, w_2)}{\eta_2(w_1, w_2)}
 &=:\kappa_0^{\frac{1}{9}}=:k_0^{\frac{1}{27}},\\
  \varphi_1(w_1, w_2)=\frac{\eta_1(w_1, w_2)}{\eta_3(w_1, w_2)}
 &=:\kappa_1^{\frac{1}{9}}=:k_1^{\frac{1}{27}},\\
  \varphi_2(w_1, w_2)=\frac{\eta_4(w_1, w_2)}{\eta_2(w_1, w_2)}
 &=:\kappa_2^{\frac{1}{9}}=:k_2^{\frac{1}{27}},\\
  \varphi_3(w_1, w_2)=\frac{\eta_4(w_1, w_2)}{\eta_5(w_1, w_2)}
 &=:\kappa_3^{\frac{1}{9}}=:k_3^{\frac{1}{27}}.
\endaligned\right.\tag 4.34$$
Here, there are four functions $\varphi_0$, $\varphi_1$,
$\varphi_2$ and $\varphi_3$, $C_4^2=6$. In fact, the dimension of
the moduli space of algebraic curves of genus three is:
$\text{dim} {\Cal M}_3=3 \times 3-3=6$. Namely, we define
$$\left\{\aligned
 K_1 &:=\int_{0}^{1} \frac{dx}{\root 3 \of{x^2 (1-x)(1-k_2 x)(1-k_3 x)}},\\
 K_2 &:=\int_{0}^{1} \frac{dx}{\root 3 \of{x^2 (1-x)(1-k_3 x)(1-k_1 x)}},\\
 K_3 &:=\int_{0}^{1} \frac{dx}{\root 3 \of{x^2 (1-x)(1-k_1 x)(1-k_2 x)}}.
\endaligned\right.\tag 4.35$$
$$\left\{\aligned
 K_4 &:=\int_{0}^{1} \frac{dx}{\root 3 \of{x^2 (1-x)(1-k_0 x)(1-k_1 x)}},\\
 K_5 &:=\int_{0}^{1} \frac{dx}{\root 3 \of{x^2 (1-x)(1-k_0 x)(1-k_2 x)}},\\
 K_6 &:=\int_{0}^{1} \frac{dx}{\root 3 \of{x^2 (1-x)(1-k_0 x)(1-k_3 x)}}.
\endaligned\right.\tag 4.36$$
Note that
$$\left\{\aligned
 K_1 &=3 \int_{0}^{1} \frac{dt}{\root 3 \of{(1-t^3)(1-k_2 t^3)(1-k_3 t^3)}},\\
 K_2 &=3 \int_{0}^{1} \frac{dt}{\root 3 \of{(1-t^3)(1-k_3 t^3)(1-k_1 t^3)}},\\
 K_3 &=3 \int_{0}^{1} \frac{dt}{\root 3 \of{(1-t^3)(1-k_1 t^3)(1-k_2 t^3)}}.
\endaligned\right.$$
$$\left\{\aligned
 K_4 &=3 \int_{0}^{1} \frac{dt}{\root 3 \of{(1-t^3)(1-k_0 t^3)(1-k_1 t^3)}},\\
 K_5 &=3 \int_{0}^{1} \frac{dt}{\root 3 \of{(1-t^3)(1-k_0 t^3)(1-k_2 t^3)}},\\
 K_6 &=3 \int_{0}^{1} \frac{dt}{\root 3 \of{(1-t^3)(1-k_0 t^3)(1-k_3 t^3)}}.
\endaligned\right.$$

\vskip 0.5 cm
\centerline{\bf 5. Transform problems and modular equations
                   associated to algebraic surfaces}
\vskip 0.5 cm

  For the group $GL(3)$, the two integrals are given by
$$\left\{\aligned
  K_3 &=\int_{0}^{1} \frac{dx}{\root 3 \of{x^2 (1-x)(1-\kappa_1^3
  x)(1-\kappa_2^3 x)}},\\
  K_6 &=\int_{0}^{1} \frac{dx}{\root 3 \of{x^2 (1-x)(1-\kappa_0^3
  x)(1-\kappa_3^3 x)}}.
\endaligned\right.$$

{\smc Theorem 5.1 (Main Theorem 3)}. {\it The following identity
of differential forms holds:
$$\frac{dw_1 \wedge dw_2}{\root 3 \of{w_1^2 w_2^2 (1-w_1)
  (1-\lambda_1^3 w_1)(1-\lambda_2^3 w_1)}}
 =-\frac{5 t_1}{t_2^2} \frac{dt_1 \wedge dt_2}{\root 3
  \of{t_1^2 t_2^2 (1-t_1)(1-\kappa_1^3 t_1)(1-\kappa_2^3 t_1)}},
  \tag 5.1$$
where
$$w_1=\frac{(\beta+\gamma) t_1+\alpha}{\beta t_1+(\alpha+\gamma)}, \quad
  w_2=\frac{t_1 [(\beta+\gamma) t_1+\alpha]^2 [\beta t_1+(\alpha+\gamma)]}
      {t_2^5},\tag 5.2$$
which is a rational transformation of order five. Here,
$$\left\{\aligned
  \alpha &=\frac{(v_1-1)(v_2-1)(v_1+v_2-2)-(u_1-1)(u_2-1)(u_1+u_2-2)}
           {2(u_1-1)(u_2-1)(v_1-1)(v_2-1)},\\
  \beta  &=\frac{-(v_1-1)(v_2-1)(2 v_1 v_2-v_1-v_2)+(u_1-1)(u_2-1)(2
           u_1 u_2-u_1-u_2)}{2(u_1-1)(u_2-1)(v_1-1)(v_2-1)},\\
  \gamma &=\frac{(v_1-1)(v_2-1)}{(u_1-1)(u_2-1)},
\endaligned\right.\tag 5.3$$
with
$$u_1=\kappa_1^3, \quad u_2=\kappa_2^3, \quad
  v_1=\lambda_1^3, \quad v_2=\lambda_2^3.$$
Moreover, the moduli $(\kappa_1, \kappa_2)$ and $(\lambda_1,
\lambda_2)$ satisfy the modular equation:
$$(\kappa_1^3-1)(\kappa_2^3-1)(\kappa_1^3-\kappa_2^3)
 =(\lambda_1^3-1)(\lambda_2^3-1)(\lambda_1^3-\lambda_2^3),\tag 5.4$$
which is an algebraic variety of dimension three. The
corresponding two algebraic surfaces are:
$$\left\{\aligned
 w_3^3 &=w_1^2 w_2^2 (1-w_1)(1-\lambda_1^3 w_1)(1-\lambda_2^3 w_1),\\
 t_3^3 &=t_1^2 t_2^2 (1-t_1)(1-\kappa_1^3 t_1)(1-\kappa_2^3 t_1),
\endaligned\right.\tag 5.5$$
which give the homogeneous algebraic equations of degree seven:
$$\left\{\aligned
 w_3^3 w_4^4 &=w_1^2 w_2^2 (w_4-w_1)(w_4-\lambda_1^3 w_1)(w_4-\lambda_2^3 w_1),\\
 t_3^3 t_4^4 &=t_1^2 t_2^2 (t_4-t_1)(t_4-\kappa_1^3 t_1)(t_4-\kappa_2^3 t_1).
\endaligned\right.\tag 5.6$$}

{\it Proof}. Let us consider the differential forms:
$$\frac{x dy \wedge dz+y dz \wedge dx+z dx \wedge dy}{\root 3 \of{f_j(x, y, z)}},$$
where $f_j(x, y, z)$ ($j=1, 2$) are homogeneous polynomials about
$x$, $y$ and $z$ of degree nine.

 Set
$$x=U(\xi, \zeta, \eta), \quad y=V(\xi, \zeta, \eta), \quad
  z=W(\xi, \zeta, \eta),$$
where $U$, $V$ and $W$ are homogeneous functions of degree $n$. By
$$\left\{\aligned
  n U&=\xi \frac{\partial U}{\partial \xi}+\zeta \frac{\partial U}{\partial \zeta}+
       \eta \frac{\partial U}{\partial \eta},\\
  n V&=\xi \frac{\partial V}{\partial \xi}+\zeta \frac{\partial V}{\partial \zeta}+
       \eta \frac{\partial V}{\partial \eta},\\
  n W&=\xi \frac{\partial W}{\partial \xi}+\zeta \frac{\partial W}{\partial \zeta}+
       \eta \frac{\partial W}{\partial \eta}
\endaligned\right.$$
and
$$\left\{\aligned
  d U&=d \xi \frac{\partial U}{\partial \xi}+d \zeta \frac{\partial U}{\partial \zeta}+
       d \eta \frac{\partial U}{\partial \eta},\\
  d V&=d \xi \frac{\partial V}{\partial \xi}+d \zeta \frac{\partial V}{\partial \zeta}+
       d \eta \frac{\partial V}{\partial \eta},\\
  d W&=d \xi \frac{\partial W}{\partial \xi}+d \zeta \frac{\partial W}{\partial \zeta}+
       d \eta \frac{\partial W}{\partial \eta},
\endaligned\right.$$
we have
$$U dV \wedge dW+V dW \wedge dU+W dU \wedge dV=H (\xi d \zeta \wedge d \eta+
  \zeta d \eta \wedge d \xi+\eta d \xi \wedge d \zeta),$$
where
$$H=\frac{1}{n} \vmatrix \format \c \quad & \c \quad & \c\\
  \frac{\partial U}{\partial \xi} & \frac{\partial V}{\partial \xi} &
  \frac{\partial W}{\partial \xi}\\
  \frac{\partial U}{\partial \zeta} & \frac{\partial V}{\partial \zeta} &
  \frac{\partial W}{\partial \zeta}\\
  \frac{\partial U}{\partial \eta} & \frac{\partial V}{\partial \eta} &
  \frac{\partial W}{\partial \eta}\\
  \endvmatrix.$$
Hence,
$$\frac{x dy \wedge dz+y dz \wedge dx+z dx \wedge dy}{\root 3 \of{f_j(x, y, z)}}
 =H \frac{\xi d \zeta \wedge d \eta+\zeta d \eta \wedge d \xi+\eta
  d \xi \wedge d \zeta}{\root 3 \of{f_j(U, V, W)}}.$$
Here, the functions $f_j(U, V, W)$ ($j=1, 2$) are of degree $9n$.
Let
$$f_j(U, V, W)=T_j^3 \varphi_j(\xi, \zeta, \eta),$$
where $T_j^3$ ($j=1, 2$) are cubic factors of degree $9n-9$. Put
$M_j=\frac{T_j}{H}$ ($j=1, 2$). Then
$$\frac{x dy \wedge dz+y dz \wedge dx+z dx \wedge dy}{\root 3 \of{f_j(x, y, z)}}
 =\frac{1}{M_j} \frac{\xi d \zeta \wedge d \eta+\zeta d \eta \wedge d \xi+\eta
  d \xi \wedge d \zeta}{\root 3 \of{\varphi_j(\xi, \zeta, \eta)}},$$
where $\varphi_j(\xi, \zeta, \eta)$ ($j=1, 2$) are homogeneous
functions of degree $9$.

  Put
$$\left\{\aligned
 f_1(x, y, z) &:=x^2 y^2 z^2 (z-x)(z-\lambda_1^3 x)(z-\lambda_2^3 x),\\
 f_2(x, y, z) &:=x^2 y^2 z^2 (z-y)(z-\lambda_0^3 y)(z-\lambda_3^3 y).
\endaligned\right.$$
Note that
$$x dy \wedge dz+y dz \wedge dx+z dx \wedge dy
 =z^3 d (\frac{x}{z}) \wedge d (\frac{y}{z}).$$
We will study the following two differential forms:
$$\omega_1:=\frac{x dy \wedge dz+y dz \wedge dx+z dx \wedge dy}
            {\root 3 \of{f_1(x, y, z)}}, \quad
  \omega_2:=\frac{x dy \wedge dz+y dz \wedge dx+z dx \wedge dy}
            {\root 3 \of{f_2(x, y, z)}}.$$
In fact,
$$\omega_1=\frac{d (\frac{x}{z}) \wedge d (\frac{y}{z})}{\root 3
  \of{(\frac{x}{z})^2 (\frac{y}{z})^2 (1-\frac{x}{z})
  (1-\lambda_1^3 \frac{x}{z})(1-\lambda_2^3 \frac{x}{z})}},$$
$$\omega_2=\frac{d (\frac{x}{z}) \wedge d (\frac{y}{z})}{\root 3
  \of{(\frac{x}{z})^2 (\frac{y}{z})^2 (1-\frac{y}{z})
  (1-\lambda_0^3 \frac{y}{z})(1-\lambda_3^3 \frac{y}{z})}}.$$
Set $\frac{x}{z}=u^3$ and $\frac{y}{z}=v^3$. Then
$$\omega_1=\frac{9 du \wedge dv}{\root 3 \of{(1-u^3)(1-\lambda_1^3
           u^3)(1-\lambda_2^3 u^3)}}, \quad
  \omega_2=\frac{9 du \wedge dv}{\root 3 \of{(1-v^3)(1-\lambda_0^3
           v^3)(1-\lambda_3^3 v^3)}}.$$

  For
$$\left\{\aligned
  f_1(x, y, z) &=x^2 y^2 z^2 (z-x)(z-\lambda_1^3 x)(z-\lambda_2^3 x),\\
  f_2(x, y, z) &=x^2 y^2 z^2 (z-y)(z-\lambda_0^3 y)(z-\lambda_3^3 y),\\
  \varphi_1(\xi, \zeta, \eta) &=\xi^2 \zeta^2 \eta^2 (\eta-\xi)(\eta-
  \kappa_1^3 \xi)(\eta-\kappa_2^3 \xi),\\
  \varphi_2(\xi, \zeta, \eta) &=\xi^2 \zeta^2 \eta^2 (\eta-\zeta)(\eta-
  \kappa_0^3 \zeta)(\eta-\kappa_3^3 \zeta),
\endaligned\right.$$
consider the following equations:
$$U^2 V^2 W^2 (W-U)(W-\lambda_1^3 U)(W-\lambda_2^3 U)
 =T_1^3 \xi^2 \zeta^2 \eta^2 (\eta-\xi)(\eta-\kappa_1^3 \xi)(\eta-\kappa_2^3 \xi)$$
and
$$U^2 V^2 W^2 (W-V)(W-\lambda_0^3 V)(W-\lambda_3^3 V)
 =T_2^3 \xi^2 \zeta^2 \eta^2 (\eta-\zeta)(\eta-\kappa_0^3 \zeta)(\eta-\kappa_3^3 \zeta).$$
The exponents of $U$, $V$, $W$ at $f_1(U, V, W)$ are $5:2:5$, the
exponents of $U$, $V$, $W$ at $f_2(U, V, W)$ are $2:5:5$. On the
other hand, $T_1^3$ and $T_2^3$ are of degree $9n-9$. This leads
us to consider the system of equations:
$$\aligned
 &U^2 V^2 W^2 (W-U)(W-\lambda_1^3 U)(W-\lambda_2^3 U)\\
=&(a \xi+b \eta)^{6} \zeta^{24} (c \xi+d \eta)^{6} \xi^2 \zeta^2
  \eta^2 (\eta-\xi)(\eta-\kappa_1^3 \xi)(\eta-\kappa_2^3 \xi)
\endaligned$$
and
$$\aligned
 &U^2 V^2 W^2 (W-V)(W-\lambda_0^3 V)(W-\lambda_3^3 V)\\
=&\xi^{24} (a \zeta+b \eta)^{6} (c \zeta+d \eta)^{6} \xi^2 \zeta^2
  \eta^2 (\eta-\zeta)(\eta-\kappa_0^3 \zeta)(\eta-\kappa_3^3 \zeta).
\endaligned$$

  For the equation
$$U^2 V^2 W^2 (W-U)(W-\lambda_1^3 U)(W-\lambda_2^3 U)
 =T_1^3 \xi^2 \zeta^2 \eta^2 (\eta-\xi)(\eta-\kappa_1^3 \xi)(\eta-\kappa_2^3 \xi),$$
put
$$\left\{\aligned
  U &=\zeta^{6} (a \xi+b \eta),\\
  V &=\xi \zeta \eta (a \xi+b \eta)^{2} (c \xi+d \eta)^{2},\\
  W &=\zeta^{6} (c \xi+d \eta).
\endaligned\right.$$
Here,
$$\left(\frac{U}{W}, \frac{V}{W}\right)=\left(\frac{a \xi+b \eta}{c \xi+d \eta},
  \frac{\xi \eta (a \xi+b \eta)^{2}(c \xi+d \eta)}{\zeta^{5}}\right).$$
This is a transform of order $5$.

  Let
$$\xi \mapsto A \xi+B \eta, \quad \zeta \mapsto \zeta, \quad
  \eta \mapsto C \xi+D \eta,$$
where
$$\left\{\aligned
 A&=-\kappa_1^6 \kappa_2^3-\kappa_1^3-\kappa_2^6+3 \kappa_1^3 \kappa_2^3,\\
 B&=\kappa_1^6+\kappa_2^6-\kappa_1^3 \kappa_2^3-\kappa_1^3-\kappa_2^3+1,\\
 C&=-\kappa_1^6 \kappa_2^6+\kappa_1^6 \kappa_2^3+\kappa_1^3 \kappa_2^6
    -\kappa_1^6-\kappa_2^6+\kappa_1^3 \kappa_2^3,\\
 D&=\kappa_1^3 \kappa_2^6+\kappa_1^6+\kappa_2^3-3 \kappa_1^3 \kappa_2^3.
\endaligned\right.$$
Then
$$\left\{\aligned
  \eta-\xi &\mapsto (C-A) \xi+(D-B) \eta=(\kappa_1^3-1)(\kappa_2^3-1)^2
                    (\eta-\kappa_1^3 \xi),\\
  \eta-\kappa_1^3 \xi &\mapsto (C-\kappa_1^3 A) \xi+(D-\kappa_1^3 B) \eta
         =(\kappa_2^3-\kappa_1^3)(\kappa_1^3-1)^2 (\eta-\kappa_2^3 \xi),\\
  \eta-\kappa_2^3 \xi &\mapsto (C-\kappa_2^3 A) \xi+(D-\kappa_2^3 B) \eta
         =(1-\kappa_2^3)(\kappa_1^3-\kappa_2^3)^2 (\eta-\xi).
\endaligned\right.$$
Note that
$$A+B-C-D=(\kappa_1^3-1)^2 (\kappa_2^3-1)^2.$$

   We are interested in the three quotients:
$$\frac{W-U}{\eta-\xi}, \quad \frac{W-\lambda_1^3 U}{\eta-\kappa_1^3 \xi},
  \quad \frac{W-\lambda_2^3 U}{\eta-\kappa_2^3 \xi}.$$
In fact, when $d-b=a-c$, i.e., set
$$b=\alpha, \quad d=\alpha+\gamma, \quad c=\beta, \quad
  a=\beta+\gamma,$$
then
$$\frac{W-U}{\eta-\xi}=\gamma \zeta^{6}.$$
Under the above action,
$$\aligned
 U \mapsto \zeta^{6} [a(A \xi+B \eta)+b (C \xi+D \eta)],\\
 W \mapsto \zeta^{6} [c(A \xi+B \eta)+d (C \xi+D \eta)].
\endaligned$$
Put
$$\left\{\aligned
 A &=(\lambda_2^3-1)^2 (c-\lambda_1^3 a) (d-\lambda_1^3 b)-(\lambda_2^3
     -\lambda_1^3) (\lambda_1^3-1) (c-\lambda_2^3 a) (d-b),\\
 B &=(\lambda_2^3-1)^2 (d-\lambda_1^3 b)^2-(\lambda_2^3-\lambda_1^3)
     (\lambda_1^3-1) (d-b) (d-\lambda_2^3 b),\\
 C &=(\lambda_2^3-\lambda_1^3) (\lambda_1^3-1) (c-a) (c-\lambda_2^3 a)-
     (\lambda_2^3-1)^2 (c-\lambda_1^3 a)^2,\\
 D &=(\lambda_2^3-\lambda_1^3) (\lambda_1^3-1) (c-a) (d-\lambda_2^3 b)-
     (\lambda_2^3-1)^2 (c-\lambda_1^3 a) (d-\lambda_1^3 b).
\endaligned\right.$$
Then
$$\left\{\aligned
 W-U &\mapsto (\lambda_1^3-1)(\lambda_2^3-1)^2 (ad-bc)(W-\lambda_1^3 U),\\
 W-\lambda_1^3 U &\mapsto (\lambda_2^3-\lambda_1^3)(\lambda_1^3-1)^2
                  (ad-bc)(W-\lambda_2^3 U),\\
 W-\lambda_2^3 U &\mapsto (1-\lambda_2^3) (\lambda_1^3-\lambda_2^3)^2
                  (ad-bc)(W-U).
\endaligned\right.$$
When $ad-bc=1$, we have
$$\left\{\aligned
 \frac{W-U}{\eta-\xi} &\mapsto \frac{(\lambda_1^3-1)(\lambda_2^3-1)^2}
 {(\kappa_1^3-1)(\kappa_2^3-1)^2} \frac{W-\lambda_1^3 U}{\eta-\kappa_1^3 \xi},\\
 \frac{W-\lambda_1^3 U}{\eta-\kappa_1^3 \xi} &\mapsto \frac{(\lambda_2^3-
 \lambda_1^3)(\lambda_1^3-1)^2}{(\kappa_2^3-\kappa_1^3)(\kappa_1^3-1)^2}
 \frac{W-\lambda_2^3 U}{\eta-\kappa_2^3 \xi},\\
 \frac{W-\lambda_2^3 U}{\eta-\kappa_2^3 \xi} &\mapsto \frac{(1-\lambda_2^3)
 (\lambda_1^3-\lambda_2^3)^2}{(1-\kappa_2^3)(\kappa_1^3-\kappa_2^3)^2} \frac{W-U}
 {\eta-\xi}.
\endaligned\right.$$
Note that
$$A+B-C-D=(\alpha+\beta+\gamma)^2 (\lambda_1^3-1)^2 (\lambda_2^3-1)^2.$$
Two different expressions about $A+B-C-D$ imply that
$$(\alpha+\beta+\gamma)^2 (\lambda_1^3-1)^2 (\lambda_2^3-1)^2
 =(\kappa_1^3-1)^2 (\kappa_2^3-1)^2.$$
Set
$$u_1=\kappa_1^3, \quad u_2=\kappa_2^3, \quad
  v_1=\lambda_1^3, \quad v_2=\lambda_2^3.$$
Then
$$\left\{\aligned
 A &=-u_1^2 u_2-u_1-u_2^2+3 u_1 u_2,\\
 B &=u_1^2+u_2^2-u_1 u_2-u_1-u_2+1,\\
 C &=-u_1^2 u_2^2+u_1^2 u_2+u_1 u_2^2-u_1^2-u_2^2+u_1 u_2,\\
 D &=u_1 u_2^2+u_1^2+u_2-3 u_1 u_2,
\endaligned\right.$$
and
$$\left\{\aligned
 A=&(v_1-1)^2 (v_2-1)^2 \alpha \beta+v_1 (v_1-1) (v_2-1)^2 \alpha \gamma
   -(v_1-1)^2 (v_2-1) \beta \gamma+\\
   &-(v_1^2 v_2+v_1+v_2^2-3 v_1 v_2) \gamma^2,\\
 B=&(v_1-1)^2 (v_2-1)^2 \alpha^2-[(v_1-1)^2 (v_2-1)+(v_1-1) (v_2-1)^2]
    \alpha \gamma+\\
   &+(v_1^2+v_2^2-v_1 v_2-v_1-v_2+1) \gamma^2,\\
 C=&-[v_1 (v_1-1) (v_2-1)^2+v_2 (v_2-1) (v_1-1)^2] \beta \gamma-(v_1-1)^2
    (v_2-1)^2 \beta^2+\\
   &-(v_1^2 v_2^2-v_1^2 v_2-v_1 v_2^2+v_1^2+v_2^2-v_1 v_2) \gamma^2,\\
 D=&-(v_1-1)^2 (v_2-1)^2 \alpha \beta+(v_1-1) (v_2-1)^2 \beta \gamma-v_2
    (v_2-1) (v_1-1)^2 \alpha \gamma+\\
   &+(v_1 v_2^2+v_1^2+v_2-3 v_1 v_2) \gamma^2.
\endaligned\right.$$

  Note that
$$ad-bc=(\beta+\gamma)(\alpha+\gamma)-\alpha \beta
 =(\alpha+\beta+\gamma) \gamma.$$
The condition $ad-bc=1$ implies that
$$(\alpha+\beta+\gamma) \gamma=1.$$
On the other hand,
$$\alpha+\beta+\gamma=\pm \frac{(\kappa_1^3-1)(\kappa_2^3-1)}
  {(\lambda_1^3-1)(\lambda_2^3-1)}.$$
Without loss of generality, we can assume that
$$\alpha+\beta+\gamma=\frac{(u_1-1)(u_2-1)}{(v_1-1)(v_2-1)}.$$
Hence,
$$\gamma=\frac{(v_1-1)(v_2-1)}{(u_1-1)(u_2-1)}.$$
Substituting the expression of $\gamma$ to the equations about $B$
and $C$, we get the following two quadratic equations about
$\alpha$ and $\beta$ respectively:
$$\aligned
  &(u_1-1)^2 (u_2-1)^2 (v_1-1)^2 (v_2-1)^2 \alpha^2+\\
 -&(u_1-1)(u_2-1)(v_1-1)(v_2-1)[(v_1-1)^2 (v_2-1)+(v_1-1)(v_2-1)^2] \alpha+\\
 +&[(v_1-1)^2 (v_2-1)^2 (v_1^2+v_2^2-v_1 v_2-v_1-v_2+1)+\\
 -&(u_1-1)^2 (u_2-1)^2 (u_1^2+u_2^2-u_1 u_2-u_1-u_2+1)]=0
\endaligned$$
and
$$\aligned
  &(u_1-1)^2 (u_2-1)^2 (v_1-1)^2 (v_2-1)^2 \beta^2+\\
 +&(u_1-1)(u_2-1)(v_1-1)(v_2-1)[v_1(v_1-1)(v_2-1)^2+v_2(v_2-1)(v_1-1)^2] \beta+\\
 +&[(v_1-1)^2 (v_2-1)^2 (v_1^2 v_2^2-v_1^2 v_2-v_1 v_2^2+v_1^2+v_2^2-v_1 v_2)+\\
 -&(u_1-1)^2 (u_2-1)^2 (u_1^2 u_2^2-u_1^2 u_2-u_1 u_2^2+u_1^2+u_2^2-u_1 u_2)]=0.
\endaligned$$
Thus,
$$\left\{\aligned
 \alpha &=\frac{(v_1-1)(v_2-1)(v_1+v_2-2) \pm \sqrt{\Delta_1}}
          {2(u_1-1)(u_2-1)(v_1-1)(v_2-1)},\\
  \beta &=\frac{-(v_1-1)(v_2-1)(2 v_1 v_2-v_1-v_2) \pm \sqrt{\Delta_2}}
          {2(u_1-1)(u_2-1)(v_1-1)(v_2-1)},
\endaligned\right.$$
where
$$\left\{\aligned
 \Delta_1= &4(u_1-1)^2 (u_2-1)^2 (u_1^2+u_2^2-u_1 u_2-u_1-u_2+1)+\\
          -&3(v_1-1)^2 (v_2-1)^2 (v_1-v_2)^2,\\
 \Delta_2= &4(u_1-1)^2 (u_2-1)^2 (u_1^2 u_2^2-u_1^2 u_2-u_1 u_2^2+u_1^2+u_2^2-u_1 u_2)+\\
          -&3(v_1-1)^2 (v_2-1)^2 (v_1-v_2)^2.
\endaligned\right.$$
Substituting the expressions of $\alpha$, $\beta$ and $\gamma$ to
the equations about $A$ and $D$, we get the next two equations
about $\Delta_1$ and $\Delta_2$:
$$\aligned
  &\pm \sqrt{\Delta_1} \cdot \pm \sqrt{\Delta_2}-(v_1-1)(v_2-1)(v_1-v_2)
   (\pm \sqrt{\Delta_1} \pm \sqrt{\Delta_2})+\\
 +&4 (u_1-1)^2 (u_2-1)^2 (u_1^2 u_2+u_1+u_2^2-3 u_1 u_2)-3 (v_1-1)^2
   (v_2-1)^2 (v_1-v_2)^2=0,
\endaligned$$
$$\aligned
  &\pm \sqrt{\Delta_1} \cdot \pm \sqrt{\Delta_2}+(v_1-1)(v_2-1)(v_1-v_2)
   (\pm \sqrt{\Delta_1} \pm \sqrt{\Delta_2})+\\
 +&4 (u_1-1)^2 (u_2-1)^2 (u_1 u_2^2+u_1^2+u_2-3 u_1 u_2)-3 (v_1-1)^2
   (v_2-1)^2 (v_1-v_2)^2=0.
\endaligned$$
Set $X=\pm \sqrt{\Delta_1}$ and $Y=\pm \sqrt{\Delta_2}$. We have
$$\left\{\aligned
  X+Y= &\frac{2(u_1-1)^3 (u_2-1)^3 (u_1-u_2)}{(v_1-1)(v_2-1)(v_1-v_2)},\\
   XY= &-2(u_1-1)^2 (u_2-1)^2 (u_1^2 u_2+u_1 u_2^2+u_1^2+u_2^2-6 u_1 u_2
        +u_1+u_2)+\\
       &+3(v_1-1)^2 (v_2-1)^2 (v_1-v_2)^2.
\endaligned\right.$$

  Note that
$$(X+Y)^2-2 XY=X^2+Y^2=\Delta_1+\Delta_2.$$
This gives a equation which connects the variables $u_1$, $u_2$
and $v_1$, $v_2$:
$$(u_1-1)^2 (u_2-1)^2 (u_1-u_2)^2=(v_1-1)^2 (v_2-1)^2 (v_1-v_2)^2.$$
Similarly, the relation
$$(XY)^2=X^2 Y^2=\Delta_1 \Delta_2$$
gives the other equation about $u_1$, $u_2$ and $v_1$, $v_2$:
$$\aligned
 &4(u_1^2+u_2^2-u_1 u_2-u_1-u_2+1)(u_1^2 u_2^2-u_1^2 u_2-u_1 u_2^2
 +u_1^2+u_2^2-u_1 u_2)+\\
 &-(u_1^2 u_2+u_1 u_2^2+u_1^2+u_2^2-6 u_1 u_2+u_1+u_2)^2\\
=&3(v_1-1)^2 (v_2-1)^2 (v_1-v_2)^2.
\endaligned$$
The identity
$$\aligned
 &4(u_1^2+u_2^2-u_1 u_2-u_1-u_2+1)(u_1^2 u_2^2-u_1^2 u_2-u_1 u_2^2
 +u_1^2+u_2^2-u_1 u_2)+\\
 &-(u_1^2 u_2+u_1 u_2^2+u_1^2+u_2^2-6 u_1 u_2+u_1+u_2)^2\\
=&3(u_1-1)^2 (u_2-1)^2 (u_1-u_2)^2
\endaligned$$
implies that the above two equations about $u_1$, $u_2$ and $v_1$,
$v_2$ are the same one:
$$(u_1-1)^2 (u_2-1)^2 (u_1-u_2)^2=(v_1-1)^2 (v_2-1)^2 (v_1-v_2)^2.$$

  In fact, the relation $\frac{W-U}{\eta-\xi}=\gamma \zeta^6$ gives that
$$\left\{\aligned
  \frac{W-\lambda_1^3 U}{\eta-\kappa_1^3 \xi}
=&\gamma \frac{(\kappa_1^3-1)(\kappa_2^3-1)^2}{(\lambda_1^3-1)
  (\lambda_2^3-1)^2} \zeta^{6},\\
  \frac{W-\lambda_2^3 U}{\eta-\kappa_2^3 \xi}
=&\gamma \frac{(\kappa_1^3-1)^3 (\kappa_2^3-1)^2 (\kappa_1^3-
  \kappa_2^3)}{(\lambda_1^3-1)^3 (\lambda_2^3-1)^2 (\lambda_1^3-
  \lambda_2^3)} \zeta^{6},\\
  \frac{W-U}{\eta-\xi}
=&\gamma \frac{(\kappa_1^3-1)^3 (\kappa_2^3-1)^3(\kappa_1^3-
  \kappa_2^3)^3}{(\lambda_1^3-1)^3 (\lambda_2^3-1)^3 (\lambda_1^3-
  \lambda_2^3)^3} \zeta^{6}.
\endaligned\right.$$
This implies that
$$(\kappa_1^3-1)^3 (\kappa_2^3-1)^3 (\kappa_1^3-\kappa_2^3)^3
 =(\lambda_1^3-1)^3 (\lambda_2^3-1)^3 (\lambda_1^3-\lambda_2^3)^3.$$
Therefore, we get the modular equation:
$$(\kappa_1^3-1)(\kappa_2^3-1)(\kappa_1^3-\kappa_2^3)
 =(\lambda_1^3-1)(\lambda_2^3-1)(\lambda_1^3-\lambda_2^3).$$

  The modular equation
$$(u_1-1)(u_2-1)(u_1-u_2)=(v_1-1)(v_2-1)(v_1-v_2)$$
implies that
$$\Delta_1=(u_1-1)^2 (u_2-1)^2 (u_1+u_2-2)^2, \quad
  \Delta_2=(u_1-1)^2 (u_2-1)^2 (2 u_1 u_2-u_1-u_2)^2.$$
Thus,
$$\alpha=\frac{(v_1-1)(v_2-1)(v_1+v_2-2) \pm (u_1-1)(u_2-1)(u_1+u_2-2)}
         {2(u_1-1)(u_2-1)(v_1-1)(v_2-1)},$$
$$\beta=\frac{-(v_1-1)(v_2-1)(2 v_1 v_2-v_1-v_2) \pm (u_1-1)(u_2-1)(2
        u_1 u_2-u_1-u_2)}{2(u_1-1)(u_2-1)(v_1-1)(v_2-1)}.$$

  We have
$$\aligned
  T_1^3
=&\frac{U^2 V^2 W^2}{\xi^2 \zeta^2 \eta^2} \cdot
  \frac{W-U}{\eta-\xi} \cdot \frac{W-\lambda_1^3 U}{\eta-\kappa_1^3 \xi}
  \cdot \frac{W-\lambda_2^3 U}{\eta-\kappa_2^3 \xi}\\
=&\zeta^{42} (a \xi+b \eta)^{6} (c \xi+d \eta)^{6} \gamma^{3}
  \frac{(\kappa_1^3-1)^{3} (\kappa_2^3-1)^{3}}{(\lambda_1^3-1)^{3}
  (\lambda_2^3-1)^{3}}.
\endaligned$$
Hence,
$$T_1=\gamma \frac{(u_1-1)(u_2-1)}{(v_1-1)(v_2-1)} \zeta^{14} (a \xi+b
      \eta)^{2} (c \xi+d \eta)^{2}
     =\zeta^{14} (a \xi+b \eta)^{2} (c \xi+d \eta)^{2} .$$
On the other hand,
$$H=\frac{1}{7} \vmatrix \format \c \quad & \c \quad & \c\\
    \frac{\partial U}{\partial \xi} & \frac{\partial V}{\partial \xi} &
    \frac{\partial W}{\partial \xi}\\
    \frac{\partial U}{\partial \zeta} & \frac{\partial V}{\partial \zeta} &
    \frac{\partial W}{\partial \zeta}\\
    \frac{\partial U}{\partial \eta} & \frac{\partial V}{\partial \eta} &
    \frac{\partial W}{\partial \eta}\\
    \endvmatrix
  =-5 \xi \eta \zeta^{12} (a \xi+b \eta)^{2} (c \xi+d \eta)^{2}.$$
Thus,
$$M_1=\frac{T_1}{H}=-\frac{1}{5} \frac{\zeta^2}{\xi \eta}.$$
Set
$$w_1=\frac{x}{z}, \quad w_2=\frac{y}{z}, \quad
  t_1=\frac{\xi}{\eta}, \quad t_2=\frac{\zeta}{\eta}.$$
Then $M=-\frac{1}{5} \frac{t_2^2}{t_1}$ and
$$w_1=\frac{U}{W}=\frac{(\beta+\gamma) t_1+\alpha}{\beta t_1+(\alpha+\gamma)}, \quad
  w_2=\frac{V}{W}=\frac{t_1 [(\beta+\gamma) t_1+\alpha]^2 [\beta t_1+(\alpha+\gamma)]}
      {t_2^5}.$$
Moreover,
$$\frac{x dy \wedge dz+y dz \wedge dx+z dx \wedge dy}{\root 3
  \of{x^2 y^2 z^2 (z-x)(z-\lambda_1^3 x)(z-\lambda_2^3 x)}}
 =\frac{dw_1 \wedge dw_2}{\root 3 \of{w_1^2 w_2^2 (1-w_1)
  (1-v_1 w_1)(1-v_2 w_1)}}.$$
$$\frac{\xi d \zeta \wedge d \eta+\zeta d \eta \wedge d \xi+\eta
  d \xi \wedge d \zeta}{\root 3 \of{\xi^2 \zeta^2 \eta^2 (\eta-\xi)
  (\eta-\kappa_1^3 \xi)(\eta-\kappa_2^3 \xi)}}
 =\frac{dt_1 \wedge dt_2}{\root 3 \of{t_1^2 t_2^2 (1-t_1)(1-u_1
  t_1)(1-u_2 t_1)}}.$$
Therefore,
$$\frac{dw_1 \wedge dw_2}{\root 3 \of{w_1^2 w_2^2 (1-w_1)
  (1-v_1 w_1)(1-v_2 w_1)}}
 =-\frac{5 t_1}{t_2^2} \frac{dt_1 \wedge dt_2}{\root 3
  \of{t_1^2 t_2^2 (1-t_1)(1-u_1 t_1)(1-u_2 t_1)}}.$$

  Note that the condition $(\alpha+\beta+\gamma) \gamma=1$ gives
that
$$\left\{\aligned
  \alpha &=\frac{(v_1-1)(v_2-1)(v_1+v_2-2)-(u_1-1)(u_2-1)(u_1+u_2-2)}
           {2(u_1-1)(u_2-1)(v_1-1)(v_2-1)},\\
  \beta  &=\frac{-(v_1-1)(v_2-1)(2 v_1 v_2-v_1-v_2)+(u_1-1)(u_2-1)(2
           u_1 u_2-u_1-u_2)}{2(u_1-1)(u_2-1)(v_1-1)(v_2-1)},\\
  \gamma &=\frac{(v_1-1)(v_2-1)}{(u_1-1)(u_2-1)}.
\endaligned\right.$$

  For
$$U^2 V^2 W^2 (W-V)(W-\lambda_0^3 V)(W-\lambda_3^3 V)
 =T_2^3 \xi^2 \zeta^2 \eta^2 (\eta-\zeta)(\eta-\kappa_0^3 \zeta)(\eta-\kappa_3^3 \zeta),$$
put
$$\left\{\aligned
  U &=\xi \zeta \eta (a \zeta+b \eta)^{2} (c \zeta+d \eta)^{2},\\
  V &=\xi^{6} (a \zeta+b \eta),\\
  W &=\xi^{6} (c \zeta+d \eta).
\endaligned\right.$$
Here,
$$\left(\frac{U}{W}, \frac{V}{W}\right)=
  \left(\frac{\zeta \eta (a \zeta+b \eta)^{2} (c \zeta+d \eta)}{\xi^{5}},
  \frac{a \zeta+b \eta}{c \zeta+d \eta}\right).$$
This is a transform of order $5$. The other argument is similar as
above. This completes the proof of Main Theorem 3.
\flushpar
$\qquad \qquad \qquad \qquad \qquad \qquad \qquad \qquad \qquad
 \qquad \qquad \qquad \qquad \qquad \qquad \qquad \qquad \qquad
 \quad \boxed{}$

{\smc Corollary 5.2}. {\it The following identity holds:
$$\frac{dy_1 \wedge dy_2}{\root 3 \of{(1-y_1^3)(1-\lambda_1^3
  y_1^3)(1-\lambda_2^3 y_1^3)}}=-\frac{5 x_1^3}{x_2^6} \frac{dx_1
  \wedge dx_2}{\root 3 \of{(1-x_1^3)(1-\kappa_1^3 x_1^3)(1-
  \kappa_2^3 x_1^3)}},\tag 5.7$$
where
$$y_1^3=\frac{(\beta+\gamma) x_1^3+\alpha}{\beta x_1^3+(\alpha+\gamma)}, \quad
  y_2^3=\frac{x_1^3 [(\beta+\gamma) x_1^3+\alpha]^2 [\beta x_1^3+(\alpha+\gamma)]}
      {x_2^{15}},\tag 5.8$$
and
$$\left\{\aligned
  \alpha &=\frac{(v_1-1)(v_2-1)(v_1+v_2-2)-(u_1-1)(u_2-1)(u_1+u_2-2)}
           {2(u_1-1)(u_2-1)(v_1-1)(v_2-1)},\\
  \beta  &=\frac{-(v_1-1)(v_2-1)(2 v_1 v_2-v_1-v_2)+(u_1-1)(u_2-1)(2
           u_1 u_2-u_1-u_2)}{2(u_1-1)(u_2-1)(v_1-1)(v_2-1)},\\
  \gamma &=\frac{(v_1-1)(v_2-1)}{(u_1-1)(u_2-1)},
\endaligned\right.$$
with
$$u_1=\kappa_1^3, \quad u_2=\kappa_2^3, \quad
  v_1=\lambda_1^3, \quad v_2=\lambda_2^3.$$
Moreover, the moduli $(\kappa_1, \kappa_2)$ and $(\lambda_1,
\lambda_2)$ satisfy the modular equation:
$$(\kappa_1^3-1)(\kappa_2^3-1)(\kappa_1^3-\kappa_2^3)
 =(\lambda_1^3-1)(\lambda_2^3-1)(\lambda_1^3-\lambda_2^3).$$}

{\it Proof}. Put $w_1=y_1^3$, $w_2=y_2^3$, $t_1=x_1^3$ and
$t_2=x_2^3$ in the Main Theorem 3.
\flushpar
$\qquad \qquad \qquad \qquad \qquad \qquad \qquad \qquad \qquad
 \qquad \qquad \qquad \qquad \qquad \qquad \qquad \qquad \qquad
 \quad \boxed{}$

\vskip 0.5 cm
\centerline{\bf 6. Triangular $s$-functions associated to
                   $U(2, 1)$ and Picard curves}
\vskip 0.5 cm

  We will study the system of partial differential equations:
$$\left\{\aligned
 \{w_1, w_2; v_1, v_2\}_{v_1} &=F_1(\gamma; v_1, v_2),\\
 \{w_1, w_2; v_1, v_2\}_{v_2} &=F_2(\gamma; v_1, v_2),\\
 [w_1, w_2; v_1, v_2]_{v_1} &=P_1(\alpha, \beta, \gamma; v_1, v_2),\\
 [w_1, w_2; v_1, v_2]_{v_2} &=P_2(\alpha, \beta, \gamma; v_1, v_2),
\endaligned\right.\tag 6.1$$
where
$$\left\{\aligned
 F_1(\gamma; v_1, v_2) &=\frac{-\gamma v_2(v_2-1)}{v_1(v_1-1)(v_1-v_2)},\\
 F_2(\gamma; v_1, v_2) &=\frac{-\gamma v_1(v_1-1)}{v_2(v_2-1)(v_2-v_1)},\\
 P_1(\alpha, \beta, \gamma; v_1, v_2) &=\frac{\alpha}{v_1}+
 \frac{\beta}{v_1-1}+\frac{\gamma}{v_1-v_2},\\
 P_2(\alpha, \beta, \gamma; v_1, v_2) &=\frac{\alpha}{v_2}+
 \frac{\beta}{v_2-1}+\frac{\gamma}{v_2-v_1}.
\endaligned\right.\tag 6.2$$
The solutions are denoted by
$$w_1=s_1(\alpha, \beta, \gamma; v_1, v_2), \quad
  w_2=s_2(\alpha, \beta, \gamma; v_1, v_2).\tag 6.3$$
Consequently,
$$\left\{\aligned
 \{s_1(\alpha, \beta, \gamma; v_1, v_2), s_2(\alpha, \beta, \gamma;
 v_1, v_2); v_1, v_2\}_{v_1} &=F_1(\gamma; v_1, v_2),\\
 \{s_1(\alpha, \beta, \gamma; v_1, v_2), s_2(\alpha, \beta, \gamma;
 v_1, v_2); v_1, v_2\}_{v_2} &=F_2(\gamma; v_1, v_2),\\
 [s_1(\alpha, \beta, \gamma; v_1, v_2), s_2(\alpha, \beta, \gamma;
 v_1, v_2); v_1, v_2]_{v_1} &=P_1(\alpha, \beta, \gamma; v_1, v_2),\\
 [s_1(\alpha, \beta, \gamma; v_1, v_2), s_2(\alpha, \beta, \gamma;
 v_1, v_2); v_1, v_2]_{v_2} &=P_2(\alpha, \beta, \gamma; v_1, v_2),
\endaligned\right.\tag 6.4$$

  By the above Main Proposition 1, for the next two systems of equations,
$$\left\{\aligned
 \{u_1, u_2; w_1, w_2\}_{w_1} &=F_1(\gamma; w_1, w_2),\\
 \{u_1, u_2; w_1, w_2\}_{w_2} &=F_2(\gamma; w_1, w_2),\\
 [u_1, u_2; w_1, w_2]_{w_1} &=P_1(\alpha, \beta, \gamma; w_1, w_2),\\
 [u_1, u_2; w_1, w_2]_{w_2} &=P_2(\alpha, \beta, \gamma; w_1, w_2),
\endaligned\right.$$
and
$$\left\{\aligned
 \{u_1, u_2; x, y\}_{x} &=F_1(\gamma^{\prime}; x, y),\\
 \{u_1, u_2; x, y\}_{y} &=F_2(\gamma^{\prime}; x, y),\\
 [u_1, u_2; x, y]_{x} &=P_1(\alpha^{\prime}, \beta^{\prime}, \gamma^{\prime}; x, y),\\
 [u_1, u_2; x, y]_{y} &=P_2(\alpha^{\prime}, \beta^{\prime}, \gamma^{\prime}; x, y),
\endaligned\right.$$
we have
$$\left(\matrix
  F_1(\gamma^{\prime}; x, y)\\
  F_2(\gamma^{\prime}; x, y)\\
  P_1(\alpha^{\prime}, \beta^{\prime}, \gamma^{\prime}; x, y)\\
  P_2(\alpha^{\prime}, \beta^{\prime}, \gamma^{\prime}; x, y)
  \endmatrix\right)
 =\frac{M(w_1, w_2; x, y)}{\frac{\partial(w_1, w_2)}{\partial(x, y)}}
  \left(\matrix
  F_1(\gamma; w_1, w_2)\\
  F_2(\gamma; w_1, w_2)\\
  P_1(\alpha, \beta, \gamma; w_1, w_2)\\
  P_2(\alpha, \beta, \gamma; w_1, w_2)
  \endmatrix\right)+
  \left(\matrix
  \{w_1, w_2; x, y\}_{x}\\
  \{w_1, w_2; x, y\}_{y}\\
    [w_1, w_2; x, y]_{x}\\
    [w_1, w_2; x, y]_{y}
  \endmatrix\right).$$

  For the integral
$$K_3=\int_{0}^{1} \frac{dx}{\root 3 \of{x^2(1-x)(1-k_1 x)(1-k_2 x)}}
     =\Gamma\left(\frac{1}{3}\right) \Gamma\left(\frac{2}{3}\right)
      F_1\left(\frac{1}{3}; \frac{1}{3}, \frac{1}{3}; 1; k_1, k_2\right),$$
corresponding to the Picard curve, we have
$$\alpha=c-b=\frac{2}{3}, \quad \beta=a+b-c+1=\frac{2}{3}, \quad
  \gamma=-b=-\frac{1}{3}.$$
Thus,
$$\left\{\aligned
 F_1 &=\frac{\frac{1}{3} v_2(v_2-1)}{v_1(v_1-1)(v_1-v_2)},\\
 F_2 &=\frac{\frac{1}{3} v_1(v_1-1)}{v_2(v_2-1)(v_2-v_1)},\\
 P_1 &=\frac{\frac{2}{3}}{v_1}+\frac{\frac{2}{3}}{v_1-1}+
       \frac{-\frac{1}{3}}{v_1-v_2},\\
 P_2 &=\frac{\frac{2}{3}}{v_2}+\frac{\frac{2}{3}}{v_2-1}+
       \frac{-\frac{1}{3}}{v_2-v_1}.
\endaligned\right.$$

  For the following two equations:
$$\left\{\aligned
 \{s_1, s_2; x, y\}_{x} &=F_1(\gamma; x, y),\\
 \{s_1, s_2; x, y\}_{y} &=F_2(\gamma; x, y),\\
 [s_1, s_2; x, y]_{x} &=P_1(\alpha, \beta, \gamma; x, y),\\
 [s_1, s_2; x, y]_{y} &=P_2(\alpha, \beta, \gamma; x, y),
\endaligned\right.$$
and
$$\left\{\aligned
 \{w_1, w_2; x, y\}_{x} &=F_1\left(-\frac{1}{3}; x, y\right),\\
 \{w_1, w_2; x, y\}_{y} &=F_2\left(-\frac{1}{3}; x, y\right),\\
 [w_1, w_2; x, y]_{x} &=P_1\left(\frac{2}{3}, \frac{2}{3}, -\frac{1}{3}; x, y\right),\\
 [w_1, w_2; x, y]_{y} &=P_2\left(\frac{2}{3}, \frac{2}{3}, -\frac{1}{3}; x, y\right),
\endaligned\right.$$
where $s_1=s_1(\alpha, \beta, \gamma; x, y)$ and $s_2=s_2(\alpha,
\beta, \gamma; x, y)$, we have
$$\frac{M(w_1, w_2; x, y)}{\frac{\partial(w_1, w_2)}{\partial(x, y)}}
  \left(\matrix \{s_1, s_2; w_1, w_2\}_{w_1}\\
                \{s_1, s_2; w_1, w_2\}_{w_2}\\
                  [s_1, s_2; w_1, w_2]_{w_1}\\
                  [s_1, s_2; w_1, w_2]_{w_2}
  \endmatrix\right)
 =\left(\matrix
  F_1(\gamma; x, y)-F_1\left(-\frac{1}{3}; x, y\right)\\
  F_2(\gamma; x, y)-F_2\left(-\frac{1}{3}; x, y\right)\\
  P_1(\alpha, \beta, \gamma; x, y)-P_1\left(\frac{2}{3},
  \frac{2}{3}, -\frac{1}{3}; x, y\right)\\
  P_2(\alpha, \beta, \gamma; x, y)-P_2\left(\frac{2}{3},
  \frac{2}{3}, -\frac{1}{3}; x, y\right)
  \endmatrix\right).$$
Thus,
$$\left(\matrix \{s_1, s_2; w_1, w_2\}_{w_1}\\
                \{s_1, s_2; w_1, w_2\}_{w_2}\\
                  [s_1, s_2; w_1, w_2]_{w_1}\\
                  [s_1, s_2; w_1, w_2]_{w_2}
  \endmatrix\right)
 =\frac{M(x, y; w_1, w_2)}{\frac{\partial(x, y)}{\partial(w_1, w_2)}}
  \left(\matrix
  F_1(\gamma; x, y)-F_1\left(-\frac{1}{3}; x, y\right)\\
  F_2(\gamma; x, y)-F_2\left(-\frac{1}{3}; x, y\right)\\
  P_1(\alpha, \beta, \gamma; x, y)-P_1\left(\frac{2}{3},
  \frac{2}{3}, -\frac{1}{3}; x, y\right)\\
  P_2(\alpha, \beta, \gamma; x, y)-P_2\left(\frac{2}{3},
  \frac{2}{3}, -\frac{1}{3}; x, y\right)
  \endmatrix\right).$$

   Let us consider the integral
$$z(x, y)=\int_{0}^{1} \frac{dt}{\root 3 \of{t(t-1)(t-x)(t-y)}},$$
which satisfies the following partial differential equations:
$$\left\{\aligned
 x(1-x)r+y(1-x)s+\left(1-\frac{5}{3}x\right)p-\frac{1}{3}yq-\frac{1}{9}z=0,\\
 y(1-y)t+x(1-y)s+\left(1-\frac{5}{3}y\right)q-\frac{1}{3}xp-\frac{1}{9}z=0,
\endaligned\right.$$
where
$$r=\frac{\partial^2 z}{\partial x^2}, \quad
  s=\frac{\partial^2 z}{\partial x \partial y}, \quad
  t=\frac{\partial^2 z}{\partial y^2}, \quad
  p=\frac{\partial z}{\partial x}, \quad
  q=\frac{\partial z}{\partial y}.$$
This gives that
$$\left\{\aligned
 &r+\left(\frac{\frac{2}{3}}{x}+\frac{\frac{2}{3}}{x-1}+
  \frac{\frac{1}{3}}{x-y}\right)p+\frac{-\frac{1}{3}y(y-1)}
  {x(x-1)(x-y)}q+\frac{\frac{1}{9}}{x(x-1)}z=0,\\
 &t+\left(\frac{\frac{2}{3}}{y}+\frac{\frac{2}{3}}{y-1}+
  \frac{\frac{1}{3}}{y-x}\right)q+\frac{-\frac{1}{3}x(x-1)}
  {y(y-1)(y-x)}p+\frac{\frac{1}{9}}{y(y-1)}z=0,\\
 &s=\frac{\frac{1}{3}}{x-y}p+\frac{-\frac{1}{3}}{x-y}q.
\endaligned\right.$$
Here, $\alpha=\frac{2}{3}$, $\beta=\frac{2}{3}$ and
$\gamma=-\frac{1}{3}$. Consequently,
$$P_1(x, y)=\frac{\frac{2}{3}}{x}+\frac{\frac{2}{3}}{x-1}+
            \frac{-\frac{1}{3}}{x-y}, \quad
  P_2(x, y)=\frac{\frac{2}{3}}{y}+\frac{\frac{2}{3}}{y-1}+
            \frac{-\frac{1}{3}}{y-x}.$$
$$F_1(x, y)=\frac{\frac{1}{3}y(y-1)}{x(x-1)(x-y)}, \quad
  F_2(x, y)=\frac{\frac{1}{3}x(x-1)}{y(y-1)(y-x)}.$$

  Let us consider the following partial differential equations:
$$\left\{\aligned
 \{ w_1, w_2; x, y \}_{x} &=\frac{\frac{1}{3} y(y-1)}{x(x-1)(x-y)},\\
 \{ w_1, w_2; x, y \}_{y} &=\frac{\frac{1}{3} x(x-1)}{y(y-1)(y-x)},\\
 [ w_1, w_2; x, y ]_{x} &=\frac{\frac{2}{3}}{x}+\frac{\frac{2}{3}}{x-1}
                         +\frac{-\frac{1}{3}}{x-y},\\
 [ w_1, w_2; x, y ]_{y} &=\frac{\frac{2}{3}}{y}+\frac{\frac{2}{3}}{y-1}
                         +\frac{-\frac{1}{3}}{y-x}.
\endaligned\right.$$
According to Main Theorem 1 and Main Theorem 2, the transformation
is given by
$$w=x^{-\frac{2}{9}} y^{-\frac{2}{9}} (x-1)^{-\frac{2}{9}} (y-1)^{-\frac{2}{9}}
    (x-y)^{-\frac{2}{9}} \root 3 \of{\Delta},\tag 6.5$$
where $\Delta=\frac{\partial(x, y)}{\partial(w_1, w_2)}$ and the
function $w$ satisfies the following partial differential
equations:
$$\left\{\aligned
 &\frac{\partial^2 w}{\partial x^2}+\left(\frac{\frac{2}{3}}{x}+
  \frac{\frac{2}{3}}{x-1}+\frac{\frac{1}{3}}{x-y}\right)
  \frac{\partial w}{\partial x}+\frac{-\frac{1}{3} y(y-1)}{x(x-1)(x-y)}
  \frac{\partial w}{\partial y}+\frac{\frac{1}{9}}{x(x-1)} w=0,\\
 &\frac{\partial^2 w}{\partial x \partial y}=\frac{\frac{1}{3}}{x-y}
  \frac{\partial w}{\partial x}+\frac{-\frac{1}{3}}{x-y}
  \frac{\partial w}{\partial y},\\
 &\frac{\partial^2 w}{\partial y^2}+\left(\frac{\frac{2}{3}}{y}+
  \frac{\frac{2}{3}}{y-1}+\frac{\frac{1}{3}}{y-x}\right)
  \frac{\partial w}{\partial y}+\frac{-\frac{1}{3} x(x-1)}{y(y-1)(y-x)}
  \frac{\partial w}{\partial x}+\frac{\frac{1}{9}}{y(y-1)} w=0.
\endaligned\right.$$
The above argument shows that a solution is given by
$$w=z(x, y)=\int_{0}^{1} \frac{dt}{\root 3 \of{t(t-1)(t-x)(t-y)}}.\tag 6.6$$

  Set
$$f_1(x, y)=\frac{y(y-1)}{x(x-1)(x-y)}, \quad
  f_2(x, y)=\frac{x(x-1)}{y(y-1)(y-x)}.\tag 6.7$$
$$P_1(\alpha, \beta, \gamma; x, y)=
  \frac{\alpha}{x}+\frac{\beta}{x-1}+\frac{\gamma}{x-y}, \quad
  P_2(\alpha, \beta, \gamma; x, y)=
  \frac{\alpha}{y}+\frac{\beta}{y-1}+\frac{\gamma}{y-x}.\tag 6.8$$
Let
$$T=\left(\matrix
    -1 &  0 & 1\\
       & -1 & 1\\
       &    & 1
    \endmatrix\right), \quad
  S_1=\left(\matrix
      1 & 0 & 0\\
      0 & 0 & 1\\
      0 & 1 & 0
      \endmatrix\right), \quad
  S_2=\left(\matrix
        &   & 1\\
        & 1 &  \\
      1 &   &
      \endmatrix\right).\tag 6.9$$
Then
$$T(x, y)=(1-x, 1-y).\tag 6.10$$
$$\left\{\aligned
  S_1(x, y) &=\left(\frac{x}{y}, \frac{1}{y}\right),\\
  S_1 T(x, y) &=\left(\frac{1-x}{1-y}, \frac{1}{1-y}\right),\\
  T S_1(x, y) &=\left(\frac{y-x}{y}, \frac{y-1}{y}\right),\\
  S_1 T S_1(x, y) &=\left(\frac{y-x}{y-1}, \frac{y}{y-1}\right).
\endaligned\right.\tag 6.11$$
$$\left\{\aligned
  S_2(x, y) &=\left(\frac{1}{x}, \frac{y}{x}\right),\\
  S_2 T(x, y) &=\left(\frac{1}{1-x}, \frac{1-y}{1-x}\right),\\
  T S_2(x, y) &=\left(\frac{x-1}{x}, \frac{x-y}{x}\right),\\
  S_2 T S_2(x, y) &=\left(\frac{x}{x-1}, \frac{x-y}{x-1}\right).
\endaligned\right.\tag 6.12$$
We have
$$\left\{\aligned
  f_1(T(x, y)) &=-f_1(x, y),\\
  f_1(S_1(x, y)) &=-f_1(x, y),\\
  f_1(S_1 T(x, y)) &=f_1(x, y),\\
  f_1(T S_1(x, y)) &=f_1(x, y),\\
  f_1(S_1 T S_1(x, y)) &=-f_1(x, y),
\endaligned\right.\tag 6.13$$
$$\left\{\aligned
  f_2(T(x, y)) &=-f_2(x, y),\\
  f_2(S_2(x, y)) &=-f_2(x, y),\\
  f_2(S_2 T(x, y)) &=f_2(x, y),\\
  f_2(T S_2(x, y)) &=f_2(x, y),\\
  f_2(S_2 T S_2(x, y)) &=-f_2(x, y),
\endaligned\right.\tag 6.14$$
and
$$\left\{\aligned
 P_1(\alpha, \beta, \gamma; T(x, y)) &=-P_1(\beta, \alpha, \gamma; x, y),\\
 P_1(\alpha, \beta, \gamma; S_1(x, y)) &=y P_1(\alpha, \gamma, \beta; x, y),\\
 P_1(\alpha, \beta, \gamma; S_1 T(x, y)) &=(y-1) P_1(\gamma, \alpha, \beta; x, y),\\
 P_1(\alpha, \beta, \gamma; T S_1(x, y)) &=-y P_1(\beta, \gamma, \alpha; x, y),\\
 P_1(\alpha, \beta, \gamma; S_1 T S_1(x, y)) &=(1-y) P_1(\gamma, \beta, \alpha; x, y),
\endaligned\right.\tag 6.15$$
$$\left\{\aligned
 P_2(\alpha, \beta, \gamma; T(x, y)) &=-P_2(\beta, \alpha, \gamma; x, y),\\
 P_2(\alpha, \beta, \gamma; S_2(x, y)) &=x P_2(\alpha, \gamma, \beta; x, y),\\
 P_2(\alpha, \beta, \gamma; S_2 T(x, y)) &=(x-1) P_2(\gamma, \alpha, \beta; x, y),\\
 P_2(\alpha, \beta, \gamma; T S_2(x, y)) &=-x P_2(\beta, \gamma, \alpha; x, y),\\
 P_2(\alpha, \beta, \gamma; S_2 T S_2(x, y)) &=(1-x) P_2(\gamma, \beta, \alpha; x, y).
\endaligned\right.\tag 6.16$$

  Hence, we find that
$$\left\{\aligned
 \{s_1(\beta, \alpha, \gamma; v_1, v_2), s_2(\beta, \alpha,
 \gamma; v_1, v_2); v_1, v_2 \}_{v_1} &=F_1(-\gamma; T(v_1, v_2)),\\
 \{s_1(\beta, \alpha, \gamma; v_1, v_2), s_2(\beta, \alpha,
 \gamma; v_1, v_2); v_1, v_2 \}_{v_2} &=F_2(-\gamma; T(v_1, v_2)),\\
 [s_1(\beta, \alpha, \gamma; v_1, v_2), s_2(\beta, \alpha,
 \gamma; v_1, v_2); v_1, v_2 ]_{v_1} &=-P_1(\alpha, \beta, \gamma; T(v_1, v_2)),\\
 [s_1(\beta, \alpha, \gamma; v_1, v_2), s_2(\beta, \alpha,
 \gamma; v_1, v_2); v_1, v_2 ]_{v_2} &=-P_2(\alpha, \beta, \gamma; T(v_1, v_2)).
\endaligned\right.\tag 6.17$$
$$\left\{\aligned
 \{s_1(\alpha, \gamma, \beta; v_1, v_2), s_2(\alpha, \gamma,
 \beta; v_1, v_2); v_1, v_2 \}_{v_1} &=F_1(-\beta; S_1(v_1, v_2)),\\
 \{s_1(\alpha, \gamma, \beta; v_1, v_2), s_2(\alpha, \gamma,
 \beta; v_1, v_2); v_1, v_2 \}_{v_2} &=F_2(-\beta; S_2(v_1, v_2)),\\
 [s_1(\alpha, \gamma, \beta; v_1, v_2), s_2(\alpha, \gamma,
 \beta; v_1, v_2); v_1, v_2 ]_{v_1} &=\frac{1}{v_2} P_1(\alpha, \beta,
 \gamma; S_1(v_1, v_2)),\\
 [s_1(\alpha, \gamma, \beta; v_1, v_2), s_2(\alpha, \gamma,
 \beta; v_1, v_2); v_1, v_2 ]_{v_2} &=\frac{1}{v_1} P_2(\alpha, \beta,
 \gamma; S_2(v_1, v_2)).
\endaligned\right.\tag 6.18$$
$$\left\{\aligned
 \{s_1(\gamma, \alpha, \beta; v_1, v_2), s_2(\gamma, \alpha,
 \beta; v_1, v_2); v_1, v_2 \}_{v_1} &=F_1(\beta; S_1 T(v_1, v_2)),\\
 \{s_1(\gamma, \alpha, \beta; v_1, v_2), s_2(\gamma, \alpha,
 \beta; v_1, v_2); v_1, v_2 \}_{v_2} &=F_2(\beta; S_2 T(v_1, v_2)),\\
 [s_1(\gamma, \alpha, \beta; v_1, v_2), s_2(\gamma, \alpha,
 \beta; v_1, v_2); v_1, v_2 ]_{v_1} &=\frac{1}{v_2-1} P_1(\alpha, \beta,
 \gamma; S_1 T(v_1, v_2)),\\
 [s_1(\gamma, \alpha, \beta; v_1, v_2), s_2(\gamma, \alpha,
 \beta; v_1, v_2); v_1, v_2 ]_{v_2} &=\frac{1}{v_1-1} P_2(\alpha, \beta,
 \gamma; S_2 T(v_1, v_2)).
\endaligned\right.\tag 6.19$$
$$\left\{\aligned
 \{s_1(\beta, \gamma, \alpha; v_1, v_2), s_2(\beta, \gamma,
 \alpha; v_1, v_2); v_1, v_2 \}_{v_1} &=F_1(\alpha; T S_1(v_1, v_2)),\\
 \{s_1(\beta, \gamma, \alpha; v_1, v_2), s_2(\beta, \gamma,
 \alpha; v_1, v_2); v_1, v_2 \}_{v_2} &=F_2(\alpha; T S_2(v_1, v_2)),\\
 [s_1(\beta, \gamma, \alpha; v_1, v_2), s_2(\beta, \gamma,
 \alpha; v_1, v_2); v_1, v_2 ]_{v_1} &=-\frac{1}{v_2} P_1(\alpha, \beta,
 \gamma; T S_1(v_1, v_2)),\\
 [s_1(\beta, \gamma, \alpha; v_1, v_2), s_2(\beta, \gamma,
 \alpha; v_1, v_2); v_1, v_2 ]_{v_2} &=-\frac{1}{v_1} P_2(\alpha, \beta,
 \gamma; T S_2(v_1, v_2)).
\endaligned\right.\tag 6.20$$
$$\left\{\aligned
 \{s_1(\gamma, \beta, \alpha; v_1, v_2), s_2(\gamma, \beta,
 \alpha; v_1, v_2); v_1, v_2 \}_{v_1} &=F_1(-\alpha; S_1 T S_1(v_1, v_2)),\\
 \{s_1(\gamma, \beta, \alpha; v_1, v_2), s_2(\gamma, \beta,
 \alpha; v_1, v_2); v_1, v_2 \}_{v_2} &=F_2(-\alpha; S_2 T S_2(v_1, v_2)),\\
 [s_1(\gamma, \beta, \alpha; v_1, v_2), s_2(\gamma, \beta,
 \alpha; v_1, v_2); v_1, v_2 ]_{v_1} &=\frac{1}{1-v_2} P_1(\alpha, \beta,
 \gamma; S_1 T S_1(v_1, v_2)),\\
 [s_1(\gamma, \beta, \alpha; v_1, v_2), s_2(\gamma, \beta,
 \alpha; v_1, v_2); v_1, v_2 ]_{v_2} &=\frac{1}{1-v_1} P_2(\alpha, \beta,
 \gamma; S_2 T S_2(v_1, v_2)).
\endaligned\right.\tag 6.21$$

{\smc Proposition 6.1 (Main Proposition 3)}. {\it The triangular
$s$-functions have the following properties of transformations:
\roster
\item
$$\left\{\aligned
  \{s_1(\alpha, \beta, \gamma; T(v_1, v_2)), s_2(\alpha, \beta, \gamma;
  T(v_1, v_2)); v_1, v_2\}_{v_1} &=F_1(\gamma; v_1, v_2),\\
  \{s_1(\alpha, \beta, \gamma; T(v_1, v_2)), s_2(\alpha, \beta, \gamma;
  T(v_1, v_2)); v_1, v_2\}_{v_2} &=F_2(\gamma; v_1, v_2),\\
  [s_1(\alpha, \beta, \gamma; T(v_1, v_2)), s_2(\alpha, \beta, \gamma;
  T(v_1, v_2)); v_1, v_2]_{v_1} &=P_1(\beta, \alpha, \gamma; v_1, v_2),\\
  [s_1(\alpha, \beta, \gamma; T(v_1, v_2)), s_2(\alpha, \beta, \gamma;
  T(v_1, v_2)); v_1, v_2]_{v_2} &=P_2(\beta, \alpha, \gamma; v_1, v_2).
\endaligned\right.\tag 6.22$$
\item
$$\left\{\aligned
 \{s_1(\alpha, \beta, \gamma; S_1(v_1, v_2)), s_2(\alpha, \beta, \gamma;
  S_1(v_1, v_2)); v_1, v_2\}_{v_1} =&F_1(\gamma; v_1, v_2),\\
 \{s_1(\alpha, \beta, \gamma; S_2(v_1, v_2)), s_2(\alpha, \beta, \gamma;
  S_2(v_1, v_2)); v_1, v_2\}_{v_2} =&F_2(\gamma; v_1, v_2),\\
  [s_1(\alpha, \beta, \gamma; S_1(v_1, v_2)), s_2(\alpha, \beta, \gamma;
  S_1(v_1, v_2)); v_1, v_2]_{v_1} =&P_1(\alpha, -2 \gamma, \beta+3 \gamma; v_1, v_2),\\
  [s_1(\alpha, \beta, \gamma; S_2(v_1, v_2)), s_2(\alpha, \beta, \gamma;
  S_2(v_1, v_2)); v_1, v_2]_{v_2} =&P_2(\alpha, -2 \gamma, \beta+3 \gamma; v_1, v_2).
\endaligned\right.\tag 6.23$$
\item
$$\left\{\aligned
 \{s_1(\alpha, \beta, \gamma; S_1 T(v_1, v_2)), s_2(\alpha, \beta, \gamma;
  S_1 T(v_1, v_2)); v_1, v_2\}_{v_1} =&F_1(\gamma; v_1, v_2),\\
 \{s_1(\alpha, \beta, \gamma; S_2 T(v_1, v_2)), s_2(\alpha, \beta, \gamma;
  S_2 T(v_1, v_2)); v_1, v_2\}_{v_2} =&F_2(\gamma; v_1, v_2),\\
  [s_1(\alpha, \beta, \gamma; S_1 T(v_1, v_2)), s_2(\alpha, \beta, \gamma;
  S_1 T(v_1, v_2)); v_1, v_2]_{v_1} =&P_1(-2 \gamma, \alpha, \beta+3 \gamma; v_1, v_2),\\
  [s_1(\alpha, \beta, \gamma; S_2 T(v_1, v_2)), s_2(\alpha, \beta, \gamma;
  S_2 T(v_1, v_2)); v_1, v_2]_{v_2} =&P_2(-2 \gamma, \alpha, \beta+3 \gamma; v_1, v_2).
\endaligned\right.\tag 6.24$$
\item
$$\left\{\aligned
 \{s_1(\alpha, \beta, \gamma; T S_1(v_1, v_2)), s_2(\alpha, \beta, \gamma;
  T S_1(v_1, v_2)); v_1, v_2\}_{v_1} =&F_1(\gamma; v_1, v_2),\\
 \{s_1(\alpha, \beta, \gamma; T S_2(v_1, v_2)), s_2(\alpha, \beta, \gamma;
  T S_2(v_1, v_2)); v_1, v_2\}_{v_2} =&F_2(\gamma; v_1, v_2),\\
  [s_1(\alpha, \beta, \gamma; T S_1(v_1, v_2)), s_2(\alpha, \beta, \gamma;
  T S_1(v_1, v_2)); v_1, v_2]_{v_1} =&P_1(\beta, -2 \gamma, \alpha+3 \gamma; v_1, v_2),\\
  [s_1(\alpha, \beta, \gamma; T S_2(v_1, v_2)), s_2(\alpha, \beta, \gamma;
  T S_2(v_1, v_2)); v_1, v_2]_{v_2} =&P_2(\beta, -2 \gamma, \alpha+3 \gamma; v_1, v_2).
\endaligned\right.\tag 6.25$$
\item
$$\left\{\aligned
 &\{s_1(\alpha, \beta, \gamma; S_1 T S_1(v_1, v_2)), s_2(\alpha, \beta, \gamma;
  S_1 T S_1(v_1, v_2)); v_1, v_2\}_{v_1}=F_1(\gamma; v_1, v_2),\\
 &\{s_1(\alpha, \beta, \gamma; S_2 T S_2(v_1, v_2)), s_2(\alpha, \beta, \gamma;
  S_2 T S_2(v_1, v_2)); v_1, v_2\}_{v_2}=F_2(\gamma; v_1, v_2),\\
 &[s_1(\alpha, \beta, \gamma; S_1 T S_1(v_1, v_2)), s_2(\alpha, \beta, \gamma;
  S_1 T S_1(v_1, v_2)); v_1, v_2]_{v_1}\\
=&P_1(-2 \gamma, \beta, \alpha+3 \gamma; v_1, v_2),\\
 &[s_1(\alpha, \beta, \gamma; S_2 T S_2(v_1, v_2)), s_2(\alpha, \beta, \gamma;
  S_2 T S_2(v_1, v_2)); v_1, v_2]_{v_2}\\
=&P_2(-2 \gamma, \beta, \alpha+3 \gamma; v_1, v_2).
\endaligned\right.\tag 6.26$$
\endroster}

{\it Proof}.

(1) Let $(u_1, u_2)=T(v_1, v_2)=(1-v_1, 1-v_2)$.

  By Main Proposition 1, we have
$$\aligned
 &\{s_1(\alpha, \beta, \gamma; T(v_1, v_2)), s_2(\alpha, \beta, \gamma;
  T(v_1, v_2)); v_1, v_2\}_{v_1}\\
=&-\{s_1(\alpha, \beta, \gamma; u_1, u_2), s_2(\alpha, \beta, \gamma;
  u_1, u_2); u_1, u_2\}_{u_1}\\
=&-F_1(\gamma; u_1, u_2)\\
=&-F_1(\gamma; T(v_1, v_2))\\
=&F_1(\gamma; v_1, v_2).
\endaligned$$
$$\aligned
 &[s_1(\alpha, \beta, \gamma; T(v_1, v_2)), s_2(\alpha, \beta, \gamma;
  T(v_1, v_2)); v_1, v_2]_{v_1}\\
=&-[s_1(\alpha, \beta, \gamma; u_1, u_2), s_2(\alpha, \beta,
  \gamma; u_1, u_2); u_1, u_2]_{u_1}\\
=&-P_1(\alpha, \beta, \gamma; u_1, u_2)\\
=&-P_1(\alpha, \beta, \gamma; T(v_1, v_2))\\
=&P_1(\beta, \alpha, \gamma; v_1, v_2).
\endaligned$$

(2) Let $(u_1, u_2)=S_1(v_1, v_2)=\left(\frac{v_1}{v_2},
\frac{1}{v_2}\right)$.

  By Main Proposition 1, we have
$$\aligned
 &\{s_1(\alpha, \beta, \gamma; S_1(v_1, v_2)), s_2(\alpha, \beta, \gamma;
  S_1(v_1, v_2)); v_1, v_2\}_{v_1}\\
=&-\{s_1(\alpha, \beta, \gamma; u_1, u_2), s_2(\alpha, \beta,
  \gamma; u_1, u_2); u_1, u_2\}_{u_1}\\
=&-F_1(\gamma; u_1, u_2)\\
=&-F_1(\gamma; S_1(v_1, v_2))\\
=&F_1(\gamma; v_1, v_2).
\endaligned$$
$$\aligned
 &[s_1(\alpha, \beta, \gamma; S_1(v_1, v_2)), s_2(\alpha,
  \beta, \gamma; S_1(v_1, v_2)); v_1, v_2]_{v_1}\\
=&3 \frac{v_1}{v_2} \{s_1(\alpha, \beta, \gamma; u_1, u_2),
  s_2(\alpha, \beta, \gamma; u_1, u_2); u_1, u_2\}_{u_1}+\\
 &+\frac{1}{v_2} [s_1(\alpha, \beta, \gamma; u_1, u_2),
  s_2(\alpha, \beta, \gamma; u_1, u_2); u_1, u_2]_{u_1}\\
=&3 \frac{v_1}{v_2} F_1(\gamma; u_1, u_2)+\frac{1}{v_2}
  P_1(\alpha, \beta, \gamma; u_1, u_2)\\
=&3 \frac{v_1}{v_2} F_1(\gamma; S_1(v_1, v_2))+\frac{1}{v_2}
  P_1(\alpha, \beta, \gamma; S_1(v_1, v_2))\\
=&-3 \frac{v_1}{v_2} F_1(\gamma; v_1, v_2)+P_1(\alpha, \gamma,
  \beta; v_1, v_2)\\
=&\frac{3 \gamma (v_2-1)}{(v_1-1)(v_1-v_2)}+\frac{\alpha}{v_1}+
  \frac{\gamma}{v_1-1}+\frac{\beta}{v_1-v_2}\\
=&\frac{\alpha}{v_1}+\frac{-2 \gamma}{v_1-1}+\frac{\beta+3 \gamma}{v_1-v_2}\\
=&P_1(\alpha, -2 \gamma, \beta+3 \gamma; v_1, v_2).
\endaligned$$

(3) Let $(u_1, u_2)=S_1 T(v_1, v_2)=\left(\frac{1-v_1}{1-v_2},
\frac{1}{1-v_2}\right)$.

  By Main Proposition 1, we have
$$\aligned
 &\{s_1(\alpha, \beta, \gamma; S_1 T(v_1, v_2)), s_2(\alpha, \beta, \gamma;
  S_1 T(v_1, v_2)); v_1, v_2\}_{v_1}\\
=&\{s_1(\alpha, \beta, \gamma; u_1, u_2), s_2(\alpha, \beta,
  \gamma; u_1, u_2); u_1, u_2\}_{u_1}\\
=&F_1(\gamma; u_1, u_2)\\
=&F_1(\gamma; S_1 T(v_1, v_2))\\
=&F_1(\gamma; v_1, v_2).
\endaligned$$
$$\aligned
 &[s_1(\alpha, \beta, \gamma; S_1 T(v_1, v_2)), s_2(\alpha,
  \beta, \gamma; S_1 T(v_1, v_2)); v_1, v_2]_{v_1}\\
=&-3 \frac{v_1-1}{v_2-1} \{s_1(\alpha, \beta, \gamma; u_1, u_2),
  s_2(\alpha, \beta, \gamma; u_1, u_2); u_1, u_2\}_{u_1}+\\
 &+\frac{1}{v_2-1} [s_1(\alpha, \beta, \gamma; u_1, u_2),
  s_2(\alpha, \beta, \gamma; u_1, u_2); u_1, u_2]_{u_1}\\
=&-3 \frac{v_1-1}{v_2-1} F_1(\gamma; u_1, u_2)+\frac{1}{v_2-1}
  P_1(\alpha, \beta, \gamma; u_1, u_2)\\
=&-3 \frac{v_1-1}{v_2-1} F_1(\gamma; S_1 T(v_1, v_2))+
  \frac{1}{v_2-1} P_1(\alpha, \beta, \gamma; S_1 T(v_1, v_2))\\
=&-3 \frac{v_1-1}{v_2-1} F_1(\gamma; v_1, v_2)+P_1(\gamma, \alpha,
  \beta; v_1, v_2)\\
=&\frac{3 \gamma v_2}{v_1(v_1-v_2)}+\frac{\gamma}{v_1}
  +\frac{\alpha}{v_1-1}+\frac{\beta}{v_1-v_2}\\
=&\frac{-2 \gamma}{v_1}+\frac{\alpha}{v_1-1}+\frac{\beta+3 \gamma}{v_1-v_2}\\
=&P_1(-2 \gamma, \alpha, \beta+3 \gamma; v_1, v_2).
\endaligned$$

(4) Let $(u_1, u_2)=T S_1(v_1, v_2)=\left(\frac{v_2-v_1}{v_2},
\frac{v_2-1}{v_2}\right)$.

  By Main Proposition 1, we have
$$\aligned
 &\{s_1(\alpha, \beta, \gamma; T S_1(v_1, v_2)), s_2(\alpha, \beta, \gamma;
  T S_1(v_1, v_2)); v_1, v_2\}_{v_1}\\
=&\{s_1(\alpha, \beta, \gamma; u_1, u_2), s_2(\alpha, \beta,
  \gamma; u_1, u_2); u_1, u_2\}_{u_1}\\
=&F_1(\gamma; u_1, u_2)\\
=&F_1(\gamma; T S_1(v_1, v_2))\\
=&F_1(\gamma; v_1, v_2).
\endaligned$$
$$\aligned
 &[s_1(\alpha, \beta, \gamma; T S_1(v_1, v_2)), s_2(\alpha,
  \beta, \gamma; T S_1(v_1, v_2)); v_1, v_2]_{v_1}\\
=&-3 \frac{v_1}{v_2} \{s_1(\alpha, \beta, \gamma; u_1, u_2),
  s_2(\alpha, \beta, \gamma; u_1, u_2); u_1, u_2\}_{u_1}+\\
 &-\frac{1}{v_2} [s_1(\alpha, \beta, \gamma; u_1, u_2),
  s_2(\alpha, \beta, \gamma; u_1, u_2); u_1, u_2]_{u_1}\\
=&-3 \frac{v_1}{v_2} F_1(\gamma; u_1, u_2)-\frac{1}{v_2}
  P_1(\alpha, \beta, \gamma; u_1, u_2)\\
=&-3 \frac{v_1}{v_2} F_1(\gamma; T S_1(v_1, v_2))-\frac{1}{v_2}
  P_1(\alpha, \beta, \gamma; T S_1(v_1, v_2))\\
=&-3 \frac{v_1}{v_2} F_1(\gamma; v_1, v_2)+P_1(\beta, \gamma,
  \alpha; v_1, v_2)\\
=&\frac{3 \gamma (v_2-1)}{(v_1-1)(v_1-v_2)}+\frac{\beta}{v_1}
  +\frac{\gamma}{v_1-1}+\frac{\alpha}{v_1-v_2}\\
=&\frac{\beta}{v_1}+\frac{-2 \gamma}{v_1-1}+\frac{\alpha+3 \gamma}{v_1-v_2}\\
=&P_1(\beta, -2 \gamma, \alpha+3 \gamma; v_1, v_2).
\endaligned$$

(5) Let $(u_1, u_2)=S_1 T S_1(v_1, v_2)=\left(\frac{v_2-v_1}
{v_2-1}, \frac{v_2}{v_2-1}\right)$.

  By Main Proposition 1, we have
$$\aligned
 &\{s_1(\alpha, \beta, \gamma; S_1 T S_1(v_1, v_2)), s_2(\alpha, \beta, \gamma;
  S_1 T S_1(v_1, v_2)); v_1, v_2\}_{v_1}\\
=&-\{s_1(\alpha, \beta, \gamma; u_1, u_2), s_2(\alpha, \beta,
  \gamma; u_1, u_2); u_1, u_2\}_{u_1}\\
=&-F_1(\gamma; u_1, u_2)\\
=&-F_1(\gamma; S_1 T S_1(v_1, v_2))\\
=&F_1(\gamma; v_1, v_2).
\endaligned$$
$$\aligned
 &[s_1(\alpha, \beta, \gamma; S_1 T S_1(v_1, v_2)), s_2(\alpha,
  \beta, \gamma; S_1 T S_1(v_1, v_2)); v_1, v_2]_{v_1}\\
=&3 \frac{v_1-1}{v_2-1} \{s_1(\alpha, \beta, \gamma; u_1, u_2),
  s_2(\alpha, \beta, \gamma; u_1, u_2); u_1, u_2\}_{u_1}+\\
 &-\frac{1}{v_2-1} [s_1(\alpha, \beta, \gamma; u_1, u_2),
  s_2(\alpha, \beta, \gamma; u_1, u_2); u_1, u_2]_{u_1}\\
=&3 \frac{v_1-1}{v_2-1} F_1(\gamma; u_1, u_2)-\frac{1}{v_2-1}
  P_1(\alpha, \beta, \gamma; u_1, u_2)\\
=&3 \frac{v_1-1}{v_2-1} F_1(\gamma; S_1 T S_1(v_1,
  v_2))-\frac{1}{v_2-1} P_1(\alpha, \beta, \gamma; S_1 T S_1(v_1, v_2))\\
=&-3 \frac{v_1-1}{v_2-1} F_1(\gamma; v_1, v_2)+P_1(\gamma, \beta,
  \alpha; v_1, v_2)\\
=&\frac{3 \gamma v_2}{v_1(v_1-v_2)}+\frac{\gamma}{v_1}
  +\frac{\beta}{v_1-1}+\frac{\alpha}{v_1-v_2}\\
=&\frac{-2 \gamma}{v_1}+\frac{\beta}{v_1-1}+\frac{\alpha+3 \gamma}{v_1-v_2}\\
=&P_1(-2 \gamma, \beta, \alpha+3 \gamma; v_1, v_2).
\endaligned$$

  By the symmetry property, we can prove the other cases.

\flushpar
$\qquad \qquad \qquad \qquad \qquad \qquad \qquad \qquad
 \qquad \qquad \qquad \qquad \qquad \qquad \qquad \qquad
 \qquad \qquad \quad \boxed{}$

{\smc Corollary 6.2}. {\it The triangular $s$-functions have the
following properties of transformations:
\roster
\item
$$\left\{\aligned
 &\{s_1(\alpha, \beta, \gamma; T(v_1, v_2)), s_2(\alpha, \beta, \gamma;
  T(v_1, v_2)); v_1, v_2\}_{v_1}\\
=&\{s_1(\alpha, \beta, \gamma; v_1, v_2), s_2(\alpha, \beta, \gamma;
  v_1, v_2); v_1, v_2\}_{v_1},\\
 &\{s_1(\alpha, \beta, \gamma; T(v_1, v_2)), s_2(\alpha, \beta, \gamma;
  T(v_1, v_2)); v_1, v_2\}_{v_2}\\
=&\{s_1(\alpha, \beta, \gamma; v_1, v_2), s_2(\alpha, \beta, \gamma;
  v_1, v_2); v_1, v_2\}_{v_2},\\
 &[s_1(\alpha, \beta, \gamma; T(v_1, v_2)), s_2(\alpha, \beta, \gamma;
  T(v_1, v_2)); v_1, v_2]_{v_1}\\
=&[s_1(\beta, \alpha, \gamma; v_1, v_2), s_2(\beta, \alpha, \gamma;
  v_1, v_2); v_1, v_2]_{v_1},\\
 &[s_1(\alpha, \beta, \gamma; T(v_1, v_2)), s_2(\alpha, \beta, \gamma;
  T(v_1, v_2)); v_1, v_2]_{v_2}\\
=&[s_1(\beta, \alpha, \gamma; v_1, v_2), s_2(\beta, \alpha,
  \gamma; v_1, v_2); v_1, v_2]_{v_2}.
\endaligned\right.\tag 6.27$$
\item
$$\left\{\aligned
 &\{s_1(\alpha, \beta, \gamma; S_1(v_1, v_2)), s_2(\alpha, \beta, \gamma;
  S_1(v_1, v_2)); v_1, v_2\}_{v_1}\\
=&\{s_1(\alpha, \beta, \gamma; v_1, v_2), s_2(\alpha, \beta,
  \gamma; v_1, v_2); v_1, v_2\}_{v_1},\\
 &\{s_1(\alpha, \beta, \gamma; S_2(v_1, v_2)), s_2(\alpha, \beta, \gamma;
  S_2(v_1, v_2)); v_1, v_2\}_{v_2}\\
=&\{s_1(\alpha, \beta, \gamma; v_1, v_2), s_2(\alpha, \beta,
  \gamma; v_1, v_2); v_1, v_2\}_{v_2},\\
 &[s_1(\alpha, \beta, \gamma; S_1(v_1, v_2)), s_2(\alpha, \beta, \gamma;
  S_1(v_1, v_2)); v_1, v_2]_{v_1}\\
=&[s_1(\alpha, -2 \gamma, \beta+3 \gamma; v_1, v_2), s_2(\alpha,
  -2 \gamma, \beta+3 \gamma; v_1, v_2); v_1, v_2]_{v_1},\\
 &[s_1(\alpha, \beta, \gamma; S_2(v_1, v_2)), s_2(\alpha, \beta, \gamma;
  S_2(v_1, v_2)); v_1, v_2]_{v_2}\\
=&[s_1(\alpha, -2 \gamma, \beta+3 \gamma; v_1, v_2), s_2(\alpha,
  -2 \gamma, \beta+3 \gamma; v_1, v_2); v_1, v_2]_{v_2}.
\endaligned\right.\tag 6.28$$
\item
$$\left\{\aligned
 &\{s_1(\alpha, \beta, \gamma; S_1 T(v_1, v_2)), s_2(\alpha, \beta, \gamma;
  S_1 T(v_1, v_2)); v_1, v_2\}_{v_1}\\
=&\{s_1(\alpha, \beta, \gamma; v_1, v_2), s_2(\alpha, \beta,
  \gamma; v_1, v_2); v_1, v_2\}_{v_1},\\
 &\{s_1(\alpha, \beta, \gamma; S_2 T(v_1, v_2)), s_2(\alpha, \beta, \gamma;
  S_2 T(v_1, v_2)); v_1, v_2\}_{v_2}\\
=&\{s_1(\alpha, \beta, \gamma; v_1, v_2), s_2(\alpha, \beta,
  \gamma; v_1, v_2); v_1, v_2\}_{v_2},\\
 &[s_1(\alpha, \beta, \gamma; S_1 T(v_1, v_2)), s_2(\alpha, \beta, \gamma;
  S_1 T(v_1, v_2)); v_1, v_2]_{v_1}\\
=&[s_1(-2 \gamma, \alpha, \beta+3 \gamma; v_1, v_2), s_2(-2 \gamma,
  \alpha, \beta+3 \gamma; v_1, v_2); v_1, v_2]_{v_1},\\
 &[s_1(\alpha, \beta, \gamma; S_2 T(v_1, v_2)), s_2(\alpha, \beta, \gamma;
  S_2 T(v_1, v_2)); v_1, v_2]_{v_2}\\
=&[s_1(-2 \gamma, \alpha, \beta+3 \gamma; v_1, v_2), s_2(-2
  \gamma, \alpha, \beta+3 \gamma; v_1, v_2); v_1, v_2]_{v_2}.
\endaligned\right.\tag 6.29$$
\item
$$\left\{\aligned
 &\{s_1(\alpha, \beta, \gamma; T S_1(v_1, v_2)), s_2(\alpha, \beta, \gamma;
  T S_1(v_1, v_2)); v_1, v_2\}_{v_1}\\
=&\{s_1(\alpha, \beta, \gamma; v_1, v_2), s_2(\alpha, \beta,
  \gamma; v_1, v_2); v_1, v_2\}_{v_1},\\
 &\{s_1(\alpha, \beta, \gamma; T S_2(v_1, v_2)), s_2(\alpha, \beta, \gamma;
  T S_2(v_1, v_2)); v_1, v_2\}_{v_2}\\
=&\{s_1(\alpha, \beta, \gamma; v_1, v_2), s_2(\alpha, \beta,
  \gamma; v_1, v_2); v_1, v_2\}_{v_2},\\
 &[s_1(\alpha, \beta, \gamma; T S_1(v_1, v_2)), s_2(\alpha, \beta, \gamma;
  T S_1(v_1, v_2)); v_1, v_2]_{v_1}\\
=&[s_1(\beta, -2 \gamma, \alpha+3 \gamma; v_1, v_2), s_2(\beta, -2
  \gamma, \alpha+3 \gamma; v_1, v_2); v_1, v_2]_{v_1},\\
 &[s_1(\alpha, \beta, \gamma; T S_2(v_1, v_2)), s_2(\alpha, \beta, \gamma;
  T S_2(v_1, v_2)); v_1, v_2]_{v_2}\\
=&[s_1(\beta, -2 \gamma, \alpha+3 \gamma; v_1, v_2), s_2(\beta, -2
  \gamma, \alpha+3 \gamma; v_1, v_2); v_1, v_2]_{v_2}.
\endaligned\right.\tag 6.30$$
\item
$$\left\{\aligned
 &\{s_1(\alpha, \beta, \gamma; S_1 T S_1(v_1, v_2)), s_2(\alpha, \beta, \gamma;
  S_1 T S_1(v_1, v_2)); v_1, v_2\}_{v_1}\\
=&\{s_1(\alpha, \beta, \gamma; v_1, v_2), s_2(\alpha, \beta,
  \gamma; v_1, v_2); v_1, v_2\}_{v_1},\\
 &\{s_1(\alpha, \beta, \gamma; S_2 T S_2(v_1, v_2)), s_2(\alpha, \beta, \gamma;
  S_2 T S_2(v_1, v_2)); v_1, v_2\}_{v_2}\\
=&\{s_1(\alpha, \beta, \gamma; v_1, v_2), s_2(\alpha, \beta,
  \gamma; v_1, v_2); v_1, v_2\}_{v_2},\\
 &[s_1(\alpha, \beta, \gamma; S_1 T S_1(v_1, v_2)), s_2(\alpha, \beta, \gamma;
  S_1 T S_1(v_1, v_2)); v_1, v_2]_{v_1}\\
=&[s_1(-2 \gamma, \beta, \alpha+3 \gamma; v_1, v_2), s_2(-2
  \gamma, \beta, \alpha+3 \gamma; v_1, v_2); v_1, v_2]_{v_1},\\
 &[s_1(\alpha, \beta, \gamma; S_2 T S_2(v_1, v_2)), s_2(\alpha, \beta, \gamma;
  S_2 T S_2(v_1, v_2)); v_1, v_2]_{v_2}\\
=&[s_1(-2 \gamma, \beta, \alpha+3 \gamma; v_1, v_2), s_2(-2
  \gamma, \beta, \alpha+3 \gamma; v_1, v_2); v_1, v_2]_{v_2}.
\endaligned\right.\tag 6.31$$
\endroster}

  For $0, 1, \lambda_1, \lambda_2, \infty \in {\Bbb P}^{1}({\Bbb
C})$, $\lambda_1, \lambda_2 \notin \{0, 1\}$ and $\lambda_1 \neq
\lambda_2$, set
$$\left\{\aligned
  J_1=J_1(\lambda_1, \lambda_2):=\frac{\lambda_2^2 (\lambda_2-1)^2}
      {\lambda_1^2 (\lambda_1-1)^2 (\lambda_1-\lambda_2)^2},\\
  J_2=J_2(\lambda_1, \lambda_2):=\frac{\lambda_1^2 (\lambda_1-1)^2}
      {\lambda_2^2 (\lambda_2-1)^2 (\lambda_2-\lambda_1)^2}.
\endaligned\right.\tag 6.32$$
The corresponding Picard curve is
$$y^3=x(x-1)(x-\lambda_1)(x-\lambda_2).\tag 6.33$$
The functions $J_1$ and $J_2$ satisfy the following relations:
$$\aligned
 &J_1(\lambda_1, \lambda_2)=J_1(T(\lambda_1, \lambda_2))
 =J_1(S_1(\lambda_1, \lambda_2))=J_1(S_1 T(\lambda_1, \lambda_2))\\
=&J_1(T S_1(\lambda_1, \lambda_2))=J_1(S_1 T S_1(\lambda_1,
  \lambda_2))=\frac{\lambda_2^2 (\lambda_2-1)^2}{\lambda_1^2
  (\lambda_1-1)^2 (\lambda_1-\lambda_2)^2}.
\endaligned\tag 6.34$$
$$\aligned
 &J_2(\lambda_1, \lambda_2)=J_2(T(\lambda_1, \lambda_2))
 =J_2(S_2(\lambda_1, \lambda_2))=J_2(S_2 T(\lambda_1, \lambda_2))\\
=&J_2(T S_2(\lambda_1, \lambda_2))=J_2(S_2 T S_2(\lambda_1,
  \lambda_2))=\frac{\lambda_1^2 (\lambda_1-1)^2}{\lambda_2^2
  (\lambda_2-1)^2 (\lambda_2-\lambda_1)^2}.
\endaligned\tag 6.35$$

\vskip 0.5 cm
\centerline{\bf 7. Four derivatives and nonlinear
                   evolution equations}
\vskip 0.5 cm

  Let
$$(w_1, w_2)=\left(\frac{a_1 u_1+a_2 u_2+a_3}{c_1 u_1+c_2 u_2+c_3},
  \frac{b_1 u_1+b_2 u_2+b_3}{c_1 u_1+c_2 u_2+c_3}\right),$$
where $\gamma=\left(\matrix a_1 & a_2 & a_3\\ b_1 & b_2 & b_3\\
c_1 & c_2 & c_3 \endmatrix\right) \in GL(3, {\Bbb C})$ and
$\Delta=\det(\gamma)$. Note that
$$\vmatrix \format \c \quad & \c\\
  \frac{\partial w_1}{\partial t} &
  \frac{\partial w_2}{\partial t}\\
  \frac{\partial w_1}{\partial x} &
  \frac{\partial w_2}{\partial x}
  \endvmatrix
 =\Delta (c_1 u_1+c_2 u_2+c_3)^{-3}
  \vmatrix \format \c \quad & \c\\
  \frac{\partial u_1}{\partial t} &
  \frac{\partial u_2}{\partial t}\\
  \frac{\partial u_1}{\partial x} &
  \frac{\partial u_2}{\partial x}
  \endvmatrix,$$
$$\vmatrix \format \c \quad & \c\\
  \frac{\partial w_1}{\partial t} &
  \frac{\partial w_2}{\partial t}\\
  \frac{\partial w_1}{\partial y} &
  \frac{\partial w_2}{\partial y}
  \endvmatrix
 =\Delta (c_1 u_1+c_2 u_2+c_3)^{-3}
  \vmatrix \format \c \quad & \c\\
  \frac{\partial u_1}{\partial t} &
  \frac{\partial u_2}{\partial t}\\
  \frac{\partial u_1}{\partial y} &
  \frac{\partial u_2}{\partial y}
  \endvmatrix,$$
$$\vmatrix \format \c \quad & \c\\
  \frac{\partial w_1}{\partial x} &
  \frac{\partial w_2}{\partial x}\\
  \frac{\partial w_1}{\partial y} &
  \frac{\partial w_2}{\partial y}
  \endvmatrix
 =\Delta (c_1 u_1+c_2 u_2+c_3)^{-3}
  \vmatrix \format \c \quad & \c\\
  \frac{\partial u_1}{\partial x} &
  \frac{\partial u_2}{\partial x}\\
  \frac{\partial u_1}{\partial y} &
  \frac{\partial u_2}{\partial y}
  \endvmatrix,$$
we have
$$\frac{\vmatrix \format \c \quad & \c\\
  \frac{\partial w_1}{\partial t} &
  \frac{\partial w_2}{\partial t}\\
  \frac{\partial w_1}{\partial x} &
  \frac{\partial w_2}{\partial x}
  \endvmatrix}{\vmatrix \format \c \quad & \c\\
  \frac{\partial w_1}{\partial x} &
  \frac{\partial w_2}{\partial x}\\
  \frac{\partial w_1}{\partial y} &
  \frac{\partial w_2}{\partial y}
  \endvmatrix}
 =\frac{\vmatrix \format \c \quad & \c\\
  \frac{\partial u_1}{\partial t} &
  \frac{\partial u_2}{\partial t}\\
  \frac{\partial u_1}{\partial x} &
  \frac{\partial u_2}{\partial x}
  \endvmatrix}{\vmatrix \format \c \quad & \c\\
  \frac{\partial u_1}{\partial x} &
  \frac{\partial u_2}{\partial x}\\
  \frac{\partial u_1}{\partial y} &
  \frac{\partial u_2}{\partial y}
  \endvmatrix}, \quad
  \frac{\vmatrix \format \c \quad & \c\\
  \frac{\partial w_1}{\partial t} &
  \frac{\partial w_2}{\partial t}\\
  \frac{\partial w_1}{\partial y} &
  \frac{\partial w_2}{\partial y}
  \endvmatrix}{\vmatrix \format \c \quad & \c\\
  \frac{\partial w_1}{\partial y} &
  \frac{\partial w_2}{\partial y}\\
  \frac{\partial w_1}{\partial x} &
  \frac{\partial w_2}{\partial x}
  \endvmatrix}
 =\frac{\vmatrix \format \c \quad & \c\\
  \frac{\partial u_1}{\partial t} &
  \frac{\partial u_2}{\partial t}\\
  \frac{\partial u_1}{\partial y} &
  \frac{\partial u_2}{\partial y}
  \endvmatrix}{\vmatrix \format \c \quad & \c\\
  \frac{\partial u_1}{\partial y} &
  \frac{\partial u_2}{\partial y}\\
  \frac{\partial u_1}{\partial x} &
  \frac{\partial u_2}{\partial x}
  \endvmatrix}.$$
On the other hand,
$$\left\{\aligned
 \{w_1, w_2; x, y\}_{x} &=\{u_1, u_2; x, y\}_{x},\\
 \{w_1, w_2; x, y\}_{y} &=\{u_1, u_2; x, y\}_{y},\\
   [w_1, w_2; x, y]_{x} &=[u_1, u_2; x, y]_{x},\\
   [w_1, w_2; x, y]_{x} &=[u_1, u_2; x, y]_{y}.
\endaligned\right.$$
This leads us to study the following four evolution equations for
$GL(3)$:
$$\frac{\vmatrix \format \c \quad & \c\\
  \frac{\partial u_1}{\partial t_1} &
  \frac{\partial u_2}{\partial t_1}\\
  \frac{\partial u_1}{\partial x} &
  \frac{\partial u_2}{\partial x}
  \endvmatrix}{\vmatrix \format \c \quad & \c\\
  \frac{\partial u_1}{\partial y} &
  \frac{\partial u_2}{\partial y}\\
  \frac{\partial u_1}{\partial x} &
  \frac{\partial u_2}{\partial x}
  \endvmatrix}
 =\{u_1, u_2; x, y\}_{x}, \quad
  \frac{\vmatrix \format \c \quad & \c\\
  \frac{\partial u_1}{\partial t_1} &
  \frac{\partial u_2}{\partial t_1}\\
  \frac{\partial u_1}{\partial y} &
  \frac{\partial u_2}{\partial y}
  \endvmatrix}{\vmatrix \format \c \quad & \c\\
  \frac{\partial u_1}{\partial y} &
  \frac{\partial u_2}{\partial y}\\
  \frac{\partial u_1}{\partial x} &
  \frac{\partial u_2}{\partial x}
  \endvmatrix}
 =[u_1, u_2; x, y]_{x},$$
$$\frac{\vmatrix \format \c \quad & \c\\
  \frac{\partial u_1}{\partial t_2} &
  \frac{\partial u_2}{\partial t_2}\\
  \frac{\partial u_1}{\partial x} &
  \frac{\partial u_2}{\partial x}
  \endvmatrix}{\vmatrix \format \c \quad & \c\\
  \frac{\partial u_1}{\partial x} &
  \frac{\partial u_2}{\partial x}\\
  \frac{\partial u_1}{\partial y} &
  \frac{\partial u_2}{\partial y}
  \endvmatrix}
 =[u_1, u_2; x, y]_{y}, \quad
  \frac{\vmatrix \format \c \quad & \c\\
  \frac{\partial u_1}{\partial t_2} &
  \frac{\partial u_2}{\partial t_2}\\
  \frac{\partial u_1}{\partial y} &
  \frac{\partial u_2}{\partial y}
  \endvmatrix}{\vmatrix \format \c \quad & \c\\
  \frac{\partial u_1}{\partial x} &
  \frac{\partial u_2}{\partial x}\\
  \frac{\partial u_1}{\partial y} &
  \frac{\partial u_2}{\partial y}
  \endvmatrix}
 =\{u_1, u_2; x, y\}_{y},$$
where the functions $u_1=u_1(x, y; t_1, t_2)$ and $u_2=u_2(x, y;
t_1, t_2)$. The most important property of the above equations is
that they are invariant with respect to the group $GL(3, {\Bbb
C})$ of linear fractional transformations:
$$(u_1, u_2) \mapsto \left(\frac{a_1 u_1+a_2 u_2+a_3}
  {c_1 u_1+c_2 u_2+c_3}, \frac{b_1 u_1+b_2 u_2+b_3}
  {c_1 u_1+c_2 u_2+c_3}\right)$$
with $\gamma=\left(\matrix a_1 & a_2 & a_3\\ b_1 & b_2 & b_3\\
c_1 & c_2 & c_3 \endmatrix\right) \in GL(3, {\Bbb C})$.

{\smc Theorem 7.1 (Main Theorem 4)}. {\it For the system of
partial differential equations:
$$\frac{\vmatrix \format \c \quad & \c\\
  \frac{\partial u_1}{\partial t_1} &
  \frac{\partial u_2}{\partial t_1}\\
  \frac{\partial u_1}{\partial x} &
  \frac{\partial u_2}{\partial x}
  \endvmatrix}{\vmatrix \format \c \quad & \c\\
  \frac{\partial u_1}{\partial y} &
  \frac{\partial u_2}{\partial y}\\
  \frac{\partial u_1}{\partial x} &
  \frac{\partial u_2}{\partial x}
  \endvmatrix}
 =\{u_1, u_2; x, y\}_{x}, \quad
  \frac{\vmatrix \format \c \quad & \c\\
  \frac{\partial u_1}{\partial t_1} &
  \frac{\partial u_2}{\partial t_1}\\
  \frac{\partial u_1}{\partial y} &
  \frac{\partial u_2}{\partial y}
  \endvmatrix}{\vmatrix \format \c \quad & \c\\
  \frac{\partial u_1}{\partial y} &
  \frac{\partial u_2}{\partial y}\\
  \frac{\partial u_1}{\partial x} &
  \frac{\partial u_2}{\partial x}
  \endvmatrix}
 =[u_1, u_2; x, y]_{x},\tag 7.1$$
$$\frac{\vmatrix \format \c \quad & \c\\
  \frac{\partial u_1}{\partial t_2} &
  \frac{\partial u_2}{\partial t_2}\\
  \frac{\partial u_1}{\partial x} &
  \frac{\partial u_2}{\partial x}
  \endvmatrix}{\vmatrix \format \c \quad & \c\\
  \frac{\partial u_1}{\partial x} &
  \frac{\partial u_2}{\partial x}\\
  \frac{\partial u_1}{\partial y} &
  \frac{\partial u_2}{\partial y}
  \endvmatrix}
 =[u_1, u_2; x, y]_{y}, \quad
  \frac{\vmatrix \format \c \quad & \c\\
  \frac{\partial u_1}{\partial t_2} &
  \frac{\partial u_2}{\partial t_2}\\
  \frac{\partial u_1}{\partial y} &
  \frac{\partial u_2}{\partial y}
  \endvmatrix}{\vmatrix \format \c \quad & \c\\
  \frac{\partial u_1}{\partial x} &
  \frac{\partial u_2}{\partial x}\\
  \frac{\partial u_1}{\partial y} &
  \frac{\partial u_2}{\partial y}
  \endvmatrix}
 =\{u_1, u_2; x, y\}_{y},\tag 7.2$$
set
$$\left\{\aligned
  v_1 &:=\{u_1, u_2; x, y\}_{x},\\
  w_1 &:=[u_1, u_2; x, y]_{x},\\
  v_2 &:=\{u_1, u_2; x, y\}_{y},\\
  w_2 &:=[u_1, u_2; x, y]_{y}.
\endaligned\right.\tag 7.3$$
Then $v_1$, $v_2$, $w_1$ and $w_2$ satisfy the following two
equations:
$$\left\{\aligned
  \frac{\partial w_1}{\partial t_2}+\frac{\partial
  v_2}{\partial t_1}-v_1 \frac{\partial v_2}{\partial y}-
  v_2 \frac{\partial w_1}{\partial x}+w_1 \frac{\partial
  v_2}{\partial x}+w_2 \frac{\partial w_1}{\partial y} &=0,\\
  \frac{\partial w_2}{\partial t_1}+\frac{\partial
  v_1}{\partial t_2}-v_1 \frac{\partial w_2}{\partial y}-
  v_2 \frac{\partial v_1}{\partial x}+w_1 \frac{\partial
  w_2}{\partial x}+w_2 \frac{\partial v_1}{\partial y} &=0.
\endaligned\right.\tag 7.4$$}

{\it Proof}. The above four equations imply that
$$\frac{\partial u_1}{\partial t_1}
 =\frac{\partial^2 u_1}{\partial x^2}-2
  \frac{\vmatrix \format \c \quad & \c\\
  \frac{\partial u_1}{\partial x} &
  \frac{\partial u_2}{\partial x}\\
  \frac{\partial^2 u_1}{\partial x \partial y} &
  \frac{\partial^2 u_2}{\partial x \partial y}
  \endvmatrix}{\vmatrix \format \c \quad & \c\\
  \frac{\partial u_1}{\partial x} &
  \frac{\partial u_2}{\partial x}\\
  \frac{\partial u_1}{\partial y} &
  \frac{\partial u_2}{\partial y}
  \endvmatrix} \frac{\partial u_1}{\partial x}, \quad
  \frac{\partial u_2}{\partial t_1}
 =\frac{\partial^2 u_2}{\partial x^2}-2
  \frac{\vmatrix \format \c \quad & \c\\
  \frac{\partial u_1}{\partial x} &
  \frac{\partial u_2}{\partial x}\\
  \frac{\partial^2 u_1}{\partial x \partial y} &
  \frac{\partial^2 u_2}{\partial x \partial y}
  \endvmatrix}{\vmatrix \format \c \quad & \c\\
  \frac{\partial u_1}{\partial x} &
  \frac{\partial u_2}{\partial x}\\
  \frac{\partial u_1}{\partial y} &
  \frac{\partial u_2}{\partial y}
  \endvmatrix} \frac{\partial u_2}{\partial x},$$
$$\frac{\partial u_1}{\partial t_2}
 =\frac{\partial^2 u_1}{\partial y^2}-2
  \frac{\vmatrix \format \c \quad & \c\\
  \frac{\partial u_1}{\partial y} &
  \frac{\partial u_2}{\partial y}\\
  \frac{\partial^2 u_1}{\partial x \partial y} &
  \frac{\partial^2 u_2}{\partial x \partial y}
  \endvmatrix}{\vmatrix \format \c \quad & \c\\
  \frac{\partial u_1}{\partial y} &
  \frac{\partial u_2}{\partial y}\\
  \frac{\partial u_1}{\partial x} &
  \frac{\partial u_2}{\partial x}
  \endvmatrix} \frac{\partial u_1}{\partial y}, \quad
  \frac{\partial u_2}{\partial t_2}
 =\frac{\partial^2 u_2}{\partial y^2}-2
  \frac{\vmatrix \format \c \quad & \c\\
  \frac{\partial u_1}{\partial y} &
  \frac{\partial u_2}{\partial y}\\
  \frac{\partial^2 u_1}{\partial x \partial y} &
  \frac{\partial^2 u_2}{\partial x \partial y}
  \endvmatrix}{\vmatrix \format \c \quad & \c\\
  \frac{\partial u_1}{\partial y} &
  \frac{\partial u_2}{\partial y}\\
  \frac{\partial u_1}{\partial x} &
  \frac{\partial u_2}{\partial x}
  \endvmatrix} \frac{\partial u_2}{\partial y}.$$

  The equations
$$v_1=\{u_1, u_2; x, y\}_{x}, \quad w_1=[u_1, u_2; x, y]_{x}$$
give that
$$\frac{\partial^2 u_1}{\partial x^2}
 =\vmatrix \format \c \quad & \c\\
  v_1 & w_1\\
  \frac{\partial u_1}{\partial x} &
  \frac{\partial u_1}{\partial y}
  \endvmatrix+2 \frac{\vmatrix \format \c \quad & \c\\
  \frac{\partial u_1}{\partial x} &
  \frac{\partial u_2}{\partial x}\\
  \frac{\partial^2 u_1}{\partial x \partial y} &
  \frac{\partial^2 u_2}{\partial x \partial y}
  \endvmatrix}{\vmatrix \format \c \quad & \c\\
  \frac{\partial u_1}{\partial x} &
  \frac{\partial u_2}{\partial x}\\
  \frac{\partial u_1}{\partial y} &
  \frac{\partial u_2}{\partial y}
  \endvmatrix} \frac{\partial u_1}{\partial x},$$
$$\frac{\partial^2 u_2}{\partial x^2}
 =\vmatrix \format \c \quad & \c\\
  v_1 & w_1\\
  \frac{\partial u_2}{\partial x} &
  \frac{\partial u_2}{\partial y}
  \endvmatrix+2 \frac{\vmatrix \format \c \quad & \c\\
  \frac{\partial u_1}{\partial x} &
  \frac{\partial u_2}{\partial x}\\
  \frac{\partial^2 u_1}{\partial x \partial y} &
  \frac{\partial^2 u_2}{\partial x \partial y}
  \endvmatrix}{\vmatrix \format \c \quad & \c\\
  \frac{\partial u_1}{\partial x} &
  \frac{\partial u_2}{\partial x}\\
  \frac{\partial u_1}{\partial y} &
  \frac{\partial u_2}{\partial y}
  \endvmatrix} \frac{\partial u_2}{\partial x}.$$

  The equations
$$v_2=\{u_1, u_2; x, y\}_{y}, \quad w_2=[u_1, u_2; x, y]_{y}$$
give that
$$\frac{\partial^2 u_1}{\partial y^2}
 =\vmatrix \format \c \quad & \c\\
  v_2 & w_2\\
  \frac{\partial u_1}{\partial y} &
  \frac{\partial u_1}{\partial x}
  \endvmatrix+2 \frac{\vmatrix \format \c \quad & \c\\
  \frac{\partial u_1}{\partial y} &
  \frac{\partial u_2}{\partial y}\\
  \frac{\partial^2 u_1}{\partial x \partial y} &
  \frac{\partial^2 u_2}{\partial x \partial y}
  \endvmatrix}{\vmatrix \format \c \quad & \c\\
  \frac{\partial u_1}{\partial y} &
  \frac{\partial u_2}{\partial y}\\
  \frac{\partial u_1}{\partial x} &
  \frac{\partial u_2}{\partial x}
  \endvmatrix} \frac{\partial u_1}{\partial y},$$
$$\frac{\partial^2 u_2}{\partial y^2}
 =\vmatrix \format \c \quad & \c\\
  v_2 & w_2\\
  \frac{\partial u_2}{\partial y} &
  \frac{\partial u_2}{\partial x}
  \endvmatrix+2 \frac{\vmatrix \format \c \quad & \c\\
  \frac{\partial u_1}{\partial y} &
  \frac{\partial u_2}{\partial y}\\
  \frac{\partial^2 u_1}{\partial x \partial y} &
  \frac{\partial^2 u_2}{\partial x \partial y}
  \endvmatrix}{\vmatrix \format \c \quad & \c\\
  \frac{\partial u_1}{\partial y} &
  \frac{\partial u_2}{\partial y}\\
  \frac{\partial u_1}{\partial x} &
  \frac{\partial u_2}{\partial x}
  \endvmatrix} \frac{\partial u_2}{\partial y}.$$
This implies that
$$\frac{\partial u_1}{\partial t_1}
 =\vmatrix \format \c \quad & \c\\
  v_1 & w_1\\
  \frac{\partial u_1}{\partial x} &
  \frac{\partial u_1}{\partial y}
  \endvmatrix, \quad
  \frac{\partial u_2}{\partial t_1}
 =\vmatrix \format \c \quad & \c\\
  v_1 & w_1\\
  \frac{\partial u_2}{\partial x} &
  \frac{\partial u_2}{\partial y}
  \endvmatrix.$$
$$\frac{\partial u_1}{\partial t_2}
 =\vmatrix \format \c \quad & \c\\
  v_2 & w_2\\
  \frac{\partial u_1}{\partial y} &
  \frac{\partial u_1}{\partial x}
  \endvmatrix, \quad
  \frac{\partial u_2}{\partial t_2}
 =\vmatrix \format \c \quad & \c\\
  v_2 & w_2\\
  \frac{\partial u_2}{\partial y} &
  \frac{\partial u_2}{\partial x}
  \endvmatrix.$$

  Note that
$$\aligned
 &\frac{\partial^2 u_1}{\partial t_1 \partial t_2}
 =\frac{\partial}{\partial t_2}\left(\frac{\partial u_1}{\partial t_1}\right)
 =\frac{\partial}{\partial t_2}\left(v_1 \frac{\partial u_1}{\partial y}-
  w_1 \frac{\partial u_1}{\partial x}\right)\\
=&\frac{\partial v_1}{\partial t_2} \frac{\partial u_1}{\partial
  y}-\frac{\partial w_1}{\partial t_2} \frac{\partial u_1}{\partial x}
  +v_1 \frac{\partial}{\partial y}\left(\frac{\partial u_1}{\partial t_2}\right)
  -w_1 \frac{\partial}{\partial x}\left(\frac{\partial u_1}{\partial t_2}\right).
\endaligned$$
Here,
$$\frac{\partial}{\partial y} \left(\frac{\partial u_1}{\partial t_2}\right)
 =\frac{\partial}{\partial y} \left(v_2 \frac{\partial u_1}{\partial x}
  -w_2 \frac{\partial u_1}{\partial y}\right)
 =\frac{\partial v_2}{\partial y} \frac{\partial u_1}{\partial x}-
  \frac{\partial w_2}{\partial y} \frac{\partial u_1}{\partial y}+
  v_2 \frac{\partial^2 u_1}{\partial x \partial y}-w_2
  \frac{\partial^2 u_1}{\partial y^2},$$
$$\frac{\partial}{\partial x} \left(\frac{\partial u_1}{\partial t_2}\right)
 =\frac{\partial}{\partial x} \left(v_2 \frac{\partial u_1}{\partial x}
  -w_2 \frac{\partial u_1}{\partial y}\right)
 =\frac{\partial v_2}{\partial x} \frac{\partial u_1}{\partial x}-
  \frac{\partial w_2}{\partial x} \frac{\partial u_1}{\partial y}+
  v_2 \frac{\partial^2 u_1}{\partial x^2}-w_2
  \frac{\partial^2 u_1}{\partial x \partial y}.$$
Hence,
$$\aligned
  \frac{\partial^2 u_1}{\partial t_1 \partial t_2}
=&\left(-\frac{\partial w_1}{\partial t_2}+v_1 \frac{\partial
  v_2}{\partial y}-w_1 \frac{\partial v_2}{\partial x}\right)
  \frac{\partial u_1}{\partial x}+\left(\frac{\partial v_1}
  {\partial t_2}-v_1 \frac{\partial w_2}{\partial y}+w_1
  \frac{\partial w_2}{\partial x}\right) \frac{\partial u_1}
  {\partial y}+\\
+&(v_1 v_2+w_1 w_2) \frac{\partial^2 u_1}{\partial x \partial y}
  -v_2 w_1 \frac{\partial^2 u_1}{\partial x^2}-v_1 w_2
  \frac{\partial^2 u_1}{\partial y^2}.
\endaligned$$
On the other hand,
$$\aligned
 &\frac{\partial^2 u_1}{\partial t_1 \partial t_2}
 =\frac{\partial}{\partial t_1}\left(\frac{\partial u_1}{\partial t_2}\right)
 =\frac{\partial}{\partial t_1}\left(v_2 \frac{\partial u_1}{\partial x}-
  w_2 \frac{\partial u_1}{\partial y}\right)\\
=&\frac{\partial v_2}{\partial t_1} \frac{\partial u_1}{\partial
  x}-\frac{\partial w_2}{\partial t_1} \frac{\partial u_1}{\partial y}
  +v_2 \frac{\partial}{\partial x}\left(\frac{\partial u_1}{\partial t_1}\right)
  -w_2 \frac{\partial}{\partial y}\left(\frac{\partial u_1}{\partial t_1}\right).
\endaligned$$
Here,
$$\frac{\partial}{\partial x} \left(\frac{\partial u_1}{\partial t_1}\right)
 =\frac{\partial}{\partial x} \left(v_1 \frac{\partial u_1}{\partial y}
  -w_1 \frac{\partial u_1}{\partial x}\right)
 =\frac{\partial v_1}{\partial x} \frac{\partial u_1}{\partial y}-
  \frac{\partial w_1}{\partial x} \frac{\partial u_1}{\partial x}+
  v_1 \frac{\partial^2 u_1}{\partial x \partial y}-w_1
  \frac{\partial^2 u_1}{\partial x^2},$$
$$\frac{\partial}{\partial y} \left(\frac{\partial u_1}{\partial t_1}\right)
 =\frac{\partial}{\partial y} \left(v_1 \frac{\partial u_1}{\partial y}
  -w_1 \frac{\partial u_1}{\partial x}\right)
 =\frac{\partial v_1}{\partial y} \frac{\partial u_1}{\partial y}-
  \frac{\partial w_1}{\partial y} \frac{\partial u_1}{\partial x}+
  v_1 \frac{\partial^2 u_1}{\partial y^2}-w_1
  \frac{\partial^2 u_1}{\partial x \partial y}.$$
Hence,
$$\aligned
  \frac{\partial^2 u_1}{\partial t_1 \partial t_2}
=&\left(\frac{\partial v_2}{\partial t_1}-v_2 \frac{\partial
  w_1}{\partial x}+w_2 \frac{\partial w_1}{\partial y}\right)
  \frac{\partial u_1}{\partial x}+\left(-\frac{\partial w_2}
  {\partial t_1}+v_2 \frac{\partial v_1}{\partial x}-w_2
  \frac{\partial v_1}{\partial y}\right) \frac{\partial u_1}
  {\partial y}+\\
+&(v_1 v_2+w_1 w_2) \frac{\partial^2 u_1}{\partial x \partial y}
  -v_2 w_1 \frac{\partial^2 u_1}{\partial x^2}-v_1 w_2
  \frac{\partial^2 u_1}{\partial y^2}.
\endaligned$$
Comparing with the different expressions of $\frac{\partial^2
u_1}{\partial t_1 \partial t_2}$, we get the following equation:
$$\aligned
 &\left(\frac{\partial w_1}{\partial t_2}+\frac{\partial
  v_2}{\partial t_1}-v_1 \frac{\partial v_2}{\partial y}-
  v_2 \frac{\partial w_1}{\partial x}+w_1 \frac{\partial
  v_2}{\partial x}+w_2 \frac{\partial w_1}{\partial y}\right)
  \frac{\partial u_1}{\partial x}\\
=&\left(\frac{\partial w_2}{\partial t_1}+\frac{\partial
  v_1}{\partial t_2}-v_1 \frac{\partial w_2}{\partial y}-
  v_2 \frac{\partial v_1}{\partial x}+w_1 \frac{\partial
  w_2}{\partial x}+w_2 \frac{\partial v_1}{\partial y}\right)
  \frac{\partial u_1}{\partial y}.
\endaligned$$
Similarly, comparing with the different expressions of
$\frac{\partial^2 u_2}{\partial t_1 \partial t_2}$, we get the
other equation:
$$\aligned
 &\left(\frac{\partial w_1}{\partial t_2}+\frac{\partial
  v_2}{\partial t_1}-v_1 \frac{\partial v_2}{\partial y}-
  v_2 \frac{\partial w_1}{\partial x}+w_1 \frac{\partial
  v_2}{\partial x}+w_2 \frac{\partial w_1}{\partial y}\right)
  \frac{\partial u_2}{\partial x}\\
=&\left(\frac{\partial w_2}{\partial t_1}+\frac{\partial
  v_1}{\partial t_2}-v_1 \frac{\partial w_2}{\partial y}-
  v_2 \frac{\partial v_1}{\partial x}+w_1 \frac{\partial
  w_2}{\partial x}+w_2 \frac{\partial v_1}{\partial y}\right)
  \frac{\partial u_2}{\partial y}.
\endaligned$$
Note that $\frac{\partial(u_1, u_2)}{\partial(x, y)} \neq 0$. The
above two equations imply that
$$\left\{\aligned
  \frac{\partial w_1}{\partial t_2}+\frac{\partial
  v_2}{\partial t_1}-v_1 \frac{\partial v_2}{\partial y}-
  v_2 \frac{\partial w_1}{\partial x}+w_1 \frac{\partial
  v_2}{\partial x}+w_2 \frac{\partial w_1}{\partial y} &=0,\\
  \frac{\partial w_2}{\partial t_1}+\frac{\partial
  v_1}{\partial t_2}-v_1 \frac{\partial w_2}{\partial y}-
  v_2 \frac{\partial v_1}{\partial x}+w_1 \frac{\partial
  w_2}{\partial x}+w_2 \frac{\partial v_1}{\partial y} &=0.
\endaligned\right.$$
$\qquad \qquad \qquad \qquad \qquad \qquad \qquad \qquad \qquad
 \qquad \qquad \qquad \qquad \qquad \qquad \qquad \qquad \qquad
 \quad \boxed{}$

  In fact, The above two equations are invariant under the
action of
$$\left\{\aligned
  v_1(x, y; t_1, t_2) &\mapsto v_1(x-a_1 t_1, y-b_1 t_2; t_1, t_2),\\
  v_2(x, y; t_1, t_2) &\mapsto v_2(x-a_2 t_1, y-b_2 t_2; t_1, t_2),\\
  w_1(x, y; t_1, t_2) &\mapsto w_1(x-a_1 t_1, y-b_1 t_2; t_1, t_2)+a_2,\\
  w_2(x, y; t_1, t_2) &\mapsto w_2(x-a_2 t_1, y-b_2 t_2; t_1,
  t_2)+b_1.
\endaligned\right.\tag 7.5$$

\vskip 2.0 cm

{\smc Department of Mathematics, Peking University}

{\smc Beijing 100871, P. R. China}

{\it E-mail address}: yanglei\@math.pku.edu.cn
\vskip 1.5 cm
\Refs

\item{[AK]} {\smc P. Appell and J. Kamp\'{e} de F\'{e}riet},
            {\it Fonctions hyperg\'{e}om\'{e}triques et
            hypersph\'{e}riques, polynomes d'Hermite},
            Gauthier-Villars, 1926.

\item{[BMS]} {\smc A. A. Beilinson, Yu. I. Manin and V. V. Schechtman},
             Sheaves of the Virasoro and Neveu-Schwarz algebras,
             In: $K$-Theory, Arithmetic and Geometry, Yu. I. Manin
             Ed., Lecture Notes in Math. {\bf 1289}, 52-66,
             Springer-Verlag, 1987.

\item{[CW1]} {\smc P. B. Cohen and J. Wolfart}, Algebraic
            Appell-Lauricella functions, Analysis {\bf 12} (1992),
            359-376.

\item{[CW2]} {\smc P. B. Cohen and J. Wolfart}, Fonctions
             hyperg\'{e}om\'{e}triques en plusieurs variables et
             espaces des modules de vari\'{e}t\'{e}s ab\'{e}liennes,
             Ann. Sci. \'{E}c. Norm. Sup. {\bf 26} (1993), 665-690.

\item{[De]} {\smc R. Dedekind}, Schreiben an Herrn Borchardt
            \"{u}ber die Theorie der elliptischen Modul-Functionen,
            J. Reine Angew. Math. {\bf 83} (1877), 265-292.

\item{[DM1]} {\smc P. Deligne and G. D. Mostow}, Monodromy of
             hypergeometric functions and non-lattice integral
             monodromy, Publ. Math. I.H.E.S. {\bf 63} (1986), 5-89.

\item{[DM2]} {\smc P. Deligne and G. D. Mostow}, {\it Commensurabilities
             among Lattices in $PU(1, n)$}, Annals of Math. Studies
             {\bf 132}, Princeton University Press, 1993.

\item{[Do]} {\smc I. Dorfman}, {\it Dirac Structures and Integrability
            of Nonlinear Evolution Equations}, John Wiley \& Sons,
            1993.

\item{[GRS]} {\smc S. Gelbart, J. Rogawski and D. Soudry}, Endoscopy,
            theta-liftings, and period integrals for the unitary
            group in three variables, Ann. of Math. {\bf 145} (1997),
            419-476.

\item{[Gir]} {\smc G. Giraud}, Sur certaines fonctions automorphes
           de deux variables, Ann. Ec. Norm. {\bf (3)}, {\bf 38}
           (1921), 43-164.

\item{[G]} {\smc J. J. Gray}, {\it Linear Differential Equations
            and Group Theory from Riemann to Poincar\'{e}}, Second
            Edition, Birkh\"{a}user, 2000.

\item{[H]} {\smc K. Heegner}, Reduzierbare Abelsche Integrale und
           transformierbare automorphe Funktionen, Math. Ann. {\bf
           131} (1956), 87-140.

\item{[He]} {\smc C. Hermite}, {\it \OE uvres de Charles Hermite},
            Vol. II, E. Picard Ed., Gauthier-Villars, 1908.

\item{[Hi]} {\smc E. Hille}, {\it Ordinary Differential Equations
             in the Complex Domain}, John Wiley \& Sons, Inc.,
             1976.

\item{[Ho1]} {\smc R. P. Holzapfel}, {\it Geometry and arithmetic
            around Euler partial differential equations}, Dt. Verl.
            d. Wiss., Berlin/Reidel, Dordrecht 1986.

\item{[Ho2]} {\smc R. P. Holzapfel}, {\it The ball and some Hilbert
            problems}, Lectures in Mathematics, ETH Z\"{u}rich,
            Birkh\"{a}user, 1995.

\item{[Kl1]} {\smc F. Klein}, \"{U}ber die Transformation der
            elliptischen Functionen und die Aufl\"{o}sung der
            Gleichungen f\"{u}nften Grades, Math. Ann. {\bf 14}
            (1879), 111-172.

\item{[Kl2]} {\smc F. Klein}, {\it Lectures on the Icosahedron and
             the Solution of Equations of the Fifth Degree},
             Translated by G. G. Morrice, second and revised edition,
             Dover Publications, Inc., 1956.

\item{[Kl3]} {\smc F. Klein}, {\it Vorlesungen \"{u}ber die
             Entwicklung der Mathematik im 19. Jahrhundert},
             Berlin, 1926.

\item{[KN]} {\smc I. M. Krichever and S. P. Novikov}, Holomorphic
            bundles on algebraic curves and nonlinear equations,
            Uspekhi Mat. Nauk {\bf 35} (1980), 47-68.

\item{[L]} {\smc R. P. Langlands}, Some contemporary problems with
            origins in the Jugendtraum, {\it Mathematical
            developments arising from Hilbert problems} (Proc.
            Sympos. Pure Math., Vol. XXVIII, Northern Illinois
            Univ., De Kalb, Ill., 1974), pp. 401-418, Amer. Math.
            Soc., Providence, R. I., 1976.

\item{[LR]} {\smc R. P. Langlands and D. Ramakrishnan}, {\it The
            Zeta Functions of Picard Modular Surfaces}, Les
            Publications CRM, Universit\'{e} de Montr\'{e}al,
            1992.

\item{[Lem]} {\smc F. Lemmermeyer}, {\it Reciprocity Laws from
            Euler to Eisenstein}, Springer Monographs in
            Mathematics, Springer-Verlag, 2000.

\item{[M]} {\smc G. D. Mostow}, Generalized Picard lattices arising
            from half-integral conditions, Publ. Math. I.H.E.S. {\bf 63}
            (1986), 91-106.

\item{[O]} {\smc A. I. Ovseevi\v{c}}, Abelian extensions of fields
            of CM type, Funkcional. Anal. i Prilo\v{z}en. {\bf 8}
            (1974), no.1, 16-24.

\item{[Pa1]} {\smc E. Papperitz}, Untersuchungen \"{u}ber die algebraische
             Transformation der hypergeometrischen Functionen,
             Math. Ann. {\bf 27} (1886), 315-356.

\item{[Pa2]} {\smc E. Papperitz}, Ueber die Darstellung der
             hypergeometrischen Transcendenten durch eindeutige
             Functionen, Math. Ann. {\bf 34} (1889), 247-296.

\item{[Pi1]} {\smc E. Picard}, Sur des fonctions de deux variables
            ind\'{e}pendantes analogues aux fonctions modulaires,
            Acta Math. {\bf 2} (1883), 114-135.

\item{[Pi2]} {\smc E. Picard}, Sur les formes quadratiques ternaires
            ind\'{e}finies \`{a} ind\'{e}termin\'{e}es conjugu\'{e}es
            et sur les fonctions hyperfuchsiennes correspondantes,
            Acta Math. {\bf 5} (1884), 121-184.

\item{[Pi3]} {\smc E. Picard}, Sur les fonctions hyperfuchsiennes
             provenant des s\'{e}ries hyperg\'{e}om\'{e}triques de
             deux variables, Ann. Ecole Norm. Sup. 3$^{e}$ s\'{e}rie
             {\bf 62} (1885), 357-384.

\item{[Pi4]} {\smc E. Picard}, M\'{e}moire sur la th\'{e}orie des
             fonctions alg\'{e}briques de deux variables, J. Math.
             Pures Appl. 4$^{e}$ s\'{e}rie {\bf 5} (1889), 135-319.

\item{[P]} {\smc H. Poincar\'{e}}, {\it Papers on Fuchsian Functions},
            Translated by J. Stillwell, Springer-Verlag, 1985.

\item{[Po]} {\smc J. Polchinski}, {\it String Theory}, Vol. I, II,
            Cambridge University Press, 1998.

\item{[R]} {\smc J. D. Rogawski}, {\it Automorphic Representations
            of Unitary Groups in Three Variables}, Annals of Math.
            Studies {\bf 123}, Princeton University Press, 1990.

\item{[Sch]} {\smc H. A. Schwarz}, \"{U}ber diejenigen F\"{a}lle, in
              welchen die Gaussische hypergeometrische Reihe eine
              algebraische Function ihres vierten Elementes
              darstellt, J. Reine Angew. Math. {\bf 75} (1872),
              292-335.

\item{[Sh1]} {\smc G. Shimura}, On purely transcendental fields of
             automorphic functions of several variables, Osaka J.
             Math. {\bf 1} (1964), 1-14.

\item{[Sh2]} {\smc G. Shimura}, Arithmetic of unitary groups,
             Ann. of Math. {\bf 79} (1964), 369-409.

\item{[Sh3]} {\smc G. Shimura}, {\it Introduction to The
             Arithmetic Theory of Automorphic Functions}, Publ.
             Math. Soc. Japan No. {\bf 11}, Iwanami Shoten and
             Princeton University Press, 1971.

\item{[Sh4]} {\smc G. Shimura}, The arithmetic of automorphic
             forms with respect to a unitary group, Ann. of Math.
             {\bf 107} (1978), 569-605.

\item{[Sh5]} {\smc G. Shimura}, {\it Abelian Varieties with
             Complex Multiplication and Modular Functions},
             Princeton Mathematical Series {\bf 46}, Princeton
             University Press, 1998.

\item{[SS]} {\smc V. V. Sokolov and A. B. Shabat}, Classification
            of integrable evolution equations, In: Soviet Scientific
            Reviews, Section C: Mathematical Physics Reviews, Vol.
            {\bf 4}, 221-280, Harwood, New York, 1984.

\item{[V]} {\smc S. G. Vl\u{a}du\c{t}}, {\it Kronecker's Jugendtraum
            and Modular Functions}, Translated from the Russian by
            M. Tsfasman, Studies in the Development of Modern
            Mathematics, Vol. {\bf 2}, Gordon and Breach Science
            Publishers, New York, 1991.

\item{[W]} {\smc H. Weber}, {\it Lehrbuch der Algebra}, Vol. III,
           Chelsea Publishing Company, New York, 1961.

\item{[Wei1]} {\smc J. Weiss}, The Painlev\'{e} property for partial
             differential equations. II: B\"{a}cklund transformation,
             Lax pairs, and the Schwarzian derivative, J. Math. Phys.
             {\bf 24} (1983), 1405-1413.

\item{[Wei2]} {\smc J. Weiss}, B\"{a}cklund transformation and the
              Painlev\'{e} property, J. Math. Phys. {\bf 27} (1986),
              1293-1305.

\item{[Wil]} {\smc G. Wilson}, On the quasi-Hamiltonian formalism
              of the KdV equation, Phys. Lett. A {\bf 132} (1988),
              445-450.

\item{[Wi]} {\smc E. Witten}, Quantum field theory, Grassmannians,
            and algebraic curves, Comm. Math. Phys. {\bf 113}
            (1988), 529-600.

\endRefs
\end{document}